\documentclass[a4paper,12pt]{article}
\usepackage[T1]{fontenc}
\usepackage[latin1]{inputenc}
\usepackage{theorem} 
\usepackage{amsmath} 
\usepackage{amsfonts}
\usepackage{amscd} 
\usepackage{amssymb} 
\usepackage{enumerate}
\usepackage[francais]{babel} 
\usepackage[all]{xy}
\usepackage{index}
\makeatletter
\@addtoreset{equation}{subsection}
\@addtoreset{equation}{section}

\makeatother
\theoremstyle{break}
\newtheorem{thm}{Th\'eor\`eme}[section]
\newtheorem{prop}[thm]{Proposition}
\newtheorem{lemme}[thm]{Lemme}
\newtheorem{cor}[thm]{Corollaire}
\newtheorem{question}[thm]{Question}
\newtheorem{hyp}[thm]{Hypoth\`ese}

\newtheorem{defi}[thm]{D\'efinition}
\newtheorem{notas}[thm]{Notations}
\newtheorem{nota}[thm]{Notation}

\newenvironment{demo}
{\par {\em D\'emonstration :}{\hspace{0.02\textwidth}}}{\hfill $\Box$\par}
\newenvironment{ex}
{\refstepcounter{thm}\noindent{\em Exemple \arabic{section}.\arabic{thm} :}{\hspace{0.02\textwidth}}}{\hfill $\Box$\par}
\newenvironment{rem}
{\refstepcounter{thm}\noindent{\em Remarque \arabic{section}.\arabic{thm} :}{\hspace{0.02\textwidth}}}{\hfill $\Box$\par}
\newenvironment{exs}
{\refstepcounter{thm}\noindent{\em Exemples \arabic{section}.\arabic{thm} :}{\hspace{0.02\textwidth}}}{\hfill $\Box$\par}
\newenvironment{rems}
{\refstepcounter{thm}\noindent{\em Remarques \arabic{section}.\arabic{thm} :}{\hspace{0.02\textwidth}}}{\hfill $\Box$\par}
\newenvironment{conv}
{\refstepcounter{thm}\noindent{\em Convention \arabic{section}.\arabic{thm} :}{\hspace{0.02\textwidth}}}{\hfill $\Box$\par}
\makeatletter
\@addtoreset{equation}{section}

\makeatother
\newcommand{\termin}[1]{{\em #1}}
\def\A{\textbf{A}}\def\C{\textbf{C}}
\def\F{\textbf{F}}\def\G{\textbf{G}}
\def\K{\textbf{K}}\def\L{\textbf{L}}
\def\N{\textbf{N}}
\def\Q{\textbf{Q}}\def\R{\textbf{R}}
 
\def\V{\textbf{V}}\def\X{\textbf{X}}
\def\Z{\textbf{Z}}

\def\bU{\boldsymbol{U}}

\def\bs{\boldsymbol{s}}\def\bz{\boldsymbol{z}}
\def\cC{{\cal C}}
\def\cG{{\cal G}}\def\cH{{\cal H}}
\def\cK{{\cal K}}\def\cL{{\cal L}}
\def\cO{{\cal O}}\def\cP{{\cal P}}
\def\cT{{\cal T}}
\def\cV{{\cal V}}

\DeclareMathAlphabet{\eulercal}{U}{eus}{m}{n}

\def\ecC{{\eulercal C}}\def\ecD{{\eulercal D}}
\def\ecF{{\eulercal F}}
\def\ecG{{\eulercal G}}\def\ecH{{\eulercal H}}

\def\ecK{{\eulercal K}}\def\ecL{{\eulercal L}}
\def\ecM{{\eulercal M}}
\def\ecO{{\eulercal O}}
\def\ecR{{\eulercal R}}
\def\ecT{{\eulercal T}}
\def\ecU{{\eulercal U}}\def\ecV{{\eulercal V}}

\def \frf {{\mathfrak f}}

\def \frh {{\mathfrak h}}
\def \frH {{\mathfrak H}}
\def \frQ {{\mathfrak Q}}

\def \longto {\longrightarrow}
\def \isom {\overset{\sim}{\to}}
\def \longisom {\overset{\sim}{\longto}}

\DeclareMathOperator{\Pic}{Pic}

\DeclareMathOperator{\Br}{Br}

\DeclareMathOperator{\Coker}{Coker}
\DeclareMathOperator{\Aut}{Aut}

\DeclareMathOperator{\rg}{rg}

\DeclareMathOperator{\PL}{PL}
\DeclareMathOperator{\Ker}{Ker}

\DeclareMathOperator{\Hom}{Hom}

\DeclareMathOperator{\Sup}{Sup}

\DeclareMathOperator{\Max}{Max}

\DeclareMathOperator{\Spec}{Spec}

\DeclareMathOperator{\Gal}{Gal}
\DeclareMathOperator{\cla}{cl}
\DeclareMathOperator{\Res}{Res}
\DeclareMathOperator{\disc}{disc}

\DeclareMathOperator{\Frob}{Fr}
\DeclareMathOperator{\NS}{NS}
\DeclareMathOperator{\Id}{Id}
\let\leq\leqslant
\let\geq\geqslant
\newcommand{\sumu}[1]{\underset{#1}{\sum}}
\newcommand{\produ}[1]{\underset{#1}{\prod}}

\newcommand{\oplusu}[1]{\underset{#1}{\oplus}}

\newcommand{\otimesu}[1]{\underset{#1}{\otimes}}

\newcommand{\Maxu}[1]{\underset{#1}{\Max}}
\def \eps {\varepsilon}
\def \vide {\varnothing}
 
\def \eqdef {\overset{\text{{\tiny{d\'ef}}}}{=}}
\def \grsym {\mathfrak{S}}

\def \ind {\mathbf{1}}

\def \wt {\widetilde}
\def \wti {\wt{\imath}}

\newcommand{\acc}[2]{\left\langle #1\, ,\,#2 \right\rangle} 
\newcommand{\norm}[1]{\left|\left| #1 \right|\right|} 
 
\newcommand{\abs}[1]{\left| #1 \right|} 
\newcommand{\card}[1]{\left[#1\right]} 
\newcommand{\adh}[1]{\overline{#1}} 
\newcommand{\map}[4]{
\begin{array}{rcl}
#1 & \longto & #2 \\
#3 & \longmapsto & #4 
\end{array} 
}

\def \beq {\begin{equation}}
\def \eeq {\end{equation}}
\def \bem {\begin{multline}}
\def \courbe {\cC}
\def \xs {X_{\Sigma}}
\def \Ps {P_{\Sigma}}

\def \ceff {C_{\text{eff}}}
\def \zetac {\zeta_{\courbe}}

\DeclareMathOperator{\intrel}{intrel}
\def \indic {\Xi}

\def \sep {\text{s}}
\newcommand{\tube}[1]{\cT\left(#1\right)} 
\newfont{\cyr}{wncyr10}
\def \cha {\text{{\cyr W}}}
\newcommand{\dualtop}[1]{\left(#1\right)^{\ast}} 
\newcommand{\dualt}[1]{#1^{\ast}} 
\def \xsl {X_{\Sigma,L}}
\def \fourier {\ecF}
\def \can {\cK}
\def \gt {{\mathfrak g}_T}

\def \KT {\ecK_T}
\def \DT   {\ecD_{T}}
\def \CT   {\ecC_{T}}
\def \DTNS {\ecD_{T_{\NS}}}
\def \CTNS {\ecC_{T_{\NS}}}
\def \DTqd {\ecD_{T_{\Ps }}}
\def \CTqd   {\ecC_{T_{\Ps }}}
\def \val {\cP}
\newcommand{\placesde}[1]{\val_{#1}}
\def \OK {\cO_K}
\def \Ov {\cO_v}
\def \V {\ecV}
\def \OV {\cO_{\V}}
\def \OS {\cO_S}
\newcommand{\qde}[1]{q_{\text{\tiny {\it #1}}}}
\newcommand{\dde}[1]{d_{\text{\tiny {\it #1}}}}
\def \qk {\qde{K}}
\def \dk {\dde{K}}
\newcommand{\ade}[1]{\A_{#1}}
\def \ak {\ade{K}}
\def \al {\ade{L}}
\def \aka {\ade{\Ka}}

\def \roa {\rho_{\alpha}}
\def \roi {\rho_{i}}
\def \chia {\chi_{\alpha}}
\def \chib {\chi_{\beta}}
\def \Ka {K_{\alpha}}
\def \da {d_{\alpha}}
\def \Da {D_{\alpha}}
\def \ia {i_{\alpha}}
\def \iaak {i_{\alpha,\ak}}
\def \na {n_{\alpha}}

\def \ga {G_{\alpha}}
\def \sa {s_{\alpha}}

\def \ya {y_{\alpha}}

\def \za {z_{\alpha}}

\def \sb {s_{\beta}}
\def \lb {l_{\beta}}
\def \eb {e_{\beta}}
\def \db {d_{\beta}}
\def \nb {n_{\beta}}

\def \wb {w_{\beta}}
\def \Db {D_{\beta}}

\def \taub {\tau_{\beta}}
\def \rob {\rho_{\beta}}
\def \zb {z_{\beta}}

\def \gb {G_{\beta}}
\newcommand{\accs}[2]{\left\langle #1\, ,\,#2 \right\rangle_{\Sigma}} 
\newcommand{\accsv}[2]{\left\langle #1\, ,\,#2 \right\rangle_{\Sigma,v}} 
\def \unit {\bU} 
\def \ca {\C^{\times}}
\def \xtg {X(T)^G}
\def \xtgd {\left(X(T)^G\right)^{\vee}}
\def \xtgu {\xtg_{\unit}}
\def \xtgc {\xtg_{\C}}
\def \xtgca {\xtg_{\ca}}
\def \extg {e\,X(T)^G}

\def \extgca {\left(\extg\right)_{\C^{\ast}}}

\def \sg {\Sigma(1)/G}
\def \csg {\C^{\,\Sigma(1)/G}}
\def \Nsg {\N^{\,\Sigma(1)/G}}
\def \casg {\left(\ca\right)^{\,\Sigma(1)/G}}
\def \pls {\PL(\Sigma)}
\def \plsg {\PL(\Sigma)^G}
\def \plsgc {\PL(\Sigma)^G_{\C}}
\def \plsgca {\plsg_{\ca}}

\def \eplsg {e\,\PL(\Sigma)^G}
\def \eplsgc {\left(e\,\PL(\Sigma)^G\right)_{\C}}
\def \eplsgca {\left(e\,\PL(\Sigma)^G\right)_{\ca}}
\def \psg {\Ps^G}
\def \psgc {\left(\Ps^G\right)_{\C}}
\def \psgca {\left(\Ps^G\right)_{\ca}}

\def \epsgca {\left(e\,\Ps^G\right)_{\ca}}
\makeindex
\newindex{not}{ndx}{nnd}{Index des notations}
\newindex{def}{ddx}{dnd}{Index des d\'efinitions}
\def\nindex{\index[not]}
\def\dindex{\index[def]}
\title{Fonction z\^eta des hauteurs des vari\'et\'es toriques non d\'eploy\'ees 
\\
Height zeta functions of nonsplit toric varieties 
}
\author{David BOURQUI}
\date{}
\begin{document}
\maketitle
{\small{\textbf{Abstract : } 
We investigate the anticanonical height zeta function of a (non necessarily split)
toric variety defined over a global field of positive characteristic,
drawing our inspiration from the method 
used by Batyrev and Tschinkel to deal with the analogous problem over a number field.
By the way, we give a detailed account of their method.

\textbf{R\'esum\'e : } 
Nous \'etudions la fonction z\^eta des hauteurs anticanonique d'une vari\'et\'e
torique (non n\'ecessairement d\'eploy\'ee) d\'efinie sur un corps global de caract\'eristique positive. 
Nous nous inspirons pour cela de la m\'ethode utilis\'ee par Batyrev et Tschinkel pour traiter la situation
analogue en caract\'eristique z\'ero, m\'ethode que nous rappelons d'ailleurs en d\'etail.
}
}


\tableofcontents

\section{Introduction}

\subsection{Position et origine du probl\`eme}

Soit $V$ une vari\'et\'e projective d\'efinie sur un corps global $K$, 
i.e. un corps de nombres ou le corps de fonctions d'une courbe projective, 
lisse et g\'eom\'etriquement int\`egre, d\'efinie sur un corps fini. 
Soit $H$ une hauteur exponentielle relative \`a un fibr\'e en droites ample. 
Alors pour tout r\'eel $B$ le nombre 
\begin{equation}
n_{V,H}(B)=\# \left\{x\in V(K),\,H(x)\leq B\right\}
\end{equation}
est fini. 
Si l'ensemble $V(K)$ est dense pour la topologie de Zariski, 
la quantit\'e $n_{V,H}(B)$ tend donc vers l'infini quand $B$ tend vers l'infini.  
Une question naturelle est alors d'essayer de d\'ecrire le comportement asymptotique de la quantit\'e $n_{V,H}(B)$, 
en d'autres termes le comportement asymptotique du nombre de points de hauteur born\'ee. 
On cherche notamment \`a interpr\'eter cette description en termes de la g\'eom\'etrie de la vari\'et\'e $V$. 
C'est l'objet d'un programme initi\'e par Manin et ses collaborateurs, 
qui s'est r\'ev\'el\'e extr\^emement riche et ouvert : 
pour la v\'erification des pr\'edictions de Manin pour des classes particuli\`eres de vari\'et\'es, 
des techniques tr\`es diverses ont pu \^etre employ\'ees. 
Ces pr\'edictions (raffin\'ees par Peyre puis Batyrev et Tschinkel) sont maintenant
\'etablies pour plusieurs classes de vari\'et\'es.
Nous renvoyons le lecteur aux textes \cite{Pey:bki} et \cite{Pey:bordeaux}
pour un \'etat g\'en\'eral de la question aux alentours de 2003 et les r\'ef\'erences 
de nombreux travaux sur le sujet.  On pourra \'egalement consulter \cite{Bro:manin_dim_2}
pour un \'etat des lieux r\'ecent concernant le cas des surfaces.

Soulignons que la tr\`es grande majorit\'e de ces travaux 
se placent dans le cas o\`u $K$ est un corps de nombres.
Ici nous nous int\'eressons au cas o\`u $K$ est de caract\'eristique non nulle, 
cas encore peu explor\'e dans la litt\'erature.
Avant toute chose, nous allons pr\'eciser l'une des pr\'edictions de Manin
concernant le comportement asymptotique de $n_{V,H}(B)$,
dans le cas o\`u le corps de base est un corps de nombres.
Elle peut s'\'enoncer de la mani\`ere suivante.

\begin{question}
\label{ques:manin:peyre:arit}
Soit $V$ une vari\'et\'e projective et lisse d\'efinie sur un corps de nombres $K$. 
On suppose que la classe du faisceau anticanonique est \`a l'int\'erieur du c\^one effectif, 
et que l'ensemble $V(K)$ des points rationnels de $V$ est dense pour la topologie de Zariski. 
Soit $t$ le rang du groupe de N\'eron-S\'everi de $V$. 
Soit $H$ une hauteur relative au faisceau anticanonique.
Existe-t-il un ouvert de Zariski non vide $U$ de $V$
et une constante $C>0$ 
tels qu'on ait 
\begin{equation}
\label{eq:formule:manin:arit}
n_{U,H}(B)\underset{B\to +\infty}{\sim}C\,B\,\log(B)^{t-1}\quad ?
\end{equation}
\end{question} 
La restriction \`a un ouvert $U$ \'eventuellement strict de $V$ est n\'ecessaire en raison de l'existence
possible de  ferm\'es acccumulateurs, 
dont un prototype est donn\'e par les diviseurs exceptionnels sur les surfaces de del Pezzo.

Soulignons que bien qu'il ait \'et\'e d\'emontr\'e que la question \ref{ques:manin:peyre:arit}
avait une r\'eponse positive pour de large classes de vari\'et\'es, un contre-exemple d\^u \`a Batyrev
et Tschinkel montre que la r\'eponse \`a cette question est n\'egative en g\'en\'eral (le contre-exemple
porte sur la puissance du logarithme apparaissant dans la formule \eqref{eq:formule:manin:arit}, 
cf. \cite{BaTs:conicbundles}).

Il existe une version fonctionnelle imm\'ediate de la question \ref{ques:manin:peyre:arit} :
il suffit de  remplacer dans l'\'enonc\'e l'hypoth\`ese <<$K$ est un corps de nombres>> par <<$K$ est un corps
global de caract\'eristique positive>>.
Cependant, la nature <<dispers\'ee>> de l'ensemble des valeurs prises par les fonctions hauteurs dans le cas 
fonctionnel entra\^\i ne
qu'une formule du type \eqref{eq:formule:manin:arit}
ne pourra jamais \^etre v\'erifi\'ee.
Plus pr\'ecis\'ement, cet ensemble de valeurs sera typiquement inclus dans $q^{\Z}$ o\`u $q$ est 
le cardinal du corps des constantes. On a donc dans ce cas
\begin{equation}
\forall n\in \N,\quad n_{V,H}\left(q^{n+\frac{1}{2}}\right)=n_{V,H}\left(q^{n}\right)
\end{equation}
et une formule du type \eqref{eq:formule:manin:arit} entra\^\i nerait alors aussit\^ot 
la contradiction $\sqrt{q}=1$.

Pour obtenir une version fonctionnelle satisfaisante de la question \ref{ques:manin:peyre:arit},
on remarque que le comportement asymptotique de $n_{U,H}(B)$ est \'etroitement li\'e, 
par des th\'eor\`emes taub\'eriens,
au comportement analytique de la s\'erie g\'en\'eratrice
\begin{equation}
\zeta_{U,H}(s)
=
\sum_{x\in U(K)}H(x)^{\,-s}
\end{equation}
($s$ d\'esignant une variable complexe),
que l'on baptise fonction z\^eta des hauteurs.
Un des moyens couramment utilis\'es pour obtenir une formule du type \eqref{eq:formule:manin:arit} est 
d'ailleurs d'\'etudier d'abord le comportement analytique de cette fonction, 
puis d'appliquer un th\'eor\`eme taub\'erien ad\'equat, tel que le r\'esultat suivant.

\begin{thm}
\label{thm:taub}
S'il existe un ouvert $U$ non vide tel que 
$
\zeta_{U,H}(s)
$
converge absolument pour $\Re(s)>1$ et
un nombre r\'eel $\varepsilon>0$ 
tels que la fonction
\begin{equation}
s\longmapsto (s-1)^t\,\zeta_{U,H}(s)
\end{equation}
se prolonge en une fonction $g$ holomorphe sur l'ouvert $\{\Re(s)>1-\eps\}$, 
et v\'erifiant $g(1)\neq 0$
alors la formule \eqref{eq:formule:manin:arit} est v\'erifi\'ee pour cet ouvert $U$
avec $C=\frac{g(1)}{(t-1)!}$.
\end{thm}

La question qui suit peut alors \^etre vue comme une version fonctionnelle de la question \ref{ques:manin:peyre:arit}.

\begin{question}
\label{ques:manin:peyre:fonc}
Soit $V$ une vari\'et\'e projective et lisse d\'efinie sur un corps global $K$
de caract\'eristique positive. 
On suppose que la classe du faisceau anticanonique est \`a l'int\'erieur du c\^one effectif, 
et que l'ensemble $V(K)$ des points rationnels de $V$ est dense pour la topologie de Zariski. 
Soit $t$ le rang du groupe de N\'eron-S\'everi de $V$. 
Soit $H$ une hauteur relative au faisceau anticanonique. 

Existe-t-il un ouvert de Zariski non vide $U$ de $V$ tel que la s\'erie
\begin{equation}
\zeta_{U,H}(s)
=
\sum_{x\in U(K)}H(x)^{\,-s}
\end{equation}
converge absolument pour $\Re(s)>1$ et, 
pour un certain $\varepsilon>0$, 
se prolonge en une fonction m\'eromorphe sur l'ouvert $\{\Re(s)>1-\varepsilon\}$, 
qui a un p\^ole d'ordre $t$ en $s=1$ ? 
\end{question} 

Naturellement, 
et en accord avec les remarques d\'ej\`a faites, 
m\^eme si cette question admet une r\'eponse positive, 
on ne pourra pas appliquer le th\'eor\`eme \ref{thm:taub}.
En caract\'eristique non nulle, 
la fonction z\^eta des hauteurs a d'autres p\^oles que $1$
sur la droite $\Re(s)=1$, ne serait-ce que ceux provenant
de la p\'eriodicit\'e de $H$.

Dans le cas des corps de nombres, 
Peyre a \'et\'e le premier dans \cite{Pey:duke} \`a proposer (moyennant quelques
hypoth\`eses suppl\'ementaires sur la vari\'et\'e $V$) une expression conjecurale
de la constante $C$ apparaissant dans la formule \eqref{eq:formule:manin:arit}.
Cette expression conjecturale d\'epend d'invariants g\'eom\'etriques et arithm\'etiques
de la vari\'et\'e $V$, ainsi que du choix de la hauteur. Elle a ensuite \'et\'e raffin\'ee
par Batyrev et Tschinkel, et adapt\'ee au cas fonctionnel par 
Peyre dans \cite{Pey:prepu:drap}. Nous rappelons la d\'efinition de la constante 
de Peyre raffin\'ee \`a la section \ref{subsec:def:height}. Nous la noterons $C^{\ast}_{V,H}$.

On a ainsi des versions raffin\'ees des questions \ref{ques:manin:peyre:arit}
et \ref{ques:manin:peyre:fonc}.
\begin{question}
\label{ques:manin:peyre:arit:raf}
Soit $V$ une vari\'et\'e projective et lisse d\'efinie sur un corps de nombres $K$. 
On suppose que la classe du faisceau anticanonique est \`a l'int\'erieur du c\^one effectif, 
et que l'ensemble $V(K)$ des points rationnels de $V$ est dense pour la topologie de Zariski. 
Soit $t$ le rang du groupe de N\'eron-S\'everi de $V$. 
Soit $H$ une hauteur relative au faisceau anticanonique.
On suppose en outre que $V$ v\'erifie les hypoth\`eses n\'ecessaires pour que 
la constante de Peyre raffin\'ee $C^{\ast}_{V,H}$ soit d\'efinie.

Existe-t-il un ouvert de Zariski non vide $U$ de $V$
tel qu'on ait 
\begin{equation}
\label{eq:formule:manin:arit:raf}
n_{U,H}(B)\underset{B\to +\infty}{\sim}C^{\ast}_{V,H}\,B\,\log(B)^{t-1}\quad ?
\end{equation}
\end{question} 
\begin{rem}
\`A la connaissance de l'auteur, dans tous les cas o\`u on sait montrer que la r\'eponse \`a la question 
\ref{ques:manin:peyre:arit} est positive, on sait \'egalement montrer
que la r\'eponse \`a la question \ref{ques:manin:peyre:arit:raf} est positive.
\end{rem}

\begin{question}
\label{ques:manin:peyre:fonc:raf}
Soit $V$ une vari\'et\'e projective et lisse d\'efinie sur un corps global $K$
de caract\'eristique positive. 
On suppose que la classe du faisceau anticanonique est \`a l'int\'erieur du c\^one effectif, 
et que l'ensemble $V(K)$ des points rationnels de $V$ est dense pour la topologie de Zariski. 
Soit $t$ le rang du groupe de N\'eron-S\'everi de $V$. 
Soit $H$ une hauteur relative au faisceau anticanonique, 
On suppose en outre que $V$ v\'erifie les hypoth\`eses n\'ecessaires pour que 
la constante de Peyre raffin\'ee $C^{\ast}_{V,H}$ soit d\'efinie.

Existe-t-il un ouvert de Zariski non vide $U$ de $V$ tel que la s\'erie
\begin{equation}
\zeta_{U,H}(s)
=
\sum_{x\in U(K)}H(x)^{\,-s}
\end{equation}
converge absolument pour $\Re(s)>1$ et, 
pour un certain $\varepsilon>0$, 
se prolonge en une fonction m\'eromorphe sur l'ouvert $\{\Re(s)>1-\varepsilon\}$, 
qui a un p\^ole d'ordre $t$ en $s=1$, 
et v\'erifiant
\begin{equation}
\lim_{s\to 1}(s-1)^{t}\zeta_{U,H}(s)=(t-1)!\,C^{\ast}_{V,H}\quad ?
\end{equation}
\end{question} 

Concernant la question \ref{ques:manin:peyre:fonc:raf}, 
le cas des espaces projectifs est trait\'e par Wan dans \cite{Wan}, 
montrant ainsi une formule figurant d\'ej\`a dans \cite{Ser:LMWT}. 
Le cas des vari\'et\'es de drapeaux, 
qui englobe le pr\'ec\'edent, 
a \'et\'e trait\'e  ind\'ependamment par Peyre dans \cite{Pey:prepu:drap},
et Lai et Yeung dans \cite{LaiYeu} 
(sans interpr\'etation de la constante dans ce dernier cas, c'est-\`a-dire 
que seule la question \ref{ques:manin:peyre:fonc} est consid\'er\'ee).

Dans ce texte, 
on \'etudie la question \ref{ques:manin:peyre:fonc:raf} 
pour une vari\'et\'e torique projective et lisse d\'efinie sur un corps 
global de caract\'eristique positive, 
non n\'ecessairement d\'eploy\'ee.

Une de motivations de ce travail est que le probl\`eme analogue 
sur les corps de nombres a d\'ej\`a \'et\'e trait\'e avec 
succ\`es\footnote{C'\'etait \'egalement le cas pour les vari\'et\'es de drapeaux.}, 
qui plus est de deux mani\`ere diff\'erentes : 
Batyrev et Tschinkel ont d\'emontr\'e dans \cite{BaTs:anis} et \cite{BaTs:manconj} 
que la r\'eponse \`a la question \ref{ques:manin:peyre:arit:raf}  
\'etait positive pour les vari\'et\'es toriques, 
en exploitant la structure de groupe du tore pour utiliser des techniques d'analyse harmonique.
Par la suite Salberger a red\'emontr\'e dans \cite{Sal:tammes} le r\'esultat dans un cadre plus restreint 
(vari\'et\'es toriques d\'eploy\'ees, d\'efinies sur $\Q$, de faisceau anticanonique globalement engendr\'e)
mais par une m\'ethode compl\`etement diff\'erente bas\'ee sur l'usage de la description explicite des
torseurs universels au-dessus des vari\'et\'es toriques.

Dans \cite{Bou:jnt} et \cite{Bou:crelle}, 
nous avons montr\'e comment, 
en s'inspirant de la m\'ethode de Salberger, 
on pouvait montrer que la r\'eponse \`a la question \ref{ques:manin:peyre:fonc:raf}
\'etait positive pour les vari\'et\'es toriques d\'eploy\'ees d\'efinies sur un corps de fonctions quelconque 
(sans hypoth\`ese sur le faisceau anticanonique).

Dans ce texte, nous adaptons au cas fonctionnel 
l'approche utilis\'ee par Batyrev et Tschinkel 
dans \cite{BaTs:anis} et  \cite{BaTs:manconj},
pour \'etendre le r\'esultat aux vari\'et\'es toriques non n\'ecessairement d\'eploy\'ees.
La sous-section suivante d\'etaille cette adaptation.

Ce texte contient \'egalement une pr\'esentation d\'etaill\'ee de la d\'emonstration du r\'esultat de Batyrev et Tschinkel,
les deux d\'emonstrations \'etant pr\'esent\'ees en parall\`ele.
La raison de ce choix est au moins double : tout d'abord, il permet de bien mettre en \'evidence les analogies et les 
diff\'erences qui existent dans le traitement du calcul de la fonction z\^eta des hauteurs des vari\'et\'es toriques
entre le cas des corps de nombres et le cas des corps de fonctions.
Ensuite, pour autant qu'il nous soit permis d'en juger, ce choix peut s'av\'erer utile 
\`a ceux qui d\'esirent comprendre en d\'etail la d\'emarche de Batyrev et Tschinkel, les articles
\cite{BaTs:anis} et  \cite{BaTs:manconj} pouvant s'av\'erer d'un abord un peu ardu
et elliptique pour le lecteur non averti.
\\~\\

\noindent\textbf{Remerciements}

Je remercie chaleureusement Antoine Chambert-Loir, Jean-Louis Colliot-Th\'el\`ene 
et Emmanuel Peyre pour leurs remarques, corrections et suggestions concernant ce
texte. J'ai envers le rapporteur une reconnaissance tout particuli\`ere
pour sa lecture minutieuse et ses innombrables corrections et suggestions.

\subsection{L'adaptation de la m\'ethode de Batyrev et Tschinkel en caract\'eristique positive}
\label{subsec:adaptation}

Dans cette section, 
nous r\'esumons bri\`evement la m\'ethode utilis\'ee dans \cite{BaTs:manconj} et \cite{BaTs:anis}, 
en expliquant quelles parties de la d\'emonstration n\'ecessitent une modification en caract\'eristique non nulle.

La premi\`ere \'etape consiste \`a d\'efinir explicitement un syst\`eme de hauteurs 
puis \`a l'\'etendre \`a l'espace ad\'elique associ\'e au tore.
La construction est strictement la m\^eme dans le cas fonctionnel. 
Elle est rappel\'ee dans les sections \ref{subsec:fct:zeta:hauteur} (nous corrigeons
au passage une erreur de Batyrev et Tschinkel dans la d\'efinition des hauteurs locales
pour les places ramifi\'ees) et \ref{subsec:formule:poisson}.

\`A ce stade, il faut d\'ej\`a noter que la topologie de 
l'espace ad\'elique associ\'e au tore a des propri\'et\'es diff\'erentes 
dans chacune des deux situation. 
Moralement, en fait, 
la situation est plus agr\'eable en caract\'eristique positive~: 
beaucoup des groupes topologiques mis en jeu sont compacts 
(notamment, le point \ref{item:3:prop:compacite} de la proposition 
\ref{prop:compacite} n'est valable qu'en caract\'eristique positive).

Disposant des fonctions hauteurs sur l'espace
ad\'elique associ\'e au tore, lequel est un groupe ab\'elien localement compact,
l'id\'ee cruciale de Batyrev et Tschinkel est d'appliquer la formule de Poisson afin d'obtenir 
une repr\'esentation int\'egrale de la fonction z\^eta des hauteurs.
Pour ce faire, 
il faut \'etablir l'int\'egrabilit\'e de la transform\'ee de Fourier de la hauteur,
laquelle se d\'ecompose en produit de transform\'ees de Fourier locales. 
On utilise dans le cas des corps de nombres une expression explicite pour 
presque toutes les transform\'ees de Fourier locales (cf. le th\'eor\`eme \ref{thm:transfo:nonram}),
et des majorations ad\'equates pour les transform\'ees de Fourier restantes.
La formule explicite d\'ecrivant presque toutes les transform\'ees de Fourier locales est 
la m\^eme dans le cas fonctionnel.
En ce qui concerne les transform\'ees de Fourier locales aux places restantes, 
leur continuit\'e suffit pour assurer la convergence dans le cas fonctionnel, 
cependant nous avons besoin de quelques renseignements 
sur la forme des fonctions obtenues (cf. la sous-section \ref{subsubsec:transfolocales:casfonc}).

Le choix d'un scindage du groupe des caract\`eres du tore permet alors de montrer
que la fonction z\^eta des hauteurs s'obtient par int\'egration 
(sur un espace vectoriel r\'eel dans le cas arithm\'etique, sur un produit de cercles dans le cas fonctionnel) 
d'une fonction qui poss\`ede une expression en terme de produit de fonctions $L$ de Hecke ; 
cf. le corollaire \ref{cor:rep:int:arit} pour le cas des corps de nombres
(cf. \'egalement \cite[Theorem 3.1.3]{BaTs:anis} et 
\cite[page 46]{BaTs:manconj}), 
et le corollaire \ref{cor:rep:int:fonc} pour le cas fonctionnel.
Il faut <<ma\^\i triser>> le comportement analytique de le fonction sous l'int\'egrale, 
et dans le cas des corps de nombres, 
on a besoin pour cela d'un contr\^ole uniforme sur les bandes verticales des fonctions $L$, 
obtenu par Rademacher via le principe de Phragmen-Lindel\"of 
(proposition \ref{prop:maj:unif:L} et \cite[Theorem 3.2.3]{BaTs:anis}). 
Dans le cas fonctionnel, 
l'holomorphie de la fonction $L\left(.,\chi\right)$ 
quand le caract\`ere $\chi$ est non trivial est suffisante.

Pour d\'eterminer les propri\'et\'es analytiques la fonction z\^eta des hauteurs, 
il s'agit maintenant 
de comprendre comment l'int\'egration modifie le comportement 
analytique de la fonction sous l'int\'egrale 
(cet \'etape n'appara\^\i t d'ailleurs pas dans le cas des tores anisotropes). 
C'est l'objet de la proposition technique de Batyrev et Tschinkel (\cite[Theorem 6.19]{BaTs:manconj}). 
La d\'emonstration proc\`ede par des applications successives du th\'eor\`eme des r\'esidus. 
Dans ce texte, nous utilisons une version raffin\'e du r\'esultat d\^ue \`a Chambert-Loir et Tschinkel
(th\'eor\`eme \ref{thm:CLTs}).

La transposition directe du lemme technique et de sa d\'emonstration en caract\'eristique non nulle 
s'av\`ere difficile \`a mettre en \oe uvre, 
car bien que la compacit\'e des espaces topologiques mis en jeu simplifie un peu les choses, 
rendant inutiles des hypoth\`eses du type contr\^ole uniforme sur les bandes verticales
(indispensables en caract\'eristique z\'ero), 
les fonctions z\^etas des hauteurs en caract\'eristique non nulle 
s'av\`erent poss\'eder plus de p\^oles 
sur la droite $\Re(s)=1$ que ceux provenant de la p\'eriodicit\'e de la hauteur.
Ce ph\'enom\`ene est m\^eme d\'ej\`a visible dans le cas des vari\'et\'es toriques d\'eploy\'ees. 
Prenons en effet l'exemple du plan projectif \'eclat\'e en un point : 
en notant $q$ le cardinal du corps des constantes, 
la formule de la page 355 de \cite{Bou:jnt} montre que la fonction z\^eta des hauteurs anticanonique 
s'\'ecrit dans ce cas
\begin{equation}
\zeta_H(s)=f_1(q^{-s})\,\zetac(3\,s-2)\,\zetac(2\,s-1)+f_2(q^{-s})\,\zetac(2\,s-1)+f_3(q^{-s}),
\end{equation}
o\`u, 
pour $i=1,2,3$, $s\mapsto f_i(q^{-s})$ est holomorphe sur le domaine $\Re(s)>\frac{1}{2}$ 
et $s\mapsto f_1(q^{-s})$ ne s'annule pas sur ce domaine. 
Ainsi l'ensemble des p\^oles situ\'es sur la droite $\Re(s)=1$ contient $\{1+\frac{2\,i\,k\,\pi}{3\,\log q}\}_{k\in\Z}$,
bien qu'on puisse v\'erifier que la fonction z\^eta des hauteurs n'est pas 
$\frac{2\,i\,\pi}{3\,\log q}$-p\'eriodique.
Rappelons que dans le cas de la caract\'eristique z\'ero, 
Batyrev et Tschinkel montrent 
que le seul p\^ole de la fonction z\^eta des hauteurs 
des vari\'et\'es toriques sur la droite $\Re(s)=1$ est $s=1$. 
La gestion des p\^oles suppl\'ementaires en caract\'eristique non nulle se r\'ev\`ele vite \^etre tr\`es d\'elicate 
(voire ing\'erable\dots) 
si on veut suivre <<au plus pr\`es>> la m\'ethode de Batyrev et Tschinkel. 

C'est pourquoi, pour aboutir \`a une version fonctionnelle du r\'esultat technique
utilis\'e dans le cas des corps de nombres, 
nous exploitons la p\'eriodicit\'e des fonctions mises en jeu pour 
les exprimer en terme de s\'eries de type combinatoire.
Les techniques utilis\'ees pour \'evaluer le comportement de ces s\'eries 
par int\'egration sont alors 
similaires \`a celles employ\'ees dans \cite{Bou:crelle}. 
Nous aboutissons ainsi aux lemmes \ref{lm:tech1},
\ref{lm:tech1:bis} et \ref{lm:tech1:bis:tilde}
qui sont le pendant en caract\'eristique positive 
du lemme technique de Batyrev et Tschinkel.

Le comportement analytique de la fonction obtenue apr\`es int\'egration
est essentiellement d\'ecrit par une suite exacte de $\Z$-modules libres de rang fini
(cf. les \'enonc\'es du th\'eor\`eme \ref{thm:CLTs} et du  lemme \ref {lm:tech0}).
La suite exacte mise en jeu n'est pas exactement la m\^eme dans le cas des corps de nombres
ou dans le cas fonctionnel.
Dans les deux cas, elle provient de la construction suivante~: 
la r\'esolution flasque du groupe des caract\`eres 
du tore par le groupe de Picard de la vari\'et\'e torique (i.e. la suite exacte \eqref{eq:resflT}). 
induit par dualit\'e une suite exacte de tores alg\'ebriques~;
on consid\`ere alors l'image de cette suite exacte par le morphisme degr\'e (d\'efinie \`a la section \ref{subsec:ledegre}).
La diff\'erence essentielle vient alors du fait que le morphisme degr\'e est surjectif dans le cas des corps de nombres mais pas dans le cas fonctionnel (cf. le lemme \ref{lm:degre:quasidep}) o\`u il est seulement de conoyau fini.

La derni\`ere \'etape de la d\'emonstration consiste \`a calculer explicitement le terme principal 
de la fonction z\^eta des hauteurs au point critique $s=1$, et \`a v\'erifier s'il 
co\"\i ncide avec la pr\'ediction de Peyre. 
On a besoin d'un th\'eor\`eme d'Ono sur le nombre de Tamagawa d'un tore alg\'ebrique (th\'eor\`eme \ref{thm:ono}), 
lequel th\'eor\`eme a \'et\'e d\'emontr\'e dans \cite{Ono:tamnum} pour tout corps global, 
mais a cependant d\^u \^etre corrig\'e par Oesterl\'e dans le cas de la caract\'eristique non nulle, 
en introduisant un facteur correctif dans la d\'efinition du nombre de Tamagawa 
(ce facteur correctif provient de la non-surjectivit\'e du degr\'e en caract\'eristique non nulle). 
On a besoin \'egalement de r\'esultats de Colliot-Th\'el\`ene et Sansuc, 
d\'emontr\'es pour tout corps global \'egalement, 
et permettant d'obtenir le lemme \ref{lm:appfaible}. 
On peut donc ici reprendre la ligne de calcul de Batyrev et Tschinkel, 
ce qui est fait dans la partie \ref{app_fonction_zeta}. 
Notons que dans le cas fonctionnel, le d\'efaut de surjectivit\'e du degr\'e fait
intervenir au cours du calcul des termes non triviaux (correspondant \`a des
cardinaux de groupes finis) qui n'apparaissent pas dans le cas des corps de nombres.
Dans une version pr\'ec\'edente de ce texte, l'auteur affirmait que l'un de ces termes
non triviaux subsistait dans l'expression finale du terme principal de la fonction z\^eta
des hauteurs (il s'agit de l'invariant $\KT$
d\'efini \`a la sous-section \ref{subsubsec:KT}), <<montrant>> ainsi que la constante pr\'edite
par Peyre et Batyrev-Tschinkel n'\'etait pas la bonne dans ce cas.
Les calculs avaient \'et\'e men\'es en supposant \`a tort que le groupe 
intervenant dans le point \ref{item:3:lm:kergammast} du lemme \ref{lm:kergammast}
\'etait trivial. Une fois ce point corrig\'e, on s'aper\c coit que tous les termes
suppl\'ementaires intervenant en caract\'eristique non nulle dans le calcul 
du terme principal de la fonction z\^eta des hauteurs se simplifient, et que 
l'expression obtenue est bien celle attendue.

\section{Rappels et notations}

\subsection{Quelques notations}
\label{subsec:not}

Nous fixons ici quelques notations
utilis\'ees dans l'ensemble du texte.

On note $\card{E}$ le cardinal d'un ensemble fini $E$.

Pour tout r\'eel $\alpha$, on note 
$\R_{>\alpha}$
l'ensemble
$
\{x\in \R, \,x>\alpha\}.
$
On d\'efinit de m\^eme de mani\`ere \'evidente les ensembles
$\R_{\geq\alpha}$, $\R_{<\alpha}$ et $\R_{\leq\alpha}$.

On note, pour tout corps $K$, 
$\overline{K}$ une cl\^oture alg\'ebrique de $K$ et 
$K^{\sep}$ la cl\^oture s\'eparable de $K$ dans $\overline{K}$.

Soit $N$ un groupe ab\'elien.
On note $\nindex{$N^{\vee}$}N^{\vee}$ le dual alg\'ebrique de $N$, 
c'est-\`a-dire le groupe ab\'elien $\Hom(N,\Z)$.
L'accouplement naturel entre $N$ et $N^{\vee}$ sera not\'e $\acc{.}{.}$. 

Soit $M$ un autre groupe ab\'elien et $f\,:\,M\to N$ un morphisme de groupes
Pour tout anneau commutatif unitaire $A$, on note $N_{A}$ le $A$-module $N\otimes_{\Z} A$
et $f_A$ le morphisme de $A$-module $f\otimes \Id\,:\,M_A\to N_A$ induit par $f$.

Pour tout anneau commutatif unitaire $A$, on d\'esigne par $A^{\times}$ le groupe des \'el\'ements inversibles de $A$.

Si $N$ est un $\Z$-module libre de rang fini, et $U$ est une partie de $N_{\R}$, on note 
\begin{equation}
\nindex{$\tube{U}$}
\tube{U}
\eqdef
\{\bs\in N_{\C},\,\Re(\bs)\in U\}
\end{equation}
le domaine tubulaire au-dessus de $U$.

\subsection{Rappels sur les corps globaux}\label{subsec:corps:globaux}

Dans tout ce texte, 
on appellera \dindex{corps de fonctions}\termin{corps de fonctions} un corps
global de caract\'eristique positive, i.e. une extension $K$ de
type fini et de degr\'e de transcendance $1$ d'un corps fini.
Le \dindex{corps des constantes}\termin{corps des constantes} de $K$ est la cl\^oture alg\'ebrique
du sous-corps premier de $K$ dans $K$. Il sera not\'e  \nindex{$\F_K$}$\F_K$.
Le choix d'une base de transcendance s\'eparable de $K$ sur $\F_K$
permet alors d'identifier $K$ \`a une extension finie s\'eparable 
du corps des fractions rationnelles 
en une ind\'etermin\'ee $\F_K(T)$. Il existe en outre une courbe projective, 
lisse et g\'eom\'etriquement int\`egre d\'efinie sur $\F_K$, unique \`a isomorphisme pr\`es,
et not\'ee \nindex{$\courbe_K$}$\courbe_K$, 
telle que $K$ s'identifie au corps de fonctions de $\courbe_K$.

On adopte dans ce texte la convention suivante :
on fixe un corps fini $k$ dont on note $q$ le cardinal. 
On supposera alors, dans tout ce texte, que les corps de fonctions consid\'er\'es ont un corps
des constantes qui contient $k$.
Pour un tel corps de fonctions $K$ on note \nindex{$\qk$}$\qk$ le cardinal du corps des constantes, 
et \nindex{$\dk$}$\dk$ son \termin{degr\'e absolu}, c'est-\`a-dire que $\dk$ v\'erifie $\qk=q^{\,\dk}$.

Soit $K$ un corps global, i.e. un corps de nombres ou un corps de fonctions. 
Dans tout ce texte, on dira \^etre dans le \dindex{cas arithm\'etique}cas arithm\'etique si le corps de base est un
corps de nombres et dans le \dindex{cas fonctionnel}cas fonctionnel si c'est un corps de fonctions.
Dans le cas arithm\'etique, 
on note \nindex{$\disc(K)$}$\disc(K)$ le discriminant absolu de $K$
et dans le cas fonctionnel
on note \nindex{$g_K$}$g_K$ le genre de $K$, c'est-\`a-dire le genre de $\courbe_K$.

On note \nindex{$\placesde{K}$}$\placesde{K}$ l'ensemble des places de $K$, et 
\nindex{$\placesde{K,f}$}$\placesde{K,f}$ (respectivement \nindex{$\placesde{K,\infty}$}$\placesde{K,\infty}$)
l'ensemble des places finies (respectivement archim\'ediennes) de $K$.
On identifiera toujours un \'el\'ement $v$ de $\placesde{K,f}$ \`a l'unique valuation normalis\'ee 
qui le repr\'esente, 
i.e. l'unique \'el\'ement de $v$ dont le groupe de valeurs est $\Z$. 
Dans le cas fonctionnel, on a $\placesde{K}=\placesde{K,f}$ 
et $\placesde{K}$ s'identifie \`a l'ensemble des points ferm\'es de la courbe~$\courbe_K$.

Pour $v\in \placesde{K}$, 
on note \nindex{$K_v$}$K_v$ le compl\'et\'e de $K$ en $v$.
Pour $v\in \placesde{K,f}$, 
on note \'egalement
\nindex{$
\Ov
$}
$
\Ov
$ 
l'anneau de valuation de $v$, 
et 
\nindex{
$
k_v
$}
$
k_v
$ 
le corps r\'esiduel de $v$. 

Pour $v\in \placesde{K,\infty}$, on d\'esignera abusivement
par $\Ov^{\times}$ le sous-groupe de $K_v^{\times}$
constitu\'e des \'el\'ements de valeur absolue $1$.

Dans le cas fonctionnel, 
on note aussi 
\nindex{$f_v$}
\begin{equation}
f_v
=
[k_v:\F_K]
\end{equation}
le degr\'e r\'esiduel <<absolu>>.

Le cardinal de $k_v$ sera not\'e \nindex{$q_v$}$q_v$.  
Dans le cas fonctionnel, il est \'egal \`a $\qk^{f_v}$.

On normalise la mesure de Haar \nindex{$dx_v$}$dx_v$ sur $K_v$ de la mani\`ere suivante :
si $v$ est finie, on choisit $dx_v$ 
de sorte qu'on ait la relation $\int_{\Ov}dx_v=1$,
si $v$ est archim\'edienne et $K_v=\R$, 
$dx_v$ est la mesure de Lebesgue usuelle 
et si $K_v$ est isomorphe \`a $\C$ on prend $dx_v=i\,dz\,d\overline{z}$.

Si $S$ est un sous-ensemble fini de $\placesde{K}$ contenant $\placesde{K,\infty}$,
on note
\nindex{$\OS$}
\begin{equation}
\OS=\bigcap_{v\notin S}\Ov
\end{equation}
l'anneau des $S$-entiers de $K$.

Soit $v\in \placesde{K,f}$.
Pour $x\in K_v$ on note
\nindex{$\abs{\,.\,}_v$}
\begin{equation}
\abs{x}_v=q_v^{\,-v(x)}.
\end{equation}
Si $v\in \placesde{K}\setminus \placesde{K,f}$,
$\abs{\,.\,}_v$ d\'esigne la valeur absolue usuelle sur $\R$ si
$K_v=\R$,
et le carr\'e de la valeur absolue usuelle 
sur $\C$ si $K_v$ est isomorphe \`a $\C$.

Avec ce choix de valeur absolues, 
on a, pour tout $v\in \placesde{K}$ et pour tout $a$ de $K_v$,
$d(ax)_v=\abs{a}_v\,dx_v$.
Pour tout $x\in K^{\times}$ on a $\abs{x}_{v}=1$ 
pour tout $v\in \placesde{K}$ sauf un nombre fini 
et la \termin{formule du produit} 
\begin{equation}
\prod_{v\in \placesde{K}} 
\abs{x}_v\,
=\,
1.
\end{equation}

La fonction z\^eta de Dedekind de $K$, 
not\'ee \nindex{$\zeta_{K}$}$\zeta_{K}$, 
est d\'efinie par la formule 
\begin{equation}
\zeta_{K}(s)
=
\prod_{v\in \placesde{K,f}}\,\left(1-q_v^{\,-s}\right)^{-1}
\end{equation}
o\`u $s$ est une variable complexe.
On dispose  des r\'esultats classiques suivants 
sur le comportement analytique de $\zeta_{K}$.
\begin{prop}\label{prop:fonc:zeta}
\begin{enumerate}
\item
Dans le cas arithm\'etique, 
la s\'erie d\'efinissant $\zeta_K(s)$ converge absolument pour $\Re(s)>1$ 
et se prolonge en une fonction m\'eromorphe sur $\C$ tout entier, 
avec un p\^ole simple en $s=1$. La fonction $(s-1)\,\zeta_K(s)$ se prolonge
en une fonction holomorphe sur $\tube{\R_{>0}}$.
\item
Dans le cas fonctionnel, 
la s\'erie d\'efinissant $\zeta_{K}(s)$ converge absolument pour $\Re(s)>1$ 
et se prolonge en une fonction m\'eromorphe sur $\C$ tout entier, 
avec un p\^ole simple en $s=1$, 
et sans z\'ero pour $\Re (s)\neq \frac{1}{2}$. 
De plus $\zeta_{K}(s)$ est une fraction rationnelle en $\qk^{\,-s}$, plus pr\'ecis\'ement on a
\begin{equation}\label{eq:fct:zeta:fonc}
\zeta_{K}(s)
=
\frac
{P(\qk^{\,-s})}
{(1-\qk^{\,-s})(1-\qk^{\,1-s})}
\end{equation}
o\`u $P$ est un polyn\^ome.
\end{enumerate}
\end{prop}
Dans le cas fonctionnel, on note \nindex{$Z_K$}$Z_K$ la fraction rationnelle v\'erifiant
\begin{equation}
\forall s\in \tube{R_{>1}},\quad Z_K(\qk^{\,-s})=\zeta_{K}(s).
\end{equation}

Nous faisons \`a pr\'esent quelques rappels sur une g\'en\'eralisation naturelle
de la fonction z\^eta de Dedekind : les fonctions $L$ d'Artin.

Soit $K$ un corps global et $M$ un $\Z$-module de rang fini 
muni d'une structure de $\Gal(K^{\sep}/K)$-module discret, 
c'est-\`a-dire une repr\'esentation
\begin{equation}
\rho\,:\,\Gal(K^{\sep}/K)\longto \Aut(M).
\end{equation} 
se factorisant \`a travers un quotient fini de $\Gal(K^{\sep}/K)$.

Soit $G$ un tel quotient, qui est donc le groupe de Galois d'une extension finie 
galoisienne $L$ de $K$.
On a une repr\'esentation 
\begin{equation}
\rho\,:\,G\longto \Aut(M).
\end{equation} 

Soit $v$ une place finie de $K$, 
$G_v$ un groupe de d\'ecomposition au-dessus de $v$ et $I_v$ le groupe d'inertie correspondant. 
Soit 
\begin{equation}
\rho_v\,:\,G_v/I_v\to \Aut(M^{\,I_v})
\end{equation}
la repr\'esentation d\'eduite de $\rho$. 
Soit $\Frob_v\in G_v/I_v$ le frobenius.
On pose, pour tout nombre complexe $s$ tel que $\Re(s)>0$,
\begin{equation}
L_v(s,M,\rho)=\frac{1}{\det(1-\rho_v(\Frob_v)\,q_v^{-s})},
\end{equation}
ce qui est bien d\'efini et ne d\'epend ni du choix de $G_v$,
ni du choix de $G$. 
Si $v$ est une place archim\'edienne, par commodit\'e d'\'ecriture,
on notera $L_v(s,M)$ la fonction constante \'egale \`a $1$.

\begin{ex}
Si $M$ est de rang $1$ et $\rho$ est la repr\'esentation triviale,
\label{ex:Mrang1}
alors pour toute place $v$ finie  on a 
$L_v(s,M,\rho)=(1-q_v^{-s})^{-1}$,
en d'autres termes $L_v(s,M,\rho)$
co\"incide avec le facteur local en $v$ du produit eul\'erien
d\'efinissant la fonction $\zeta_K$.
\end{ex}

\begin{prop}
\begin{enumerate}
\item
Le produit eul\'erien
\begin{equation}
s\mapsto \prod_{v\in \placesde{K}
} L_v(s,M,\rho)
\end{equation}
converge normalement sur tout compact de $\tube{\R_{>1}}$, et d\'efinit
donc une fonction holomorphe sur $\tube{\R_>1}$, not\'ee $L(s,M,\rho)$.
\item
La fonction $L(\,.\,,M,\rho)$ se prolonge en une fonction m\'eromorphe sur $\C$, ayant un p\^ole d'ordre
$\rg\left(M^{\rho\left(\Gal(K^{\sep}/K)\right)}\right)$ en $s=1$.
\end{enumerate}
\end{prop}
\begin{demo}
Le r\'esultat est vrai si $\rho$ est la repr\'esentation triviale, d'apr\`es l'exemple \ref{ex:Mrang1}
et la proposition \ref{prop:fonc:zeta}.
Si $\rho$ est irr\'eductible et $M$ est de rang sup\'erieur \`a $2$,
d'apr\`es \cite[Satz 3]{Artin:uber} le r\'esultat est encore vrai et $L(\,.\,,M,\rho)$
se prolonge en fait en une fonction m\'eromorphe sur tout le plan complexe
qui est holomorphe et inversible au voisinage de $1$.

Dans le cas g\'en\'eral, 
soit $G$ un quotient fini de $\Gal(K^{\sep}/K)$ \`a travers lequel $\rho$ se factorise.
On a une d\'ecomposition
\begin{equation}
M_{\C}=M^{\rho(G)}_{\C}\bigoplus \oplusu{i\in I} M_i
\end{equation}
o\`u, pour $i\in I$, $M_i$ est un sous-espace $\rho(G)$-stable de dimension sup\'erieur \`a 2  et la repr\'esentation $\rho_i\,:\,G\to \Aut(M_i)$ est irr\'eductible. On a alors pour toute place finie $v$
\begin{equation}
L_v(s,M,\rho)=\left(\frac{1}{1-q_v^{-s}}\right)^{\rg(M^{\rho(G)})}
\, \prod_{i\in I} L_v(s,M_i,\rho_i).
\end{equation}
Le r\'esultat en d\'ecoule.
\end{demo}

On pose 
\begin{equation}
\ell(M,\rho)=\lim_{s\to 1} (s-1)^{\rg\left(M^{\rho(\Gal(K^{\sep}/K))}\right)}\,L(s,M,\rho)
\end{equation}
Par la suite on notera
tr\`es souvent \nindex{$L_v(s,M)$}$L_v(s,M)$ (respectivement \nindex{$L(s,M)$}$L(s,M)$, 
respectivement \nindex{$\ell(M)$}$\ell(M)$) 
en lieu et place de $L_v(s,M,\rho)$ (respectivement $L(s,M,\rho)$, respectivement $\ell(M,\rho)$) 
lorsque la repr\'esentation $\rho$ sera clairement indiqu\'ee par le contexte.

On \'enonce \`a pr\'esent deux lemmes \'el\'ementaires sur ces fonctions $L$.
\begin{lemme}\label{lm:comp:L:exseq}
Soit 
\begin{equation}
0
\longto 
M_1
\longto 
M_2
\longto 
M_3
\longto
0
\end{equation}
une suite exacte de $\Gal(K^{\sep}/K)$-modules discrets
qui sont libres de rang fini en tant que $\Z$-modules.
On a alors pour toute place $v$ de $K$
\begin{equation}
\forall s\in \tube{\R_{>0}},\quad L_v(s,M_2)=L_v(s,M_1)\,L_v(s,M_3)
\end{equation}
En particulier, on a
\begin{equation}
\forall s\in \tube{\R_{>1}},\quad L(s,M_2)=L(s,M_1)\,L(s,M_3)
\end{equation}
et 
\begin{equation}
\ell(M_2)=\ell(M_1)\,\ell(M_3).
\end{equation}
\end{lemme}

\begin{lemme}\label{lm:lvzgh}
Soit $L/K$ une extension finie galoisienne de groupe $G$, $H$ un sous-groupe
de $G$, $K'$ le corps global $K^{H}$.
Alors, pour toute place $v$ de $K$ et tout $s\in \tube{\R_{>0}}$, 
on a
\begin{equation}
L_v(s,\Z[G/H])=
\prod_{
\substack{
w\in \placesde{K'}
\\
w|v
}
}
\frac{1}{1-q_w^{-s}}
\end{equation}
En particulier, on a 
\begin{equation}
L(\,.\,,\Z[G/H])=\zeta_{K'}
\end{equation}
et
\begin{equation}
\ell(\Z[G,H])=\Res_{s=1} \zeta_{K'}(s).
\end{equation}
\end{lemme}

\subsection{Hauteurs d'Arakelov, mesure de Tamagawa et constante de Peyre}
\label{subsec:def:height}

Dans cette section, 
nous rappelons bri\`evement la construction des hauteurs d'Arakelov sur une vari\'et\'e 
projective d\'efinie sur un corps global \`a partir de m\'etriques ad\'eliques sur les fibr\'es en droites, 
renvoyant \`a \cite[chapitre 2]{Pey:ecoledete} pour plus de d\'etails, notamment les preuves omises.
Suivant Peyre (\cite{Pey:duke} et \cite{Pey:prepu:drap}), 
nous expliquons ensuite comment une m\'etrique ad\'elique 
sur le faisceau anticanonique d'une vari\'et\'e projective lisse $V$ 
d\'efinie sur un corps global $K$
induit une mesure sur l'espace ad\'elique $V(\ak)$ associ\'e \`a $V$.
La d\'efinition g\'en\'erale de l'espace ad\'elique associ\'e \`a une vari\'et\'e alg\'ebrique est trait\'ee dans 
\cite[I.2]{Wei:AAG} et \cite[I.3]{Oes:invent}. 
Ici, comme $V$ est projective, $V(\ak)$ est l'espace
\begin{equation}
\prod_{v\in \placesde{K}}V(K_v)
\end{equation}
muni de la topologie produit.

Enfin nous rappelons comment la mesure ainsi construite sur $V(\ak)$ permet alors, si $V$ v\'erifie
des hypoth\`eses suppl\'ementaires, de d\'efinir la constante $C^{\ast}_{V,H}$ apparaissant dans les questions 
\ref{ques:manin:peyre:arit:raf}
et \ref{ques:manin:peyre:fonc:raf}. 
L\`a encore, 
nous renvoyons \`a \cite{Pey:duke}  et  \cite{Pey:prepu:drap} pour plus de d\'etails.

Soit $V$ une vari\'et\'e projective, lisse et g\'eom\'etriquement int\`egre d\'efinie sur un corps global $K$, 
et soit $\cL$ un fibr\'e en droites sur $V$. 
Si $v$ est une place de $K$, une \dindex{m\'etrique $v$-adique}\termin{m\'etrique $v$-adique} sur $\cL$ est la donn\'ee, 
pour tout point $K_v$-rationnel 
$x\,:\,\Spec(K_v)\to V$, 
d'une norme $v$-adique $\norm{.}_{v,x}$ sur $\cL(x)=x^{\ast}\cL$
(rappelons que l'on a en particulier $\norm{as}_{v,x}=\abs{a}_v \norm{s}_{v,x}$
pour tout $s$ de $\cL(x)$ et tout $a$ de $K_v$), 
telle que pour tout ouvert de Zariski $U$ de $V$ et toute section $s$ de $\cL$ sur $U$, 
l'application
\begin{equation}
x\mapsto \norm{s(x)}_{v,x}
\end{equation}
est continue sur $U(K_v)$ pour la topologie $v$-adique. 

\begin{rem}
Si $v$ est une place finie de $K$, la donn\'ee d'un mod\`ele projectif $\cV$ de $V$ sur $\Spec(\Ov)$ 
et d'un mod\`ele $\ecL$ de $\cL$ sur $\cV$ d\'efinit de mani\`ere naturelle une telle m\'etrique.
Soit en effet $x\in V(K_v)$. 
Par le crit\`ere valuatif de propret\'e, $x$ d\'efinit un point $\wt{x}\,:\,\Spec(\Ov)\to \cV$.
Le $\Ov$-module $\ecL_{\wt{x}}=\wt{x}^{\ast}\ecL$ est alors libre de rang $1$. Soit $s_0$
un g\'en\'erateur. Le changement de base $\Ov\to K_v$ 
induit un isomorphisme naturel de $\cL(x)$ sur $\ecL_{\wt{x}}\otimes_{\Ov} K_v$,
qui permet de d\'efinir une norme $v$-adique sur $\cL(x)$ par la formule
\begin{equation}
\forall s\in \ecL(x),\quad \norm{s}_v=\abs{\frac{s}{s_0}}_v.
\end{equation}
Comme deux g\'en\'erateurs de $\ecL_{\wt{x}}$ se d\'eduisent l'un de l'autre par
multiplication par un \'el\'ement inversible de $\Ov$, cette d\'efinition est bien
ind\'ependante du choix de $s_0$.
\end{rem}

Une \dindex{m\'etrique ad\'elique}\termin{m\'etrique ad\'elique} sur $\cL$ est alors la donn\'ee d'une famille 
$(\norm{.}_v)_{v\in \placesde{K}}$ 
de m\'etriques $v$-adiques sur $\cL$ v\'erifiant la condition suivante : 
il existe un ensemble fini $S$ de places de $K$ contenant les places archim\'ediennes, 
un mod\`ele projectif $\ecV$ de $V$ sur $\Spec(\cO_S)$, 
et un mod\`ele $\ecL$ de $\cL$ sur $\cV$, 
tels que pour presque tout $v$ de $\placesde{K}\setminus S$ la m\'etrique $\norm{.}_v$ est d\'efinie par le couple 
$(\ecV\times_{\Spec(\cO_S)} \Spec(\Ov), \ecL\times_{\Spec(\cO_S)} \Spec(\Ov))$. 

Donnons un exemple d'une m\'etrique ad\'elique sur $\cL$ lorsque $\cL$ est engendr\'e par ses sections
globales.
On fixe une base $\{s_0,\dots,s_r\}$ du $K$-espace vectoriel $H^0(X,\cL)$.
On pose alors, pour tout $v\in \placesde{K}$, pour tout $x\in V(K_v)$
et toute section locale $s$ de $\cL$ ne s'annulant pas en $x$
\begin{equation}
\norm{s(x)}_{v,x}^{-1}=
\Maxu{i=0,\dots,r} 
\left(
\abs{
\frac 
{s_i(x)}
{s(x)}
}_v
\right).
\end{equation}
Il est \`a noter que cette d\'efinition ne d\'epend pas du choix de la base de $H^0(X,\cL)$.
\begin{lemme}
La famille de m\'etriques $v$-adiques ainsi d\'efinie est une m\'etrique ad\'elique sur $\cL$.
\end{lemme}
\begin{defi}\label{def:met:adel:st}
Cette m\'etrique ad\'elique sera appel\'ee \dindex{m\'etrique ad\'elique standard}\termin{m\'etrique ad\'elique standard}
sur $\cL$.
\end{defi}

Soit \`a pr\'esent $\cL$ et $\cL'$ deux fibr\'es en droites sur $V$, 
et $(\norm{.}_v)_{v\in \placesde{K}}$ (respectivement
$(\norm{.}'_v)_{v\in \placesde{K}}$) 
une m\'etrique ad\'elique sur $\cL$ (respectivement $\cL'$).

Pour toute place $v$ de $K$, 
on d\'efinit ainsi la m\'etrique $v$-adique produit 
$\norm{.}_v\otimes \norm{.}'_v$
sur $\cL\otimes \cL'$ :
pour tout $x\in V(K_v)$, tout $s\in \cL(x)$ et tout $s'\in \cL'(x)$, on pose
\begin{equation}
(\norm{.}_v\otimes \norm{.}'_v) (s\otimes s')=\norm{s}_{v,x}\,\norm{s'}'_{v,x}
\end{equation}
On d\'efinit \'egalement la m\'etrique $v$-adique duale $\norm{.}^{\vee}_v$ sur $\cL^{-1}$ :
pour tout $x\in V(K_v)$ et tout $s^{\vee}\in \cL^{-1}(x)$, on pose
\begin{equation}
\norm{s^{\vee}}^{\vee}_{v,x}=\frac{\abs{\acc{s^{\vee}}{s}}_v}{\norm{s}_{v,x}}
\end{equation}
o\`u $s$ est un \'el\'ement quelconque de $\ecL(x)\setminus\{0\}$.
\begin{lemme}\label{lm:tensomet}
$\left(\norm{.}_v\otimes \norm{.}'_v\right)_{v\in\placesde{K}}$ est une m\'etrique ad\'elique sur $\cL\otimes \cL'$,
et
$\left(\norm{.}^{\vee}_v\right)_{v\in\placesde{K}}$ est une m\'etrique ad\'elique sur $\cL^{-1}$.
\end{lemme}

Soit \`a pr\'esent $L$ une extension finie de $K$.
On suppose donn\'ee une m\'etrique ad\'elique  
$(\norm{.}_\V)_{\V\in \placesde{L}}$ sur $\cL_{L}$.
Pour tout $v\in \placesde{K}$, 
soit $\V$ une place de $L$ divisant $v$.
On identifie $V(K_v)$ \`a un sous-ensemble de $V(L_{\V})$
Soit $x$ un \'el\'ement de $V(K_v)$.
On pose alors, pour tout \'el\'ement $s$ de $\cL(x)$,
\begin{equation}
\norm{s}_{v,x}=\norm{s}_{\V,x}^{\frac{1}{[L_{\V}:K_v]}}.
\end{equation}
\begin{lemme}\label{lm:restr:metrique}
La m\'etrique $(\norm{.}_v)_{v\in \placesde{K}}$ ainsi d\'efinie 
est une m\'etrique ad\'elique sur $\cL$.
\end{lemme}

Une \dindex{hauteur d'Arakelov}\termin{hauteur d'Arakelov} sur $V$ est un couple $(\ecL,(\norm{.}_v))$ o\`u $\cL$ 
est un fibr\'e en droites sur $V$ et $(\norm{.}_v)$ une m\'etrique ad\'elique sur $\cL$. 

Soit $(\ecL,(\norm{.}_v))$ une hauteur d'Arakelov sur $V$.
Soit $x\in V(K)$ et $s$ une section locale de $\cL$ qui ne s'annule pas en $x$.
Presque tous les facteurs du produit infini 
\begin{equation}\label{eq:def:height}
\produ{v\in \placesde{K}}\norm{s(x)}_{v,x}^{-1}
\end{equation}
sont \'egaux \`a $1$ ; le produit infini converge donc et sa valeur ne d\'epend pas du choix de $s$,
d'apr\`es la formule du produit. On la note \nindex{$H_{\cL}$}$H_{\cL}(x)$ : c'est la \termin{hauteur} (exponentielle)\label{page:def:height}
du point $x$ associ\'ee \`a la hauteur d'Arakelov $(\ecL,(\norm{.}_v))$.

Soit $V$ une vari\'et\'e alg\'ebrique projective et lisse de dimension $d$ d\'efinie sur $K$, 
et $\left(\norm{.}_v\right)_{v\in \placesde{K}}$ une m\'etrique ad\'elique sur le faisceau anticanonique. 
On note $H$ la hauteur associ\'ee.

Pour tout $v\in \placesde{K}$, une telle m\'etrique d\'efinit une mesure $\omega_{V,v}$ sur $V(K_v)$ 
(on rappelle qu'on a choisi une normalisation de la mesure de Haar sur $K_v$, cf. la section \ref{subsec:corps:globaux}).

Nous supposons \`a pr\'esent que  
que la vari\'et\'e $V$ v\'erifie les hypoth\`eses suivantes 
(ce sont les hypoth\'eses \'enonc\'ees dans \cite[2.1]{Pey:prepu:drap}) :
\begin{enumerate}
\item
la classe du faisceau anticanonique de $V$ appartient \`a l'int\'erieur du c\^one effectif de $V$;
\item
le groupe $\Pic(V^{\sep})$ (o\`u $V^{\sep}=V\times_{K}K^{\sep}$) 
est un $\Z$-module libre de rang fini, et co\"\i ncide avec $\Pic(\adh{V})$ 
(o\`u $\adh{V}=V\times_{K}\adh{K}$)~;
\item 
les groupes de cohomologie $H^1(V,\ecO_V)$ et $H^2(V,\ecO_V)$ sont nuls~;
\item
si $\ell$ est un nombre premier distinct de la caract\'eristique, la partie $p$-primaire
de $\Br(\adh{V})$ est finie.
\end{enumerate}
De  telles hypoth\`eses sont v\'erifi\'ees en particulier par les vari\'et\'es toriques projectives
et lisses (cf. \cite[Remarque 2.1.1.]{Pey:prepu:drap} et \cite[Exemple 2.1.4]{Pey:circle}).

Pour des vari\'et\'es v\'erifiant ces hypoth\`eses, 
la mesure de Tamagawa sur l'espace ad\'elique $V(\ak)$ est donn\'ee par la formule\footnote{
Un des points d\'elicats de la construction est de montrer la convergence du produit~;
ceci n\'ecessite d'une part la formule de Weil
reliant volume $v$-adique et nombre de points sur le corps r\'esiduel,
d'autre part les conjectures de Weil prouv\'ees par Deligne.}
\nindex{$\omega_V$}
\begin{equation}
\omega_{V}
=
c_{K,\dim(V)}
\,
\ell\left(\Pic(V^{\sep})\right)
\,
\prod_{v\in \placesde{K}}L_v(1,\Pic(V^{\sep}))^{-1}\,\omega_{V,v},
\end{equation}
o\`u l'on a pos\'e, pour tout entier $d\geq 0$, 
\nindex{$c_{K,d}$}
\begin{equation}\label{eq:def:ckv}
c_{K,d}
=
\left\{
\begin{array}{ll}
\disc(K)^{-\frac{d}{2}}  & \text{dans le cas arithm\'etique,} \\
\qde{K}^{\,(1-g_{K})\,d} & \text{dans le cas fonctionnel.}
\end{array}
\right.
\end{equation}

Nous sommes \`a pr\'esent en mesure de d\'efinir la constante $C^{\ast}_{V,H}$ 
apparaissant dans les questions \ref{ques:manin:peyre:arit:raf} et \ref{ques:manin:peyre:fonc:raf}. 

On d\'efinit en fait trois constantes \`a partir des donn\'ees pr\'ec\'edentes.

L'invariant \nindex{$\alpha^{\ast}(V)$}$\alpha^{\ast}(V)$ est d\'efini comme 
\begin{equation}
\label{eq:defalpha}
\alpha^{\ast}(V)
=
\int_{\ceff(V)^{\vee}}
e^{\,-\acc{y}{\omega_{V}^{-1}}}\,dy,
\end{equation}
o\`u $\ceff(V)\subset \Pic(V)_{\R}$ est le c\^one effectif de $V$, 
\begin{equation}
\ceff(V)^{\vee}=\left\{y\in \Pic(V)^{\vee}_{\R},\quad \forall x\in \ceff(V),\quad \acc{y}{x}\geq 0\right\}
\end{equation}
et $dy$
est la mesure de Lebesgue sur $\Pic(V)^{\vee}_{\R}$ normalis\'ee
par le r\'eseau $\Pic(V)^{\vee}$.

L'invariant \nindex{$\beta(V)$}$\beta(V)$ est d\'efini comme 
\begin{equation}
\beta(V)=\card{H^1(K,\Pic(V^{\sep}))}.
\end{equation} 

Enfin on d\'efinit
\nindex{$\tau_H\left(V\right)$}
\begin{equation}
\tau_H\left(V\right)=\omega_{V}\left(\adh{V(K)}\right),
\end{equation}
en d'autre termes $\tau_{H}\left(V\right)$ est le volume pour la mesure $\omega_{V}$ 
de l'adh\'erence de l'ensemble des points rationnels de $V$ dans l'espace ad\'elique $V(\ak)$.

La constante \nindex{$C^{\ast}_{V,H}$}$C^{\ast}_{V,H}$ est alors \'egale par d\'efinition \`a
$\frac{1}{\left(\rg(\Pic(V))-1\right)!}\,C_{V,H}$
o\`u
\nindex{$C_{V,H}$}
\begin{equation}
C_{V,H}=\alpha^{\ast}(V)\,\beta(V)\,\tau_H\left(V\right).
\end{equation}
Soulignons que la n\'ecessit\'e d'introduire la constante $\beta(V)$ a \'et\'e mise en \'evidence
par le r\'esultat de Batyrev et Tschinkel sur les vari\'et\'es toriques.

\section{Tores alg\'ebriques}

\subsection{Quelques rappels}

Soit $K$ un corps.
Un \dindex{tore alg\'ebrique}\termin{tore alg\'ebrique} d\'efini sur $K$ (de dimension~$d$) 
est un groupe alg\'ebrique $T$ d\'efini sur $K$ 
tel qu'il existe un isomorphisme de $\adh{K}$-groupes alg\'ebriques
\begin{equation}
T_{\adh{K}}\longisom (\G_{m,\adh{K}})^d.
\end{equation}
Si $T$ est un tore alg\'ebrique d\'efini sur $K$, on dit qu'une extension $L$ de $K$
d\'eploie $T$ s'il existe un isomorphisme de $L$-groupes alg\'ebriques
\begin{equation}
T_L\longisom (\G_{m,L})^d.
\end{equation}
Par \cite[Proposition 1.2.1]{Ono:algtor},
si $T$ est un tore alg\'ebrique d\'efini sur $K$
il existe une extension s\'eparable finie $L$ de $K$ qui d\'eploie $T$. 
En particulier il existe une extension galoisienne finie $L$ de $K$ qui d\'eploie $T$.

Soit $T$ un tore alg\'ebrique d\'efini sur $K$ de dimension $d$. 
On note $\nindex{$X(T)$}X(T)$ le groupe des caract\`eres de $T$, 
i.e. le groupe des morphismes de $K^{\sep}$-groupes alg\'ebriques
de $T_{K^{\sep}}$ vers $\G_{m,K^{\sep}}$.
C'est un $\Z$-module libre de rang $d$, 
sur lequel le groupe $\Gal(K^{\sep}/K)$ agit contin\^ument, 
et qui d\'epend fonctoriellement de $T$.
Si $L$ est une extension galoisienne de $K$ qui d\'eploie $T$,
l'action de $\Gal(K^{\sep}/K)$ sur $X(T)$ se factorise \`a travers $\Gal(L/K)$.
En outre le foncteur qui \`a un tore alg\'ebrique $T$ associe $X(T)$
d\'efinit une \'equivalence entre la cat\'egorie des tores alg\'ebriques d\'efinis sur $K$ 
et la cat\'egorie des $\Z$-modules libres de rang fini muni d'une action continue de $\Gal(K^{\sep}/K)$.
Si $L/K$ est un extension finie galoisienne, 
cette \'equivalence induit une \'equivalence
entre la cat\'egorie des tores alg\'ebriques d\'efini sur $K$ et d\'eploy\'es par $L$ 
et la cat\'egorie des $\Gal(L/K)$-modules qui sont libres de rang fini comme $\Z$-modules.

Soit $T$ un tore alg\'ebrique d\'efini sur $K$ et 
d\'eploy\'e par une extension finie galoisienne $L$ de $K$, de groupe 
de Galois $G$. Alors pour toute $K$-alg\`ebre $K'$ le groupe $T(K')$ des $K'$-points de $T$
s'identifie canoniquement \`a
\begin{equation}\label{eq:fonc:points:T}
\Hom_G(X(T),(L\otimes_K K')^{\times})=\left(X(T)^{\vee}\otimes (L\otimes_K K')^{\times}\right)^G.
\end{equation}
En particulier, le groupe $T(K)$ des points $K$-rationnels de $T$ s'identifie canoniquement \`a
\begin{equation}\label{eq:points:rat:T}
\Hom_G(X(T),L^{\times})=\left(X(T)^{\vee}\otimes L^{\times}\right)^G.
\end{equation}

\begin{exs}
Un exemple imm\'ediat de tore alg\'ebrique est fourni par les tores d\'eploy\'es, 
c'est-\`a-dire les groupes alg\'ebriques $K$-isomorphes \`a un produit de copies de $\G_{m,K}$.

Un autre exemple, 
important pour la suite, 
est donn\'e par la situation suivante : 
soit $K_0/K$ une extension finie s\'eparable, 
et $L/K$ une extension finie galoisienne de groupe $G$ contenant $K_0$. 
Soit $G_0$ le groupe de Galois de $L/K_0$. 
Au $G$-module 
$\Z[G/G_0]$ correspond, 
par l'\'equivalence de cat\'egories ci-dessus, 
la \termin{restriction \`a la Weil}  de $K_0$ \`a $K$ de $\G_m$, 
not\'ee \nindex{$\Res_{K_0/K}\G_m$}$\Res_{K_0/K}\G_m$. 
En particulier on a la propri\'et\'e
\begin{equation}
\left(\Res_{K_0/K}\G_m\right)(K)=\left(K_0\right)^{\times}.
\end{equation}

Un tore alg\'ebrique sur $K$ est dit 
\dindex{tore alg\'ebrique quasi-d\'eploy\'e}\termin{quasi-d\'eploy\'e} s'il est isomorphe 
sur $K$ \`a un produit de tores du type $\Res_{K_0/K}\G_m$. 
Un tore alg\'ebrique sur $K$
(respectivement un tore alg\'ebrique sur $K$ d\'eploy\'e par une extension galoisienne
finie de groupe $G$)
est quasi-d\'eploy\'e 
si et seulement si 
son groupe des caract\`eres est un $\Gal(K^{\sep}/K)$-module de permutation (respectivement
un $G$-module de permutation), 
c'est-\`a-dire poss\`ede une $\Z$-base stable sous l'action de $\Gal(K^{\sep}/K)$
respectivement $G$).

En fait dans la situation ci-dessus, pour tout tore alg\'ebrique $T$ d\'efini sur $K_0$, 
on peut d\'efinir la restriction \`a la Weil de $K_0$ \`a $K$ de $T$ (cf. \cite[\S 1.4]{Ono:algtor}), 
qui est un tore alg\'ebrique d\'efini sur $K$ v\'erifiant en particulier
\begin{equation}
\left(\Res_{K_0/K}T\right)(K)
=
T(K_0).
\end{equation}
Par la suite seul le cas $T=\G_m$ nous sera utile. 
\end{exs}

\begin{notas}
Soit $f\,:\,T_1\to T_2$ un morphisme de $K$-tores alg\'ebriques.
Le morphisme de $\Gal(K^{\sep}/K)$-modules $X(T_2)\to X(T_1)$
induit fonctoriellement par $f$ sera alors encore not\'e $f$.
Le morphisme de groupes $T_1(K)\to T_2(K)$  induit par $f$ sera not\'e 
\nindex{$f_K$, o\`u $f$ est un morphisme de tores alg\'ebriques}$f_K$. 
\end{notas}

\begin{lemme}\label{lm:res:qd}
Soit $G$ un groupe fini et $M$ un $\Z$-module libre de rang fini,
muni d'une action de $G$. Il existe une suite exacte de $G$-modules
\begin{equation}
0\to M \to P \to Q \to 0
\end{equation}
o\`u $P$ et $Q$ sont libres de rang fini en tant que $\Z$-modules,
et $P$ est un $G$-module de permutation.
\end{lemme}
\begin{rem}
Ce lemme \'el\'ementaire nous sera utile dans la preuve des propositions \ref{prop:conoyau:fini} 
et \ref{prop:compacite}. 
Une version plus fine sera utilis\'ee 
dans la section \ref{subsec:flasque}.
\end{rem}
\begin{demo}
Le dual $P^{\vee}$ d'un $G$-module de permutation $P$ \'etant encore un $G$-module de permutation,
il suffit de construire 
un morphisme de $G$-modules surjectif $P'\to M^{\vee}$,
avec $P'$ un $G$-module de permutation. Ceci peut se faire
de la mani\`ere suivante : soit $(e_1,\dots,e_d)$ une $\Z$-base de $M^{\vee}$
et pour $i=1,\dots,d$, $G_i$ le stabilisateur de $e_i$ pour l'action de $G$.
Pour $i=1,\dots,d$, il existe alors un unique morphisme de $G$-module $\Z[G/G_i]\to M^{\vee}$
envoyant $G_i$ sur $e_i$. Le morphisme somme
$\oplus_{i=1}^d \Z[G/G_i]\to M^{\vee}$ r\'epond \`a la question.
\end{demo}

\subsection{L'espace ad\'elique associ\'e \`a un tore alg\'ebrique}

Soit $K$ un corps global. 
On note \nindex{$\G_m(\ak)$}$\G_m(\ak)$ le groupe des id\`eles de $K$, 
muni de la topologie ad\'elique classique, 
qui en fait un groupe topologique ab\'elien localement compact. 
L'injection diagonale 
$\G_m(K)\hookrightarrow \G_m(\ak)$ 
identifie $\G_m(K)$ \`a un sous-groupe discret de $\G_m(\ak)$. 

Pour tout $v\in \placesde{K}$, 
on note 
$\G_m(\Ov)=\Ov^{\times}$ le sous-groupe compact
maximal de $\G_m(K_v)=K_v^{\times}$. 
Soit
\nindex{$\K(\G_m)$}
\begin{equation}
\K(\G_m)=\prod_{v\in \placesde{K}}\G_m(\Ov),
\end{equation}
c'est le sous-groupe compact maximal de $\G_m(\ak)$.

Si $L/K$ est une extension finie,
$\G_m(\ak)$ s'injecte naturellement dans $\G_m(\al)$. 
Si de plus $L/K$ est galoisienne de groupe $G$, 
on a une action naturelle de $G$ sur $\G_m(\al)$ et alors 
\begin{equation}
\G_m(\al)^{\,G}=\G_m(\ak).
\end{equation}

On note \nindex{$C_K$}$C_K=\G_m(\ak)/\G_m(K)$ le groupe des classes d'id\`eles de $K$.

Nous d\'ecrivons maintenant la g\'en\'eralisation de ces notions 
\`a un tore alg\'ebrique $T$ quelconque. 
Bien entendu pour $T=\G_m$ on retrouvera les d\'efinitions pr\'ec\'edentes. 
La construction de l'espace ad\'elique associ\'e \`a un tore alg\'ebrique $T$ est en fait un cas particulier 
de la construction g\'en\'erale de l'espace ad\'elique associ\'e \`a une vari\'et\'e alg\'ebrique d\'efinie sur $K$ 
(cf. \cite[I.2]{Wei:AAG} et \cite[I.3]{Oes:invent}). 
Elle peut se faire de la fa\c con suivante.

Pour toute place $v$ de $K$, 
$T(K_v)$ est muni naturellement d'une structure de groupe topologique ab\'elien localement compact. 
On d\'esigne alors (abusivement) par \nindex{$T(\Ov)$}$T(\Ov)$ le sous-groupe compact maximal de $T(K_v)$. 
En fait, si $S$ d\'esigne l'ensemble des places de $K$ archim\'ediennes ou ramifi\'ees dans une extension
galoisienne $L$ de groupe $G$ d\'eployant $T$ 
et $S_L$ 
l'ensemble des places de $L$ divisant une place de $S$, 
le sch\'ema en groupes
\begin{equation}
\ecT=\Spec\left(\cO_{S_L}\otimes X(T)\right)^{G}
\end{equation}
est un mod\`ele de $T$ sur $\Spec\left(\cO_{S}\right)$, 
et pour toutes les places $v$ en dehors de $S$ on a $T(\Ov)=\ecT(\Ov)$, 
d'o\`u la notation adopt\'ee.

L'espace ad\'elique associ\'e \`a $T$ 
est alors le sous-groupe du groupe produit $\produ{v} T(K_v)$
d\'ecrit par
\nindex{$T(\ak)$}
\begin{equation}
T(\ak)
=
\left\{
(t_v)\in \prod_v T(K_v),\,\,t_v\in T(\Ov)\,\,\text{pour presque tout }
v\in\placesde{K}
\right\}
.
\end{equation}
Pour un sous-ensemble fini $S$ de $\placesde{K}$, on consid\'erera aussi
\nindex{$T(\ak)_S$}
$T(\ak)_S$ le sous-groupe de $\produ{v} T(K_v)$ d\'ecrit par
\begin{equation}
T(\ak)_S
=
\left\{
(t_v)\in \prod_v T(K_v),\,\,\forall v\notin S,\,\,t_v\in T(\Ov)
\right\}
,
\end{equation}
de sorte que $T(\ak)$ est la r\'eunion des $T(\ak)_S$ pour $S$ d\'ecrivant
l'ensemble des parties finies de $\placesde{K}$.

Comme sous-groupe de $\produ{v} T(K_v)$, $T(\ak)$ est un groupe ab\'elien.
On le munit de la topologie dont une base 
d'ouverts est donn\'ee par les sous-ensembles du type
\begin{equation}
\prod_{v\in S}U_v\times\prod_{v\notin S} T(\Ov)
\end{equation}
o\`u $S$ est un ensemble fini de places de $K$ et, 
pour $v\in S$, $U_v$ est un ouvert de $T(K_v)$. 
Cette topologie, qui est plus fine que la topologie issue de la topologie produit sur $\prod T(K_v)$, 
fait de $T(\ak)$ un groupe topologique localement compact. 
On peut alors identifier $T(K)$ \`a un sous-groupe discret de $T(\ak)$.

On note
\nindex{$\K(T)$}
\begin{equation}
\K(T)=T(\ak)_{\vide}=\prod_{v\in\placesde{K}} T(\Ov),
\end{equation}
c'est le sous-groupe compact maximal de $T(\ak)$.

Si $L$ est une extension finie galoisienne de $K$ d\'eployant $T$, 
de groupe de Galois $G$, 
on dispose \'egalement, 
comme pour $T(K)$, 
d'une description simple de tous ces groupes en termes du $G$-module $X(T)$, 
pratique pour manipuler des suites exactes. 
Pour toute place $v$ de $K$, on a d'apr\`es \eqref{eq:fonc:points:T}
une identification canonique
\begin{equation}\label{eq:points:T:Kv}
T(K_v)\longisom \Hom_G\left(X(T),\left(L\otimes_K K_v\right)^{\times}\right).
\end{equation}
Rappelons que $\left(L\otimes_K K_v\right)^{\times}$ s'identifie \`a $\underset{\V|v}{\prod}L_{\V}^{\times}$. 
L'identification \eqref{eq:points:T:Kv}
induit une identification
\begin{equation}\label{eq:points:T:Ov}
T(\Ov)\longisom \Hom_G\left(X(T),\prod_{\V|v}\OV^{\times}\right).
\end{equation}
Si on choisit une place $\V$ divisant $v$ et 
si on note $G_v$ son groupe de d\'ecomposition, 
les identifications \eqref{eq:points:T:Kv} et \eqref{eq:points:T:Ov}
induisent respectivement des identifications
\begin{equation}\label{eq:points:T:Kv:Gv}
T(K_v)\longisom \Hom_{G_v}\left(X(T),L_{\V}^{\times}\right)
\end{equation}
et
\begin{equation}\label{eq:points:T:Ov:Gv}
T(\Ov)\longisom \Hom_{G_v}\left(X(T),\OV^{\times}\right).
\end{equation}
On peut identifier $T(\ak)$ au groupe
\begin{equation}
\Hom_G(X(T),\G_m(\al))
\end{equation}
muni de la topologie induite par celle de $\G_m(\al)$, 
et $\K(T)$ au groupe
\begin{equation}
\Hom_G(X(T),\K(\G_{m,L})).
\end{equation}
\begin{ex} 
Soit $K_0$ une extension finie s\'eparable de $K$, et $T=\Res_{K_0/K}\G_m$. 
Alors $T(\ak)$ s'identifie canoniquement \`a $\G_m\left(\A_{K_0}\right)$.  
\end{ex}
\begin{notas}
Soit $f\,:\,T_1\to T_2$ un morphisme $K$-de tores alg\'ebriques, et $v$ une place
de $K$. Le morphisme continu de groupes topologiques $T_2(\ak)\to T_1(\ak)$ 
(respectivement $T_2(K_v)\to T_1(K_v)$) induit fonctoriellement
par $f$ sera alors not\'e \nindex{$f_{\ak}$, o\`u $f$ est un morphisme de tores alg\'ebriques}
$f_{\ak}$ (respectivement \nindex{$f_v$, o\`u $f$ est un morphisme de tores alg\'ebriques}$f_v$). 
\end{notas}

\subsection{Le degr\'e}
\label{subsec:ledegre}
\subsubsection{D\'efinitions}\label{subsec:def:ledegre}

Soit $K$ un corps global. Pour toute place $v$ de $K$, 
on a un morphisme de valuation 
d\'efini par 
\nindex{$\deg_{K,v}$}
\begin{equation}
\deg_{K,v}\,:\,\map{\G_m(K_v)}{\Z}
{x_v}{v(x_v)}
\end{equation}
si $v$ est finie et
\begin{equation}\label{eq:valumarchi}
\deg_{K,v}\,:\,
\map{\G_m(K_v)}{\R}
{x_v}{\log \abs{x_v}_v}
\end{equation}
si $v$ est archim\'edienne.
Dans les deux cas, ce morphisme a pour noyau $\G_m(\Ov)$.

On en d\'eduit un morphisme <<degr\'e>>
d\'efini par
\nindex{$\deg_{K}$}
\begin{equation}
\deg_K\,:\,
\map{\G_m(\ak)}{\R}
{(x_v)}{\sumu{v\in \placesde{K}}\log(q_v)\,\deg_{K,v}(x_v)}
\end{equation}
dans le cas arithm\'etique
et
\begin{equation}
\deg_K\,:\,\map{\G_m(\ak)}{\Z}
{(x_v)}{\sumu{v\in \placesde{K}}f_v\,\deg_{K,v}(x_v)}.
\end{equation}
Dans le cas arithm\'etique, le morphisme $\deg_K$ est surjectif, 
car le morphisme $\deg_{K,v}$ est surjectif
si $v$ est archim\'edienne.
Dans le cas fonctionnel,
$\deg_K$ est surjectif d'apr\`es \cite[VII \S$\,$5, Cor 6]{Wei:BNT}.

Dans les deux cas, 
le noyau de $\deg_K$ contient
$\G_m(K)$ (par la formule du produit) et $\K(\G_m)$.
Il sera not\'e \nindex{$\G_m(\ak)^1$}$\G_m(\ak)^1$.

Soit $L$ une extension finie 
de $K$. Dans le cas arithm\'etique, 
on a la relation
\begin{equation}\label{eq:formule_degre:arit}
{\deg_L}_{|_{\G_m(\ak)}}=[L:K]\,\deg_K.
\end{equation}
Dans le cas fonctionnel, 
on a la relation
\begin{equation}\label{eq:formule_degre}
\dde{L}\,{\deg_L}_{|_{\G_m(\ak)}}=\dk\,[L:K]\,\deg_K,
\end{equation}
o\`u l'on rappelle que $q^{\,\dde{L}}$ et $q^{\,\dk}$ 
repr\'esentent les cardinaux des corps des constantes 
des corps de fonctions $L$ et $K$ respectivement.
\\~\\
\indent
Soit \`a pr\'esent $T$ un tore alg\'ebrique d\'efini sur $K$, 
d\'eploy\'e par une extension finie galoisienne $L$ de groupe de Galois $G$. 
Soit $m$ un \'el\'ement de $X(T)^{\,G}$, c'est-\`a-dire un morphisme 
de $K$-groupes $T\to \G_m$
(le morphisme dual est alors le morphisme de $G$-modules $\Z\to X(T)$ qui envoie $1$ sur $m$).

Par composition du morphisme continu $m_{\ak}\,:\,T(\ak)\to \G_m(\ak)$ 
avec 
$\deg_K$ 
on obtient un morphisme continu 
$
T(\ak)\to \R
$
dans le cas arithm\'etique et
$
T(\ak)\to \Z
$
dans le cas fonctionnel.
On note \nindex{$\deg_T$}$\deg_T$
le morphisme d\'efini par 
\begin{equation}
\begin{array}{rl}
T(\ak)&\longto \Hom\left(X(T)^{\,G},\R\right)\\
t&\longmapsto \left[m\mapsto(\deg_K\circ m_{\ak})(t)\right]
\end{array}
\end{equation}
dans le cas arithm\'etique,
et
\begin{equation}
\begin{array}{rl}
T(\ak)&\longto \Hom\left(X(T)^{\,G},\Z\right)\\
t&\longmapsto \left[m\mapsto(\deg_K\circ m_{\ak})(t)\right]
\end{array}
\end{equation}
dans le cas fonctionnel. 
Dans ce dernier cas, 
le morphisme $\deg_T$ 
n'est autre que le morphisme $\theta$ d\'efini par Oesterl\'e dans \cite[I.5.5]{Oes:invent}, 
compos\'e avec $\log_{\qde{K}}$. 

Dans les deux cas, on a $\deg_{\G_m}=\deg_K$.
Le r\'esultat suivant est imm\'ediat.
\begin{lemme}\label{lm:fonc:degre}
Le morphisme $\deg_T$ est fonctoriel dans le sens suivant : 
soit
\begin{equation}
f\,:\,T_1\longto T_2
\end{equation}
est un morphisme de $K$-tores alg\'ebriques d\'eploy\'es par $L$.
Le morphisme
\begin{equation}
f\,:\,X(T_2)\longto X(T_1)
\end{equation}
de $G$-module associ\'e
induit par dualit\'e des morphismes 
\begin{equation}
f^{\vee}\,:\,\left(X(T_1)^{G}\right)^{\vee}\longto \left(X(T_2)^{G}\right)^{\vee},
\end{equation}
et
\begin{equation}
f^{\vee}_{\R}\,:\,\left(X(T_1)^{G}\right)_{\R}^{\vee}\longto \left(X(T_2)^{G}\right)_{\R}^{\vee}.
\end{equation}
On a alors
\begin{equation}
\deg_{T_2}\circ f_{\ak}=f^{\vee}_{\R}\circ \deg_{T_1}
\end{equation}
dans le cas arithm\'etique
et
\begin{equation}
\deg_{T_2}\circ f_{\ak}=f^{\vee}\circ \deg_{T_1}
\end{equation}
dans le cas fonctionnel.
\end{lemme}

Pour un tore alg\'ebrique $T$ quelconque, \`a l'instar de $\deg_K$, 
$\deg_T$ se d\'ecompose en une somme de degr\'es locaux,
que nous d\'ecrivons \`a pr\'esent.
Soit $v$ une place de $K$, $G_v$ un groupe
de d\'ecomposition au-dessus de $v$ et $m$ un \'el\'ement de $X(T)^{G_v}$,
c'est-\`a-dire un morphisme de $K_v$-tores $T_{K_v}\to \G_{m,K_v}$.
Par composition du morphisme continu $m_v\,:\,T(K_v)\to \G_m(K_v)$ avec $\deg_{K,v}$ on obtient un morphisme continu 
$
T(K_v)\to \Z
$
si $v$ est finie et
$
T(K_v)\to \R
$
si $v$ est archim\'edienne. 
On note \nindex{$\deg_{T,v}$}$\deg_{T,v}$
le morphisme d\'efini par 
\begin{equation}
\begin{array}{rl}
T(K_v)&\longto \Hom\left(X(T)^{\,G_v},\Z\right)\\
t&\longmapsto \left[m\mapsto(\deg_{K,v}\circ m_{v})(t)\right]
\end{array}
\end{equation}
si $v$ est finie et
\begin{equation}
\begin{array}{rl}
T(K_v)&\longto \Hom\left(X(T)^{\,G_v},\R\right)\\
t&\longmapsto \left[m\mapsto(\deg_{K,v}\circ m_{v})(t)\right]
\end{array}
\end{equation}
si $v$ est archim\'edienne. 
Dans les deux cas, le noyau de $\deg_{T,v}$ est $T(\Ov)$.

Notons \nindex{$i_{T,v}$}$i_{T,v}$ l'injection continue naturelle de groupes topologiques $T(K_v)\to T(\ak)$.

Dans le cas arithm\'etique, pour toute place $v$, le diagramme 
\begin{equation}
\xymatrix{
T(K_v)\ar[r]^{i_{T,v}}\ar[d]^{\log(q_v)\,\deg_{T,v}}&T(\ak)\ar[d]^{\deg_{T}} \\
\Hom\left(X(T)^{\,G_v},\R\right)\ar[r] & \Hom(X(T)^G,\R)
}
\end{equation}
(o\`u la fl\`eche horizontale du bas est le morphisme naturel de restriction) 
est commutatif, et on a
\begin{equation}\label{eq:dec:deg_T:arit}
\deg_T=\sum_{v\in \placesde{K}}\log(q_v)\,\deg_{T,v}.
\end{equation}

Dans le cas fonctionnel, 
pour toute place $v$ le diagramme 
\begin{equation}
\xymatrix{
T(K_v)\ar[r]^{i_{T,v}}\ar[d]^{f_v\,\deg_{T,v}}&T(\ak)\ar[d]^{\deg_{T}} \\
\Hom\left(X(T)^{\,G_v},\Z\right)\ar[r] & \Hom(X(T)^G,\Z)
}
\end{equation}
(o\`u la fl\`eche horizontale du bas est le morphisme naturel de restriction) 
est commutatif, et on a 
\begin{equation}\label{eq:dec:deg_T:fonc}
\deg_T=\sum_{v\in \placesde{K}}f_v\,\deg_{T,v}.
\end{equation}

Le noyau de $\deg_T$ sera not\'e \nindex{$T(\ak)^1$}$T(\ak)^1$. 
Par la formule du produit, 
$T(K)$ est contenu dans $T(\ak)^1 $. 
Par ailleurs, comme chaque morphisme $\deg_{T,v}$ 
a pour noyau $T(\Ov)$, 
$T(\ak)^1$ contient $\K(T)$. 
En outre, 
$T(\ak)^1$  s'identifie au groupe
\begin{equation}
T(\ak)^1=\Hom_G(X(T),\G_m(\al)^1).
\end{equation}
Pour toute partie finie $S$ de $\placesde{K}$, on note
\nindex{$T(\ak)^1_S$}
\begin{equation}
T(\ak)^1_S=T(\ak)_S\cap T(\ak)^1.
\end{equation}

On d\'efinit \`a pr\'esent une variante du morphisme $\deg_T$
et des degr\'es locaux $\deg_{T,v}$ qui nous sera utile par la suite.
Contrairement \`a $\deg_T$, la d\'efinition de cette variante
d\'epend du choix de l'extension $L$ d\'eployant $T$.
Pour comparer les deux notions, nous aurons besoin du lemme \'el\'ementaire suivant.
\begin{lemme}
\label{lm:anis}
Soit $M$ un $G$-module qui est un $\Z$-module libre de rang fini. 
Alors la fl\`eche naturelle
\begin{equation}
\left(M^{\vee}\right)^G
\longto
\left(M^G\right)^{\vee}
\end{equation}
est une injection de conoyau fini.
\end{lemme}
\begin{demo}
Le quotient $M/M^G$ \'etant sans torsion, on a une suite exacte
\begin{equation}
0\to 
(M/M^G)^{\vee}
\to M^{\vee} 
\to \left(M^G\right)^{\vee}
\to 0
\end{equation}
En prenant les $G$-invariants, 
on obtient la suite exacte
\begin{equation}
0\to 
\left(\left(M/M^G\right)^{\vee}\right)^G
\to \left(M^{\vee} \right)^G
\to \left(M^G\right)^{\vee}
\to H^1\left(G,(M/M^G)^{\vee}\right)
\end{equation}
dont l'avant-derni\`ere fl\`eche est la fl\`eche de l'\'enonc\'e.
Son conoyau est donc fini.
Soit $\phi$ un \'el\'ement de $\left(M^{\vee} \right)^G=\Hom_G(M,\Z)$ 
dont la restriction \`a $M^G$ est nulle.
Pour tout $x\in M$, $N_G\,x\eqdef \sumu{g\in G} g.x$ est un \'el\'ement de $M^G$ 
et  on a 
\begin{equation}
0=\phi(N_G\,x)=\card{G}\phi(x)
\end{equation}
donc $\phi(x)=0$.
Ainsi $\phi$ est nulle, 
d'o\`u l'injectivit\'e 
(en d'autres termes,
le dual d'un $G$-module anisotrope est anisotrope).
\end{demo}

On se place dans le cas arithm\'etique.
On consid\`ere la suite exacte
\begin{equation}
0\longto \G_m(\al)^1 \longto \G_m(\al) \overset{\deg_L}{\longto} \R\longto 0.
\end{equation}
Tensorisons par $X(T)^{\vee}$ et prenons les $G$-invariants. 
On obtient la suite exacte
\begin{equation}
0 \longto T(\ak)^1  \longto T(\ak)\longto \left(X(T)^{\vee}\right)^{\,G}_{\R}
\end{equation}
et donc un morphisme 
\begin{equation}
\nindex{$\deg_{T,L}$ (cas arithm\'etique)}\deg_{T,L}\,:\,T(\ak)\longto \left(X(T)^{\vee}\right)^{\,G},
\end{equation}
de noyau $T(\ak)^1$. 
On notera encore $\deg_{T,L}$ le morphisme obtenu par composition avec le
morphisme de restriction
$\left(X(T)^{\vee}\right)^{\,G}_{\R}
\hookrightarrow 
\left(X(T)^{\,G}\right)^{\vee}_{\R}$  (qui est un isomorphisme d'apr\`es le lemme \ref{lm:anis}),
On a donc, pour tout $t\in T(\ak)$ et tout $m\in X(T)^G$,
\begin{equation}
\acc{\deg_{T,L}(x)}{m}=\deg_L(m_{\ak}(t)).
\end{equation}
D'apr\`es \eqref{eq:formule_degre:arit}
on a donc 
\begin{equation}\label{eq:formule_degre:arit:degTL}
\deg_{T,L}=[L:K]\deg_T,
\end{equation}
ce qui montre en particulier que
contrairement au morphisme $\deg_T$, 
le morphisme $\deg_{T,L}$ d\'epend du choix de l'extension d\'eployant $T$.

Remarquons que le diagramme
\begin{equation}
\xymatrix{
T(\ak)\ar[rr]^<<<<<<<<<<{\deg_{T,L}}\ar[d]&&\Hom_G(X(T),\R) \ar[d] \\
T_L(\al) \ar[rr]^<<<<<<<<<<{\deg_{T_L}} &&\Hom(X(T),\R)
}
\end{equation}
(o\`u les fl\`eches verticales sont les inclusions naturelles)
est commutatif.

Pla\`c;ons nous \`a pr\'esent dans le cas fonctionnel.
On consid\`ere la suite exacte
\begin{equation}
0\longto \G_m(\al)^1 \longto \G_m(\al) \overset{\deg_L}{\longto} \Z\longto 0.
\end{equation}
Tensorisons par $X(T)^{\vee}$ et prenons les $G$-invariants. 
On obtient la suite exacte
\begin{equation}
0 \longto T(\ak)^1  \longto T(\ak)\longto \left(X(T)^{\vee}\right)^{\,G}
\end{equation}
et donc un morphisme 
\begin{equation}
\nindex{$\deg_{T,L}$ (cas fonctionnel)}\deg_{T,L}\,:\,T(\ak)\longto \left(X(T)^{\vee}\right)^{\,G},
\end{equation}
de noyau $T(\ak)^1$. 
On notera encore $\deg_{T,L}$ le morphisme obtenu par composition avec le
morphisme de restriction
$\left(X(T)^{\vee}\right)^{\,G} 
\hookrightarrow 
\left(X(T)^{\,G}\right)^{\vee}$  (qui est injectif d'apr\`es le lemme \ref{lm:anis}),
On a donc, pour tout $t\in T(\ak)$ et tout $m\in X(T)^G$,
\begin{equation}
\acc{\deg_{T,L}(t)}{m}=\deg_L(m_{\ak}(t)).
\end{equation}
D'apr\`es \eqref{eq:formule_degre},
on a donc la relation
\begin{equation}
d_L\deg_{T,L}=d_K\,[L:K]\deg_T
\end{equation}
ce qui montre que contrairement au morphisme $\deg_T$, 
le morphisme $\deg_{T,L}$ d\'epend du choix de l'extension d\'eployant $T$.

Remarquons que le diagramme
\begin{equation}
\xymatrix{
T(\ak)\ar[rr]^<<<<<<<<<<{\deg_{T,L}}\ar[d]&&\Hom_G(X(T),\Z) \ar[d] \\
T_L(\al) \ar[rr]^<<<<<<<<<<{\deg_{T_L}} &&\Hom(X(T),\Z)
}
\end{equation}
(o\`u les fl\`eches verticales sont les inclusions naturelles)
est commutatif.

Tout comme le morphisme $\deg_{T}$, 
le morphisme $\deg_{T,L}$ peut se d\'ecomposer (de mani\`ere non canonique) 
en une somme de degr\'es locaux,
que nous d\'ecrivons \`a pr\'esent.

Soit $v$ une place finie de $K$, $\V$ une place de $L$ divisant $v$
et $G_v$ son groupe de d\'ecomposition.
On consid\`ere la suite exacte de $G_v$-modules.
\begin{equation}
0\longto \OV^{\times} \longto K_{\V}^{\times} \overset{\V}{\longto} \Z\longto 0.
\end{equation}
Tensorisons par $X(T)^{\vee}$ et prenons les $G_v$-invariants. 
On obtient la suite exacte
\begin{equation}
0 \longto T(\Ov)  \longto T(K_v)\longto \left(X(T)^{\vee}\right)^{\,G_v}
\end{equation}
et donc un morphisme 
\begin{equation}\label{eq:def:degTLV}
\nindex{$\deg_{T,L,\V}$, $\V$ finie}\deg_{T,L,\V}\,:\,T(K_v)\longto \left(X(T)^{\vee}\right)^{\,G_v}
\end{equation} 
de noyau $T(\Ov)$. 
Remarquons que le diagramme
\begin{equation}\label{eq:diag:degtlvtkv}
\xymatrix{
T(K_v)\ar[rr]^<<<<<<<<<<{\deg_{T,L,\V}}\ar[d]&&\Hom_{G_v}(X(T),\Z) \ar[d] \\
T_L(L_{\V}) \ar[rr]^<<<<<<<<<<{\deg_{T_L,\V}} &&\Hom(X(T),\Z)
}
\end{equation}
(o\`u les fl\`eches verticales sont les inclusions naturelles)
est commutatif.

Soit \`a pr\'esent $v$ une place archim\'edienne, $\V$ une place de $L$ divisant $v$
et $G_v$ son groupe de d\'ecomposition.
On consid\`ere la suite exacte de $G_v$-modules.
\begin{equation}
0\longto \OV^{\times} \longto K_{\V}^{\times} \overset{\log \abs{\,.\,}_{\V}}{\longto} \R\longto 0.
\end{equation}
Tensorisons par $X(T)^{\vee}$ et prenons les $G_v$-invariants. 
On obtient la suite exacte
\begin{equation}
0 \longto T(\Ov) \longto T(K_v)\longto \left(X(T)^{\vee}\right)^{\,G_v}_{\R}
\end{equation}
et donc un morphisme 
\begin{equation}
\nindex{$\deg_{T,L,\V}$, $\V$ archim\'edienne}\deg_{T,L,\V}\,:\,T(K_v)\longto \left(X(T)^{\vee}\right)_{\R}^{\,G_v}
\end{equation} 
de noyau $T(\Ov)$.
Remarquons que le diagramme
\begin{equation}\label{eq:diag:degtlvtkv:archi}
\xymatrix{
T(K_v)\ar[rr]^<<<<<<<<<<{\deg_{T,L,\V}}\ar[d]&&\Hom_{G_v}(X(T),\R) \ar[d] \\
T_L(L_{\V}) \ar[rr]^<<<<<<<<<<{\deg_{T_L,\V}} &&\Hom(X(T),\R)
}
\end{equation}
(o\`u les fl\`eches verticales sont les inclusions naturelles)
est commutatif.

Dans le cas arithm\'etique, pour toute place $v$ le diagramme 
\begin{equation}
\xymatrix{
T(K_v)\ar[r]^{i_{T,v}}\ar[d]^{\log(q_{\V})\,\deg_{T,L,\V}}&T(\ak)\ar[d]^{\deg_{T,L}} \\
\Hom_{G_v}(X(T),\R)  \ar[r] & \Hom(X(T)^G,\R)
}
\end{equation}
(o\`u la fl\`eche horizontale du bas est le morphisme de restriction)
est  commutatif,
et on a 
\begin{equation}\label{eq:dec:degTL:arit}
\deg_{T,L}=\sum_{v\in \placesde{K}} \log(q_{\V})\,\deg_{T,L,\V}.
\end{equation}

Dans le cas fonctionnel, pour toute place $v$ le diagramme 
\begin{equation}
\xymatrix{
T(K_v)\ar[r]^{i_{T,v}}\ar[d]^{f_{\V}\,\deg_{T,L,\V}}&T(\ak)\ar[d]^{\deg_{T,L}} \\
\Hom_{G_v}(X(T),\Z)  \ar[r] & \Hom(X(T)^G,\Z)
}
\end{equation}
(o\`u la fl\`eche horizontale du bas est le morphisme de restriction)
est  commutatif
et on a
\begin{equation}\label{eq:dec:degTL:fonc}
\deg_{T,L}=\sum_{v\in \placesde{K}} f_{\V}\,\deg_{T,L,\V}.
\end{equation}

\subsubsection{Propri\'et\'es du degr\'e local dans le cas d'une place finie 
}\label{subsubsec:ledegreafini}
On consid\`ere toujours un tore alg\'ebrique $T$ d\'efini sur un corps global $K$
et d\'eploy\'e par une extension galoisienne $L$ de groupe $G$
\begin{lemme}
\label{lemme:draxl}
Soit $v$ une place finie de $K$ non ramifi\'ee dans $L$, $\V$ une place
de $L$ divisant $v$ et $G_v$ le groupe de d\'ecomposition correspondant.
Alors le morphisme
\begin{equation}\label{eq:degTLV}
\deg_{T,L,\V}\,:\,T(K_v)\longto \left(X(T)^{\vee}\right)^{\,G_v}
\end{equation}
est surjectif, et induit donc un isomorphisme
\begin{equation}
T(K_v)/T(\Ov)\longisom (X(T)^{\vee})^{G_v}.
\end{equation}
\end{lemme}
\begin{rem}
On verra plus loin (cf. la proposition \ref{prop:draxl})
que pour une place finie quelconque le morphisme \eqref{eq:degTLV}
est de conoyau fini.
\end{rem}
\begin{demo}
L'argument qui suit est repris de la page 449 de \cite{Dr}.
Au vu de la construction de $\deg_{T,L,\V}$, le conoyau du morphisme
\eqref{eq:degTLV} est le groupe $H^1(G_v,X(T)^{\vee}\otimes \OV^{\times})$.
Or, 
si  $v$ est non ramifi\'ee, 
$\OV^{\times}$ est cohomologiquement trivial, 
et comme $X\left(T\right)^{\vee}$ est sans torsion, 
d'apr\`es \cite[IX, $\S\,5$, Corollaire]{Ser:corps}, 
$X\left(T\right)^{\vee}\otimes \OV^{\times}$ est encore
cohomologiquement trivial, 
d'o\`u le r\'esultat. 
\end{demo}
\begin{lemme}\label{lm:degTv:degTLV:ari}\label{lm:degTv:degTLV:fonc}
Soit $v$ une place finie de $K$,
$\V$ une place de $L$ divisant $v$, $G_v$
son groupe de d\'ecomposition et $e_v$ l'indice de ramification de $v$ dans $L$.

On a alors
\begin{equation}
\deg_{T,L,\V}=e_v\,\deg_{T,v}.
\end{equation}
\end{lemme}
\begin{demo}
On a en effet, pour $t\in T(K_v)$,
et pour $m\in X(T)^G$,
\begin{align}
\acc{\deg_{T,L,\V}(t)}{m}
&=
\V(m_v(t))
\\
&=
e_v \,v(\acc{x}{t})
\\
&
=
e_v\,\acc{\deg_{T,v}(t)}{m}.
\end{align}
\end{demo}

\subsubsection{Propri\'et\'es du degr\'e dans le cas arithm\'etique}\label{subsubsec:ledegrearit}

\begin{lemme}\label{lm:degTv:degTLV:ari:archi}
Soit $v$ une place archim\'edienne de $K$,
$\V$ une place de $L$ divisant $v$, $G_v$
son groupe de d\'ecomposition.

On a alors
\begin{equation}
\deg_{T,L,\V}={[L_{\V}:K_v]}\,\deg_{T,v}
\end{equation}
\end{lemme}
\begin{demo}
Soit $t\in T(K_v)$ et $m\in X(T)^G$. On a
\begin{align}
\acc{\deg_{T,L,\V}(x)}{t}
&=
\log\left(\abs{m_v(t)}_{\V}\right)
\\
&=
{[L_{\V}:K_v]}\,\log\left(\abs{m_v(t)}_{v}\right)
\\
&
=
{[L_{\V}:K_v]}\,\acc{\deg_{T,v}(t)}{m}.
\end{align}
\end{demo}
\begin{lemme}\label{lm:varchitkvtov}
Soit $v$ une place archim\'edienne de $K$,
$\V$ une place de $L$ divisant $v$ et $G_v$ son 
groupe de d\'ecomposition. 
Le morphisme 
\begin{equation}\label{eq:varchi:tkvto}
\deg_{T,L,\V}\,:\,T(K_v)\longto \Hom_{G_v}(X(T),\R)
\end{equation}
est surjectif.

En particulier le morphisme
\begin{equation}
\deg_{T,L,\V}\,:\,T(K_v)\longto \Hom(X(T)^G,\R)
\end{equation}
est surjectif.
\end{lemme}
\begin{demo}
Le morphisme \eqref{eq:varchi:tkvto}
est surjectif car il admet pour section le morphisme
\begin{equation}
\Hom_{G_v}(X(T),\R) \longto \Hom_{G_v}(X(T),\R^{\times})\subset \Hom_{G_v}(X(T),L_{\V}^{\times})\isom T(K_v) 
\end{equation}
obtenu par composition soit avec l'exponentielle si $L_{\V}=\R$,
soit avec le carr\'e de l'exponentielle si $L_V\isom \C$.

La deuxi\`eme assertion vient du fait 
que le morphisme 
de $\R$-espaces vectoriels
\begin{equation}
\Hom_{G_v}(X(T),\R)\to \Hom(X(T)^{G_v},\R)
\end{equation}
est un isomorphisme, 
et que le morphisme de restriction
\begin{equation}
\Hom(X(T)^{G_v},\R)\to \Hom(X(T)^{G},\R)
\end{equation}
est surjectif.
\end{demo}

\begin{cor}\label{cor:ari:surj:degT}
Dans le cas arithm\'etique, $\deg_T$ et $\deg_{T,L}$  sont surjectifs.
\end{cor}
\begin{demo}
Pour toute place $v$ archim\'edienne, on a d'apr\`es le lemme
\ref{lm:varchitkvtov}
\begin{equation}
\deg_{T,L,\V}(T(K_v))=\Hom(X(T)^G,\R).
\end{equation}
ce qui montre la surjectivit\'e de $\deg_{T,L}$
d'apr\`es \eqref{eq:dec:degTL:arit}.

Par ailleurs, pour toute place $v$ archim\'edienne, 
$T(K_v)$ est un groupe divisible,
et $\deg_{T,v}$ et $\deg_{T,L,\V}$ sont proportionnels.
On en d\'eduit qu'on a 
\begin{equation}
\deg_{T,v}(T(K_v))=\Hom(X(T)^G,\R).
\end{equation}
ce qui montre la surjectivit\'e de $\deg_{T}$
d'apr\`es \eqref{eq:dec:deg_T:arit}.
\end{demo}

\subsubsection{Propri\'et\'es du degr\'e dans le cas fonctionnel}

Contrairement \`a ce qui se passe dans le cas 
arithm\'etique,
$\deg_T$ n'est plus n\'ecessairement surjectif dans le cas fonctionnel, 
comme le montre le lemme suivant.

\begin{lemme}
\label{lm:degre:quasidep}
Soit $K$ un corps de fonctions et $L/K$ une extension finie galoisienne,
de groupe de Galois $G$.
Soit $K_0$ une extension de $K$ contenue dans $L$, 
$T$ le tore $\Res_{K_0/K}\G_m$ et $G_0=\Gal(L/K_0)$. 
Soit $d_0$ tel que $q_{\!\text{\tiny {\it K}}_0}=\qk^{d_0}$, 
de sorte que $d_0=\frac{{d_{\!\text{\tiny {\it K}}_0}}}{\dk}$. 
Alors, via les identifications naturelles $T(\ak)=\G_m\left(\A_{K_0}\right)$ et $(X(T)^{\,G})^{\vee}=\Z$, 
le morphisme $\deg_T$ n'est autre que le morphisme $d_0\,\deg_{K_0}$, et son image est $d_0\,\Z$. 
En particulier $\deg_T$ n'est pas n\'ecessairement surjectif.
\end{lemme}

\begin{demo} 
Soit $t\in \G_m(\ade{K_0})$. 
Via l'identification 
\begin{equation}
T(\ak)\longisom \G_m(\ade{K_0}),
\end{equation}
il lui correspond l'\'el\'ement de 
\begin{equation}
T(\ak)=\Hom_G(\Z[G/G_0],\G_m(\al))
\end{equation}
qui envoie $G_0$ sur $t$. 
Le $\Z$-module $\Z[G/G_0]^{\,G}$ est de rang 1 engendr\'e par $\underset{g\in G/G_0}{\sum} g\,G_0$. 
Le morphisme $\deg_T(t)$ envoie alors  $\underset{g\in G/G_0}{\sum} g\,G_0$ sur
\begin{align}
\deg_K
\left(
\prod_{g\in G/G_0}g\,t\right)
&=
\frac{\dde{L}}{\dk\,[L:K]}\sum_{g\in G/G_0}\deg_L(g\,t)\\
&=
\frac{\dde{L}\,[G]}{\dk\,[L:K]\,[G_0]}\deg_L(t)\\
&=
\frac{\dde{L}}{\dk\,[L:K_0]}\deg_L(t)\\
&=
d_0\,\frac{\dde{L}}{d_{\!\text{\tiny {\it K}}_0}\,[L:K_0]}\deg_L(t)\\
&=
d_0\,\deg_{K_0}(t),
\end{align}
la premi\`ere et la derni\`ere \'egalit\'e venant de la formule \eqref{eq:formule_degre}. 
Le fait que l'image est $d_0\,\Z$ d\'ecoule alors de l'existence de diviseurs de degr\'e 1 (\cite[VII, \S $\,$5,Cor 6]{Wei:BNT}).
\end{demo}

\begin{nota}
Dans le cas fonctionnel, on note \nindex{$\DT$}$\DT$
l'image du morphisme $\deg_T$.
\end{nota}

\begin{prop}\label{prop:conoyau:fini}
On se place dans le cas fonctionnel.
Le morphisme
\begin{equation}
\deg_{T,L}\,:\,T(\ak)\longto \left(X(T)^{\vee}\right)^{\,G}
\end{equation}
est de conoyau fini, et il en est de m\^eme du morphisme
\begin{equation}
\deg_T\,:\,T(\ak)\longto \left(X(T)^{G}\right)^{\vee}.
\end{equation}
\end{prop}
\begin{demo} 
C'est un cas particulier de \cite[I.5.6.b]{Oes:invent}, o\`u la preuve
est donn\'ee pour tout groupe lin\'eaire alg\'ebrique. 
Nous donnons  ici une preuve pour les tores alg\'ebriques.

D'apr\`es le lemme \ref{lm:res:qd}, il existe une suite exacte de $G$-modules
\begin{equation}
0\to X(T) \to P \to Q\to 0.
\end{equation}
o\`u $P$ et $Q$ sont libres de rang fini en tant que $\Z$-modules,
et $P$ est un $G$-module de permutation.
Notons $T_P$ (respectivement $T_Q)$ le $K$-tore de module de caract\`ere $P$
(respectivement $Q$).

Le conoyau de la fl\`eche 
$\left(P^{\vee}\right)^{\,G}\to \left(X(T)^{\vee}\right)^{\,G}$ 
est le groupe $H^1\left(G,Q^{\vee}\right)$ qui est fini.
Comme $T_P$ est quasi-d\'eploy\'e,
d'apr\`es le lemme \ref{lm:degre:quasidep} 
le morphisme $\deg_{T_P,L}$ est de conoyau fini. 
Le diagramme commutatif
\begin{equation}
\xymatrix{
T_P(\ak)\ar[r]\ar[d]^{\deg_{T_P,L}}&T(\ak)\ar[d]^{\deg_{T,L}} \\
\left(P^{\vee}\right)^{\,G}\ar[r]& \left(X(T)^{\vee}\right)^{\,G}}
\end{equation}
permet de conclure pour $\deg_{T,L}$. 
Or on a la formule
\begin{equation}
d_L\deg_{T,L}=d_K\,[L:K]\deg_T,
\end{equation}
et $\left(X(T)^{\vee}\right)^{\,G}$ 
est un sous-module d'indice fini 
de $\left(X(T)^{G}\right)^{\vee}$.  
On en d\'eduit le r\'esultat pour $\deg_T$.
\end{demo}

\begin{nota}
Dans le cas fonctionnel, on note 
\nindex{$\CT$}$\CT$ le conoyau de $\deg_T$.
\end{nota}

On consid\`ere toujours $T$ un tore alg\'ebrique d\'efini sur le corps de fonctions $K$.  
Les r\'esultats de la fin de cette section donnent des renseignements sur les groupes $\DT$ et $\CT$.
Soit $\cG$ le groupe de Galois de $K^{\sep}/K$, il contient un sous-groupe distingu\'e $\cH$ 
tel que le quotient $\cG/\cH$ s'identifie \`a ${\mathfrak G}$ le groupe de Galois absolu de $k$. 
La repr\'esentation continue de $\cG$ dans $\Aut(X(T))$ induit une repr\'esentation continue
\begin{equation}
\varrho\,:\,{\mathfrak G}\longto\Aut\left(X(T)^{\,\cH}\right).
\end{equation}
Soit $\gt={\mathfrak G}/\Ker(\varrho)$ et $d_T=\card{\gt}$. 
Ainsi le corps fini \`a $\qde{K}^{\,d_T}$ \'el\'ements est le corps des constantes minimal 
d'une extension galoisienne $L/K$ d\'eployant $T$.

Pour tout $G$-module $M$, 
on note $\nindex{$N_G$}N_G$ la norme sur $M$,
i.e. le morphisme qui \`a $m\in M$ associe l'\'el\'ement de $M^{\,G}$
\begin{equation}
N_G(m)=\sum_{g\in G} g\,m.
\end{equation}

\begin{lemme}
\label{lm:dt}
L'image de $N_G T(\al)$ par $\deg_T$ est $d_T\,\left(X(T)^{\,G}\right)^{\vee}$.
\end{lemme}

\begin{demo}
Soit 
\begin{equation}
\phi\in T(\al)=\Hom(X(T),\G_m(\al)).
\end{equation}
Le morphisme $\deg_T(N_G \phi)$ envoie $m\in X(T)^{\,G}$ sur
\begin{align}
\deg_K \left(\prod_{g\in G} (g\,\phi)(m)\right)
&
=
\deg_K \left(\prod_{g\in G} g.\left(\phi(g^{-1}\,m)\right)\right)\\
&
=
\deg_K \left(\prod_{g\in G} g.\left(\phi(m)\right)\right)\\
&
=
\frac{d_T}{[L:K]}\deg_L \left(\prod_{g\in G} g.\left(\phi(m)\right)\right)\\
&
=
d_T\,\deg_L(\phi(m)).
\end{align}
La troisi\`eme \'egalit\'e provient de la formule \eqref{eq:formule_degre},
compte tenu du fait qu'on a $d_T=d_L/d_K$. 
Ainsi l'image de $N_G T(\al)$ par $\deg_T$ est incluse dans $d_T \,\left(X(T)^{\,G}\right)^{\vee}$.

Montrons \`a pr\'esent que la fl\`eche
\begin{equation}
\deg_T\,:\,N_G T(\al)\longto d_T\,\left(X(T)^{\,G}\right)^{\vee}
\end{equation}
est surjective. 

Soit $\psi$ un \'el\'ement de $d_T\,\left(X(T)^{\,G}\right)^{\vee}$.
On peut construire un morphisme 
\begin{equation}
\phi\in \Hom(X(T)^G,\G_m(\al))
\end{equation}
tel que pour $m\in X(T)^{\,G}$ on ait
\begin{equation}\label{eq:phi}
\forall m\in X(T)^{\,G},\quad\deg_L(\phi(m))=\frac{\psi(m)}{d_T}.
\end{equation}
Soit en effet $t\in \G_m(\al)$ un id\`ele tel que $\deg_L(t)=1$
Choisissons en outre une base de $X(T)^{\,G}$ $(m_1,\dots,m_r)$ .
Pour $i=1,\dots,r$, $\psi(m_i)$ s'\'ecrit $d_T\,n_i$ avec $n_i\in \Z$.
On pose alors $\phi(m_i)=t^{n_i}$.

Comme $X(T)^G$ est en tant que $\Z$-module
un facteur direct de $X(T)$, 
un tel morphisme $\phi$ s'\'etend en un morphisme
\begin{equation}
\phi\in \Hom(X(T),\G_m(\al)).
\end{equation} 

On v\'erifie que pour un tel $\phi$, on a, compte tenu
de \eqref{eq:phi}, $\deg_T(N_G\,\phi)=\psi$.

\end{demo}

\begin{cor}
$\DT$ contient $d_T\,\left(X(T)^{\,G}\right)^{\vee}$.
\end{cor}

\begin{cor}
\label{cor:ct0}
Si $k$ est alg\'ebriquement clos dans $L$, on a $\CT=0$.
\end{cor}

\begin{demo}
En effet dans ce cas on a $d_T=1$.
\end{demo}

\subsection{Groupe de classes}

Rappelons que pour tout corps global $K$, 
$C_K$ d\'esigne le groupe des classes d'id\`eles de $K$.

On consid\`ere un tore alg\'ebrique d\'efini sur corps global $K$,
d\'eploy\'e par une extension finie galoisienne $L$ de groupe $G$.

Le groupe
\begin{equation}
T(C_L)\eqdef \Hom(X(T),C_L)
\end{equation}
est isomorphe \`a $T(\al)/T(L)$.
Il est par ailleurs muni naturellement d'une action de $G$.
Nous posons
\begin{equation}
\nindex{$T(C_K)$}T(C_K)
=
\Hom_G(X(T),C_L)=T(C_L)^{\,G}.
\end{equation}
Cette d\'efinition est ind\'ependante du choix de l'extension d\'eployant $T$.

Notons que la suite exacte longue de cohomologie associ\'ee \`a la suite exacte de $G$-modules
\begin{equation}
0\longto 
X(T)^{\vee}\otimes \G_m(L)
\longto 
X(T)^{\vee}\otimes\G_m(\al)
\longto X(T)^{\vee}\otimes C_L 
\longto 0
\end{equation}
fournit la suite exacte
\begin{equation}
0\longto T(K) \longto T(\ak) \longto T(C_K)\longto H^1(G,T(L))\longto H^1(G,T(\al))
\end{equation}
et donc la suite exacte
\begin{equation}
0\longto T(\ak)/T(K)\longto T(C_K)\longto \cha(T) \longto 0,
\end{equation}
o\`u $\cha(T)$ est le \termin{groupe de Tate-Shafarevich} de $T$, d\'efini par   
\begin{equation}\label{eq:def:cha}
\nindex{$\cha(T)$}\cha(T)=\Ker\left(H^1(G,T(L))\to H^1(G,T(\al))\right),
\end{equation}
cette d\'efinition ne d\'ependant pas du choix de l'extension d\'eployant $T$.

Ainsi le groupe $T(\ak)/T(K)$ s'injecte dans $T(C_K)$ mais ne lui est  en g\'en\'eral pas \'egal
(pour un exemple avec $\cha(T)\neq 0$ cf. par exemple \cite[p.224, \S G.]{CTS:Requiv}).

Un cas important d'\'egalit\'e se produit 
quand $T$ est la restriction \`a la Weil de $K_0$ \`a $K$ de $\G_m$, 
pour $K_0$ extension s\'eparable de $K$. 
Dans ce cas, $T(C_K)$ s'identifie \`a $C_{K_0}$. 
En effet, on a alors, par le th\'eor\`eme de Hilbert 90, 
\begin{equation}
C_{K_0}=C_{L}^{\,\Gal(L/K_0)}.
\end{equation}
Bien entendu l'isomorphisme 
\begin{equation}
T(C_K)\longisom T(\ak)/T(K)
\end{equation}
est encore valable si $T$ est quasi-d\'eploy\'e.

\subsection{La dualit\'e de Nakayama}

\subsubsection{Les groupes de cohomologie \`a la Tate}
Nous effectuons quelques rappels sur les propri\'et\'es des groupes de cohomologie \`a la Tate, 
renvoyant \`a \cite[Chapitre VIII]{Ser:corps} pour plus de d\'etails. 
Soit $G$ un groupe fini. 
Pour tout $n\in \Z$, 
on peut d\'efinir sur la cat\'egorie des $G$-modules un foncteur $\widehat{H}^n(G,\,.\,)$, 
\`a valeurs dans la cat\'egorie des $\Z$-modules, 
v\'erifiant entre autres les propri\'et\'es suivantes : 
\begin{itemize}
\item  Pour $n\geq 1$, ce foncteur co\"\i ncide avec le foncteur classique $H^n(G,\,.\,)$, 
$n$-\`eme foncteur d\'eriv\'e droit du foncteur <<points fixes sous $G$>> : $M\mapsto M^{\,G}$.
\item  
Si $M$ est un $G$-module et 
\begin{equation}
N_G\,:\,m\mapsto \sum_{g\in G}g.m
\end{equation}
est la norme, 
on a $\widehat{H}^0(G,M)=M^{\,G}/N_G\,M$.
\item  Si 
\begin{equation}
0\longto M'\longto M\longto M''\longto 0
\end{equation}
est une suite exacte de $G$-modules, on a une suite exacte longue
\begin{gather}
\dots 
\longto 
\widehat{H}^{n-1}(G,M'')
\overset{\delta}{\longto}
\widehat{H}^{n}(G,M')
\longto
\widehat{H}^{n}(G,M)\\ 
\longto\widehat{H}^{n}(G,M'')
\overset{\delta}{\longto}
\widehat{H}^{n+1}(G,M')
\longto
\widehat{H}^{n+1}(G,M)\dots
\end{gather}
\item  Pour tout $n\in \Z$, 
$\widehat{H}^{n}(G,M)$ est tu\'e par la multiplication par $[G]$. 
\item
Si $M$ est de type fini, pour tout $n\in \Z$, 
$\widehat{H}^{n}(G,M)$ est fini.
\item  Si  $M$ et $M'$ sont des $G$-modules et $m$ et $n$ sont dans $\Z$, il existe
une application $\Z$-bilin\'eaire 
\begin{equation}\label{eq:cup_produit}
\widehat{H}^m(G,M)\otimes \widehat{H}^n(G,M')\longto\widehat{H}^{m+n}(G,M\otimes M').
\end{equation}
fonctorielle en $M$ et $M'$.
\item
Pour tout $G$-module $M$, le morphisme naturel de $G$-module 
\begin{equation}
M\otimes M^{\vee}\to \Z
\end{equation}
induit par fonctorialit\'e un morphisme trace
\begin{equation}
\widehat{H}^0(G,M\otimes M^{\vee})\to \widehat{H}^0(G,\Z)\isom \Z/\card{G}.
\end{equation}
En composant l'application bilin\'eaire \eqref{eq:cup_produit} pour $M'=M^{\vee}=\Hom(M,\Z)$ 
et $n=-n$ avec la trace, on obtient une dualit\'e parfaite
\begin{equation}
\widehat{H}^{-n}(G,M)\otimes \widehat{H}^{n}(G,M^{\vee})\longto \Z/\card{G}
\end{equation}
qui permet en particulier d'identifier $\widehat{H}^{-n}(G,M)$ 
au dual $\Hom\left(\widehat{H}^n(G,M^{\vee}),\Q/\Z\right)$ de 
$\widehat{H}^n(G,M^{\vee}$ (qui co\"\i ncide avec son dual topologique 
si $\widehat{H}^n\left(G,M^{\vee}\right)$ est fini), et inversement.
\end{itemize}

Dans toute la suite, 
pour tout $n\in \Z$, 
nous noterons $H^n(G,M)$ le groupe $\widehat{H}^n(G,M)$, 
sauf dans le cas o\`u $n$ vaut explicitement $0$, 
o\`u nous conserverons la notation $\widehat{H}^0(G,M)$ pour \'eviter toute confusion. 

\subsubsection{\'Enonc\'e de la dualit\'e de Nakayama}
Soit $K$ un corps global et 
$T$ un tore alg\'ebrique sur $K$, 
d\'eploy\'e par
une extension galoisienne $L$ de groupe $G$.
La dualit\'e de Nakayama donne un moyen simple de calculer la cohomologie du $G$-module $T(C_L)$. 
Le r\'esultat suivant est une cons\'equence de \cite[Theorem 3]{Na}. 
La th\'eorie du corps de classes permet de d\'efinir un g\'en\'erateur canonique $\alpha$ du groupe $H^2(G,C_L)$, 
appel\'e classe fondamentale (cf. \cite[p.359]{HoNa:cohomology_cft}).
\begin{thm}[Nakayama]
\label{thm:dual:naka}
Le cup-produit par $\alpha$ induit pour tout $n\in \Z$ un isomorphisme
\begin{equation}
H^{n}\left(G,X(T)^{\vee}\right)\longisom H^{n+2}(G,T(C_L)).
\end{equation}
En particulier les groupes de cohomologie $H^{n}(G,T(C_L))$ sont finis pour~tout~$n$.
\end{thm}
Rappelons que pour tout $n$ le groupe fini 
$H^{n}\left(G,X(T)^{\vee}\right)$ 
s'identifie canoniquement au dual 
$\Hom(H^{-n}(G,X(T)),\Q/\Z)=\dualt{H^{-n}(G,X(T))}$
du groupe fini $H^{-n}(G,X(T))$.

\subsection{Cocompacit\'e}

Soit $K$ un corps global.
Nous notons \nindex{$C^1_K$}$C^1_K=\G_m(\ak)^1/\G_m(K)$. 
C'est un groupe compact par \cite[IV\S4, Theorem 6]{Wei:BNT}. 
Soit \`a pr\'esent $T$ un tore alg\'ebrique sur $K$ d\'eploy\'e par une extension
galoisienne finie $L$ de groupe $G$.
Posons
\nindex{$T\left(C_K\right)^1$}
\begin{equation}
T\left(C_K\right)^1=\Hom_G(X(T),C^1_L),
\end{equation}
ce qui ne d\'epend pas du choix de l'extension $L$ d\'eployant $K$.

\begin{lemme}\label{lm:h1gtcl_fini}
Le groupe $H^1\left(G,T(C_L)^1\right)$ est fini.
\end{lemme}

\begin{demo}
Dans le cas fonctionnel, on consid\`ere la suite exacte
\begin{equation}
0
\longto 
X(T)^{\vee}\otimes C^1_L 
\longto 
X(T)^{\vee}\otimes C_L
\overset{\text{Id}\otimes \deg_L}{\longto} 
X(T)^{\vee} 
\longto 
0.
\end{equation}
La suite exacte de cohomologie sous $G$ associ\'ee fournit la suite exacte
\begin{equation}
T(C_K)\longto \left(X(T)^{\vee}\right)^{\,G} \longto H^1\left(G,T(C^1_L)\right) \longto H^1(G,T(C_L)).
\end{equation}
Par dualit\'e de Nakayama, 
$H^1\left(G,T(C_L)\right)$ est fini 
et la fl\`eche 
\begin{equation}
T(C_K)\longto \left(X(T)^{\vee}\right)^{\,G}
\end{equation}
n'est autre que la factorisation de $\deg_{T,L}$ par $T(C_K)$, 
son conoyau est donc fini.

Dans le cas arithm\'etique, on utilise la suite exacte
\begin{equation}
0
\longto 
X(T)^{\vee}\otimes C^1_L 
\longto 
X(T)^{\vee}\otimes C_L
\overset{\text{Id}\otimes \deg_L}{\longto} 
X(T)^{\vee}_{\R} 
\longto 
0
\end{equation}
et un argument strictement similaire (ici $\deg_{T,L}$
est d'ailleurs surjectif d'apr\`es le corollaire \ref{cor:ari:surj:degT}).
\end{demo}

\begin{prop}\label{prop:compacite}
\begin{enumerate}
\item Le 
\label{item:1:prop:compacite}
quotient $T(\ak)^1/T(K)$ est compact. 
\item Dans 
\label{item:2:prop:compacite}
le cas arithm\'etique, $T(\ak)/\K(T).T(K).T(\ak)_{\placesde{K,\infty}}$ est fini.
\item Dans 
\label{item:3:prop:compacite}
le cas fonctionnel,  $T(\ak)^1/\K(T).T(K)$ est fini.
\end{enumerate}
\end{prop}

\begin{demo} 
Le point  \ref{item:1:prop:compacite} est un cas particulier de \cite[IV.1.3.]{Oes:invent}, 
o\`u il est montr\'e plus g\'en\'eralement pour un groupe r\'esoluble.
Nous donnons une preuve de ce r\'esultat pour un tore alg\'ebrique.

Si $T$ est quasi-d\'eploy\'e, le r\'esultat d\'ecoule
imm\'ediatemment de la compacit\'e de $\G_m(\al)^1/\G_m(L)$
pour tout corps global $L$.

Par ailleurs, si $T$ est quasi-d\'eploy\'e.
on a un isomorphisme
\begin{equation}
T\left(C_K\right)^1\longisom T(\ak)^1 /T(K)
\end{equation}
En effet, il suffit de le montrer pour $T=\Res_{K_0/K}\G_m$
o\`u $K_0$ est une extension de $K$ contenue dans $L$.
Or, par le th\'eor\`eme de Hilbert 90, 
on a 
\begin{equation}
\left(C^1_L\right)^{\,\Gal(L/K_0)}=C^1_{K_0}.
\end{equation}
Ainsi, si $T$ est quasi-d\'eploy\'e, $T\left(C_K\right)^1$ 
est un compact.

Pour $T$ quelconque, d'apr\`es le lemme \ref{lm:res:qd}, il existe une suite exacte de $G$-modules
\begin{equation}
0\to X(T) \to P \to Q\to 0.
\end{equation}
o\`u $P$ et $Q$ libres de rang fini en tant que $\Z$-modules,
et $P$ est un $G$-module de permutation.
 Notons $T_P$ (respectivement $T_Q)$ le $K$-tore de module de caract\`ere $P$
(respectivement $Q$).

On a le diagramme commutatif suivant
\begin{equation}
\xymatrix{
T_P(\ak)^1/T_P(K)\ar[r]\ar[d]^{\wr}&T(\ak)^1 /T(K) \ar@{^{(}->}[d]\\
T_P\left(C_K\right)^1\ar[r]&T\left(C_K\right)^1}
\end{equation}
d'apr\`es lequel il suffit, pour montrer la compacit\'e de 
$T(\ak)^1 /T(K)$,
 de montrer que la fl\`eche 
$T_P\left(C_K\right)^1\to T\left(C_K\right)^1$ est de conoyau fini.

Or la suite exacte
\begin{equation}
0\longto Q^{\vee}\longto P^{\vee} \longto X(T)^{\vee} \longto 0
\end{equation}
fournit la suite exacte
\begin{equation}
T_P\left(C_K\right)^1
\longto 
T\left(C_K\right)^1
\longto 
H^1\left(G,T_Q(C_L)^1\right).
\end{equation}
D'apr\`es le lemme \ref{lm:h1gtcl_fini}, $H^1\left(G,T_Q(C_L)^1\right)$ est fini,
d'o\`u le r\'esultat.

Montrons le point \ref{item:2:prop:compacite}, qui, dans le cas o\`u $T=\G_m$,
est le th\'eor\`eme de finitude du nombre de classes d'id\'eaux d'un corps de nombres.
D'apr\`es le lemme \ref{lm:varchitkvtov},
on a 
\begin{equation}
T(\ak)=T(\ak)^1\,.T(\ak)_{\placesde{K,\infty}}
\end{equation}
et il suffit donc de montrer que le quotient
\begin{equation}
T(\ak)^1/\K(T).T(K).T(\ak)^1_{\placesde{K,\infty}}
\end{equation}
est fini. Mais ce dernier groupe est discret,
et compact car $T(\ak)^1/T(K)$ est compact d'apr\`es le point \ref{item:1:prop:compacite}.

Montrons \`a pr\'esent le point \ref{item:3:prop:compacite}.
De fa\`c;on semblable au point \ref{item:2:prop:compacite}, il peut se d\'eduire du point \ref{item:1:prop:compacite} :
la compacit\'e de $T(\ak)^1/T(K)$ entra\^\i ne celle de $T(\ak)^1/\K(T).T(K)$, 
et donc la finitude de ce dernier groupe dans le cas fonctionnel car il est 
alors \'egalement discret. 

Cependant  le point \ref{item:3:prop:compacite} peut se retrouver directement, 
en remarquant qu'on a 
\begin{equation}
\left(T(\al)^1\right)^{\,G}=T(\ak)^1,
\end{equation}
d'o\`u une injection
\begin{equation}
T(\ak)^1 /\left(\K(T_L).T(L)\right)^{\,G} \hookrightarrow \left(T(\al)^1/\K(T_L).T(L)\right)^{\,G}.
\end{equation}
Or $T(\al)^1/\K(T_L).T(L)$ n'est autre que $\left(\Pic^0(\courbe_L)\right)^{\dim(T)}$. 
C'est donc un groupe fini.
On en d\'eduit que le groupe 
\begin{equation}
T(\ak)^1 /\left(\K(T_L).T(L)\right)^{\,G}
\end{equation}
est fini. 

Par ailleurs, $\K(T_L)\cap T(L)$ s'identifie \`a un produit de $\dim(T)$ copies du groupe
des fonctions r\'eguli\`eres inversibles sur la courbe projective $\courbe_L$, et est donc fini.
Il en en donc de m\^eme pour $H^1\left(G,\K(T_L)\cap T(L)\right)$.
Maintenant la suite exacte
\begin{equation}
0\longto \K(T).T(K) \longto \left(\K(T_L).T(L)\right)^{\,G} \longto H^1\left(G,\K(T_L)\cap T(L)\right)
\end{equation}
montre que $\K(T).T(K)$ est d'indice fini dans $\left(\K(T_L).T(L)\right)^{\,G}$, 
d'o\`u la finitude de $T(\ak)^1 /\K(T).T(K)$. 
\end{demo}

\subsection{R\'esultats locaux}\label{subsec:calcullocal}

\begin{lemme}[Ono]
\label{lm:ono}
Soit $v$ une place finie de $K$.
Soit $1\to T''\to T'\to T \to 1$ une suite exacte de tores alg\'ebriques sur $K$.
Soit $L$ une extension galoisienne finie d\'eployant $T$, $T'$
et $T"$.
Si $v$ est non ramifi\'ee dans $L/K$, 
le morphisme
\begin{equation}
T'(\Ov)\longto T(\Ov)
\end{equation}
est surjectif.
\end{lemme}
\begin{demo} 
La preuve figure dans \cite[Lemma 4.2.1]{Ono_aritag}, 
nous la rappelons. 
Soit $\V$ une place de $L$ au-dessus de $v$, 
et $G_v$ le groupe de d\'ecomposition de $\V$.
On a  
\begin{gather}
T(\Ov)\isom \Hom_{G_v}(X(T),\OV^{\times}),\\
T'(\Ov)\isom \Hom_{G_v}(X(T'),\OV^{\times})
\end{gather}
et une suite exacte
\begin{equation}
T'(\Ov)
\longto 
T(\Ov)
\longto H^1(G_v,X\left(T''\right)^{\vee}\otimes \OV^{\times}).
\end{equation}

Or, 
comme  $v$ est non ramifi\'ee, 
$\OV^{\times}$ est cohomologiquement trivial, 
et comme $X\left(T''\right)^{\vee}$ est sans torsion, 
d'apr\`es \cite[IX, $\S\,5$, Corollaire]{Ser:corps}, 
$X\left(T''\right)^{\vee}\otimes \OV^{\times}$ est encore
cohomologiquement trivial, 
d'o\`u le r\'esultat. 
\end{demo}

Soit $T$ un tore alg\'ebrique d\'efini sur un corps global $K$,
d\'eploy\'e par une extension finie galoisienne $L$ de groupe de Galois 
$G$.
On suppose donn\'ee une suite exacte de $G$-modules
\begin{equation}\label{eq:exseq:calcullocal}
0\longto X(T) \overset{\gamma}{\longto P} \longto Q \longto 0
\end{equation}
o\`u $P$ et $Q$ sont des $\Z$-modules libres de rang fini, 
et $P$ est un $G$-module de permutation.
Le but de ce qui suit est de pr\'eciser le comportement du degr\'e local
vis-\`a-vis de cette suite exacte. Comme cons\'equence, on obtient la finitude
du conoyau du degr\'e local.

Choisissons une base $G$-stable $(n_i)_{i\in I}$ du $G$-module de permutation $P$,
dont on note $(n_i^{\vee})_{i\in I}$ la base duale.
Pour $\wti\in I/G$,  choisissons en outre
un \'el\'ement $n_{\wti}$ de l'orbite $\wti$. 
Soit $G_{\wti}$ son stabilisateur.
Ces choix permettent d'identifier  
$P$ (et $P^{\vee}$) au $G$-module de permutation
$\oplusu{{\wti}\in I/G}\Z[G/G_{\wti}]$. 
Si on d\'esigne par $K_{\wti}$ le corps $L^{G_{\wti}}$,
cette identification iduit un isomorphisme
\begin{equation}\label{eq:iso:tp}
T_P\longisom \produ{{\wti}\in I/G}\Res_{K_{\wti},K}\G_m.
\end{equation}

Soit $v$ une place finie de $K$, $\V$ une place de $L$ divisant $v$
et $G_v$ son groupe de d\'ecomposition.
Rappelons qu'on note encore $\V$ la valuation normalis\'ee repr\'esentant $\V$.

L'isomorphisme \eqref{eq:iso:tp} induit au niveau des $K_v$-points un isomorphisme
\begin{equation}\label{eq:iso:tpkv}
T_P(K_v)
\longisom 
\prod_{{\wti}\in I/G}\,\G_m(K_{\wti}\otimes K_v)
=
\prod_{{\wti}\in I/G} 
\,
\prod_{
\substack{w\in  \placesde{K_{\wti}}\\ w|v}
}
\,
\G_m(K_w).
\end{equation}

Pour $\wti\in I/G$, consid\'erons l'application
\begin{equation}
\begin{array}{rcl}
G/G_{\wti}&\longto&\{w\in \placesde{K_{\wti}},\,v|w\}\\
g&\longmapsto&g^{-1}.\V|_{K_{\wti}} \end{array}
\end{equation}
Ce n'est autre que le passage au quotient par l'action de $G_v$, 
d'o\`u une correspondance entre $\wti/G_v$ et les places de $K_{\wti}$ au-dessus de $v$. 

Pour $j\in \wti/G_v$, nous notons $w_j$ la place de $K_{\wti}$ au-dessus de $v$
donn\'ee par cette correspondance.
On choisit en outre $n_j$ un \'el\'ement quelconque de $j$. 

On pose  alors
\begin{equation}
\tau_j=\sum_{n \in G_v.n_j}\,\gamma^{\vee}(n^{\vee}),
\end{equation}
de sorte que $\tau_j$ est un \'el\'ement de $(X(T)^{\vee})^{\,G_v}$.

Pour $\wti\in I/G$ et $j\in \wti/G_v$, on identifie, via l'isomorphisme
\eqref{eq:iso:tpkv},
$\G_m(K_{w_j})$ \`a un sous-groupe de $T_P(K_v)$.
Soit $\pi_{w_j}$ une uniformisante de $K_{w_j}$
et $e_j$ l'indice de ramification de $w_j$ dans $L$.

\begin{lemme}\label{lm:jvgammavpiwj}
On a la relation
\begin{equation}\label{eq:jvgammavpiwj}
\deg_{T,L,\V}\left[\gamma_v(\pi_{w_j})\right]=e_j\,\tau_j.
\end{equation}
\end{lemme}

\begin{demo}
Le morphisme $\gamma$ et la d\'ecomposition
\begin{equation}
P\isom \bigoplus_{\wti\in I/G} \bigoplus_{j\in \wti/G_v} \oplusu{n\in G_v.n_j} \Z\,n
\end{equation}
induisent des  morphismes de $G_v$-module
\begin{equation}
X(T)\overset{\gamma_j}{\longto} \bigoplus_{n\in G_v.n_j} \Z\,n
\end{equation}
qui \`a leur tour induisent par dualit\'e 
des morphismes de $K_v$-tores 
\begin{equation}
\Res_{K_{w_j}/K_v}\G_m\longto T.
\end{equation}
Pour tout $j$, ce dernier morphisme induit au niveau des $K_v$-points
un morphisme 
\begin{equation}
\Res_{K_{w_j}/K_v}\G_m(K_v)\isom \G_m(K_{w_j})\longto T(K_v)
\end{equation}
qui n'est autre que la restriction de $\gamma_v$ \`a $\G_m(K_{w_j})$.

Il peut donc se d\'ecrire comme le morphisme
\begin{equation}\label{eq:descr}
\Hom_{G_v}
\left(
\oplusu{n\in G_v.n_j} \Z\,n, L_{\V}^{\times}
\right) 
\overset{\circ\, \gamma_j}{\longto}
\Hom_{G_v}(X(T), L_{\V}^{\times}).
\end{equation}
On calcule alors $\gamma_v(\pi_{w_j})$ en utilisant la description
\eqref{eq:descr}.

Via l'identification
\begin{equation}
\G_m(K_{w_j})\isom \Hom_{G_v}(\oplusu{n\in G_v.n_j} \Z\,n, L_{\V}) 
\end{equation}
l'\'el\'ement $\pi_{w_j}$ de $\G_m(K_{w_j})$ correspond
au morphisme  
qui envoie $n_j$ sur $\pi_{w_j}$

Son image par $\gamma_v$ est donc l'\'el\'ement de $\Hom_{G_v}(X(T),L_{\V}^{\times})$
donn\'e par 
\begin{equation}
m\mapsto \prod_{n\in G_v.n_j}(g_n\,\pi_{w_j})^{\acc{m}{\gamma^{\vee}(n^{\vee})}},
\end{equation}
o\`u pour $n\in G_v.n_j$ on note $g_n$ un \'el\'ement de $G_v$ tel que $n=g_n.n_j$.

On en d\'eduit que $\deg_{T,L,\V}(\gamma_v(\pi_{w_j}))$
est l'\'el\'ement de $(X(T)^{\vee})^{G_v}$ qui \`a $m\in X(T)$
associe
\begin{align}
\sum_{n\in G_v.n_j} \acc{m}{\gamma^{\vee}(n^{\vee})}\,\V(g_n\,\pi_{w_j})
&=\sum_{n\in G_v.n_j} \acc{m}{\gamma^{\vee}(n^{\vee})}\,\V(\pi_{w_j})\\
&=\acc{m}{e_j\,\sum_{n\in G_v.n_j}\gamma^{\vee}(n^{\vee})}\\
&=\acc{m}{e_j\,\tau_j}
\end{align}
d'o\`u le r\'esultat.
\end{demo}

Comme cons\'equence du lemme \ref{lm:jvgammavpiwj}, 
on obtient le r\'esultat suivant, qui  est signal\'e \`a la page 449 de  \cite{Dr}.
\begin{prop}
\label{prop:draxl}
Le morphisme 
\begin{equation}\label{eq:degTLV:2}
\deg_{T,L,\V}\,:\,T(K_v)\longto \left(X(T)^{\vee}\right)^{\,G_v}
\end{equation}
est de conoyau fini.
\end{prop}
\begin{rem}
Si en outre $v$ n'est pas ramifi\'ee dans $L$, on a d\'ej\`a vu 
que le conoyau est en fait trivial (cf. le lemme \ref{lemme:draxl}).
\end{rem}
\begin{demo}
Comme $H^1(G_v,Q^{\vee})$ est fini, la fl\`eche 
$\left(P^{\vee}\right)^{\,G_v}\to \left(X(T)^{\vee}\right)^{\,G_v}$
(induite par la restriction de $\gamma^{\vee}$ \`a $\left(P^{\vee}\right)^{\,G_v}$)
est de conoyau fini. 
Comme $\left(P^{\vee}\right)^{\,G_v}$ est engendr\'e par les \'el\'ements
\begin{equation}
\sum_{n\in G_v.n_j} n^{\vee}
\end{equation}
pour $\wti$ d\'ecrivant $I/G$ et $j$ d\'ecrivant $\wti/G_v$, 
on voit que le sous-$\Z$-module 
de $\left(X(T)^{\vee}\right)^{\,G_v}$ engendr\'e par les 
$\tau_j$ pour $j\in \wti/G_v$ et $\wti\in I/G$ est d'indice fini dans  $\left(X(T)^{\vee}\right)^{\,G_v}$.
Il en est donc de m\^eme du sous-$\Z$-module engendr\'e par les $e_j\,\tau_j$.
Or, le lemme \ref{lm:jvgammavpiwj} montre que $\deg_{T,L,\V}(\gamma_v(T_P(K_v)))$ 
contient ce dernier sous-module.

\end{demo}

\begin{rem}
Le raisonnement pr\'ec\'edent pourrait permettre en fait de montrer  
en une seule fois la proposition \ref{prop:draxl}
et le lemme \ref{lemme:draxl}, 
en choisissant la suite exacte
\eqref{eq:exseq:calcullocal}
de sorte que $Q$ soit un $G$-module flasque (cf. la sous-section \ref{subsubsec:introflasque}).
En effet, dans ce cas,
la fl\`eche $\left(P^{\vee}\right)^{\,G_v}\to \left(X(T)^{\vee}\right)^{\,G_v}$
est surjective, et non plus seulement de conoyau fini. 
Or, si $v$ n'est pas ramifi\'ee dans $L$, tous les $e_j$ sont \'egaux \`a $1$, et le raisonnement
pr\'ec\'edent montre que $\deg_{T,L,\V}$ est surjective. 
\end{rem}

\subsection{R\'esolution flasque d'un tore alg\'ebrique et applications}\label{subsec:flasque}

\subsubsection{Rappels et notations}\label{subsubsec:introflasque}

La notion de r\'esolution flasque d'un tore alg\'ebrique a \'et\'e introduite par Colliot-Th\'el\`ene et Sansuc 
dans \cite{CTS:Requiv} en vue de l'\'etude de la R-\'equivalence sur les tores.

Soit $T$ un tore alg\'ebrique d\'efini sur un corps $K$, 
d\'eploy\'e par une extension finie galoisienne $L$ de groupe de Galois $G$. 
Rappelons qu'un $G$-module $M$ est dit \dindex{$G$-module flasque}\termin{flasque} 
si pour tout sous-groupe $H$ de $G$ on a 
$H^{-1}(H,M)=0$. 
Par \cite[lemme 3, page 181]{CTS:Requiv} il existe une suite exacte de $G$-modules 
\begin{equation}\label{eq:res:fl}
0\longto X(T) \overset{\gamma}{\longto} P \longto Q \longto 0 
\end{equation}
o\`u $P$ et $Q$ sont libres de rang fini comme $\Z$-modules, $P$ est de permutation et $Q$ est flasque
(ceci g\'en\'eralise le lemme \ref{lm:res:qd}).
Un telle suite exacte est appel\'ee une \dindex{r\'esolution flasque}\termin{r\'esolution flasque} de $X(T)$.

On fixe pour toute la section  \ref{subsec:flasque} un tore alg\'ebrique $T$ d\'efini sur un corps 
global $K$, 
d\'eploy\'e par une extension finie galoisienne $L$ de groupe de Galois $G$. 
et une r\'esolution flasque \eqref{eq:res:fl} de $X(T)$.
Soit $T_P$ le tore alg\'ebrique associ\'e au $G$-module $P$. 
C'est un tore quasi-d\'eploy\'e. 
Soit $T_Q$ le tore alg\'ebrique associ\'e au $G$-module $Q$.

\subsubsection{Un r\'esultat local}

On conserve les hypoth\`eses et notations introduites dans la sous-section \ref{subsubsec:introflasque}.

On consid\`ere dans cette partie une place finie $v$ de $K$.
On note $\V$ une place de $L$ au-dessus de $v$
et $G_v$ le groupe de d\'ecomposition de $\V$.

\begin{prop}\label{prop:tpkvtkv}
Si $v$ est non ramifi\'ee dans $L/K$, 
le morphisme
\begin{equation}
\gamma_v\,:\,T_P(K_v)\longto T(K_v)
\end{equation}
est surjectif.
\end{prop}

\begin{demo}
La suite exacte de cohomologie associ\'ee \`a la suite exacte
de $G_v$-modules
\begin{equation}
1\to T_Q(L_{\V})\to T_P(L_{\V}) \to T(L_{\V}) \to 1
\end{equation}
et le fait que $T_P$ soit d\'eploy\'e, donc que $H^1(G_v,T_P(L_{\V}))$
soit nul d'apr\`es Hilbert 90, montre
que le conoyau de $\gamma_v$
est 
\begin{equation}
H^1(G_v,T_Q(L_{\V}))
=
H^1(G_v,Q^{\vee}\otimes L_{\V}^{\times}).
\end{equation}
La nullit\'e de ce groupe 
d\'ecoule de l'application de la dualit\'e de Nakayama 
aux corps locaux : 
celle-ci fournit un isomorphisme
\begin{equation}
H^1(G_v,Q^{\vee}\otimes L_{\V}^{\times})\longisom
H^{-1}(G_v,Q^{\vee}).
\end{equation}
Par ailleurs, 
$v$ \'etant non ramifi\'ee, 
$G_v$ est cyclique,
d'o\`u un isomorphisme
\begin{equation}
H^{-1}(G_v,Q^{\vee})
\longisom
H^{1}(G_v,Q^{\vee}).
\end{equation}
Or le groupe $H^{1}(G_v,Q^{\vee})$ est nul car $Q$ est flasque.

\end{demo}
\begin{cor}\label{cor:LKnonram:tpaktaksurj}
On se place dans le cas fonctionnel et 
on suppose que le tore $T$ est d\'eploy\'e par une extension $L/K$ non ramifi\'ee.
Alors le morphisme $T_P(\ak)\to T(\ak)$ est surjectif.
\end{cor}
\begin{demo}
Ceci d\'ecoule de la proposition \ref{prop:tpkvtkv} 
et du lemme \ref{lm:ono}.
\end{demo}
\subsubsection{Approximation faible}

On conserve les hypoth\`eses et notations introduites dans la sous-section \ref{subsubsec:introflasque}.

Soit $\adh{T(K)}$ l'adh\'erence de $T(K)$ dans 
$\underset{v\in\placesde{K}}{\prod}T(K_v)$ 
muni de la topologie produit 
et soit
\begin{equation}\label{eq:def:at}
\nindex{$A(T)$}A(T)=\left(\prod_{v\in\placesde{K}}T(K_v)\right)/\adh{T(K)},
\end{equation}
c'est le groupe d'obstruction \`a l'approximation faible.
Il est nul si $T=\G_m$ et plus g\'en\'eralement si $T$
est quasi-d\'eploy\'e, d'apr\`es \cite[p. 334]{Has}.

Soit $p$ le morphisme $T(\ak)/T(K)\to A(T)$ induit par la projection $T(\ak)\to A(T)$. 

Pour tout ensemble fini $S$ de places de $K$ 
contenant les places archim\'ediennes, 
nous notons 
\begin{equation}
\nindex{$T(\ak)^S$}T(\ak)^S=T(\ak)\,\bigcap \,\prod_{v\notin S} T(K_v).
\end{equation}
et 
\begin{equation}
\adh{T(K)}^{\,\,S}=\adh{T(K)}\,\bigcap \,\produ{v\in S} T(K_v)
\end{equation}
(en d'autres termes, $\adh{T(K)}^{\,\,S}$ est l'adh\'erence de l'image de $T(K)$
dans $\produ{v\in S} T(K_v)$).

\begin{prop}
\label{lm:scindage}\label{lm:appfaible}
\begin{enumerate}
\item\label{item:lm:scindage:1}
On suppose que $S$ contient les places ramifi\'ees dans $L/K$. 
On a alors un scindage
\begin{equation}\label{eq:scindage1}
\adh{T(K)}
=
\left(\prod_{v\notin S}T(K_v)\right)\times \adh{T(K)}^{\,\,S}.
\end{equation}
En particulier on a un scindage
\begin{equation}\label{eq:scindage2}
\adh{T(K)}\cap T(\ak)=T(\ak)^S\times \adh{T(K)}^{\,\,S}
\end{equation}
et les \'egalit\'es 
\begin{equation}\label{eq:scindage3}
A(T)=\left(\prod_{v\in S}T(K_v)\right)/\adh{T(K)}^{\,\,S}=T(\ak)/\left(\adh{T(K)}\cap T(\ak)\right).
\end{equation}
\item\label{item:lm:scindage:2}
La suite 
\begin{equation}\label{eq:exsqtpaktakat}
T_P(\ak)/T_P(K)
\longto
T(\ak)/T(K)
\overset{p}{\longto}
A(T)
\longto 0
\end{equation}
est exacte.
\item\label{item:lm:scindage:3}
Il existe une suite exacte
\begin{equation}
\label{eq:exseqcts}
0\longto A(T) \longto H^1(G,Q)^{\ast}\longto \cha(T) \longto 0.
\end{equation}
En particulier, $A(T)$ est fini et on a  
\begin{equation}\label{eq:A=H1/cha}
\card{A(T)}
=
\frac
{\card{H^1(G,Q)}}
{\card{\cha(T)}}.
\end{equation}
\end{enumerate}
\end{prop}
\begin{rems}
\begin{enumerate}
\item
Le point \ref{item:lm:scindage:1} est d\^u \`a Voskresenskii.
\item
Dans le cas des corps de nombres, le point \ref{item:lm:scindage:2} 
est indiqu\'e dans \cite[Theorem 3.1.1]{BaTs:anis}. 
\item
Le point \ref{item:lm:scindage:3}  a \'et\'e initialement d\'emontr\'e par Voskresenskii
 en caract\'eristique z\'ero et sous une forme l\'eg\`erement
diff\'erente (\cite[Thm. 6]{Vos:bir}). 
Il est d\^u sous la forme donn\'ee ici \`a Colliot-Th\'el\`ene
et Sansuc (\cite[Proposition 19 (iB)]{CTS:Requiv}).
La d\'emonstration qui en est faite ci-dessous 
reprend, \`a des d\'etails de pr\'esentation pr\`es, celle des auteurs de \cite{CTS:Requiv}. 
Elle figure dans ce texte d'une part par souci de compl\'etude, d'autre
part parce que les arguments utilis\'es permettent aussi d'obtenir le point 
\ref{item:lm:scindage:2}.
\item
Si on ne suppose plus $Q$ flasque dans la suite
exacte \ref{eq:res:fl},
Draxl montre dans \cite{Dr} que  le conoyau du morphisme
\begin{equation}
T_P(\ak)/T_P(K)
\longto
T(\ak)/T(K)
\end{equation}
est de cardinal
\begin{equation}
\frac{\card{H^1(G,Q)}}{\cha(T)}.
\end{equation}
Ce r\'esultat
permet de retrouver le point 
\ref{item:lm:scindage:2}
\`a partir du point \ref{item:lm:scindage:3}.

Il permet aussi de montrer que les points \ref{item:lm:scindage:2}
et \ref{item:lm:scindage:3}
ne sont pas n\'ecessairement v\'erifi\'es si on ne suppose plus $Q$ flasque 
dans la suite exacte \ref{eq:res:fl} (contrairement \`a ce qui est affirm\'e
en haut de la page 3231 de \cite{BaTs:hzf}).

Consid\'erons par exemple le cas o\`u $G$ est cyclique, et o\`u la suite exacte
\begin{equation}
0\longto X(T) \longto P \longto Q\longto 0
\end{equation}
est la suite duale de la suite exacte
\begin{equation}
0\longto I_G \longto \Z[G] \overset{\eps}{\longto} \Z\longto 0
\end{equation}
o\`u $\eps$ est l'augmentation $\sum a_g\,g\mapsto \sum a_g$.

Comme $T$ est d\'eploy\'e, on a $A(T)=1$, et d'apr\`es le r\'esultat 
de Draxl le conoyau de la fl\`eche
$
T_P(\ak)/T_P(K)
\to
T(\ak)/T(K)$
est de cardinal $\card{G}$.
\item
La relation \eqref{eq:A=H1/cha} nous servira pour 
le calcul du terme principal de la fonction z\^eta des hauteurs dans la section \ref{app_fonction_zeta}.
\end{enumerate}
\end{rems}
\begin{demo}

Comme $T_P$ est quasi-d\'eploy\'e, on a $A(T_P)=0$
soit en d'autres termes
$\produ{v}\,T_P(K_v)=\adh{T_P(K)}$.
Ceci ajout\'e au fait que 
$\produ{v} \gamma_v$ 
est continu pour la topologie produit 
montre les inclusions
\begin{equation}\label{eq:inclu:gamma}
\left(\produ{v} \gamma_v\right)\left(\prod_{v} T_P(K_v)\right)
\subset 
\adh{\left(\produ{v} \gamma_v\right) (T_P(K))}
\subset 
\adh{T(K)}.
\end{equation}

Or, d'apr\`es la proposition \ref{prop:tpkvtkv},
l'image de
$\produ{v\notin S}T_P(K_v)$ 
par 
$\produ{v} \gamma_v$
est $\produ{v\notin S}T(K_v)$.
On en d\'eduit que $\adh{T(K)}$ contient $\produ{v\notin S}T(K_v)$,
ce qui, compte tenu du fait que $\adh{T(K)}$ est un sous-groupe 
de $\produ{v}\, T(K_v)$,
montre l'\'egalit\'e \eqref{eq:scindage1}.
Les \'egalit\'es \eqref{eq:scindage2} et \eqref{eq:scindage3}
en d\'ecoule aussit\^ot.

Pour la suite de la preuve, on consid\`ere le diagramme commutatif et exact
\begin{equation}\label{eq:grodiag}
\xymatrix{
&
0
&
&
\\
0\ar[r]
&
\cha(T)\ar[u]
&
&
\\
T_P(C_K)\ar[r]\ar[u]
&
T(C_K) 
\ar[r]\ar[u] 
&
H^1(G,T_Q(C_L)) 
\ar[r]
&
0
\\
T_P(\ak)\ar[r]^{\gamma_{\,\ak}}\ar[u]
&
T(\ak) 
\ar[r]^<<<<<{\partial} \ar[u]
&
H^1(G,T_Q(\al)) 
\ar[r]\ar[u]
&
0
\\
T_P(K)
\ar[r]\ar[u]
&
T(K) 
\ar[r] 
\ar[u] 
&
H^1(G,T_Q(L))
\ar[u]_{\,\nu} 
\ar[r]
&
0
\\
0\ar[u]
&
0\ar[u]
&
}
\end{equation}
La deuxi\`eme (respectivement troisi\`eme, respectivement quatri\`eme) ligne
s'obtient \`a partir de la suite exacte longue de cohomologie tir\'ee de la suite 
obtenue en tensorisant le dual de la suite exacte \eqref{eq:res:fl} 
avec $C_L$ (respectivement $\G_m(\al)$, respectivement $\G_m(L)$), 
en prenant les $G$-invariants, et en remarquant que :
\begin{itemize}
\item
par dualit\'e de Nakayama,
\begin{equation}
H^1(G,P^{\vee}\otimes C_L)\longisom H^1(G,P)^{\,\ast}
\end{equation}
et $P$ \'etant de permutation, $H^1(G,P)=0$ ;
\item
$H^1(G,T_P(\al))$ et $H^1(G,T_P(L))$ sont nuls ;
ceci d\'ecoule du fait
que $T_P$ est quasi-d\'eploy\'e et de Hilbert 90 
(si $k'/k$ est une extension finie galoisienne de corps,
$H^1(\Gal(k'/k),\G_m(k))=0$~; si en outre $k$ est global,
$H^1(\Gal(k'/k),\G_m(\ade{k'}))=0$).
\end{itemize}
La premi\`ere (respectivement deuxi\`eme, respectivement troisi\`eme) colonne
s'obtient \`a partir de la suite exacte longue de cohomologie tir\'ee de la suite 
obtenue en tensorisant la suite exacte 
\begin{equation}
1\longto \G_m(L)\longto \G_m(\al)\longto C_L
\end{equation}
par $P^{\vee}$ (respectivement $X(T)^{\vee}$), respectivement $Q^{\vee}$),
de la d\'efinition de $\cha(T)$ et du fait que $T_P$ \'etant quasi-d\'eploy\'e, on a $\cha(T_P)=0$.

Nous allons montrer que le sous-groupe $\partial^{-1}(\nu(H^1(G,T_Q(L))$
est \'egal \`a $\adh{T(K)}\cap T(\ak)$, soit, en d'autre termes,
que le groupe $A(T)$ s'identifie au conoyau de $\nu$.

De \eqref{eq:inclu:gamma}
on d\'eduit aussit\^ot l'inclusion
\begin{equation}\label{eq:inclu:gamma:ak}
\Ker(\partial)=\gamma_{\,\ak}\left(T_P(\ak)\right)
\subset 
\adh{T(K)}\cap T(\ak)
\end{equation}
Du diagramme \eqref{eq:grodiag} on d\'eduit l'\'egalit\'e
\begin{equation}
\partial^{-1}(\nu(H^1(G,T_Q(L))=T(K).\Ker(\partial)
\end{equation}
soit
\begin{equation}\label{eq:pnuinadhtk}
\partial^{-1}(\nu(H^1(G,T_Q(L))\subset \adh{T(K)}\cap T(\ak)
\end{equation}

Par ailleurs on a vu que l'image
de $\produ{v}\,T_P(K_v)$
par $\produ{v}\gamma_v$
\'etait \'egale \`a
\begin{equation}
\produ{v\notin S}T(K_v)\times \produ{v\in S}\gamma_v(T_P(K_v)).
\end{equation}
Or, pour tout $v$, $\gamma_v$ est ouverte
d'apr\`es \cite[Proposition 2.7(a)]{Sal:tammes}.

Ainsi l'image
de $\produ{v}T_P(K_v)$
par $\produ{v}\gamma_v$
est ouverte
dans $\produ{v}T(K_v)$.

Mais, d'apr\`es le lemme \ref{lm:ono}, on a
\begin{equation}
\gamma_{\,\ak} (T_P(\ak))=T(\ak) \bigcap \left(\produ{v}\gamma_v\right)\left(\produ{v}T_P(K_v)\right).
\end{equation}

En particulier, l'image de $T_P(\ak)$ par $\gamma_{\,\ak}$ 
est ouverte dans $T(\ak)$ pour la topologie produit.
Or cette image n'est autre que $\Ker(\partial)$.
Ainsi le sous-groupe $\partial^{-1}(\nu(H^1(G,T_Q(L))$
est ouvert dans $T(\ak)$ (toujours pour la topologie produit),
donc ferm\'e. Comme il contient $T(K)$, on a 
\begin{equation}\label{eq:adhtkinpnu}
\adh{T(K)}\cap T(\ak)\subset \partial^{-1}(\nu(H^1(G,T_Q(L))
\end{equation}
De \eqref{eq:pnuinadhtk} et \eqref{eq:adhtkinpnu} on d\'eduit l'\'egalit\'e
cherch\'ee
\begin{equation}\label{eq:adhtk=pnu}
\partial^{-1}(\nu(H^1(G,T_Q(L))=\adh{T(K)}\cap T(\ak)
\end{equation}

Le fait que le groupe $A(T)$ s'identifie au conoyau de $\nu$
a deux cons\'equences~: d'une part 
le lemme du serpent et le diagramme commutatif exact
\begin{equation}
\xymatrix{
0\ar[d]
&
0\ar[d]
&
&
\\
T_P(K)
\ar[r]\ar[d]
&
T(K) 
\ar[r] 
\ar[d] 
&
H^1(G,T_Q(L))
\ar[d]
\ar[r]
&
0
\\
T_P(\ak)\ar[r]\ar[d]
&
T(\ak) 
\ar[r]
\ar[d]
&
H^1(G,T_Q(\al)) 
\ar[r]
\ar[d]
&
0
\\
T_P(\ak)/T_P(K)\ar[d]
&
T(\ak)/T(K)
\ar[d]
&
A(T)
\ar[d]
&
\\
0
&
0
&
0
&
}
\end{equation}
montrent le point 3.

D'autre part, une chasse au diagramme standard 
dans 
\eqref{eq:grodiag}
montre l'existence d'une suite exacte
\begin{equation}
0\longto A(T) \longto H^1(G,T_Q(C_L))\longto \cha(T) \longto 0.
\end{equation}

Par dualit\'e de Nakayama, on a un isomorphisme
\begin{equation}
H^1(G,T_Q(C_L))\longisom H^1(G,Q)^{\,\ast}
\end{equation}
d'o\`u la suite exacte \eqref{eq:exseqcts}. 
\end{demo}
\begin{cor}\label{cor:nonramAT=0}
On se place dans le cas fonctionnel et 
on suppose que le tore $T$ est d\'eploy\'e par une extension $L/K$ non ramifi\'ee.
Alors $A(T)$ est trivial.
\end{cor}
\begin{demo}
Ceci d\'ecoule du corollaire \ref{cor:LKnonram:tpaktaksurj}
et du point \ref{item:lm:scindage:2} de la proposition \ref{lm:scindage}. 
\end{demo}

\subsubsection{Un invariant des tores alg\'ebriques d\'efinis sur les corps de fonctions}
\label{subsubsec:KT}

On introduit dans cette sous-section un invariant des tores alg\'ebriques d\'efinis un corps de fonctions.
Cet invariant intervient naturellement dans le calcul du terme principal
de la fonction z\^eta des hauteurs d'une vari\'et\'es torique en caract\'eristique positive. 
La question de savoir si cet invariant
est trivial ou non n'est nullement \'evidente, et a \'et\'e r\'esolue par Colliot-Th\'el\`ene et Suresh.
Comme d\'ej\`a indiqu\'e, dans une version pr\'ec\'edente de ce texte, il \'etait affirm\'e \`a tort que cet invariant
subsistait dans l'expression finale du terme principal de la fonction z\^eta des hauteurs des vari\'et\'es 
toriques.

On conserve les hypoth\`ese et notations introduites dans la sous-section \ref{subsubsec:introflasque}.
et on suppose en outre dans cette sous-section
que  $K$ est un corps de fonctions, 
de corps des constantes $k$.
De la r\'esolution flasque \eqref{eq:res:fl}
on tire
au niveau des espaces ad\'eliques la suite exacte
\begin{equation}
0\longto T_Q(\ak) \longto T_P(\ak)\longto T(\ak) .
\end{equation}
Enfin, en prenant l'image de cette suite exacte
par l'application $\deg_T$ (cf. section \ref{subsec:ledegre}), 
on obtient un complexe
\begin{equation}
\ecD_{T_Q} \longto \ecD_{T_P}\longto \DT.
\end{equation}
\begin{lemme}\label{lm:cokerdtfini}
Le conoyau de la fl\`eche
$
\ecD_{T_P}\longto \DT
$
est fini, et il ne d\'epend pas du choix de la r\'esolution flasque.
\end{lemme}
\begin{demo}
On a le diagramme commutatif suivant
\begin{equation}
\xymatrix{
\ecD_{T_Q} \ar[r]\ar[d]
&
\ecD_{T_P} \ar[r]\ar[d] 
&
\DT  \ar[d]
\\
\left(Q^{\,G}\right)^{\vee} \ar[r]
&
\left(P^{\,G}\right)^{\vee} \ar[r]  
&
\left(X(T)^{\,G}\right)^{\vee}  
}
\end{equation}
o\`u les fl\`eches verticales sont de conoyau fini.
Par ailleurs le quotient $P^G/X(T)^G$ s'injecte dans $Q^G$, 
et est donc sans torsion.
Ainsi le morphisme $\left(P^{\,G}\right)^{\vee}\to \left(X(T)^{\,G}\right)^{\vee}$ est surjectif, 
et donc la fl\`eche
\begin{equation}
\ecD_{T_P}\longto \DT
\end{equation}
est de conoyau fini.

Montrons que ce conoyau ne d\'epend pas du choix de la r\'esolution flasque. En effet 
soit $P_1$ un $G$-module de permutation, on obtient \`a partir
de \eqref{eq:res:fl} 
une nouvelle r\'esolution flasque 
\begin{equation}\label{eq:res:fl:1} 
\ecR_1\,:\,0\longto X(T) \longto P\oplus P_1 \longto Q\oplus P_1 \longto 0,
\end{equation}
qui induit un morphisme $T_P\times T_{P_1} \to T$ lequel se factorise en
\begin{equation}\label{eq:facto}
T_P\times T_{P_1}\longto T_P \longto T,
\end{equation}
o\`u la premi\`ere fl\`eche est la projection naturelle et la deuxi\`eme le morphisme induit par \eqref{eq:res:fl}.
On a un isomorphisme naturel $\ecD_{T_P\times T_{P_1}}\isom \ecD_{T_P}\times \ecD_{T_{P_1}}$
et le morphisme $\ecD_{T_P}\times \ecD_{T_{P_1}}\to \ecD_{T_P}$ induit par 
la premi\`ere fl\`eche de \eqref{eq:facto} n'est autre que la projection naturelle.
On a donc
\begin{equation}
\Coker(\ecD_{T_P\times T_{P_1}}\to \DT)=\Coker(\ecD_{T_P}\to \DT).
\end{equation}
Par ailleurs si
\begin{equation}\label{eq:res:fl:2} 
0\longto X(T) \longto P'\longto Q' \longto 0.
\end{equation}
est une autre r\'esolution flasque de $X(T)$, 
par \cite[lemme 5]{CTS:Requiv}, 
les $G$-modules $Q$ et $Q'$ sont isomorphes 
apr\`es addition de $G$-modules de permutation convenables.
\end{demo}
On note $\nindex{$\KT$}\KT$ le cardinal du conoyau
de la fl\`eche $
\ecD_{T_P}\longto \DT
$.

\begin{prop}
\label{prop:kt=1}
On a $\KT=1$ dans les cas suivants : 
\begin{itemize}
\item $T$ v\'erifie l'approximation faible (i.e. $A(T)=0$)
\item $T$ est anisotrope, 
\item $T$ est d\'eploy\'e par une extension dans laquelle $k$ est alg\'ebriquement clos.
\end{itemize}
\end{prop}

\begin{demo}
Rappelons que l'on d\'esigne par $L$ une extension galoisienne finie
d\'eployant $T$, et que son groupe de Galois est not\'e $G$.

Si $A(T)=0$, le lemme \ref{lm:appfaible} montre que le morphisme
\begin{equation}
T_P(\ak)\longto T(\ak)
\end{equation}
est surjectif. 
Ainsi $\KT=1$.

Si $T$ est anisotrope, $X(T)^{\,G}=0$, donc $\DT=0$, d'o\`u le r\'esultat.

Supposons $T$ d\'eploy\'e par une extension dans laquelle
$k$ est alg\'ebriquement clos. Alors $k$ est encore
alg\'ebriquement clos dans la cl\^oture galoisienne d'une telle extension.
On peut ainsi  supposer que $k$
est alg\'ebriquement clos dans $L$.
D'apr\`es le corollaire \ref{cor:ct0} 
on a $\CT=0$, 
i.e. 
\begin{equation}
\DT=\left(X(T)^{\,G}\right)^{\vee},
\end{equation}
et de m\^eme
\begin{equation}
\ecD_{T_P}
=
\left(P^{\,G}\right)^{\vee}.
\end{equation}
Mais, comme d\'ej\`a vu dans la preuve du lemme \ref{lm:cokerdtfini},  
la fl\`eche
\begin{equation}
\left(P^{\,G}\right)^{\vee} \to \left(X(T)^{\,G}\right)^{\vee}
\end{equation}
est surjective, d'o\`u le r\'esultat.
\end{demo}

\begin{rem}
Le cas o\`u $A(T)=0$ 
se produit en particulier quand $T$ est quasi-d\'eploy\'e.
Plus g\'en\'eralement, 
d'apr\`es la suite exacte \eqref{eq:exseqcts}, 
il se produit si  $H^1(G,Q)=0$. 
D'apr\`es \cite[Corollaire 2]{CTS:Requiv} et la dualit\'e de Nakayama, 
ceci est v\'erifi\'e en particulier 
quand $T$ est d\'eploy\'e par une extension 
m\'etacyclique\footnote{dont le groupe de Galois est \`a Sylow cycliques} 
de $K$.
\end{rem}

Colliot-Th\'el\`ene et Suresh ont exhib\'e dans \cite{CTSu} un exemple de tore alg\'ebrique $T$
ne v\'erifiant pas $\KT=1$, que  nous d\'ecrivons \`a pr\'esent.
Soit $G=\Z/2\times \Z/2$ le groupe de Klein, $\tau$ et $\sigma$ des \'el\'ements de $G$
tels que $G$ est engendr\'e par $\tau$ et $\sigma$.
On note $I_G$ le noyau de l'augmentation $\Z[G]\to \Z$
et $N_G$ l'\'el\'ement $\sumu{g\in G} g$ de $\Z[G]$.
Soit $N$ le sous-$G$-module de $\Z[G]$
d\'efini par 
\begin{equation}
N=\{n\in \Z[G],\quad \exists m\in I_G,\quad \sigma.n-n=m+\tau\,m\text{ et }\tau.n-n=m+\sigma\,m\}
\end{equation}
Soit $L/K$ est une extension galoisienne de groupe $G$ et $T$ le tore alg\'ebrique
sur $K$ dont le module des caract\`eres est $N^{\vee}$. 
Alors $T$ s'identifie \`a un sous-tore de $\Res_{L/K}\G_m$,
et $T(K)$ s'identifie \`a un sous-groupe de $\Res_{L/K}\G_m(K)=L^{\times}$, 
plus pr\'ecis\'ement
\begin{equation}\label{eq:descr:TK}
T(K)=\{y\in 
L^{\times},\quad\exists x\in L^{\times},\quad 
N_{L/K}(x)=1,\,\,
^\sigma y\,y^{-1}=x\,^\tau x\text{ et }^\tau y\,y^{-1}=x\,^\sigma x\}.
\end{equation}
\begin{thm}[Colliot-Th\'el\`ene, Suresh]
\label{thm:ctsu}
Soit $k$ un corps fini de caract\'eristique diff\'erente de $2$ dans lequel $-1$ est un carr\'e et $K=k(t)$ le corps
des fractions rationnelles en une ind\'etermin\'ee sur $k$. Soit $u$ un \'el\'ement de $k$ qui n'est pas un carr\'e.
Soit $L$ l'extension $K\left(\sqrt{u},\sqrt{t}\right)$. C'est une extension galoisienne de $K$ de groupe 
$G$. On note $\sigma$ (respectivement $\tau$) le g\'en\'erateur du groupe de Galois de  
$K\left(\sqrt{u}\right)/K$ (respectivement $K\left(\sqrt{t}\right)/K$). Soit $T$ le tore alg\'ebrique
sur $K$ de module de cocaract\`eres $N$. Alors $\KT$ est diff\'erent de $1$.
\end{thm}
La d\'emonstration utilise une r\'esolution flasque explicite de $T$,
et la proposition suivante, que nous utiliserons \`a la sous-section \ref{subsec:rem:casfonc} pour montrer que sur certaines
vari\'et\'es toriques, les hauteurs logarithmique canoniques locales peuvent ne pas \^etre \`a valeurs enti\`eres.
\begin{prop}\label{prop:ctsu}
On conserve les notations du th\'eor\`eme \ref{thm:ctsu}.
\begin{enumerate}
\item
On a $\sqrt{u}\sqrt{t}\in T(K)$.
\item
Soit $v$ la place de $K$ d'uniformisante $t$, $\V$ l'unique place de
$L$ divisant $v$ et $G_v$ son groupe de d\'ecomposition. 
Alors $G_v=G$, $v$ est ramifi\'ee, et l'indice de ramification est $2$.
\item
On a 
\begin{equation}
\deg_{T,L,\V}\left(\sqrt{u}\sqrt{t}\right)= N_G
\end{equation}
et
\begin{equation}
\left(X(T)^{\vee}\right)^{G_v}=\Z\,N_G.
\end{equation}
\end{enumerate}
\end{prop}
\begin{demo}
Soit $i\in k$ tel que $i^2=-1$.
Alors $N_{L/K}i=i^4=1$. Or on a 
\begin{equation}
^\sigma \left(\sqrt{u}\sqrt{t}\right)\,\left(\sqrt{u}\sqrt{t}\right)^{-1}
=
-\left(\sqrt{u}\sqrt{t}\right)\,\left(\sqrt{u}\sqrt{t}\right)^{-1}
=
-1=i^2=i\,^\tau i
\end{equation}
et
\begin{equation}
^\tau \left(\sqrt{u}\sqrt{t}\right)\,\left(\sqrt{u}\sqrt{t}\right)^{-1}
=
-\left(\sqrt{u}\sqrt{t}\right)\,\left(\sqrt{u}\sqrt{t}\right)^{-1}
=
-1=i^2=i\,^\sigma i.
\end{equation}
Ainsi, d'apr\`es \eqref{eq:descr:TK}, $\sqrt{u}\sqrt{t}$ est dans $T(K)$, ce qui montre le premier point.

Le deuxi\`eme point est imm\'ediat.

Pour le troisi\`eme point, on note qu'on a le diagramme commutatif suivant
\begin{equation}
\xymatrix{
T(K_v)\ar[d]^{\deg_{T,L,\V}}\ar[r]& \Res_{L/K}\G_m(K_v)\isom L_{\V}\ar[d]^{\deg_{\Res_{L/K}\G_m,L,\V}}
\\
\left(X(T)^{\vee}\right)^{\,G_v}\ar[r]
&
\Z[G]^{G_v}=\Z\,N_G
}
\end{equation}
dont les fl\`eches horizontales sont injectives.

Pour tout $y\in L_{\V}^{\times}$, on a
\begin{equation}
\deg_{\Res_{L/K}\G_m,L,\V}(y)=\V(y)\,N_G.
\end{equation}
Ainsi
\begin{equation}
\deg_{T,L,\V}\left(\sqrt{u}\sqrt{t}\right)
=
\deg_{\Res_{L/K}\G_m,L,\V}\left(\sqrt{u}\sqrt{t}\right)
=
\V\left(\sqrt{u}\sqrt{t}\right)\,N_G
=
N_G.
\end{equation}
On en d\'eduit aussit\^ot que $\left(X(T)^{\vee}\right)^{G_v}=\Z\,N_G$.
\end{demo}

Comme le $G$-module $P$ apparaissant dans \eqref{eq:res:fl} est de permutation, 
le d\'ebut de la suite exacte
longue de cohomologie tir\'ee de \eqref{eq:res:fl}
s'\'ecrit
\begin{equation}
0\longto X(T)^{\,G} \longto P^{\,G}\longto Q^{\,G} \longto H^1(G,X(T)) \longto 0.
\end{equation}
L'homologie du complexe 
\begin{equation}
\left(Q^{\,G}\right)^{\vee}\to \left(P^{\,G}\right)^{\vee} 
\to 
\left(X(T)^{\,G}\right)^{\vee}
\end{equation}
est donc de cardinal $\card{H^1(G,X(T))}$.

Le lemme suivant est utilis\'e dans le calcul du terme principal de la fonction z\^eta des hauteurs
dans le cas fonctionnel.

\begin{lemme}
\label{prop:relationkt}
L'homologie du complexe
\begin{equation}
\ecD_{T_Q} \longto \ecD_{T_P}\longto \DT
\end{equation}
est de cardinal
\begin{equation}
\frac
{\card{H^1(G,X(T))}\,\card{\ecC_{T_Q}}\,\card{\CT}\,\KT}
{\card{\ecC_{T_P}}}
.
\end{equation}
\end{lemme}

\begin{demo}
Notons $\ecH$ ce cardinal.
On a le diagramme commutatif suivant
\begin{equation}
\xymatrix{
&
0 \ar[d]
&
0 \ar[d] 
&
0 \ar[d]
&
\\
0 \ar[r]
&
\ecD_{T_Q} \ar[r]\ar[d]
&
\ecD_{T_P} \ar[r]\ar[d] 
&
\DT  \ar[d] \ar[r]
&
0
\\
0 \ar[r]
&
\left(Q^{\,G}\right)^{\vee} \ar[r]\ar[d]
&
\left(P^{\,G}\right)^{\vee} \ar[r]\ar[d]
&
\left(X(T)^{\,G}\right)^{\vee}  \ar[r]\ar[d]
&
0
\\
0 \ar[r]
&
\ecC_{T_Q} \ar[r]\ar[d]
&
\ecC_{T_P} \ar[r]\ar[d] 
&
\CT \ar[r]\ar[d]
&
0
\\
&
0
&
0 
&
0 
&
}
\end{equation}
o\`u les verticales sont exactes et les horizontales sont des complexes, 
en d'autres termes on a une suite exacte de complexes. 
Par ailleurs les complexes sont exacts en : 
$\ecD_{T_Q}$, 
$\left(Q^{\,G}\right)^{\vee}$,
$\left(X(T)^{\,G}\right)^{\vee}$
et $\CT$.
Les caract\'eristiques d'Euler-Poincar\'e de ces complexes sont alors, 
de haut en bas :
$
\frac
{\KT}
{\ecH}
$
,
$
\frac
{1}
{\card{H^1(G,X(T))}}
$
et
$
\frac
{\card{\ecC_{T_Q}}\,\card{\CT}}
{\card{\ecC_{T_P}}}
$
.
Comme on a une suite exacte de complexes,
on en d\'eduit la relation
\begin{equation}
\frac
{\ecH}
{\KT}
\,
\frac
{1}
{\card{H^1(G,X(T))}}
\,
\frac
{\card{\ecC_{T_P}}}
{\card{\ecC_{T_Q}}\,\card{\CT}}
=
1,
\end{equation}
d'o\`u le r\'esultat.
\end{demo}

\subsection{Mesure ad\'elique et nombre de Tamagawa d'un tore alg\'ebrique}\label{sec:tamagawa:tore}

Si $X$ est une vari\'et\'e alg\'ebrique lisse d\'efinie sur un corps global $K$, 
on a d\'ej\`a rappel\'e qu'une m\'etrisation du faisceau anticanonique de $X$ permet, 
pour tout $v\in\placesde{K}$, 
de construire une mesure $\omega_{X,v}$ sur l'espace analytique $X(K_v)$. 
Rappelons aussi que tout choix d'une section globale $\omega$ partout 
non nulle du faisceau anticanonique en fournit une m\'etrisation par la formule
\begin{equation}
\forall\,x\in X(K_v),\quad
\forall\,s\in \omega_X^{-1}(x),\quad 
\norm{s}_v=\frac{\abs{s(x)}_v}{\abs{\omega(x)}_v},
\end{equation}
le choix de la valeur absolue $|\,.\,|_v$ sur $\omega_X^{-1}(x)$ \'etant arbitraire.  
Si $X=G$ est un groupe alg\'ebrique et $\omega$ est de plus choisie $G$-invariante \`a gauche, 
les mesures locales obtenues sont des mesures de Haar \`a gauche.

Ono, 
dans l'article \cite{Ono:algtor}, 
d\'efinit le nombre de Tamagawa d'un tore alg\'ebrique 
d\'efini sur un corps global. 
Dans \cite{Ono:tamnum}, 
il \'etablit une relation simple entre ce nombre de Tamagawa et certains invariants de type cohomologique du tore 
(cf. le th\'eor\`eme \ref{thm:ono} ci-dessous), 
montrant en particulier la rationalit\'e du nombre de Tamagawa 
(chose nullement \'evidente sur la d\'efinition initiale). 
Cette relation joue un r\^ole important dans l'interpr\'etation du terme principal 
des fonctions z\^eta des hauteurs des vari\'et\'es toriques, 
dans le cas arithm\'etique comme dans le cas fonctionnel. 
Cependant, dans le cas des corps de fonctions, il s'est av\'er\'e que la d\'efinition du nombre de Tamagawa 
d'un tore alg\'ebrique donn\'ee dans \cite{Ono:algtor} \'etait incompatible avec la relation du th\'eor\`eme \ref{thm:ono}. 
Par la suite Oesterl\'e a montr\'e dans \cite{Oes:invent} qu'en introduisant dans la d\'efinition d'Ono un facteur correctif 
\'egal au cardinal du conoyau du degr\'e (ce que nous avons not\'e $\CT$), 
la relation du th\'eor\`eme \ref{thm:ono} devenait correcte. 
Nous rappelons dans cette section la construction du nombre de Tamagawa d'un tore alg\'ebrique, 
et le r\'esultat principal de \cite{Ono:tamnum}.

Soit $T$ un tore alg\'ebrique de dimension $d$, d\'efini sur un corps global $K$. 
Soit $\Omega_T$ une $d$-forme diff\'erentielle $K$-rationnelle $T$-invariante sur $T$ 
(une telle forme est uniquement d\'etermin\'ee \`a multiplication par un \'el\'ement de $K^{\times}$ pr\`es). 
Cette forme induit donc pour tout $v\in\placesde{K}$ une mesure de Haar $\omega_{T,v}$ sur $T(K_v)$.

On a alors pour presque toute place finie $v$ la relation
\begin{equation}
\int\limits_{T(\Ov)}\omega_{T,v}
=
\frac
{1}
{L_v(1,X(T))}.
\end{equation}
Posons, si $v$ est finie, 
\begin{equation}\label{eq:def:dmuv}
d\mu_v
=\frac
{1}
{L_v(1,X(T))}
\,\omega_{T,v},
\end{equation}
et si $v$ est archim\'edienne
\begin{equation}
d\mu_v=\omega_{T,v}.
\end{equation}
On aura alors
\begin{equation}
\int\limits_{T(\Ov)}d\mu_v=1
\end{equation}
pour presque tout $v$.

On peut alors d\'efinir une mesure de Haar $\omega_{T}$ sur $T(\ak)$ en posant
\begin{equation}\label{eq:defomegaT}
\omega_{T}
=
c_{K,\dim(T)}\,\prod_{v\in \placesde{K}}d\mu_v
\end{equation}
(cf \eqref{eq:def:ckv} pour la d\'efinition de $c_{K,\dim(T)}$).

Si la forme $\Omega_T$ est chang\'ee en $\lambda\,\Omega_T$,
o\`u $\lambda\in K^{\times}$, $d\mu_v$ est chang\'ee en $\abs{\lambda}_vd\mu_v$
pour tout $v$.
Ainsi, par la formule du produit, 
$\omega_T$ ne d\'epend pas du choix de la forme $\Omega_T$. 

Notons qu'en particulier
\begin{equation}
\int\limits_{\K(T)}\omega_{T}
=
c_{K,\dim(T)}\,
\prod_{v}\int\limits_{T(\Ov)}d\mu_v
\end{equation}
est non nul.

\`A partir de \nindex{$\omega_T$}$\omega_T$, 
on construit une mesure de Haar sur $T(\ak)^1/T(K)$ 
de la mani\`ere suivante.

Dans le cas arithm\'etique, 
soit $dt$ la mesure de Lebesgue sur $(X(T)^G)^{\vee}_{\R}$,
normalis\'ee par le r\'eseau $(X(T)^G)^{\vee}$. 
Rappelons que $\deg_T$ induit
un isomorphisme de groupes topologiques
\begin{equation}
\deg_T\,:\,T(\ak)/T(\ak)^1\longisom \left(X(T)^G\right)^{\vee}_{\R}.
\end{equation}
Soit $\wt{\omega_T}$ la mesure quotient 
sur $T(\ak)/T(K)$ induite par $\omega_T$,
$T(K)$ \'etant muni de la mesure discr\`ete.
Soit \nindex{$\omega_T^1$}$\omega_T^1$ la mesure sur $T(\ak)^1/T(K)$
d\'efinie par la relation
\begin{equation}
\wt{\omega_T}=\omega_T^1\,.\,\left(\deg_T^{-1}\right)_{\ast}(dt).
\end{equation}

Dans le cas fonctionnel, 
comme $T(\ak)^1 $ est ouvert dans $T(\ak)$, 
la restriction de $\omega_T$ \`a $T(\ak)^1$ fournit 
une mesure de Haar sur $T(\ak)^1$. 
Soit $\omega_T^1$ la mesure quotient sur $T(\ak)^1/T(K)$.

On pose, 
dans le cas arithm\'etique,
\begin{equation}\label{eq:defbTarit}
b(T)=\int\limits_{T(\ak)^1 /T(K)}\!\!\omega_T^1
\end{equation}
et dans le cas fonctionnel
\begin{equation}\label{eq:defbTfonc}
b(T)=\log(\qk)^{-rg\left(X(T)^{\,G}\right)}\,\int\limits_{T(\ak)^1 /T(K)}\!\!\omega_T^1.
\end{equation}
(rappelons que $T(\ak)^1/T(K)$ est compact).

On pose alors, suivant Ono,
\begin{equation}\label{eq:def:taut:arit}
\nindex{$\tau(T)$}\tau(T)=\frac{b(T)}{\ell(X(T))}
\end{equation}
dans le cas arithm\'etique et, 
suivant Oesterl\'e (\cite[I.5.9 et 5.12.]{Oes:invent}),
\begin{equation}\label{eq:def:taut:fonc}
\tau(T)=\frac{1}{\card{\CT}}\,\frac{b(T)}{\ell(X(T))}.
\end{equation}
dans le cas fonctionnel.
Le nombre $\tau(T)$ est appel\'e \termin{nombre de Tamagawa} du tore alg\'ebrique $T$. 

L'objet de l'article \cite{Ono:tamnum}, 
corrig\'e par Oesterl\'e dans le cas fonctionnel, 
est la d\'emonstration du 
\begin{thm}[Ono, Oesterl\'e]\label{thm:ono}
On a la relation 
\begin{equation}
\tau(T)=\frac{\card{H^1(G,X(T))}}{\card{\cha(T)}}.
\end{equation}
\end{thm}
Ce r\'esultat, 
comme d\'ej\`a indiqu\'e, 
sera utile lors de l'interpr\'etation du terme principal de la fonction z\^eta des hauteurs.

\section{Hauteurs sur une vari\'et\'e torique et fonction z\^eta associ\'ee}

\subsection{G\'eom\'etrie des vari\'et\'es toriques}
\label{subsec:vte:tor}

\subsubsection{Vari\'et\'es toriques d\'eploy\'ees}
Nous rappelons la construction 
des vari\'et\'es toriques d\'eploy\'ees. 
On renvoie aux r\'ef\'erences classiques sur les vari\'et\'es toriques, 
comme \cite{Ful:toric} et \cite{Oda:conv} pour plus de d\'etails et
la preuve des assertions qui suivent.

Soit $M$ un $\Z$-module libre de rang fini et $\Lambda$ un c\^one de $M_{\R}$.
L'int\'erieur relatif de $\Lambda$ sera not\'e \nindex{$\intrel$}$\intrel \left(\Lambda\right)$.
Le \dindex{c\^one dual}\termin{c\^one dual} de $\Lambda$ est not\'e $\Lambda^{\vee}$ et est d\'efini par 
\begin{equation}
\Lambda^{\vee}\eqdef\{y\in M^{\vee}_{\R},\quad \forall x\in \Lambda,\quad \acc{y}{x}\geq 0\}
\end{equation}
Le c\^one $\Lambda$ 
dit \dindex{c\^one poly\'edral rationnel}\termin{poly\'edral rationnel} 
s'il est engendr\'e par un ensemble fini d'\'el\'ements de $M$.
Si $\Lambda$ est un c\^one poly\'edral rationnel, son dual 
est encore poly\'edral rationnel. 
Le c\^one $\Lambda$ est dit \termin{strictement convexe} si 
\begin{equation}
\Lambda\cap -\Lambda=\{0\}.
\end{equation}
Remarquons que $\Lambda$ est strictement convexe si et seulement si $\Lambda^{\vee}$
est d'int\'erieur non vide.

On fixe un corps de base $L$.
Soit $M$ un $\Z$-module libre de rang fini et $T_L=\Spec(L[M])$, 
c'est-\`a-dire que $T_L$ est le tore alg\'ebrique d\'efini sur $L$ ayant pour groupe de caract\`eres $M$. 

Soit $\sigma$ un c\^one poly\'edral rationnel strictement convexe de $M^{\vee}_{\R}$.
\`A un tel c\^one est associ\'e une vari\'et\'e affine normale 
d\'efinie sur $L$ 
\begin{equation}
X_{\sigma,L}=\Spec(L[\sigma^{\vee}\cap M]),
\end{equation}
munie naturellement d'une action de $T_L$.

Un \dindex{\'eventail}\termin{\'eventail} de $M^{\vee}$ est un ensemble fini $\Sigma$ de c\^ones poly\'edraux rationnels 
strictement convexes de $M^{\vee}_{\R}$, v\'erifiant les conditions suivantes~: 
\begin{itemize}
\item toute face d'un c\^one de $\Sigma$ est un c\^one de $\Sigma$ ; 
\item l'intersection de deux c\^ones de $\Sigma$ est une face de chacun des deux c\^ones. 
\end{itemize}

Si $\sigma$ et $\sigma'$ sont deux c\^ones de $\Sigma$, les inclusions 
$\sigma\cap \sigma' \subset \sigma$ et $\sigma\cap \sigma' \subset \sigma'$ 
induisent des immersions ouvertes $X_{\sigma\cap \sigma',L}\hookrightarrow X_{\sigma,L}$
et $X_{\sigma\cap \sigma',L}\hookrightarrow X_{\sigma',L}$, qui sont compatibles
avec les actions de $T_L$.
Ceci donne un proc\'ed\'e de recollement des vari\'et\'es $X_{\sigma,L}$ pour $\sigma$ d\'ecrivant
les c\^ones de $\Sigma$, lequel proc\'ed\'e est compatible aux
actions de $T_L$, et permet de construire une vari\'et\'e normale $X_{\Sigma,L}$, 
munie d'une action de $T_L$ : c'est la vari\'et\'e torique (d\'efinie sur $L$) associ\'ee \`a l'\'eventail $\Sigma$.

Un \'eventail $\Sigma$ est dit \termin{non d\'eg\'en\'er\'e} si ses c\^ones ne sont pas
inclus dans un sous-espace strict de $M^{\vee}_{\R}$.

Un \'eventail $\Sigma$ est dit \termin{r\'egulier} si tout c\^one de $\Sigma$ 
est engendr\'e par une partie d'une $\Z$-base de $M^{\vee}$, 
et \termin{complet} si les c\^ones de $\Sigma$ recouvrent $M^{\vee}_{\R}$.
Un \'eventail $\Sigma$ est r\'egulier (respectivement complet) si et seulement si la vari\'et\'e $X_{\Sigma,L}$
est lisse (respectivement compl\`ete).

Un \'eventail $\Sigma$ est dit \termin{projectif} si la vari\'et\'e $X_{\Sigma,L}$ est projective.

Les \dindex{rayons d'un \'eventail}\termin{rayons} de $\Sigma$ sont les c\^ones de $\Sigma$ de dimension 1.
On note \nindex{$\Sigma(1)$}$\Sigma(1)$ l'ensemble des rayons de $\Sigma$. 
Pour $\alpha\in \Sigma(1)$, on note $\roa$ l'\'el\'ement de $M^{\vee}$ qui engendre 
le mono\"\i de $\alpha\cap M^{\vee}$.
Pour tout c\^one $\sigma$ de $\Sigma$ nous notons 
\nindex{$\sigma(1)$}
\begin{equation}
\sigma(1)=\{\alpha\in\Sigma(1),\, \alpha\subset \sigma\}
\end{equation}
(ainsi $\{0\}(1)=\vide$). 

L'application qui \`a un \'el\'ement $\alpha$ de $\Sigma(1)$ associe l'adh\'erence dans 
$\xsl$ de la $T_L$-orbite ferm\'ee de $X_{\alpha,L}$
d\'efinit une bijection de $\Sigma(1)$ sur l'ensemble des diviseurs irr\'eductibles de $\xsl$ 
contenu dans le bord $\xsl\setminus T_L$. 
Pour tout rayon $\alpha$, on note $\Da$ le diviseur irr\'eductible ainsi associ\'e \`a $\alpha$ ;
c'est un diviseur $T_L$-invariant.
On note \nindex{$\Ps$}$\Ps$ le $\Z$-module libre de base $(\Da)_{\alpha\in \Sigma(1)}$ ;
ce n'est autre que le groupe des diviseurs de Weil $T_L$-invariants sur $\xsl$.

On note \'egalement \nindex{$\PL(\Sigma)$}$\PL(\Sigma)$ le groupe des applications $\Sigma$-lin\'eaires par morceaux sur 
$M^{\vee}$, c'est-\`a-dire les applications $\varphi\,:\,M^{\vee}\to \Z$ telles que 
la restriction de $\varphi$ \`a $\sigma \cap M^{\vee}$ est lin\'eaire pour tout c\^one $\sigma$ de $\Sigma$.

On suppose \`a pr\'esent $\Sigma$ r\'egulier et non d\'eg\'en\'er\'e.
L'application 
\begin{equation}\label{eq:iso:pls:ps}
\varphi\mapsto \sum_{\alpha\in \Sigma(1)}\varphi(\roa)\,\Da
\end{equation}
est alors un isomorphisme de groupes qui permet d'identifier $\PL(\Sigma)$ \`a $\Ps$. Par la suite nous utiliserons 
souvent cette identification. 

Comme le groupe de Picard de $T_L$ est trivial, l'application qui \`a un \'el\'ement de
$\Ps$ associe sa classe dans $\Pic(\xsl)$, 
induit une suite exacte
\begin{equation}
0\longto L[T_L]^{\times}/L^{\times} \longto \Ps  \longto \Pic(\xsl)\longto 0,
\end{equation}
o\`u la fl\`eche $L[T_L]^{\times}/L^{\times} \longto \Ps$ est induite par l'application
qui \`a une fonction rationnelle associe son diviseur.
D'apr\`es le lemme de Rosenlicht, $L[T_M]^{\times}/L^{\times}$ est isomorphe
au groupe des caract\`eres de $T_M$, c'est-\`a-dire $M$, 
d'o\`u la suite exacte
\begin{equation}\label{eq:definir}
0\longto M \overset{\gamma}{\longto} \Ps  \longto \Pic(\xsl)\longto 0.
\end{equation}
Pour $m\in M$, on a 
\begin{equation}\label{eq:gammam:dep}
\gamma(m)=\sum_{\alpha\in \Sigma(1)} \acc{m}{\roa} \Da.
\end{equation}

En outre $\Pic(\xsl)$ est un $\Z$-module libre de rang fini.

On a le r\'esultat suivant :
\begin{prop}\label{prop:antican}
La classe dans $\Pic(\xsl)$ du faisceau anticanonique de $\xsl$
co\"\i ncide avec la classe du diviseur $\sumu{\alpha\in \Sigma(1)}\Da$.
\end{prop}

\subsubsection{Vari\'et\'es toriques non d\'eploy\'ees}\label{subsubsec:vte:tor:nondep}
Nous rappelons \`a pr\'esent la construction des vari\'et\'es toriques non n\'ecessairement d\'eploy\'ees.
On se limitera au cas des vari\'et\'es toriques projectives et lisses.
On se donne un tore alg\'ebrique $T$ d\'efini sur un corps $K$. 
Soit $L/K$ une extension finie galoisienne de groupe $G$ d\'eployant $T$. 
Soit $\Sigma\subset X(T)^{\vee}_{\R}$ un \'eventail projectif et lisse
et $\xsl$ la vari\'et\'e projective et lisse associ\'ee.

On suppose en outre que l'action de $G$ sur $X(T)^{\vee}_{\R}$ pr\'eserve les c\^ones de $\Sigma$ 
(on dira alors que $\Sigma$ est un \dindex{$G$-\'eventail}$G$-\'eventail\footnote{
On peut montrer qu'un $G$-\'eventail projectif et lisse de $X(T)^{\vee}_{\R}$
existe toujours, cf. \cite{CTHaSk:comp}.}). 
Alors $G$ agit aussi sur le sch\'ema  projectif $\xsl$, 
et on peut donc consid\'erer le quotient de $\xsl$ par $G$. 
Ce quotient est une vari\'et\'e $\xs$ d\'efinie sur $K$, 
qui est une compactification \'equivariante projective et lisse de $T$ (cf. \cite[\S 1]{Vos:proj}).

Le $\Z$-module $\Ps$ est alors muni naturellement d'une action de $G$, et le $G$-module
r\'esultant est un $G$-module de permutation. Par ailleurs, l'action de $G$ sur les c\^ones de $\sigma$
induit naturellement une action de $G$ sur $\PL(\Sigma)$ et l'isomorphisme
\eqref{eq:iso:pls:ps} est un isomorphisme de $G$-modules.

Nous notons \nindex{$\sg$}$\sg$ l'ensemble des orbites de $\Sigma(1)$ sous l'action de $G$. 
Pour chaque $\alpha\in \sg$ nous choisissons arbitrairement 
un \'el\'ement de $\alpha$, 
nous notons \nindex{$\roa$}$\roa$ le g\'en\'erateur de cet \'el\'ement 
ainsi que \nindex{$\ga$}$\ga$ le stabilisateur de $\roa$, de sorte que le $G$-ensemble 
$\alpha$ s'identifie \`a $G/\ga$.
On d\'eduit de ces choix un isomorphisme
\begin{equation}\label{eq:idps}
\Ps\longisom \underset{\alpha\in\Sigma(1)/G}{\bigoplus}\Z\left[G/\ga\right]
\end{equation}
d'o\`u un isomorphismes de $K$-tores
\begin{equation}\label{eq:idtps}
T_{\Ps}\longisom \underset{\alpha\in\Sigma(1)/G}{\prod}\Res_{L^{\,\ga}/K} \G_m.
\end{equation}
Pour $\alpha\in \sg$, on notera \nindex{$\Ka$}$\Ka$ le corps $L^{\,\ga}$, et \nindex{$\ia$}$\ia$ 
la projection $G$-\'equivariante
de $\Ps$ sur $\Z\left[G/\ga\right]$ induite par l'isomorphisme \eqref{eq:idps}.
Le morphisme de $K$-tores associ\'e est l'injection 
$\Res_{\Ka/K} \G_m\to T_{\Ps}$ induite par l'isomorphisme \eqref{eq:idtps}. 

Pour $\alpha\in \sg$, on note
$\Da=\sumu{\beta\in \alpha} D_{\beta}$. Ainsi $(\Da)_{\alpha\in \sg}$
est une base de $\Ps^G$. On note $(\Da^{\vee})_{\alpha\in \sg}$ sa base duale.

La suite exacte de $\Z$-modules libres de rang fini
\begin{equation}\label{eq:resflT}
0\longto X(T) \overset{\gamma}{\longto} \Ps  \longto \Pic(\xsl)\longto 0
\end{equation}
est une suite exacte de $G$-modules.
Il en r\'esulte par dualit\'e une suite exacte de tores alg\'ebriques
\begin{equation}\label{eq:resflT2}
0\longto T_{\NS} \longto T_{\Ps } \longto T\longto 0.
\end{equation}
Comme $\Ps$ est un $G$-module de permutation, on a $H^1(G,\Ps)=0$.
Ainsi, en prenant les $G$-invariants dans \eqref{eq:resflT}, on obtient la suite exacte
\begin{equation}\label{eq:resflT:G}
0\longto X(T)^G \overset{\gamma}{\longto} \Ps^G  \longto \Pic(\xs)\longto H^1(G,X(T))\longto 0.
\end{equation}
D'apr\`es \eqref{eq:gammam:dep}, on a 
\begin{equation}\label{eq:gammam:nondep}
\forall m\in \xtg,\quad \gamma(m)
=
\sumu{\alpha\in\sg} \acc{m}{\roa} \Da.
\end{equation}
On en d\'eduit qu'on a
\begin{equation}\label{eq:gammaveeDaroa}
\forall \alpha\in \sg,\quad \gamma^{\vee}\left(\Da^{\vee}\right)=\roa.
\end{equation}

On a en outre le r\'esultat suivant :
\begin{lemme}\label{lm:picxsflasque}
$\Pic(\xsl)$ est un $G$-module flasque. 
\end{lemme}
En d'autres termes, la suite exacte \eqref{eq:resflT}
est une r\'esolution flasque de $X(T)$.
\begin{demo}
Dans le cas o\`u $K$ est de caract\'eristique z\'ero, ce r\'esultat
est contenu dans la preuve de la proposition 6 (page 189) de \cite{CTS:Requiv}.
Le cas o\`u $K$ est de caract\'eristique non nulle en d\'ecoule.
Choisissons en effet un corps $K'$ de caract\'eristique 0 
tel qu'il existe une extension galoisienne $L'$ de $K'$ de groupe $G$
(rappelons qu'un proc\'ed\'e classique pour construire une telle extension $L'/K'$
est de plonger $G$ dans le groupe sym\'etrique $\grsym_n$ pour $n$ convenable,
et de consid\'erer l'extension $\Q(T_1,\dots,T_n)/\Q(\sigma_1,\,\dots,\sigma_n)$
o\`u les $T_i$ sont des ind\'etermin\'ees et les $\sigma_i$ les polyn\^omes sym\'etriques
\'el\'ementaires en ces ind\'etermin\'ees ; cette extension est galoisienne de groupe
$\grsym_n$, et l'extension $\Q(T_1,\dots,T_n)/\Q(T_1,\dots,T_n)^G$ fournit l'extension 
cherch\'ee).

Il existe donc un tore alg\'ebrique $T_P$ d\'efini sur $K'$ et d\'eploy\'e par $L'$, 
et une compactification lisse et projective $X_{\Sigma,K'}$ de $T_P$, 
telle que les $G$-modules $\Pic(X_{\Sigma,L'})$ et $\Pic(X_{\Sigma,L})$ sont isomorphe, 
et d'apr\`es \cite{CTS:Requiv} $\Pic(X_{\Sigma,L'})$ est flasque.

\end{demo}

\subsection{Hauteurs sur une vari\'et\'e torique}
\label{subsec:fct:zeta:hauteur} 

\subsubsection{Hauteurs locales}
\label{subsubsec:haut:loc} 

On conserve les objets et notations introduits \`a la section pr\'ec\'edente,
et on suppose d\'esormais que $K$ un corps global.

On a un accouplement naturel
\begin{equation}\label{eq:accs:naturel}
\accs{}{}\,:\,\PL(\Sigma)\times X(T)^{\vee}\longto \Z
\end{equation}
lin\'eaire en le premier facteur, et qui prolonge l'accouplement naturel
\begin{equation}
\acc{}{}\,:\,X(T)\times X(T)^{\vee}\longto \Z,
\end{equation}
c'est-\`a-dire que le diagramme
\begin{equation}\label{eq:diag:xtpls}
\xymatrix{
X(T)\times X(T)^{\vee}\ar[rr]^<<<<<<<<<<{\acc{}{}}
\ar[d]^{\gamma\times \Id}
&&\Z\ar[d]^{\Id}\\
\PL(\Sigma)\times X(T)^{\vee}\ar[rr]^<<<<<<<<<<{\accs{}{}}&&\Z
}
\end{equation}
est commutatif.

Soit $v$ une place de $K$. On va d\'efinir un syst\`eme de hauteurs locales en $v$.  
On note $\V$ une place de $L$ divisant $v$ et $G_v$ le groupe de d\'ecomposition 
correspondant.

Supposons tout d'abord $v$ finie. 
Notons $e_v$ l'indice de ramification de $v$ dans l'extension $L/K$.
On d\'efinit un <<produit d'intersection local>>
\begin{equation}
\accsv{}{}\,:\,\PL(\Sigma)^{G_v}\times T(K_v)\longto \frac{1}{e_v}\,\Z
\end{equation}
par la formule
\begin{equation}\label{eq:defaccsv:vfinie}
\accsv{\varphi}{t}=\frac{1}{e_v}\accs{\varphi}{\deg_{T,L,\V}(t)}.
\end{equation}
On v\'erifie que cette d\'efinition ne d\'epend ni du choix de $\V$,
ni du choix de l'extension $L$ d\'eployant $T$.
Par lin\'earit\'e en le premier facteur, on \'etend $\accs{}{}$ en un accouplement
\begin{equation}
\accsv{}{}\,:\,\PL(\Sigma)^{G_v}_{\C}\times T(K_v)\longto \C.
\end{equation}

\begin{lemme}\label{lm:diag:xtgvplsgv}
Le diagramme
\begin{equation}\label{eq:diag:xtgvplsgv}
\xymatrix{
X(T)^{G_v}\times T(K_v)\ar[rr]^<<<<<<<<<<{\acc{\,.\,}{\deg_{T,v}(.)}}
\ar[d]^{\gamma\times \Id}
&&\Z\ar[d]^{\Id}\\
\PL(\Sigma)^{G_v}\times T(K_v)\ar[rr]^<<<<<<<<<<{\accsv{}{}} &&\frac{1}{e_v}\,\Z
}
\end{equation}
est commutatif
\end{lemme}
\begin{demo}
Soit $m\in X(T)^{G_v}$ et $t\in T(K_v)$. On a
\begin{align}
\accsv{\gamma(m)}{t}
&
=
\frac{1}{e_v}\,\accs{\gamma(m)}{\deg_{T,L,\V}(t)}
\\
&
=\label{eq:frac1ev}
\frac{1}{e_v}\,\acc{m}{\deg_{T,L,\V}(t)}
\\
&
=\label{eq:frac1ev2}
\frac{1}{e_v}\,\acc{m}{e_v\,\deg_{T,v}(t)}
\\
&
=
\acc{m}{\deg_{T,v}(t)}.
\end{align}
L'\'egalit\'e \eqref{eq:frac1ev} vient de la commutativit\'e du diagramme \eqref{eq:diag:xtpls},
et l'\'egalit\'e \eqref{eq:frac1ev}
du lemme \ref{lm:degTv:degTLV:fonc}.
\end{demo}

Supposons \`a pr\'esent $v$ archim\'edienne.
On d\'efinit un <<produit d'intersection local>>
\begin{equation}
\accsv{}{}\,:\,\PL(\Sigma)^{G_v}\times T(K_v)\longto \R
\end{equation}
par la formule
\begin{equation}\label{eq:hauteur:locale:archi}
\accsv{\varphi}{t}
=
\frac{1}{[L_{\V}:K_v]}\,\accs{\varphi}{\deg_{T,L,\V}(t)}.
\end{equation}
On v\'erifie que cette d\'efinition ne d\'epend ni du choix de $\V$,
ni du choix de l'extension $L$ d\'eployant $T$.
Par lin\'earit\'e en le premier facteur, on \'etend $\accsv{}{}$ en un accouplement
\begin{equation}
\accsv{}{}\,:\,\PL(\Sigma)^{G_v}_{\C}\times T(K_v)\longto \C.
\end{equation}

\begin{lemme}\label{lm:diag:xtgvplsgv:archi}
Le diagramme
\begin{equation}\label{eq:diag:xtgvplsgv:archi}
\xymatrix{
X(T)^{G_v}\times T(K_v)\ar[rr]^<<<<<<<<<<{\acc{\,.\,}{\deg_{T,v}(.)}}
\ar[d]^{\gamma \times \Id}
&&\R\ar[d]^{\Id}\\
\PL(\Sigma)^{G_v}\times T(K_v)\ar[rr]^<<<<<<<<<<{\accsv{}{}} &&\R
}
\end{equation}
est commutatif.
\end{lemme}
\begin{demo}
Soit $m\in X(T)^G_v$ et $t\in T(K_v)$. On a
\begin{align}
\accsv{\gamma (m)}{t}
&
=
\frac{1}{[L_{\V}:K_v]}\,\accs{\gamma (m)}{\deg_{T,L,\V}(t)}
\\
&
=
\frac{1}{[L_{\V}:K_v]}
\,\acc{m}{\deg_{T,L,\V}(t)}
\\
&
=
\frac{1}{[L_{\V}:K_v]}\,
\acc{m}{[L_{\V}:K_v]\,\deg_{T,v}(t)}
\\
&
=
\acc{m}{\deg_{T,v}(t)}.
\end{align}
La deuxi\`eme \'egalit\'e vient de la commutativit\'e du diagramme \eqref{eq:diag:xtpls},
et la troisi\`eme du lemme \ref{lm:degTv:degTLV:ari:archi}.
\end{demo}

Pour toute place $v$ de $K$, on d\'efinit alors un syst\`eme de hauteurs (exponentielles) locales
\begin{equation}\label{eq:hauteur:locale}
H_{v}\,:
\begin{array}{rcl}
\PL(\Sigma)^{\,G_v}_{\C}\times\,T(K_v)&\longto& \C\\
(\varphi,t)&\longmapsto&\exp\left[\log(q_v)\,\accsv{\varphi}{t}\right].
\end{array}
\end{equation}

\begin{rem}
\`A la pr\'esentation pr\`es, le syst\`eme de hauteurs
utilis\'e est le m\^eme que celui d\'ecrit dans 
\cite[Definition 2.1.5]{BaTs:anis}.
Notons toutefois que les auteurs de \cite{BaTs:anis}
omettent le facteur $\frac{1}{e_v}$
apparaissant dans 
la d\'efinition \eqref{eq:defaccsv:vfinie}
ce qui, pour les places ramifi\'ees, rend leur d\'efinition incorrecte 
au sens o\`u les points 
\ref{item:3:lm:metrique}
et
\ref{item:4:lm:metrique}
du lemme \ref{lm:metrique} ne sont plus v\'erifi\'es.
\end{rem}

\begin{lemme}\label{lm:metrique}
Soit $v$ une place de $K$, $\V$ une place de $L$ divisant $v$,
et $G_v$ le groupe de d\'ecomposition correspondant.
\begin{enumerate}
\item\label{item:1:lm:metrique}
Pour tout
$\varphi\in\PL(\Sigma)^{G_v}_{\C}$,
la fonction $H_v(\varphi,.)$ est invariante sous l'action de $T(\Ov)$. 
\item\label{item:2:lm:metrique}
Pour  tous
$\varphi_1,\varphi_2\in\PL(\Sigma)^{G}_{\C}$,
on a
\begin{equation}\label{eq:Hphi1pphi2}
H_v(\varphi_1+\varphi_2,.)=H_v(\varphi_1,.)\,H_v(\varphi_2,.).
\end{equation}
\item\label{item:3:lm:metrique}
On identifie $T(K_v)$ \`a un sous-groupe de $T(L_\V)$.
On a alors pour tout $\varphi\in \PL(\Sigma)^{G_v}$
et tout $t\in T(K_v)$
\begin{equation}\label{eq:comp:Hphi:ext}
H_{\V}(\varphi,t)=H_v(\varphi,t)^{[L_\V:K_v]}.
\end{equation}
\item\label{item:4:lm:metrique}
Soit $\varphi\in PL(\Sigma)^{\,G}$. 
Il existe
une unique m\'etrique $v$-adiques  $\norm{.}^{\varphi}_v$
sur $\ecO_{\xs}(\varphi)$ v\'erifiant
\begin{equation}
\forall t\in T(K_v),\,\quad H_{v}(\varphi,t)=\left(\norm{\ind(t)}^{\varphi}_v\right)^{-1}
\end{equation}
o\`u $\ind$ est la fonction r\'eguli\`ere sur $\xs$ constante \'egale \`a $1$ 
(identifi\'ee \`a la section rationnelle canonique du fibr\'e $\ecO_{\xs}(\varphi)$).

La m\'etrique $(\norm{.}^{\varphi}_v)_{v\in \placesde{K}}$ est une m\'etrique ad\'elique
sur le fibr\'e en droites $\ecO_{\xs}(\varphi)$.
\end{enumerate}
\end{lemme}
\begin{demo}
(cf. \'egalement \cite[Thm 2.1.6]{BaTs:anis})
Les points \ref{item:1:lm:metrique} et 
\ref{item:2:lm:metrique}
 d\'ecoulent imm\'ediatemment de la d\'efinition.

Pour le point \ref{item:3:lm:metrique}, 
supposons d'abord $v$ finie.
D'apr\`es la commutativit\'e du diagramme \eqref{eq:diag:degtlvtkv},
on peut \'ecrire,
pour $t\in T(K_v)$ et $\varphi\in \PL(\Sigma)^{G_v}$,  
\begin{align}
\log(q_{\V})\, \accs{\varphi}{\deg_{T_L,\V}(t)}
&
=
\log(q_{\V})\, \accs{\varphi}{\deg_{T,L,\V}(t)}
\\
&
=
\log\left(q_v^{\frac{[L_\V:K_v]}{e_v}}\right) \accs{\varphi}{\deg_{T,L,\V}(t)}
\\
&
=
[L_{\V}:K_v]\,\log(q_v)\,\accsv{\varphi}{t}.
\end{align}
d'o\`u le r\'esultat.

Supposons \`a pr\'esent $v$ archim\'edienne.  
D'apr\`es la commutativit\'e du diagramme 
\ref{eq:diag:degtlvtkv:archi},
on peut \'ecrire
pour $t\in T(K_v)$ et $\varphi\in \PL(\Sigma)^{G_v}$,  
\begin{align}
\accs{\varphi}{\deg_{T_L,\V}(t)}
&
=
\accs{\varphi}{\deg_{T,L,\V}(t)}
\\
&
=
[L_{\V}:K_v] \,\frac{1}{[L_{\V}:K_v]}\,\accs{\varphi}{\deg_{T,L,\V}(t)}
\\
&
=
[L_{\V}:K_v]\,\accsv{\varphi}{t}
\end{align}
d'o\`u le r\'esultat.

Passons \`a la d\'emonstration du point \ref{item:4:lm:metrique}.

L'unicit\'e de la m\'etrique vient de la densit\'e de $T(K_v)$ 
dans $\xs(K_v)$
et de la continuit\'e de l'application $x\mapsto \norm{s(x)}_{v,x}$
pour toute m\'etrique $v$-adique $\norm{.}_v$ sur $\ecO_{\xs}(\varphi)$ 
et toute section locale 
$s$ de $\ecO_{\xs}(\varphi)$.

D'apr\`es le lemme \ref{lm:restr:metrique}
et la relation \eqref{eq:comp:Hphi:ext}, il suffit de montrer l'existence
de la m\'etrique dans le cas o\`u le tore $T$ est d\'eploy\'e. 

Comme tout fibr\'e en droites sur $\xs$ est le quotient de deux fibr\'es 
en droites engendr\'es par leurs sections, il suffit d'apr\`es le lemme
\ref{lm:tensomet} et \eqref{eq:Hphi1pphi2}
de montrer l'existence de la m\'etrique pour les $\varphi$ tels que
$\ecO_{\xs}(\varphi)$ est engendr\'e par ses sections.

On note alors $(\norm{.}^{\text{st}}_v)_{v\in \placesde{K}}$ 
la m\'etrique ad\'elique standard
sur $\ecO_{\xs}(\varphi)$ (cf. d\'efinition \ref{def:met:adel:st}) 
et on va montrer 
\begin{equation}
\forall v\in \placesde{K},\,
\quad\forall t\in T(K_v),\,\quad H_{v}(\varphi,t)=\left(\norm{\ind(t)}_v^{\text{st}}\right)^{-1}.
\end{equation}

Soit $\{m_0,\dots,m_r\}$ une base de $H^0(\xs,\ecO(\varphi))$
constitu\'ee de caract\`eres de $T$ (cf. \cite[Lemma, p. 66]{Ful:toric}).

Pour $n\in X(T)^{\vee}$, $\ecO_{\xs}(\varphi)$ \'etant engendr\'e par ses sections, 
on a alors d'apr\`es \cite[Proposition, p. 68]{Ful:toric}
\begin{equation}
\varphi(n)=\Maxu{i=0,\dots, r}{\acc{m_i}{n}}.
\end{equation}

Soit d'abord $v$ une place finie de $K$.
Pour $t\in T(K_v)$, on a 
\begin{equation}
\accs{\varphi}{\deg_{T,K,v}(t)}=\Maxu{i=0,\dots, r} \acc{m_i}{\deg_{T,K,v}(t)}
=
\Maxu{i=0,\dots, r} v(m_i(t))
\end{equation}
soit
\begin{equation}
H_v(\varphi,x)
=
q_v^{
\,\accs{\varphi}{\deg_{T,K,v}(t)}}
=\Maxu{i=0,\dots,r} q_v^{\,v(m_i(t))}
=
\Maxu{i=0,\dots,r} \abs{m_i(t)}_v.
\end{equation}

Soit \`a pr\'esent $v$ une place archim\'edienne.
Pour $t\in T(K_v)$, on a 
\begin{equation}
\varphi(\deg_{T,K,v}(t))
=\Maxu{i=0,\dots, r} \acc{m_i}{\deg_{T,K,v}(t)}
=\Maxu{i=0,\dots, r} \log \abs{m_i(t)}_{v}
\end{equation}
soit
\begin{equation}
H_v(\varphi,t)
=
\exp(\accs{\varphi}{\deg_{T,K,v}(t)})
=
\Maxu{i=0,\dots,r} \abs{m_i(t)}_{v}.
\end{equation}

Mais par d\'efinition de la m\'etrique ad\'elique standard, 
on a dans tous les cas 
\begin{equation}
\left(\norm{\ind(t)}_v^{\text{st}}\right)^{-1}=\Maxu{i=0,\dots,r} \abs{m_i(t)}_v
\end{equation}
d'o\`u le r\'esultat.\end{demo}
\begin{rem}\label{rem:metrique:canonique}
Pour tout entier $d\geq 2$, notons $[d]$ l'endomorphisme de $\xs$ induit par l'endomorphisme $t\mapsto t^d$
de $T$. Pour tout $\varphi\in PL(\Sigma)^{\,G}$, on a alors un isomorphisme canonique
\begin{equation}\label{eq:dastphiisomOphid}
[d]^{\ast} \ecO_{\xs}(\varphi)\isom \ecO_{\xs}(\varphi)^d.
\end{equation}
Par ailleurs toute m\'etrique $v$-adique sur ${\varphi}_v$ induit naturellement 
des m\'etriques sur $[d]^{\ast} \ecO_{\xs}(\varphi)$ et 
$\ecO_{\xs}(\varphi)^d$ respectivement.
La m\'etrique $v$-adique $\norm{.}^{\varphi}_v$ d\'efinie ci-dessus
est alors la m\'etrique $v$-adique canonique sur $\ecO_{\xs}(\varphi)$, 
au sens o\`u c'est l'unique m\'etrique $v$-adique sur $\ecO_{\xs}(\varphi)$
tel que l'isomorphisme \ref{eq:dastphiisomOphid} soit une isom\'etrie.
\end{rem}
\subsubsection{Remarques sur le cas fonctionnel}\label{subsec:rem:casfonc}

Dans le cas fonctionnel, une d\'efinition alternative des
hauteurs locales nous sera utile.
Soit $v$ une place de $K$. L'accouplement $H_v$ se r\'e\'ecrit
\begin{equation}\label{eq:hauteur:locale:bis}
H_{v}\,:
\begin{array}{rcl}
\PL(\Sigma)^{\,G_v}_{\C}\times\,T(K_v)&\longto& \C\\
(\varphi,t)&\longmapsto&q^{\,f_v\,\accsv{\varphi}{t}}.
\end{array}
\end{equation}
Supposons \`a pr\'esent que $v$ v\'erifie l'hypoth\`ese suivante  :
\begin{hyp}\label{hyp:hinz}
Pour tout $\varphi\in \PL(\Sigma)^{\,G}$
et pour tout $t\in T(K_v)$, $f_v\,\accsv{\varphi}{t}$ est entier 
\end{hyp}
Pour $t\in T(K_v)$ on peut alors \'etendre par lin\'earit\'e la fonction
\begin{equation}
\map{\plsg}{\Z}{\varphi}{f_v\,\accsv{\varphi}{t}}
\end{equation}
en une application
\begin{equation}
\map{\plsgca}{\Z\otimes \ca=\ca}{\psi}{\accsv{\psi}{t}^{\,f_v}}.
\end{equation}
On obtient finalement un accouplement
\nindex{$\frH_v$}
\begin{equation}\label{eq:hauteur:locale:log:bis}
\frH_v\,:\,
\map{\plsgca\times T(K_v)}{\ca}{(\psi,t)}
{\accsv{\psi}{t}^{\,f_v}}.
\end{equation}
Par ailleurs le morphisme du groupe $\C$ vers le groupe $\C^{\times}$ 
donn\'e par $s\mapsto q^s$
induit un morphisme 
\begin{equation}
\map{\PL(\Sigma)^{\,G}_{\C}}{\PL(\Sigma)^{\,G}_{\C^{\times}}}{\varphi}{q^{\,\varphi}}
\end{equation}
et on a 
\begin{equation}\label{eq:hvphi=frhvqphi}
\forall \varphi\in \PL(\Sigma)^{\,G}_{\C},
\quad
\forall t\in T(K_v),
\quad
H_v(\varphi,t)=\frH_v\left(q^{\,\varphi},t\right).
\end{equation}

L'hypoth\`ese \ref{hyp:hinz} est v\'erifi\'ee par exemple si $v$
n'est pas ramifi\'ee.
En g\'en\'eral, l'hypoth\`ese \ref{hyp:hinz} n'est pas toujours v\'erifi\'ee, 
comme le montrent le lemme \ref{lm:1:ex:hnotinz} ci-dessous, combin\'e 
\`a la proposition \ref{prop:ctsu}.
Par la suite, pour l'\'etude de la fonction z\^eta des 
hauteurs, on se placera tr\`es souvent, par souci de simplification, 
dans le cas o\`u l'hypoth\`ese \ref{hyp:hinz} est v\'erifi\'ee pour toute place $v$ (ce qui est le cas
par exemple si $T$ est d\'eploy\'e par une extension non ramifi\'ee). 
On expliquera en appendice les modifications techniques n\'ecessaires
pour adapter le calcul au cas g\'en\'eral.

\begin{lemme}\label{lm:1:ex:hnotinz}
Soit $T$ un tore alg\'ebrique d\'efini sur un corps de fonctions $K$
et $L$ une extension galoisienne de $K$ de groupe $G$ d\'eployant $T$.
Soit $v$ une place de $K$ v\'erifiant les conditions suivantes :
\begin{enumerate}
\item $f_v=1$ ;
\item $v$ est ramifi\'ee dans $L$ ;
\item $G_v=G$ ;
\item $\left(X(T)^{\vee}\right)^{G_v}$ est de rang 1 ;
\item $\deg_{T,L,\V}\,:\,T(K_v)\to \left(X(T)^{\vee}\right)^{G_v}$ est surjective.
\end{enumerate}
Alors il existe un $G$-\'eventail $\Sigma$ projectif et lisse de $X(T)^{\vee}$
tel que $v$ ne v\'erifie pas l'hypoth\`ese \ref{hyp:hinz}.
\end{lemme}
\begin{demo}
Soit $m$ un g\'en\'erateur de $\left(X(T)^{\vee}\right)^{G_v}$.
D'apr\`es \cite{CTHaSk:comp}, on peut construire un $G$-\'eventail $\Sigma$
projectif et lisse dont $\R_{\geq 0}\,m$ est l'un des rayons. 
Comme $\R_{\geq 0}\,m \cap X(T)^{\vee}=m$ et que $m$ est $G$-invariant, l'\'el\'ement
$\varphi$ de $\pls$ qui envoie $m$ sur $1$ et tous les autres rayons sur $0$ est dans
$\plsg$. Soit $t$ un \'el\'ement de $T(K_v)$ v\'erifiant $\deg_{T,L,\V}(t)=m$.
Alors on a $\varphi(\deg_{T,L,\V}(t))=\varphi(m)=1$.
Ainsi 
\begin{equation}
f_v\,\accsv{\varphi}{t}=\frac{1}{e_v}\,\varphi(\deg_{T,L,\V}(t))=\frac{1}{e_v}.
\end{equation}
\end{demo}

\subsubsection{Hauteurs globales et fonction z\^eta des hauteurs}
Pour $t\in T(K)$ 
et $\varphi\in \PL(\Sigma)^G_{\C}$
on pose
\begin{equation}
H(\varphi,t)
=
\prod_{v\in \placesde{K}}\,H_{v}(\varphi,t).
\end{equation}
Si $\varphi$ est un \'el\'ement de $PL(\Sigma)^{\,G}$, 
d'apr\`es le point \ref{item:4:lm:metrique} du lemme \ref{lm:metrique} et la d\'efinition \eqref{eq:def:height}
$H(\varphi,.)$ est la restriction \`a $T(K)$ d'une hauteur d'Arakelov 
associ\'ee au faisceau inversible $\ecO_{\xs}(-\varphi)$.

On pose alors 
\begin{equation}
\zeta_{H(\varphi,\,.\,)}(s)
=
\sum_{t\in T(K)}\,H(\varphi,t)^{-s}
\end{equation}
pour tout $s\in \C$ tel que la s\'erie converge, 

Plus g\'en\'eralement, on pose
\nindex{$\zeta_H$}
\begin{equation}
\zeta_H(\varphi)=\sum_{t\in T(K)}\,H(-\varphi,t)
\end{equation}
pour tout $\varphi\in PL(\Sigma)^{\,G}_{\C}$ tel que la s\'erie converge.

\subsection{Mesure et nombre de Tamagawa d'une vari\'et\'e torique}

On note \nindex{$\varphi_0$}$\varphi_0$ l'\'el\'ement de $PL(\Sigma)^{\,G}$ 
correspondant au diviseur $\sum_{\alpha\in \sg} \Da$.
D'apr\`es la proposition \ref{prop:antican},
$\ecO_{\xs}(\varphi_0)$
est isomorphe
au faisceau anticanonique de $\xs$. 
Ainsi la m\'etrique ad\'elique $(\norm{.}_v^{\varphi_0})$ dont provient la hauteur $H(\varphi_0,\,.\,)$ 
induit pour tout $v$ une mesure $\omega_{\xs,v}$ sur $\xs(K_v)$.

La restriction \`a $T$ de l'inverse de la section rationnelle canonique
de  $\ecO_{\xs}(\varphi_0)$ est  
une $d$-forme diff\'erentielle $K$-rationnelle $T$-invariante sur $T$.
Cette restriction s'\'ecrit donc $\alpha\,\Omega_T$, o\`u $\alpha\in K^{\ast}$
et $\Omega_T$ est la forme utilis\'ee dans la section \ref{sec:tamagawa:tore} 
pour la construction des mesures de Haar $\omega_{T,v}$ sur $T(K_v)$.

D'apr\`es le point \ref{item:4:lm:metrique} du lemme \ref{lm:metrique}, on a donc pour tout $v\in \placesde{K}$
\begin{equation}\label{eq:normeqH}
\left(\norm{\Omega_T(\,.\,)^{-1}}^{\varphi_0}_v\right)^{-1}=\abs{\alpha}_v\,H_{v}(\varphi_0,\,.\,).
\end{equation}
En vue du calcul du terme principal de la fonction z\^eta des hauteurs, 
il faut comparer les mesures $\omega_{\xs,v}$ 
et 
$H_{v}(\varphi_0,\,.\,)\,\omega_{T,v}$ sur $T(K_v)$. 
C'est l'objet de la proposition 3.4.4. de \cite{BaTs:anis}, dont nous rappelons l'\'enonc\'e et la d\'emonstration.
\begin{lemme}[Batyrev,Tschinkel]
\label{lm:restriction:mesure}
Il existe une famille $(\alpha_v)\in \left(\R_{>0}\right)^{\placesde{K}}$ avec $\alpha_v=1$ pour presque tout $v$ et 
\begin{equation}
\prod_{v\in \placesde{K}} \alpha_v=1
\end{equation}
telle que
\begin{equation}
\forall\,v\in \placesde{K},\quad\omega_{\xs,v}|_{T(K_v)}
=
\alpha_v\,H_{v}(\varphi_0,.)\,\omega_{T,v}.
\end{equation}
\end{lemme}
\begin{demo}
Soient $t\in T(K_v)$ et $U$ un ouvert de Zariski de $T$ contenant $t$ 
tel qu'il existe un $K$-morphisme \'etale
\begin{equation}
f\,:\,U\longto \A^d_K,
\end{equation}
v\'erifiant (en notant $(y_1,\dots,y_d)$ les coordonn\'ees sur $\A^d_K$)
\begin{equation}
f^{\ast}
\left(
dy_1\wedge \dots \wedge dy_d
\right)
=
\beta^{-1}\,\Omega_{U},
\end{equation}
o\`u $\beta$ est une fonction r\'eguli\`ere inversible sur $U$. 
Ce morphisme induit pour un ouvert analytique $W$ de $U(K_v)$ contenant $t$
et assez petit 
un isomorphisme analytique
\begin{equation}
f\,:\,W\isom f(W)\subset K_v^d.
\end{equation}
Sur $W$, la mesure $\omega_{T,v}$ s'exprime alors comme
\begin{equation}
\omega_{T,v}=|\beta|_{v}\,\left(f^{-1}\right)_{\ast}(dy_1\dots dy_d).
\end{equation}
Par ailleurs, sur $W$ on a
\begin{equation}
\omega_{\xs,v}=|\beta|_{v}\left(f^{-1}\right)_{\ast}\left(\norm{\Omega_T(f^{-1}(y))^{-1}}^{\varphi_0}_v 
\,\,dy_1\dots dy_d\right)
\end{equation}
soit d'apr\`es \eqref{eq:normeqH}
\begin{equation}
\omega_{\xs,v}=\abs{\beta}_{v}
\left(f^{-1}\right)_{\ast}
\left(\abs{\alpha}_v\,H_v(\varphi_0,f^{-1}(y)) dy_1\dots dy_d\right)
\end{equation}
d'o\`u le r\'esultat.
\end{demo}

Rappelons (cf. section \ref{subsec:def:height}) 
que le \termin{nombre de Tamagawa} de $\xs$ est alors par d\'efinition
\begin{equation}\label{eq:defgammahxs}
\gamma_H\left(\xs\right)=\omega_{\xs}\left(\overline{\xs(K)}\right),
\end{equation}
o\`u
\begin{equation}\label{eq:defomegaxs}
\omega_{\xs}
=
c_{K,\dim(\xs)}\,
\ell(\Pic(\xsl))
\prod_{v\in \placesde{K}}L_v(1,\Pic(\xsl))^{-1}\,\omega_{\xs,v}.
\end{equation}

\begin{lemme}[Batyrev, Tschinkel]
\label{lm:int_T=int_X}
On a 
\begin{equation}
\int\limits_{\adh{T(K)}\cap T(\ak)}\!\!\!\!\!\omega_{\xs}=\int\limits_{\overline{\xs(K)}}\!\omega_{\xs}.
\end{equation}
\end{lemme}
\begin{demo}
C'est la proposition 3.4.5 de \cite{BaTs:anis}, 
\'enonc\'ee uniquement dans le cas arithm\'etique,
mais la preuve marche aussi dans le cas fonctionnel.
Nous la rappelons.

Soit $\adh{\xs(K)}^{\,\,S}$ l'adh\'erence de l'image de $\xs(K)$ 
dans $\produ{v\in S}\xs(K_x)$.
Rappelons que nous avions not\'e $\adh{T(K)}^{\,\,S}$ l'adh\'erence de l'image de $T(K)$ 
dans $\produ{v\in S}T(K_v)$. Comme $\produ{v\in S}T(K_v)$ est ouvert dans 
 $\produ{v\in S}\xs(K_v)$, on a 
\begin{equation}
\adh{\xs(K)}^{\,\,S} \cap \prod_{v\in S} T(K_v)
=
\adh{\left(\prod_{v\in S} T(K_v)\right)\cap \xs(K)}
=
\adh{T(K)}^{\,\,S}
\end{equation}

On a \'evidemment
\begin{equation}
\adh{\xs(K)} \subset \adh{\xs(K)}^{\,\,S} \times \produ{v\notin S} \xs(K_v).
\end{equation}
Montrons l'inclusion inverse. 
Soit  $W_S$ un ouvert de $\adh{\xs(K)}^{\,\,S}$, $T$ un sous-ensemble fini
de $\placesde{K}\setminus S$ et pour $v\in T$, $U_v$  un ouvert de $\xs(K_v)$.
Il suffit de montrer 
que l'ouvert de $\adh{\xs(K)}^{\,\,S} \times \produ{v\notin S} \xs(K_v)$ d\'efini par
\begin{equation}
W\eqdef W_S\times \produ{v\notin\,S\,\cup\, T } \xs(K_v)\times \produ{v\in T}U_v
\end{equation}
rencontre $\xs(K)$.

Or, pour tout $v$, d'apr\`es \cite[Lemma 3.2]{PlRa:alggp}, $T(K_v)$ est dense dans $\xs(K_v)$.
On en d\'eduit que pour $v\in T$,  $U_v$ rencontre $T(K_v)$.
 et que $W_S$ rencontre $\prod_{v\in S} T(K_v)$.
Comme $W_S$ est inclus dans $\adh{\xs(K)}^{\,\,S}$, $W_S$ rencontre  
en fait  
\begin{equation}
\adh{\xs(K)}^{\,\,S} \cap \prod_{v\in S} T(K_v)=\adh{T(K)}^{\,\,S}.
\end{equation}

Ainsi $W$ rencontre $\adh{T(K)}^{\,\,S}\times \produ{v\notin S} T(K_v)$. 
Comme ce dernier ensemble est inclus dans $\adh{\xs(K)}^{\,\,S} \times \produ{v\notin S} \xs(K_v)$,
$W\cap \adh{T(K)}^{\,\,S}\times \produ{v\notin S} T(K_v)$
est ouvert dans $\adh{T(K)}^{\,\,S}\times \produ{v\notin S} T(K_v)$.
D'apr\`es \eqref{eq:scindage1}, ce dernier ensemble co\"\i ncide avec $\adh{T(K)}$.
On en d\'eduit que $W$ rencontre $T(K)$.

On a donc montr\'e
\begin{equation}
\adh{\xs(K)}=\adh{\xs(K)}^{\,\,S} \times \produ{v\notin S} \xs(K_v).
\end{equation}

Ainsi 
\begin{equation}
\int\limits_{\overline{\xs(K)}}\omega_{\xs}
=\int\limits_{\adh{\xs(K)}^{\,\,S}} \otimesu{v\in S} \omega_{v,\xs}
\times \prod_{v\notin S} \omega_{v,\xs}(\xs(K_v)).
\end{equation}

Comme $\xs\setminus T$ est un ferm\'e alg\'ebrique propre de $\xs$,
l'ensemble $(\xs\setminus T)(K_v)$ est de mesure nulle pour 
$\omega_{\xs,v}$ d'apr\`es 
\cite[10.1.3 Exemple b) et 10.1.6]{Bki:diff}.
On a donc 
\begin{equation}
\int\limits_{\overline{\xs(K)}}\omega_{\xs}
=\int\limits_{\adh{\xs(K)}^{\,\,S}} \otimesu{v\in S} \omega_{v,\xs}
\times \prod_{v\notin S} \omega_{v,\xs}(T(K_v)).
\end{equation}
Par ailleurs on a 
\begin{equation}
\int\limits_{\adh{T(K)}\cap T(\ak)}\!\!\!\!\omega_{\xs}
=\int\limits_{\adh{T(K)}^{\,\,S}} \otimesu{v\in S} \omega_{v,\xs}
\times \prod_{v\notin S} \omega_{v,\xs}(T(K_v))
\end{equation}
et il suffit pour conclure de montrer que 
$\adh{\xs(K)}^{\,\,S}\setminus \adh{T(K)}^{\,\,S}$
est de mesure nulle pour $\otimesu{v\in S} \omega_{v,\xs}$.

Comme $\adh{\xs(K)}^{\,\,S}\cap \produ{v\in S}T(K_v)=\adh{T(K)}^{\,\,S}$,
on a 
\begin{equation}\label{eq:adhsubprod}
\left(\adh{\xs(K)}^{\,\,S}\setminus \adh{T(K)}^{\,\,S}\right)
\subset \left(\produ{v\in S}\xs(K_v)\setminus \produ{v\in S}T(K_v)\right)
\end{equation}
Pour tout $v$, on a $\omega_{\xs,v}(\xs(K_v)\setminus T(K_v))=0$,  donc
$\produ{v\in S}\xs(K_v)\setminus \produ{v\in S}T(K_v)$ est de mesure nulle pour $\otimesu{v\in S}\omega_{\xs,v}$,
d'o\`u le r\'esultat d'apr\`es \eqref{eq:adhsubprod}.

\end{demo}

\subsection{Le r\'esultat}\label{subsec:resultat}

Nous sommes \`a pr\'esent en mesure d'\'enoncer le r\'esultat 
obtenu par Batyrev et Tschinkel dans \cite{BaTs:manconj},
ainsi que notre r\'esultat qui en est un analogue dans le cas fonctionnel.

\begin{thm}[Batyrev, Tschinkel]\label{thm:theo:ba:ts}
Soit $K$ un corps de nombres, 
$L/K$ une extension finie galoisienne de groupe $G$, 
et $\Sigma$ un $G$-\'eventail projectif et lisse. 
Soit $\xs$ la vari\'et\'e torique projective et lisse, 
d\'efinie sur $K$, 
qui lui est associ\'ee. 
C'est une compactification d'un tore alg\'ebrique $T$.

Alors la s\'erie $\zeta_H\left(s\,\varphi_0\right)$ converge absolument pour $s\in \tube{\R_{>1}}$
et, pour un certain $\varepsilon>0$, 
la fonction
\begin{equation}
f\,:\,s\longmapsto (s-1)^{\rg(\Pic(\xs))}\,\zeta_H\left(s\,\varphi_0\right)
\end{equation}
se prolonge en une fonction holomorphe 
sur le domaine $\tube{\R_{>1-\varepsilon}}$,
v\'erifiant
\begin{equation}
f(1)
=
\alpha^{\ast}(\xs)\,\card{H^1\left(G,\Pic(\xsl)\right)}\,\gamma_H\left(\xs\right).
\end{equation}
\end{thm}

Nous obtenons pour notre part le th\'eor\`eme suivant.

\begin{thm}\label{thm:theoprinc}
Soit $K$ un corps de fonctions, 
$L/K$ une extension finie galoisienne de groupe $G$, 
et $\Sigma$ un $G$-\'eventail projectif et lisse. 
Soit $\xs$ la vari\'et\'e torique projective et lisse, 
d\'efinie sur $K$, 
qui lui est associ\'ee. 
C'est une compactification d'un tore alg\'ebrique $T$.

Alors la s\'erie $\zeta_H\left(s\,\varphi_0\right)$ converge absolument pour $s\in \tube{\R_{>1}}$
et pour un certain $\varepsilon>0$ se prolonge en une fonction m\'eromorphe 
sur $\tube{\R_{>1-\varepsilon}}$. 
Ce prolongement a un p\^ole d'ordre le rang du groupe de Picard de $\xs$ en $s=1$, et on a 
\begin{equation}
\lim_{s\to 1}(s-1)^{\rg(\Pic(\xs))}\,
\zeta_H\left(s\,\phi_0\right)
=
\alpha^{\ast}(\xs)\,\card{H^1\left(G,\Pic(\xsl)\right)}\,\gamma_H\left(\xs\right).
\end{equation}
\end{thm}

\begin{rems}
\begin{enumerate}
\item
Notre preuve, comme d\'ej\`a indiqu\'e, 
est fortement inspir\'ee des preuves de \cite{BaTs:anis} et \cite{BaTs:manconj}.
\item
Dans le cas arithm\'etique comme dans le cas fonctionnel, l'ensemble $\xs(K)$
est bien Zariski dense dans $\xs$ d'apr\`es \cite[Corollary 7.12]{BoSp:rationality}.
\item
Contrairement au cas arithm\'etique, on ne peut esp\'erer dans le cas fonctionnel prolonger la fonction
$s\mapsto (s-1)^{\rg(\Pic(\xs))}\,\zeta_H\left(s\,\phi_0\right)$ en une fonction holomorphe
sur un domaine du type $\tube{\R_{>1-\varepsilon}}$. En effet dans ce cas la fonction
z\^eta des hauteurs a d'autres p\^oles sur la droite $\{\Re(s)=1\}$, ne serait-ce que ceux d\^u au fait
qu'elle admet des p\'eriodes imaginaires pures. Par exemple, au moins dans le cas o\`u
la vari\'et\'e $\xs$ v\'erifie l'hypoth\`ese \ref{hyp:hinz} pour toute place $v$, la fonction $\zeta_H\left(s\,\phi_0\right)$
est $\frac{2\,i\,\pi}{\log(q)}$-p\'eriodique.
\end{enumerate}
\end{rems}

\subsection{Strat\'egie de Batyrev et Tschinkel}
\label{subsec:formule:poisson}

Une des id\'ees essentielles de Batyrev et Tschinkel pour \'etudier la fonction z\^eta des hauteurs
des vari\'et\'es toriques 
est d'appliquer la formule de Poisson pour obtenir une repr\'esentation int\'egrale de $\zeta_H(\varphi)$. 
Avant de pr\'eciser ce qui pr\'ec\`ede,
rappelons quelques faits \'el\'ementaires d'analyse harmonique abstraite. 

\subsubsection{Un peu d'analyse harmonique}
On note \nindex{$\unit$}$\unit$ le groupe des nombres complexes de module $1$.
Soit $\ecG$ un groupe ab\'elien localement compact.
Son \termin{dual topologique} est le groupe 
topologique, not\'e $\dualt{\ecG}$,
form\'e de l'ensemble des morphismes continus de 
$\ecG$ dans $\unit$.
C'est encore un groupe ab\'elien localement compact.

\begin{exs}\label{conv:dualtopMR}
Soit $M$ un $\Z$-module libre de rang fini.
\begin{enumerate}
\item
Le morphisme qui 
\`a $m\in M^{\vee}_{\R}$
associe le caract\`ere 
\begin{equation}
m\mapsto\exp\left(i\acc{m^{\vee}}{m}\right).
\end{equation}
est un isomorphisme du groupe topologique $M^{\vee}_{\R}$
sur le groupe topologique $\dualt{M_{\R}}$. 
Par la suite, on identifiera toujours
$\dualt{M_{\R}}$ \`a 
$M^{\vee}_{\R}$ au moyen de cet isomorphisme.
\item
Le morphisme qui \`a $m^{\vee}\otimes z\in M^{\vee}_{\unit}$
associe
\begin{equation}
m \mapsto z^{\acc{m^{\vee}}{m}}
\end{equation}
est un isomorphisme du groupe topologique $M^{\vee}_{\unit}$
sur le dual topologique de $M$.
Par la suite, on identifiera toujours
$\dualt{M}$ \`a 
$M^{\vee}_{\unit}$ au moyen de cet isomorphisme.
\end{enumerate}
\end{exs}

Soit $\ecG$ un groupe ab\'elien localement compact
muni d'une mesure de Haar $dg$.
Soit $F\,:\,\ecG\to \C$ une fonction de classe $\L^1$.
Sa \termin{transform\'ee de Fourier par rapport \`a $dg$} est la fonction 
$\fourier F\,:\,\dualt{\ecG}\to \C$ d\'efinie par 
\begin{equation}
\forall \chi\in \dualt{\ecG},\quad
\fourier F (\chi)=\int\limits_{\ecG} F(g)\,\chi(g)\,dg.
\end{equation}
Il existe alors (cf. \cite[D\'efinition 4, p.118 et Th\'eor\`eme 3, p. 123]{Bki:spec})
une unique mesure de Haar $dg^{\ast}$
sur $\dualt{\ecG}$ v\'erifiant la formule d'inversion de Fourier,
c'est-\`a-dire la propri\'et\'e suivante :
soit $F\,:\,\ecG\to \C$ une fonction de classe $\L^1$ 
telle que  $\fourier F$ est de classe $\L^1$ sur $\dualt{\ecG}$ ;
alors, pour presque tout $g$ de $\ecG$,  
on a la formule
\begin{equation}\label{eq:inv:fourier}
F(g)
=
\int\limits_{\dualt{\ecG}}\overline{\chi(g)}\,(\fourier F)(\chi)\dualt{dg}(\chi).
\end{equation}
La mesure de Haar $\dualt{dg}$ sera appel\'ee \termin{mesure duale} de la mesure de Haar $dg$.
Moyennant l'identification canonique de $\dualtop{\dualt{\ecG}}$ \`a $\ecG$, on a $\dualtop{\dualt{dg}}=dg$.

Le lemme suivant donne deux exemples standards de mesure duale.
\begin{lemme}\label{lm:ex:mesure:duale}
Soit $\ecG$ un groupe topologique ab\'elien localement compact,
et $dg$ une mesure de Haar sur $\ecG$.
\begin{enumerate}
\item
On suppose que $\ecG$ est compact, de sorte que $\dualt{\ecG}$
est discret.  Alors $\dualt{dg}$
est la mesure de Haar sur $\dualt{\ecG}$ pour laquelle chaque 
point a pour masse $\frac{1}{\int\limits_{\ecG} dg}$.
\item
On suppose que $\ecG=N_{\R}$ o\`u $N$ est un $\Z$-module libre de rang fini et
que $dg$ est la mesure de Lebesgue sur $N_{\R}$ normalis\'ee par le r\'eseau $N$.
Alors $\dualt{dg}$ est la mesure de Lebesgue sur $N^{\vee}_{\R}$
normalis\'ee par le r\'eseau $N^{\vee}$, divis\'ee par $(2\,\pi)^{\rg(N)}$.
\end{enumerate}
\end{lemme}
\begin{demo}
On applique \eqref{eq:inv:fourier} en prenant pour $F$ l'indicatrice de $G$ dans le premier
cas, et l'application $(x_1,\dots,x_n)\mapsto \exp\left(-\pi(x_1^2+\dots+x_n^2)\right)$ dans le second cas.
\end{demo}

On peut \`a pr\'esent rappeler un \'enonc\'e de la formule de Poisson (cf. \cite[Proposition 8, p. 127]{Bki:spec}).
\begin{thm}[Formule de Poisson]\label{thm:formule:poisson}
Soit $\ecG$ un groupe ab\'elien localement compact 
muni d'une mesure de Haar $dg$,
et $\ecH$ un sous-groupe ferm\'e de $\ecG$ muni d'une mesure de Haar $dh$. 
Soit $dx$ la mesure de Haar sur le quotient $\ecG/\ecH$ 
normalis\'ee par la relation $dg=dx\,dh$. 
On munit le groupe $\dualtop{\ecG/\ecH}$ de la mesure de Haar $\dualt{dx}$ duale de la mesure $dx$,
et on identifie ce groupe au sous-groupe de $\dualt{\ecG}$
constitu\'e des caract\`eres de $\ecG$ triviaux sur $\ecH$.

Soit $F\,:\,\ecG\to \C$ une fonction de classe $\L^1$ et $\fourier F$ sa transform\'ee de Fourier par rapport \`a $dg$. 
On suppose que $\fourier F$ est $\L^1$ sur $\dualtop{\ecG/\ecH}$.

Alors, pour presque tout $g$ de $\ecG$, la propri\'et\'e suivante est v\'erifi\'ee :
la fonction $h\mapsto F(g\,h)$ est $\L^1$ sur $\ecH$ et on 
a la formule
\begin{equation}
\int\limits_{\ecH}F(g\,h)dh
=
\int\limits_{\dualtop{\ecG/\ecH}}\overline{\chi(g)}\,\fourier F(\chi)\dualt{dx}(\chi).
\end{equation}
\end{thm}
\begin{cor}\label{cor:formule:poisson}
On conserve les notations et hypoth\`eses du th\'eor\`eme 
\ref{thm:formule:poisson}. On se donne en outre un sous-groupe compact $\ecK$
de $\ecG$, de volume non nul, et on suppose que la fonction $F$ est $\ecK$-invariante.
Alors la fonction $F$ est $\L^1$ sur $\ecH$ et on 
a la formule
\begin{equation}
\int\limits_{\ecH}F(h)dh
=
\int\limits_{\dualtop{\ecG/\ecK.\ecH}}\fourier F(\chi)\dualt{dx}(\chi).
\end{equation}
\end{cor}
\begin{rem}
Comme $\ecK$ est compact, $\dualtop{\ecG/\ecK.\ecH}$
est un sous-groupe ouvert de $\dualtop{\ecG/\ecH}$, et la restriction de 
$\dualt{dx}$ \`a $\dualtop{\ecG/\ecK.\ecH}$ est une mesure de Haar sur $\dualtop{\ecG/\ecK.\ecH}$.

Comme $F$ est $\ecK$-invariante, $\fourier F(\chi)$
est nulle si $\chi$ n'est pas trivial sur $\ecK$. Ainsi l'hypoth\`ese
que $\fourier F$ est $\L^1$ sur $\dualtop{\ecG/\ecH}$ \'equivaut \`a l'hypoth\`ese
que $\fourier F$ est $\L^1$ sur $\dualtop{\ecG/\ecK.\ecH}$
\end{rem}
\begin{demo}
D'apr\`es le th\'eor\`eme \ref{thm:formule:poisson}, comme $\ecK$
est de volume non nul, il existe un \'el\'ement $k$ de $\ecK$ tel
que $h\mapsto F(k.h)$ est de classe $\L^1$ sur $\ecH$ et tel
qu'on ait la formule
\begin{equation}
\int\limits_{\ecH}F(k.h)dh
=
\int\limits_{\dualtop{\ecG/\ecH}}\overline{\chi(k)}\,(\fourier F)(\chi)\dualt{dx}(\chi).
\end{equation}
Comme $F$ est $\ecK$-invariante, $\fourier F(\chi)$
est nulle si $\chi$ n'est pas trivial sur $\ecK$.
Le corollaire s'en d\'eduit aussit\^ot.
\end{demo}

\subsubsection{Application \`a la fonction z\^eta des hauteurs}

Pour appliquer ce qui pr\'ec\`ede \`a l'\'etude de la fonction z\^eta des hauteurs, 
on commence par \'etendre la hauteur 
en une fonction continue sur $T(\ak)$ en posant pour tout $(t_v)\in T(\ak)$
et tout $\varphi$ de $\PL(\Sigma)_{\C}^G$
\begin{equation}
H(\varphi,(t_v))=\prod_{v\in \placesde{K}} H_{v}(\varphi,t_v).
\end{equation}
Ainsi $H(\varphi,\,.\,)$ est invariante sous l'action de $\K(T)$ 
(ceci sera utile pour annuler beaucoup de transform\'ees de Fourier). 

Par la suite, pour tout $\varphi$ de $\PL(\Sigma)_{\C}^G$, 
nous noterons pour all\'eger l'\'ecriture 
$H(\varphi)$ la fonction
$H(\varphi,\,.\,)$.

On va appliquer le corollaire \ref{cor:formule:poisson} \`a $\ecG=T(\ak)$ 
muni de la mesure de Haar 
$\omega_{T}$ 
d\'efinie \`a la section \ref{sec:tamagawa:tore}, 
$\cH=T(K)$, 
identifi\'e \`a un sous-groupe discret de $T(\ak)$ et muni de la mesure discr\`ete,
et $\ecK=\K(T)$ (notons qu'on a bien 
$\int\limits_{\K(T)}\!\!\!\omega_{T}\neq 0$,
cf. la section \ref{sec:tamagawa:tore}).

La mesure de Haar sur $\dualtop{T(\ak)/T(K)}$ duale de la mesure quotient sur $T(\ak)/T(K)$, 
sera not\'ee $d\chi$. 

Compte tenu du fait que  $H(-\varphi)$ est $\K(T)$-invariante
pour tout \'el\'ement $\varphi$ de $\PL(\Sigma)^{\,G}$,
on d\'eduit imm\'ediatement du corollaire \ref{cor:formule:poisson}
le lemme suivant.

\begin{lemme}\label{lm:formule:poisson}
Soit $\varphi$ un \'el\'ement de $\PL(\Sigma)^{\,G}$
tel que $H(-\varphi)$ est $\L^1$ sur $T(\ak)$. 
Soit \nindex{$\fourier H(-\varphi)$}$\fourier H(-\varphi)$ 
la transform\'ee de Fourier de $H(-\varphi)$ par rapport \`a la mesure 
$\omega_{T}$.
Supposons en outre que $\fourier H(-\varphi)$  est
$\L^1$ sur $\dualtop{T(\ak)/\K(T).T(K)}$.

Alors $H(-\varphi)$ est $\L^1$ sur $T(K)$ et on a la formule
\begin{equation}\label{eq:formule:poisson:zeta}
\sum_{t\in T(K)}\,H(-\varphi,t)
=
\!\!\!\!
\int\limits_{\dualtop{T(\ak)/\K(T)\,T(K)}}\!\!\!\!\!\!\!\!\fourier\,H(-\varphi)(\chi)\,d\chi.
\end{equation}
\end{lemme}

\begin{rem}
Si $H(-\varphi)$ est int\'egrable, 
\label{rem:lm:formule:poisson}
la fonction $\fourier H(-\varphi)$ est continue
sur $\dualt{T(\ak)}$, et donc sa restriction \`a
$\dualtop{T(\ak)/T(K)}$ est int\'egrable sur tout compact.

Or, dans le cas fonctionnel,   $\dualtop{T(\ak)/\K(T).T(K)}$
est compact. Ainsi dans ce cas l'int\'egrabilit\'e de $H(-\varphi)$ entra\^\i ne automatiquement
l'int\'egrabilit\'e de  $\fourier H(-\varphi)$.
\end{rem}

Le but de la partie suivante est de montrer,
via un calcul explicite des transform\'ees de Fourier locales
en presque toutes les places et des estim\'ees ad hoc
aux places restantes, 
que l'expression sous l'int\'egrale 
dans la formule \eqref{eq:formule:poisson:zeta} 
est une fonction 
m\'eromorphe de $\varphi$ dont on contr\^ole les p\^oles.

Les lemmes techniques rappel\'es
dans la partie \ref{sec:eval:int:arit} pour le cas arithm\'etique
et d\'evelopp\'es dans la partie
\ref{sec:eval:int:fonc} dans le cas fonctionnel permettront ensuite de montrer 
que l'int\'egrale elle-m\^eme est une fonction m\'eromorphe
dont on contr\^ole le comportement analytique.

\section{Calcul des transform\'ees de Fourier 
et expression int\'egrale de la fonction z\^eta des hauteurs}

\subsection{Caract\`eres du groupe des id\`eles}\label{subsec:carac:ideles}

Soit $E$ un corps global. 

Si $v$ est une place non archim\'edienne de $E$,
tout caract\`ere de $\G_m(K_v)$ trivial sur $\G_m(\Ov)$ est de la forme
\begin{equation}
x\mapsto z^{\,v(x)}
\end{equation}
o\`u $z$ est un \'el\'ement de $\unit$.

Si $v$ est une place archim\'edienne de $E$,
tout caract\`ere de $\G_m(K_v)$ trivial sur $\G_m(\Ov)$ est de la forme
\begin{equation}
x\mapsto e^{\,i\,t\,\log\abs{x}_v}
\end{equation}
o\`u $t$ est un \'el\'ement de $\R$.

Soit \`a pr\'esent $\chi$ un caract\`ere de $\G_m(\ade{E})$.
Pour $v\in\placesde{E}$ on note $\chi_v$ le caract\`ere induit sur $\G_m(E_v)$ via
l'injection naturelle de $\G_m(E_v)$ dans $\G_m(\ade{E})$.
Si $\chi$ est trivial sur $\K(\G_m)$, il s'\'ecrit
\begin{equation}
\chi 
\,:\, 
(x_v)
\longmapsto 
\prod_{v\in\placesde{E,f}}\,z_v^{\,v(x_v)}
\prod_{v\in\placesde{E,\infty}}\,e^{i\,t_v\,\log\abs{x_v}_v}
\end{equation}
avec $(z_v)\in \unit^{\,\placesde{E,f}}$ et
$(t_v)\in \R^{\,\placesde{E,\infty}}$. 
Pour tout $v\in \placesde{E,f}$ on a donc $z_v=\chi_v(\pi_v)$. 

Si $\chi$ est de plus trivial 
sur $\G_m(\ade{E})^1$, il s'\'ecrit
$\chi' \circ \deg_E$, o\`u $\chi'$ est
un caract\`ere de $\G_m(\ade{E})/\G_m(\ade{E})^1$.
Dans le cas arithm\'etique, ce dernier groupe est isomorphe \`a $\R$ 
et il existe donc un $y\in \R$ tel que $\chi$ s'\'ecrit
\begin{equation}
x\longmapsto e^{i\,y\,\deg_E(x)}.
\end{equation}

Dans le cas fonctionnel, 
$\G_m(\ade{E})/\G_m(\ade{E})^1$ est isomorphe \`a $\Z$ 
et il existe donc un  $z\in \unit$ tel que $\chi$ s'\'ecrit
\begin{equation}
x\longmapsto z^{\deg_E(x)}.
\end{equation}

Revenons plus g\'en\'eralement \`a un caract\`ere de $\G_m(\ade{E})$
trivial sur $\K(\G_m)$.
On d\'efinit la fonction $L$ associ\'ee :
\nindex{$L_E(\chi,s)$}
\begin{equation}
L_E(\chi,s)=\prod_{v\in\placesde{E,f}}\frac{1}{1-\chi_v(\pi_v)\,q_v^{\,-s}}.
\end{equation}
Ce produit eul\'erien converge absolument pour $\Re(s)>1$.

D'apr\`es \cite[VII, $\S\,7$, Thm 6]{Wei:BNT}, 
on a
\begin{lemme}\label{lm:chi:non:triv}
Si $\chi$ n'est pas trivial sur $\G_m(\ade{E})^1$, 
la fonction $s\mapsto L_{E}(s,\chi)$ est holomorphe sur $\C$.
Dans le cas fonctionnel, 
elle s'exprime comme un polyn\^ome en $\qde{E}^{-s}$.
\end{lemme}
Dans le cas fonctionnel on notera
\nindex{$\ecL_E(\chi,\,.\,)$}$\ecL_E(\chi,\,.\,)$ le polyn\^ome v\'erifiant
\begin{equation}
\ecL_E(\chi,\qde{E}^{-s})
=
L_E(\chi,s).
\end{equation}

\subsection{Caract\`eres de $T(\ak)$}
\label{subsec:carac:T(ak)}
Soit $T$ un tore alg\'ebrique d\'efini sur un corps global $K$.
Nous reprenons les notations de la partie \ref{subsubsec:vte:tor:nondep},
le corps de base \'etant bien s\^ur le corps global $K$. 
Dans le cas fonctionnel,
le corps des constantes de $K$ est d\'esormais suppos\'e de cardinal $q$.

Rappelons en particulier 
qu'on a une suite exacte de $G$-modules
\begin{equation}\label{eq:exsqgammapi}
0
\longto X(T) 
\overset{\gamma}{\longto} 
\Ps 
{\longto}  
\Pic(\xsl)
\longto 
0
\end{equation}
et que $\Pic(\xsl)$ est un $G$-module flasque.

Pour $\chi\in T(\ak)^{\ast}$, on note $\chia$ 
le caract\`ere $\chi\circ\gamma_{\ak}\circ \iaak$, 
de sorte que 
$\chia$ est un \'el\'ement de $\G_m(\aka)^{\,\ast}$.

Le but de ce qui suit est de d\'ecrire le morphisme de caract\`eres
\begin{equation}
\dualtop{T(\ak)/T(\ak)^1}
\longto 
\dualtop{T_{\Ps}(\ak)/T_{\Ps}(\ak)^1}
\end{equation}
induit par la composition avec $\gamma_{\,\ak}$,
\`a l'aide des morphisme de caract\`eres
\begin{equation}
\map
{\dualtop{T(\ak)/T(\ak)^1}}
{\dualtop{\G_m(\aka)/\G_m(\aka)^1}}
{\chi}
{\chia}.
\end{equation}

\subsubsection{Cas arithm\'etique}\label{subsub:arit}

Commen\`c;ons par d\'ecrire les caract\`eres de $T(\ak)$ triviaux sur $T(\ak)^1$.
On a un morphisme surjectif
\begin{equation}
\deg_T\,:\,T(\ak)\longto \Hom(X(T)^G,\R)
\end{equation}
de noyau $T(\ak)^1$, qui induit donc un isomorphisme
\begin{equation}\label{eq:isotaktak1xtgrdual}
\deg_T\,:\,T(\ak)/T(\ak)^1 \longisom \Hom\left(X(T)^G,\R\right).
\end{equation}

Pour tout $y\in X(T)^G_{\R}$, on note \nindex{$\chi_y$}$\chi_y$ l'\'el\'ement de 
$\dualtop{T(\ak)/T(\ak)^1}$ qui lui correspond via l'isomorphisme dual de 
\eqref{eq:isotaktak1xtgrdual} (cf. la convention \ref{conv:dualtopMR}). 
On a donc 
\begin{equation}
\forall t\in T(\ak),\quad \chi_y(t)=\exp(\,i\,\acc{y}{\deg_T(t)}).
\end{equation}
En particulier, si $v$ est une place de $K$, on a
\begin{align}\label{eq:chiyv}
\forall t\in T(K_v),\quad (\chi_y)_v(t)&
=\exp(\,i\,\acc{y}{\deg_{T}\circ i_{T,v}(t)})
\\
&=\exp(\,i\,\log(q_v)\,\acc{y}{\deg_{T,v}(t)}).
\end{align}

On notera aussi $\ya$ le r\'eel $\acc{\roa}{y}$.

\begin{lemme}\label{lm:carac:arit}
Soit $y\in X(T)^G_{\R}$. 
Soit $\alpha\in \sg$. 
Alors pour $x\in \G_m(\aka)$ on a 
\begin{equation}
\left(\chi_y\right)_{\alpha}(x)=\exp(\,i\,\acc{\roa}{y} \,\deg_{\Ka}(x)).
\end{equation}
\end{lemme}
\begin{demo}
On a
\begin{align}
\left(\chi_y\right)_{\alpha}(x)
&
=
\chi_y\left(\gamma_{\ak}\circ \iaak(x)\right)
\\
&
=
\exp(\,i\,\acc{y}{\deg_T \circ \gamma_{\ak}\circ \iaak (x)})
\label{eq:ut_fonct_deg_1}
\\
&
=
\label{eq:ut_fonct_deg_2}
\exp\left(\,i\,\acc{y}{\gamma_{\R}^{\vee}\circ i_{\alpha,\R}^{\vee}\left[\deg_{\Res_{\Ka/K}\G_m}(x)\right]}\right)
\\
&
=
\exp\left(\,i\,\acc{y}{\gamma_{\R}^{\vee}\circ i_{\alpha,\R}^{\vee}\left[\deg_{\Ka}(x)\right]}\right)
\\
&
=
\label{eq:ut_fonct_deg_3}
\exp\left(\,i\,\acc{y}{\gamma_{\R}^{\vee}\left[\deg_{\Ka}(x)\,\Da^{\vee}\right]}\right)
\\
&
=
\label{eq:ut_fonct_deg_4}
\exp\left(\,i\,\deg_{\Ka}(x)\acc{y}{\roa}\right).
\end{align}
Le passage de \eqref{eq:ut_fonct_deg_1} \`a \eqref{eq:ut_fonct_deg_2} d\'ecoule
du lemme \ref{lm:fonc:degre}, et le passage de \eqref{eq:ut_fonct_deg_3}
\`a \eqref{eq:ut_fonct_deg_4} vient de \eqref{eq:gammaveeDaroa}.
\end{demo}

\subsubsection{Cas fonctionnel}\label{subsub:fonc}

L\`a encore, commen\`c;ons par d\'ecrire les caract\`eres de $T(\ak)$ triviaux sur $T(\ak)^1$.
On a un morphisme surjectif
\begin{equation}
\deg_T\,:\,T(\ak)\longto \DT
\end{equation}
de noyau $T(\ak)^1$, 
qui induit donc un isomorphisme
\begin{equation}\label{eq:isotaktak1DT}
\deg_T\,:\,T(\ak)/T(\ak)^1 \longisom \DT
\end{equation}
et par dualit\'e un isomorphisme entre $\dualtop{T(\ak)/T(\ak)^1}$ et $\left(\DT^{\vee}\right)_{\unit}$.
Pour tout $\bz\in \xtgu\subset \left(\DT^{\vee}\right)_{\unit}$ 
on note \nindex{$\chi_{\bz}$}$\chi_{\bz}$ le caract\`ere de 
$\dualtop{T(\ak)/T(\ak)^1}$
correspondant via cet isomorphisme. On a ainsi 
\begin{equation}
\forall t\in T(\ak),\quad \chi_{\bz}(t)=\acc{\bz}{\deg_T(t)}.
\end{equation}
En particulier, si $v$ est une place de $K$, on a
\begin{equation}\label{eq:chizv}
\forall \bz\in \xtgu,\quad \forall x\in T(K_v),\quad (\chi_{\bz})_v(x)=\acc{\bz}{\deg_{T,v}}^{f_v}.
\end{equation}

Nous noterons \nindex{$\da$}$\da$ le degr\'e absolu de $\Ka$. 
On a donc $\card{\F_{\Ka}}=q^{\da}$.
En outre, d'apr\`es le lemme \ref{lm:degre:quasidep}, 
on a 
\begin{equation}
\card{\CTqd}=\prod_{\alpha\in\Sigma(1)/G}\da.
\end{equation}

Rappelons le d\'ebut de la suite exacte \eqref{eq:resflT:G}
\begin{equation}
0\longto \xtg \overset{\gamma}{\longto} \Ps^G \longto  \Pic(\xs)
\end{equation}
qui montre que $\Ps^G/\xtg$ est sans torsion.

Ainsi le morphisme dual de $\gamma$
\begin{equation}
\left(\gamma\right)^{\vee}\,:\,
\left(\Ps^G\right)^{\vee}
\longto 
\Hom(X(T)^{\,G},\Z)
\end{equation}
est surjectif. 
D'apr\`es \eqref{eq:gammaveeDaroa} ce morphisme envoie $\Da^{\vee}$
sur $\roa$.

On obtient un diagramme
\begin{equation}\label{eq:diag}
\xymatrix{
\G_m(\aka)\ar[r]^{\iaak}\ar[d]_{\deg_{\Res_{\Ka/K}\G_m}}&T_{\Ps}(\ak)
\ar[r]^<<<<<<<<<<{\gamma_{\,\ak}} \ar[d]_{\deg_{T_{\Ps }}} & T(\ak)\ar[d]_{\deg_T}\\
 \Z\ar[r]^{i_{\alpha}^{\vee}}&\left(\Ps^G\right)^{\vee}\ar[r]^<<<<<<{\gamma^{\vee}}& \Hom(X(T)^{\,G},\Z)
}
\end{equation}
Le lemme \ref{lm:fonc:degre}
montre que
ce diagramme est commutatif.

On note \nindex{$\DT^0$}$\DT^0$ l'image de $\DTqd$ par $\gamma^{\vee}$ dans $\DT$.

\begin{lemme}\label{lm:carac:fonc}
Soit $\alpha\in \Sigma(1)/G$. 
\begin{enumerate}
\item
$\da\,\roa$ est dans $\DT^0$.
\item
On a 
\begin{equation}
\forall x\in \G_m(\aka),\quad
\forall \bz\in \xtgu,\quad 
\left(\chi_{\bz}\right)_{\alpha}(x)=\acc{\bz}{\deg_{\Ka}(x)\,\da\,\roa}.
\end{equation}
\end{enumerate}
\end{lemme}
\begin{demo}
Montrons qu'on a, pour tout $x\in \G_m(\aka)$,
\begin{equation}\label{eq:gammagvee}
\gamma^{\vee}\left(\deg_{T_{\Ps}}\left[\iaak(x)\right]\right)
=
\da\,\roa\,\deg_{\Ka}(x).
\end{equation}
En effet, on a, par commutativit\'e de \eqref{eq:diag},
et le lemme \ref{lm:degre:quasidep},
\begin{align}
\gamma^{\vee}\left(\deg_{T_{\Ps}}\left[\iaak(x)\right]\right)
&
=
\gamma^{\vee}\left(i_{\alpha}^{\vee}\left[\deg_{\Res_{\Ka/K}\G_m}(x)\right]\right)
\\
&
=
\gamma^{\vee}\left(i_{\alpha}^{\vee}\left[\da\,\deg_{\Ka}(x)\right]\right)
\\
&
=
\deg_{\Ka}(x)\,\da\,\roa.
\end{align}
Le fait que $\da\roa$ soit dans $\DT^0$ d\'ecoule alors de 
de \eqref{eq:gammagvee}
et de l'existence
d'\'el\'ements $x$ de $\G_m(\aka)$ v\'erifiant $\deg_{\Ka}(x)=1$.

Toujours par commutativit\'e de \eqref{eq:diag}, on a, pour tout $x\in \G_m(\aka)$,
\begin{equation}
\deg_T\left(\gamma_{\ak}\circ \iaak(x)\right)
=
\gamma^{\vee}\left(i_{\alpha}\left[\da\,\deg_{\Ka}(x)\right]\right)
=
\da\,\roa\,\deg_{\Ka}(x).
\end{equation}

Pour $\bz\in \xtgu$, on a alors 
\begin{align}
\left(\chi_{\bz}\right)_{\alpha}(x)
&
=
\chi_{\bz}\left(\gamma_{\ak}\circ \iaak(x)\right)
\\
&
=
\acc{\bz}{\deg_T\left(\gamma_{\ak}\circ \iaak(x)\right)}
\\
&
=\acc{\bz}{\deg_{\Ka}(x)\,\da\,\roa}.
\end{align}
\end{demo}

\begin{rem}
La commutativit\'e du diagramme \eqref{eq:diag}
et la surjectivit\'e du morphisme 
\begin{equation}
\gamma^{\vee}\,:\,\left(\Ps^G\right)^{\vee}\longto \xtgd
\end{equation}
permettent de retrouver le fait que le conoyau de $\deg_T$ est fini, c'est-\`a-dire
de red\'emontrer la proposition \ref{prop:conoyau:fini}.
On notera que si $L$ a le m\^eme corps des constantes que $K$, 
les $\da$ sont tous \'egaux \`a un, 
$\deg_{T_{\Ps }}$ est surjectif et donc $\deg_T$ \'egalement, 
ce qui red\'emontre le corollaire \ref{cor:ct0}.
\end{rem}

\subsection{Pr\'eliminaires au calcul des transform\'ees de Fourier}
\label{subsec:prelim}

Soit $v$ une place finie de $K$. 
Nous choisissons une place $\V$ de $L$ au-dessus de $v$, nous notons $G_v$ son groupe de d\'ecomposition. 
On notera encore $\V$ la valuation normalis\'ee de $L$ qui repr\'esente $\V$.

Nous notons $\Sigma(1)/G_v$ l'ensemble des orbites de $\Sigma(1)$ sous l'action de $G_v$, 
et pour $\alpha\in \Sigma(1)/G$ nous notons $\alpha/G_v$ 
le sous-ensemble de $\Sigma(1)/G_v$ des orbites incluses dans $\alpha$. 
Soit $\Sigma^{\,G_v}$ l'ensemble des c\^ones de $\Sigma$ globalement invariants 
sous l'action de $G_v$, en d'autres termes les c\^ones $\sigma$ 
dont l'ensemble des rayons $\sigma(1)$ est stable sous l'action de $G_v$. 
Pour un c\^one $\sigma\in \Sigma^{G_v}$, 
nous notons $\sigma(1)/G_v$ l'ensemble des orbites de $\sigma(1)$ sous l'action de $G_v$, 
et pour $\alpha\in \Sigma(1)/G$ nous notons $\alpha/G_v$ 
le sous-ensemble de $\sigma(1)/G_v$ des orbites incluses dans $\alpha$.

Dans toute la suite de cet texte, nous adoptons pour all\'eger l'\'ecriture la convention suivante :
si $(X_{\alpha})_{\alpha\in \Sigma(1)/G}$ est une famille d'objets index\'ee par $\Sigma(1)/G$,
pour tout $\beta\in \Sigma(1)/G_v$, $X_{\beta}$ d\'esigne $X_{\alpha}$ o\`u 
$\alpha$ est l'unique \'el\'ement de $\Sigma(1)/G$ contenant $\beta$.

Pour $\beta\in \Sigma(1)/G_v$ nous notons $\lb $ le cardinal de $\beta$
et  $\rob'$ un g\'en\'erateur d'un \'el\'ement quelconque de $\beta$.
Si $g_{\beta}\in G$ est tel que $\rob'=g_{\beta}.\rob$, on a donc
\begin{equation}\label{eq:lbeta}
\lb =
\frac
{\card{G_v}}
{\card{\left(g_{\beta}\,\gb\,g_{\beta}^{-1}\cap G_v\right)}}.
\end{equation}

On pose 
\begin{equation}
\taub =\sum_{\rho\in G_v.\rob'}\,\rho,
\end{equation}
de sorte que $\taub $ est un \'el\'ement de $(X(T)^{\vee})^{\,G_v}$.

Rappelons que l'on note $\intrel(\sigma)$ l'int\'erieur relatif d'un c\^one $\sigma$.

\begin{lemme}\label{lm:xtgv}
Pour tout $x$ \'el\'ement de $(X(T)^{\,\vee})^{\,G_v}$, 
il existe un \'el\'ement $\sigma\in \Sigma^{\,G_v}$ tel que $x$ s'\'ecrit
\begin{equation}
x=\sum_{\beta\in \sigma(1)/G_v} \nb\,\taub 
\end{equation}
o\`u les $\nb$ sont des entiers strictement positifs.

En particulier, $(X(T)^{\,\vee})^{\,G_v}$ est recouvert par les ensembles 
\begin{equation}
\intrel(\sigma)\cap (X(T)^{\,\vee})^{\,G_v}
\end{equation}
pour $\sigma$ d\'ecrivant $\Sigma^{G_v}$. 
\end{lemme}
\begin{demo}
Soit $x\in(X(T)^{\,\vee})^{\,G_v}$. Alors $x$ est dans l'int\'erieur relatif d'un unique c\^one $\sigma$ de $\Sigma$. 
Par ailleurs, pour tout $g\in G_v$, $x=g\,x$ est dans l'int\'erieur relatif du c\^one $g\,\sigma$ de $\Sigma$, 
et donc $g\,\sigma=\sigma$. 
Ainsi $\sigma\in \Sigma^{\,G_v}$. 
Le reste du lemme en d\'ecoule ais\'ement, compte tenu du fait que l'\'eventail est r\'egulier.
\end{demo}

Pour $\alpha\in \Sigma(1)/G$, consid\'erons l'application
\begin{equation}
\begin{array}{rcl}
G/\ga&\longto&\{w\in \placesde{\Ka},\,v|w\}\\
g&\longmapsto&g^{-1}.\V|_{\Ka} \end{array}
\end{equation}
Ce n'est autre que le passage au quotient par l'action de $G_v$, 
d'o\`u une correspondance entre $\alpha/G_v$ et les places de $\Ka$ au-dessus de $v$.

Soit $\beta\in \alpha/G_v$ et $\wb$ la place correspondante.
Soit $g_{\beta}\in G$ tel que $\rob'=g_{\beta}.\roa$.

On a 
\begin{equation}
[L_{\V}:K_v]
=
\card{G_v}
\end{equation}
et 
\begin{equation}
[L_{\V}:(\Ka)_{\wb}]
=
\card{ g_{\beta}\,\ga\,g_{\beta}^{-1}\cap G_v}
\end{equation}

Ainsi, si $v$ est non ramifi\'ee dans $L$, 
on a 
\begin{equation}
[k_{\V}:k_v]
=
\card{G_v}
\end{equation}
et 
\begin{equation}
[k_{\V}:k_{\wb}]
=
\card{ g_{\beta}\,\ga\,g_{\beta}^{-1}\cap G_v}
\end{equation}
d'o\`u, en combinant les deux \'egalit\'es pr\'ec\'edentes
et la formule \eqref{eq:lbeta},
\begin{equation}
\lb =[k_{\wb}:k_v].
\end{equation}
Ainsi on a 
\begin{equation}\label{eq:form:qwb}
q_{\wb}=q_v^{\,\lb }.
\end{equation}

Dans le cas fonctionnel (et toujours en supposant $v$ non ramifi\'ee dans $L$), 
on a de plus 
\begin{equation}
q_{\wb}=q_{\Ka}^{f_{\wb}}=q^{\,\da\,f_{\wb}},
\end{equation}
soit 
\begin{equation}\label{eq:form:fvlb}
f_v\,\lb =\da\,f_{\wb}.
\end{equation} 

\begin{rem}
Bien que nous n'en ayons pas besoin,
on peut noter que si $v$ est  ramifi\'ee dans $L$,
en notant $e_v$ l'indice de ramification de $v$ dans $L$ 
et $\eb$ l'indice de ramification
de $\wb$ dans $L$, un calcul similaire montre qu'on a
\begin{equation}
\lb =\frac{e_v}{\eb }[k_{\wb}:k_v].
\end{equation}
soit dans le cas fonctionnel
\begin{equation}\label{eq:form:fvlb:ram}
f_v\,\eb \,\lb =e_v\,\da\,f_{\wb}.
\end{equation} 
\end{rem}

\subsection{Les transform\'ees de Fourier locales}\label{subsec:transfo:four:locales}
Rappelons que $\plsgc$ s'identifie canoniquement \`a  $\left(\Ps^G\right)_{\C}$, et donc
\`a $\C^{\,\sg}$ via le choix de la base $(\Da)_{\alpha\in \sg}$. 
Nous noterons d\'esormais $\bs$ ou $(\sa)_{\alpha\in \sg}$ un \'el\'ement de $\plsgc$.

\subsubsection{Cas d'une place finie quelconque}\label{subsubsec:quelc}

\begin{prop}\label{prop:fv:place:finie:quelc}
Soit $v$ une place finie de $K$.
\begin{enumerate}
\item
Pour tout $\bs\in\tube{\R_{>0}^{\Sigma(1)/G}}$
\label{item:1:prop:fv:place:finie:quelc}
la fonction $H_{v}(-\bs,\,.\,)$ est int\'egrable
sur $T(K_v)$.
\item
Soit $\chi$ un \'el\'ement de $\dualtop{T(\ak)/\K(T)}$.
\label{item:2:prop:fv:place:finie:quelc}
La fonction
\begin{equation}
\bs
\mapsto
\int\limits_{T(K_v)}\,
H_{v}(-\bs,t)
\chi_v(t)\,
d\mu_v(t)
\end{equation}
est holomorphe
sur $\tube{\R_{>0}^{\Sigma(1)/G}}$.
On la note $f_{v}(\chi,\,.\,)$.
\item
\label{item:lm:majfvchi:vnonar}
Pour tout compact $\K$ de $\R^{\,\Sigma(1)/G}_{>0}$,
il existe une constante $C>0$
telle qu'on ait pour tout $\chi\in \dualtop{T(\ak)/\K(T)}$
et tout $\bs\in \tube{\K}$ la majoration
\begin{equation}
\abs{f_v(\chi,\bs)}
\leq
C.
\end{equation}
\end{enumerate}
\end{prop}
\begin{demo}
Pour tout $\bs\in \plsgc$ on a
\begin{equation}
\abs{H_v(\bs,\,.\,)}
=
H_v(\Re(\bs),\,.\,).
\end{equation}
Comme la fonction $H_{v}(-\bs,\,.\,)$ est $T(\Ov)$-invariante, 
elle induit une fonction sur $T(K_v)/T(\Ov)$, 
qui sera not\'ee de fa\c con identique. 
Pour $\bs\in \plsgc$
on a 
\begin{align} 
&\quad\int\limits_{T(K_v)}\,\abs{H_{v}(-\bs,t)}d\mu_v (t)\\
&=\int\limits_{T(K_v)}\,H_{v}(-\Re(\bs),t)d\mu_v (t)\\
&
=
\left(\int\limits_{T(\Ov)}d\mu_v \right)\,\sum_{t\in T(K_v)/T(\Ov)}\,H_{v}(-\Re(\bs),t).
\end{align}
Soit $\V$ une place de $L$ divisant $v$ et $G_v$
son groupe de d\'ecomposition.
Soit $t\in T(K_v)/T(\Ov)$.
Alors $\deg_{T,L,\V}(t)$ est un \'el\'ement de $\left(X(T)^{\vee}\right)^{G_v}$.
D'apr\`es le lemme \ref{lm:xtgv}
il existe un unique \'el\'ement $\sigma\in \Sigma^{\,G_v}$ tel que $\deg_{T,L,\V}(t)$ s'\'ecrit
\begin{equation}
\deg_{T,L,\V}(t)=\sum_{\beta\in \sigma(1)/G_v} \nb\,\taub 
\end{equation}
o\`u les $\nb$ sont des entiers strictement positifs.
On a alors 
\begin{equation}
H_{v}
\left(
-\bs
,
t
\right)
=
q_v^{
\,\,
-\frac{1}{e_v}\underset{\beta\in \sigma(1)/G_v}{\sum}\,\nb \,\lb \,\sb
}.
\end{equation}

Ainsi, compte tenu du fait que $\deg_{T,L,\V}\,:\,T(K_v)/T(\Ov)\to \left(X(T)^{\vee}\right)^{G_v}$
est injectif
et du lemme \ref{lm:xtgv}, on a
\begin{equation}
\sum_{t\in T(K_v)/T(\Ov)}\,H_{v}(-\Re(\bs),t)
\leq
\sum_{\sigma\in\Sigma^{G_v}}\,\,
\prod_{\beta\in \sigma(1)/G_v}
\sum_{(\nb )\in \Z_{>0}^{\sigma(1)/G_v}}
q_v^{
\,\,
-\frac{1}{e_v}\underset{\beta\in \sigma(1)/G_v}{\sum}\,\nb \,\lb \,\sb
}.
\end{equation}
Le membre de gauche est visiblement fini d\`es qu'on a  $\Re(\sa)>0$ pour tout $\alpha$.
La proposition en d\'ecoule facilement.
\end{demo}
\begin{rem}\label{rem:fv=}
Pour $\chi\in \dualtop{T(\ak)/\K(T)}$, 
la d\'emonstration montre que pour tout $\bs\in \tube{\R_{>0}^{\sg}}$, la s\'erie
\begin{equation}
\sum_{t\in T(K_v)/T(\Ov)}\,\chi(t)\,H_{v}(-\bs,t)
\end{equation}
est absolument convergente et qu'on a 
\begin{equation}
f_v(\chi,\bs)=\left(\int\limits_{T(\Ov)}\!d\mu_v\right)\sum_{t\in T(K_v)/T(\Ov)}\,\chi(t)\,H_{v}(-\bs,t).
\end{equation}
\end{rem}

\subsubsection{Cas arithm\'etique}

\begin{lemme}\label{lm:rel:hauteur:carac:arit}
On se place dans le cas arithm\'etique.
Soit $v$ une place de $K$. On a 
\begin{equation}
\forall\,y\in  X(T)^G_{\R},
\quad
\forall t\in T(K_v),
\quad
\left(\chi_y\right)_v(t)
=
H_{v}
\left(i\,\gamma_{\R}(y),t\right).
\end{equation}
\end{lemme}
\begin{demo}
Soit $\V$ une place de $L$ divisant $v$. 
Pour $y\in  X(T)^G_{\R}$
et $t\in T(K_v)$ on a,
compte tenu de \eqref{eq:chiyv},
de la d\'efinition \eqref{eq:hauteur:locale} de $H_v$
et du lemme \ref{lm:diag:xtgvplsgv} 
(respectivement  \ref{lm:diag:xtgvplsgv:archi})
si $v$ est finie (respectivement archim\'edienne)~:
\begin{align}
(\chi_y)_v(t)
&=\exp(\,i\,\log(q_v)\,\acc{y}{\deg_{T,v}(t)})\\
&=\exp(\,i\,\log(q_v)\,\accsv{\gamma_{\R}(y)}{t})\\
&=\exp(\log(q_v)\,\accsv{i\,\gamma_{\R}(y)}{t})\\
&=H_v(i\,\gamma_{\R} (y),t).
\end{align}
d'o\`u le r\'esultat.
\end{demo}

\begin{cor}
On a 
\begin{multline}
\forall\,y\in  X(T)^G_{\R},
\quad
\forall\,\bs\in \plsgc,
\quad
\forall t\in T(K_v),
\\
H_{v}(\bs,t)
\left(\chi_y\right)_v(t)
=
H_{v}
\left(\bs+i\,\gamma_{\R}(y),t\right).
\end{multline}
\end{cor}

\begin{cor}
\label{cor:transfo:ram:arit}
Soit $v$ une place de $K$ et $\chi$ un \'el\'ement de $\dualtop{T(\ak)/\K(T)}$.
Alors on a
\begin{multline}
\forall\,y\in  X(T)^G_{\R},
\quad
\forall\,\bs\in \tube{\R^{\,\Sigma(1)/G}_{>0}},\\
f_v(\chi.\chi_y,\bs)
=
f_{v}
\left(
\chi,
\bs-i\,\gamma_{\R}(y)
\right).
\end{multline}
\end{cor}

Dans le reste de cette section, on s'int\'eresse aux propri\'et\'es des
transform\'ees de Fourier locales aux places archim\'ediennes.
Ces propri\'et\'es seront cruciales pour appliquer la formule de Poisson
dans le cas arithm\'etique.

Soit $v\in \placesde{K,\infty}$ et $\V$ une place de $L$
divisant $v$, de groupe de d\'ecomposition $G_v$. 
D'apr\`es le lemme \ref{lm:varchitkvtov}, 
$T(K_v)/T(\Ov)$ 
s'identifie naturellement \`a $\left(X(T)_{\R}^{\vee}\right)^{G_v}$.

On note $\Sigma^{G_v}_{\max}$
l'ensemble des c\^ones de $\Sigma^{G_v}$ de dimension maximale.
Soit $\sigma\in \Sigma^{G_v}_{\max}$. 
Pour tout $j\in \sigma(1)/G_v$,
soit $l_j$ le cardinal de $j$ (ainsi $l_j=1$ ou $2$), 
$\rho_j$ un g\'en\'erateur d'un \'el\'ement de $j$,
et
\begin{equation}
\tau_j=\sum_{\rho\in G_v\,\rho_j} \rho.
\end{equation}
Si $\alpha$ est l'\'el\'ement de $\Sigma(1)/G$ tel que $j\subset \alpha$,
la notation $s_j$ d\'esigne $\sa$.

La proposition suivante donne les propri\'et\'es et estimations n\'ecessaires 
pour les transform\'ees
de Fourier locales aux places archim\'ediennes.

\begin{prop}\label{prop:estim:transfo:locale:archi}
Soit $v$ un \'el\'ement de $\placesde{K,\infty}$.
\begin{enumerate}
\item
Pour tout $\bs\in\tube{\R_{>0}^{\Sigma(1)/G}}$
\label{item:1:prop:estim:transfo:locale:archi}
la fonction $H_{v}(-\bs,\,.\,)$ est int\'egrable
sur $T(K_v)$.
\item
Soit $\chi$ un \'el\'ement de $\dualtop{T(\ak)/\K(T)}$.
\label{item:2:prop:estim:transfo:locale:archi}
La fonction
\begin{equation}
\bs
\mapsto
\int\limits_{T(K_v)}\,
H_{v}(-\bs,t)
\chi_v(x)\,
d\mu_v (x)
\end{equation}
est holomorphe
sur $\tube{\R_{>0}^{\Sigma(1)/G}}$.
On la note $f_{v}(\chi,\,.\,)$.
\item
On a, pour 
\label{item:3:prop:estim:transfo:locale:archi}
tout \'el\'ement $\chi$ de $\dualtop{T(\ak)/\K(T)}$,
\begin{equation}
\forall \bs\in \tube{\R^{\Sigma(1)/G}_{>0}},\quad
f_v(\chi,\bs)
=
\sum_{\sigma\in \Sigma^{G_v}_{\max}}
\frac{1}
{\produ{j\in \sigma(1)/G_v}
\frac{l_j}{[L_{\V}:K_v]}
\,
(s_j+i\,\acc{\rho_j}{\chi_v})}.
\end{equation}
\item
Soit $\eps$ v\'erifiant $0<\eps<1$, $\K$ un compact de $\R^{\Sigma(1)/G}_{>1-\eps}$
\label{item:maj:prop:estim:transfo:locale:archi}
et $\norm{.}$ une norme sur $\dualtop{T(K_v)/T(\Ov)}$.
Il existe alors une constante $C>0$ telle qu'on ait,
pour tout $v$ de $\placesde{K,\infty}$,
et tout $\chi\in \dualtop{T(\ak)/\K(T)}$,
\begin{equation}\label{eq:maj:fchiv:varch}
\forall \bs\in \tube{\K},
\quad
\abs{f_v(\chi,\bs)}
\leq
\frac{C}
{1+\norm{\chi_v}}
\,
\sum_{\sigma\in \Sigma^{G_v}_{\max}}
\frac
{1+\sumu{j\in \sigma(1)/G_v}\abs{\Im(s_j)}
}
{
\produ{j\in \sigma(1)/G_v}
(1+\abs{\acc{\rho_j}{\chi_v}+\Im(s_j)})
}
\end{equation}
\end{enumerate}
\end{prop}

\begin{demo}
Dans tout ce qui suit, on identifie 
$T(K_v)/T(\Ov)$ \`a $\left(X(T)^{\vee}\right)^{G_v}_{\R}$
et donc $\dualt{T(K_v)/T(\Ov)}$
\`a $X(T)^{G_v}_{\R}$.

On commence par remarquer la chose suivante.
Soit $t\in T(K_v)/T(\Ov)$.
Alors 
il existe un \'el\'ement $\sigma$  de  $\Sigma^{\,G_v}_{\max}$
et des r\'eels positifs $(t_j)_{j\in \sigma(1)/G_v}$ 
tels qu'on ait
\begin{equation}
t=\sum_{j\in \sigma(1)/G_v} t_j\,\tau_j.
\end{equation}

D'apr\`es \eqref{eq:hauteur:locale:archi}
et \eqref{eq:hauteur:locale}, on a alors, 
pour tout $(\sa)\in \C^{\Sigma(1)}$,
\begin{equation}
H_{v}
\left(
\bs
,
t
\right)
=
\exp\left(\,\frac{1}{[L_{\V}:K_v]}\sumu{j\in \sigma(1)/G_v}\,t_j \,l_j\,s_j\right)
\end{equation}
On a \'egalement, pour tout $\chi_v\in \dualt{T(K_v)/T(\Ov)}$,
\begin{equation}
\chi_{v}(t)
=
\exp\left(\,
i
\,
\acc{t}{\chi_v}
\right)
=
\exp\left(\,i\,\sumu{j\in \sigma(1)/G_v}t_j\,l_j\,\acc{\rho_j}{\chi_v}\right).
\end{equation}

Pour tout $\bs\in \C^{\Sigma(1)/G}$ on a
$
\abs{H_v(\bs,\,.\,)}
=
H_v(\Re(\bs),\,.\,).
$
Comme la fonction $H_{v}(-\bs,\,.\,)$ est $T(\Ov)$-invariante, 
elle induit une fonction sur $T(K_v)/T(\Ov)$, 
qui sera not\'ee de fa\c con identique.

Pour $\bs\in \C^{\Sigma(1)/G}$, on a alors
\begin{align} 
\int\limits_{T(K_v)}\abs{H_{v}(-\bs,t)}d\mu_v (t)
&=\int\limits_{T(K_v)}H_{v}(-\Re(\bs),t)d\mu_v (t)\\
&
=
\int\limits_{T(K_v)/T(\Ov)}\!\!\!\!H_{v}(-\Re(\bs),t)dt,
\end{align}
o\`u $dt$ est la mesure de Lebesgue sur $T(K_v)/T(\Ov)\isom \left(X(T)_{\R}^{\vee}\right)^{G_v}$,
normalis\'ee par le r\'eseau $\left(X(T)^{\vee}\right)^{G_v}$

Comme 
$\left(X(T)_{\R}^{\vee}\right)^{G_v}$
est recouvert par les \'el\'ements de $\Sigma^{G_v}_{\max}$
et que l'intersection de deux \'el\'ements distincts
de $\Sigma^{G_v}_{\max}$ est de mesure de Lebesgue nulle, on a
\begin{equation}
\int\limits_{T(K_v)/T(\Ov)}\!\!\!\!H_{v}(-\Re(\bs),t)dt
=
\sum_{\sigma\in \Sigma^{G_v}_{\max}}
\int_{\sigma}H_{v}(-\Re(\bs),t)dt
\end{equation}
Pour $\sigma\in \Sigma^{G_v}_{\max}$, les $(\tau_j)_{j\in \sigma(1)/G_v}$
engendrent $\sigma$ et forment une base de $(X(T)^{\vee})^{G_v}$.
Ainsi, si on identifie $\sigma$ \`a $\R_{\geq 0}^{\sigma(1)/G_v}$ au moyen
de la base $(\tau_j)$,
la  restriction de la mesure $dx$ \`a $\sigma$ s'identifie \`a la mesure produit
$\otimesu{j\in \sigma(1)/G_v}dx_j$ sur $\R_{\geq 0}^{\sigma(1)/G_v}$.
On a donc
\begin{equation}
\int_{\sigma}H_{v}(-\Re(\bs),t)dt
=
\!\!\!\!
\int\limits_{\R_{\geq 0}^{\sigma(1)/G_v}}\!\!\!\!H_{v}\left(-\Re(\bs),\sum t_j\,\tau_j\right)
\otimesu{j\in \sigma(1)/G_v}dt_j
\end{equation}
soit
\begin{equation}
\int\limits_{T(K_v)/T(\Ov)}\!\!\!\!\!\!H_{v}(-\Re(\bs),t)dt
=
\sum_{\sigma\in \Sigma^{G_v}_{\max}}
\,\,\,
\int\limits_{\R_{\geq 0}^{\sigma(1)/G_v}}
\!\!\!\!\exp\left(\frac{\sum -t_j\,l_j\,\Re(s_j)}{[L_{\V}:K_v]} \right)\otimesu{j\in \sigma(1)/G_v}dt_j
\end{equation}
Cette derni\`ere expression est finie d\`es que $\Re(s_j)>0$ pour tout $j\in \sigma(1)/G_v$
et tout $\sigma\in \Sigma^{G_v}_{\max}$, c'est-\`a-dire d\`es que  $\Re(\sa)>0$
pour tout $\alpha\in \Sigma(1)/G$.
Elle vaut alors
\begin{equation}
\sum_{\sigma\in \Sigma^{G_v}_{\max}}
\prod_{j\in \sigma(1)/G_v} \frac{1}{\frac{l_j}{[L_{\V}:K_v]}\,\Re(s_j)}.
\end{equation}
On en d\'eduit les points \ref{item:1:prop:estim:transfo:locale:archi}
et \ref{item:2:prop:estim:transfo:locale:archi}.

Montrons le point 
\ref{item:3:prop:estim:transfo:locale:archi}.
Soit $\bs\in \tube{\R^{\Sigma(1)/G}_{>0}}$.
On a 
\begin{align} 
&
\quad\int\limits_{T(K_v)}\,\chi_v(x)\,H_{v}(-\bs,t)d\mu_v (t)
\\
&=\int\limits_{T(K_v)/T(\Ov)}\,\chi_v(t)\,H_{v}(-\bs,t)dt
\\
&
=
\sum_{\sigma\in \Sigma^{G_v}_{\max}}
\int_{\sigma}\chi_v(t)\,H_{v}\left(-\bs,t\right)dt
\\
&
=
\sum_{\sigma\in \Sigma^{G_v}_{\max}}
\,\,\,
\int\limits_{\R_{\geq 0}^{\sigma(1)/G_v}}
\!\!\!\!
\chi_v\left(\sum t_j\,\tau_j\right)
\,
H_{v}\left(-\bs,\sum t_j\,\tau_j\right)\otimesu{j\in \sigma(1)/G_v}dt_j
\\
&
=
\sum_{\sigma\in \Sigma^{G_v}_{\max}}
\,\,\,\,
\int\limits_{\R_{\geq 0}^{\sigma(1)/G_v}}
\exp\left(\,\frac{\sum t_j\,l_j\,(s_j+i\,\acc{\rho_j}{\chi_v})}{[L_{\V}:K_v]}\right)
\otimesu{j\in \sigma(1)/G_v}dt_j
\\
&
=
\sum_{\sigma\in \Sigma^{G_v}_{\max}}
\prod_{j\in \sigma(1)/G_v} \frac{1}{\frac{l_j}{[L_{\V}:K_v]}\,(s_j+i\,\acc{\rho_j}{\chi_v})}.
\label{eq:res_fvchisa}
\end{align}

Montrons le point \ref{item:maj:prop:estim:transfo:locale:archi}.
Soit $(e_l)_{l\in L}$ une base de $(X(T)^{\vee})^{G_v}$
et $(e_l^{\vee})$ sa base duale.
Soit $\bs\in \tube{\R^{\Sigma(1)/G}_{>0}}$ et $\chi\in \dualtop{T(\ak)/\K(T)}$.
Soit $g$ (respectivement $h$) la fonction de $\R^L$ dans $\C$ donn\'ee par 
\begin{equation}
\forall (t_l)\in \R^L,\quad g(t_l)=H_v\left(-\bs,\sumu{l\in L} t_l\,e_l\right), 
\end{equation}
respectivement
\begin{equation}
\forall (t_l)\in \R^L,\quad 
h(t_l)
=\chi_v\left(\sumu{l\in L} t_l\,e_l\right)
=
\exp\left(i \sumu{l\in L} t_l \acc{e_l}{\chi_v}\right).
\end{equation}
On a ainsi , pour $l\in L$,
\begin{equation}
\frac{\partial }{\partial t_l}h(t)=i\,\acc{e_l}{\chi_v} h(t).
\end{equation}
Pour tout $l$ v\'erifiant $\acc{e_l}{\chi_v}\neq 0$ on obtient 
en int\'egrant par parties
\begin{align}
f_v(\chi,\bs)
&=\int\limits_{\R^L}g(t)\,h(t)dt
\\
&=-\frac{1}{i\,\acc{e_l}{\chi_v}}\int\limits_{\R^L} h(t) \frac{\partial }{\partial t_l}g(t)dt
\\
&=
-\sum_{\sigma\in \Sigma^{G_v}_{\max}}
\frac{1}{i\,\acc{e_l}{\chi_v}} 
\int_{\sigma}
h(t) \frac{\partial }{\partial t_l}g(t)dt.
\end{align}
La restriction de $g$ \`a $\sigma\in \Sigma^{G_v}_{\max}$ s'\'ecrit 
(en notant $(\tau^{\vee}_j)_{j\in \sigma(1)/G_v}$ la base duale
de $(\tau_j)_{j\in \sigma(1)/G_v}$)
\begin{align}
g(t)&=\exp\left(\,-\frac{1}{[L_{\V}:K_v]}\,\sumu{j\in \sigma(1)/G_v}\,\acc{\sumu{l\in L} t_l\,e_l}{\tau_j^{\vee}} \,l_j\,s_j\right)\\
&=\exp\left(\,-\frac{1}{[L_{\V}:K_v]}\sumu{l\in L} t_l \sumu{j\in \sigma(1)/G_v} \acc{e_l}{\tau_j^{\vee}} \,l_j\,s_j\right)
\end{align}
et finalement sur $\sigma$ on a
\begin{equation}
\frac{\partial }{\partial t_l}g(t)
=
-\left(\frac{1}{[L_{\V}:K_v]}\,\sumu{j\in \sigma(1)/G_v} \acc{e_l}{\tau_j^{\vee}} \,l_j\,s_j\right) g(t)
\end{equation} 
On en d\'eduit
\begin{equation}
i\,\acc{e_l}{\chi_v}\,f_v(\chi,\bs)
=
\sum_{\sigma\in \Sigma^{G_v}_{\max}}
\left(\sumu{j\in \sigma(1)/G_v} \acc{e_l}{\tau_j^{\vee}} \,\frac{l_j}{[L_{\V}:K_v]}\,s_j\right)
\int_{\sigma}
h(t) g(t)dt
\end{equation}
soit
\begin{multline}\label{eq:iaccelchivetc}
i\,\acc{e_l}{\chi_v}\,f_v(\chi,\bs)
\\
=
\sum_{\sigma\in \Sigma^{G_v}_{\max}}
\left(\sum_{j\in \sigma(1)/G_v} \acc{e_l}{\tau_j^{\vee}}\,\frac{l_j}{[L_{\V}:K_v]}\,s_j \right) 
\prod_{j\in \sigma(1)/G_v} \frac{1}{\frac{l_j}{[L_{\V}:K_v]}\,(s_j+i\,\acc{\rho_j}{\chi_v})}.
\end{multline}
En  prenant les modules dans \eqref{eq:res_fvchisa} et \eqref{eq:iaccelchivetc}
pour tous les $l$ tels que $\acc{e_l}{\chi_v}\neq 0$ et en ajoutant le tout,
on obtient apr\`es une majoration \'evidente
\begin{multline}\label{eq:fvchisa}
\left(1+\sum_{l\in L} \abs{\acc{e_l}{\chi_v}}\right)
\abs{f_v(\chi,\bs)}
\\
\leq
\sum_{\sigma\in \Sigma^{G_v}_{\max}}
\left(1+\sum_{l\in L} \abs{\sum_{j\in \sigma(1)/G_v} \acc{e_l}{\tau_j^{\vee}}\,\frac{l_j}{[L_{\V}:K_v]}\,s_j} \right) 
\prod_{j\in \sigma(1)/G_v} \frac{1}{\abs{\frac{l_j}{[L_{\V}:K_v]}\,(s_j+i\,\acc{\rho_j}{\chi_v})}}.
\end{multline}
Soit $\norm{.}_{e}$ 
la norme sur 
$X(T)^{G_v}_{\C}$
d\'efinie par 
\begin{equation}
\forall x\in X(T)^{G_v}_{\C},\quad \norm{x}_{e}=\sum_{l\in L} \abs{\acc{e_l}{x}}.
\end{equation}
De \eqref{eq:fvchisa}, on tire pour tout $\bs\in \tube{\R_{>0}^{\Sigma(1)/G}}$ 
et tout $\chi_v\in \dualt{T(K_v)/T(\Ov)}$
la majoration
\begin{equation}
\abs{f_v(\chi_v,\bs)}
\leq
\frac{2^{\rg\left(X(T)\right)}}
{1+\norm{\chi_v}_{e}}
\sum_{\sigma\in \Sigma^{G_v}_{\max}}
\prod_{j\in \sigma(1)/G_v} 
\frac{1+\norm{\sumu{j\in \sigma(1)/G_v} \frac{l_j}{[L_{\V}:K_v]}\,s_j\,\tau_j^{\vee}}_{e}}
{\abs{s_j+i\,\acc{\rho_j}{\chi_v})}}.
\end{equation}
L'\'equivalence des normes sur $X(T)^{G_v}_{\C}$ montre l'existence d'une constante
$C_1$ v\'erifiant, pour tout $\bs\in \C^{\Sigma(1)/G}$,
\begin{equation}
\norm{\sumu{j\in \sigma(1)/G_v}  \frac{l_j}{[L_{\V}:K_v]} s_j
\,\tau_j^{\vee}}_{e}
\leq C_1 \sumu{j\in \sigma(1)/G_v} \abs{\frac{l_j}{[L_{\V}:K_v]}s_j}
\leq C_1 \sumu{j\in \sigma(1)/G_v} \abs{s_j}
.
\end{equation}
Il existe par ailleurs une constante $C_2>0$
telle qu'on ait pour tout 
$\bs\in \tube{\K}$, pour tout $\sigma\in \Sigma^{G_v}_{\max}$ 
et tout $j\in \sigma(1)/G_v$
\begin{equation}
\abs{s_j+i\,\acc{\rho_j}{\chi_v})}\geq C_2 (1+\abs{\acc{\rho_j}{\chi_v}+\Im(s_j)}).
\end{equation}
Le point \ref{item:maj:prop:estim:transfo:locale:archi} s'en d\'eduit ais\'ement.
\end{demo}

On notera
\nindex{$f_{\infty}(\chi,\,.\,)$}
\begin{equation}
f_{\infty}(\chi,\,.\,)
=
\produ{v\in\placesde{K,\infty}}
f_{v}(\chi,\,.\,).
\end{equation}

On pose
\nindex{$X(T)_{\R,\infty}$}
\begin{equation}
X(T)_{\R,\infty}=\bigoplus_{v\in \placesde{K,\infty}}X(T)^{G_v}_{\R}.
\end{equation}

Soit $\chi$ un \'el\'ement de $\dualtop{T(\ak)/\K(T)}$.
Pour toute place $v$ de $K$, $\chi_v$ 
est alors un \'el\'ement de $\dualtop{T(K_v)/T(\Ov)}$.
Si $v$ est archim\'edienne, on a un isomorphisme
\begin{equation}
\dualt{\deg_{T,v}}\,:\,
\X(T)^{G_v}_\R
\longisom
\dualtop{T(K_v)/T(\Ov)}
\end{equation}
\dindex{morphisme <<type \`a l'infini>>}
\begin{defi}
Le morphisme <<type \`a l'infini>> est l'application
qui \`a $\chi\in \dualtop{T(\ak)/\K(T)}$ associe
\begin{equation}
\chi_{\infty}\eqdef\left(\left(\dualt{\deg_{T,v}}\right)^{-1}\chi_v\right)\in X(T)_{\R,\infty}.
\end{equation}
\end{defi}

Soit
$
\Sigma^{\infty}
$ 
l'\'eventail produit des 
$\Sigma^{G_v}$ pour
$v\in \placesde{K,\infty}$ : c'est l'\'eventail de
$X(T)_{\R,\infty}$ dont les c\^ones sont des produits
des c\^ones des \'eventails $\Sigma^{G_v}$
pour $v\in \placesde{K,\infty}$.
L'ensemble $\Sigma^{\infty}_{\max}$ des c\^ones de 
$\Sigma^{\infty}$ de dimension maximale est
donc l'ensemble des produits d'\'el\'ements 
de $\Sigma^{G_v}_{\max}$ pour
$v\in \placesde{K,\infty}$.

\begin{cor}\label{cor:maj:fchiv:varch:2}
Soit $\eps$ v\'erifiant $0<\eps<1$ et $\K$ un compact de $\R^{\Sigma(1)/G}_{>1-\eps}$.
Il existe une constante $C>0$ telle qu'on ait
\begin{multline}\label{eq:maj:fchiv:varch:2}
\forall \bs\in \tube{\K},\quad \forall \chi\in \dualtop{T(\ak)/\K(T)},\\
\abs{f_{\infty}(\chi,\bs)}
\leq
\frac{C}
{1+\norm{\chi_{\infty}}}
\,
\sum_{\wt{\sigma}\in \Sigma^{\infty}_{\max}}
\frac
{
1+\sumu{i\in \wt{\sigma}(1)}\abs{\Im(s_i)}
}
{
\produ{i\in \wt{\sigma}(1)}
(1+\abs{\acc{\roi}{\chi_\infty}+\Im(s_i)})
}
\end{multline}
\end{cor}
\begin{demo}
On applique le point \ref{item:maj:prop:estim:transfo:locale:archi} 
de la proposition \ref{prop:estim:transfo:locale:archi} en prenant pour \'eventail $\Sigma^{\infty}$.
\end{demo}

\subsubsection{Cas fonctionnel}\label{subsubsec:transfolocales:casfonc}

On introduit d'abord quelques d\'efinitions.

On rappelle qu'au moyen de la base $(\Da)_{\alpha\in \sg}$ on identifie
$\psg$ \`a $\Z^{\sg}$.
Ceci permet d'identifier $\psgc$ \`a $\csg$
et $\psgca$ \`a $\casg$ (donc \`a un sous-ensemble de $\csg$).

Dans toute la suite de ce texte, le terme \dindex{s\'erie formelle}\termin{s\'erie formelle} 
d\'esignera un \'el\'ement de $\C[[\za]]_{\alpha\in \sg}$
et le terme \dindex{mon\^ome}\termin{mon\^ome} un mon\^ome de $\C[\za]_{\alpha\in \sg}$. 
On identifiera un mon\^ome \`a la fonction $\csg\to \C$ qu'il induit.

Soit 
\begin{equation}
P=\sum_{(\na)\in \N^{\sg}} a_{(\na)}\prod \za^{\na}
\end{equation}
une s\'erie formelle. 
Soit $\bz$ un \'el\'ement de $\casg$ tel que la s\'erie d\'efinissant $P(\bz)$ converge
absolument. Si on voit $\bz$ comme un \'el\'ement de $\psgca$, on a donc
\begin{equation}\label{eq:P(bz)=}
P(\bz)=\sum_{(\na)\in \N^{\sg}}a_{(\na)} \prod \acc{\bz}{\Da^{\vee}}^{\na}
\end{equation}

\begin{lemme}\label{lm:sgmon}
Soit $(\na)$ un \'el\'ement de $\N^{\sg}$.
Alors l'application qui \`a $\bz\in \psgca=\casg$ associe $\acc{\bz}{\sum \na \Da^{\vee}}$ s'\'etend en un
unique mon\^ome sur $\csg$.
\end{lemme}
\begin{lemme}\label{lm:accgammaubzza=accbzroa}
Soit 
\begin{equation}
P=\sum_{(\na)\in \N^{\sg}} a_{(\na)}\prod \za^{\na}
\end{equation}
une s\'erie formelle. Alors pour tout $\bz\in \xtgca$ tel que la s\'erie
d\'efinissant $P(\gamma_{\ca}(\bz))$ converge absolument, on a 
\begin{equation}
P\left(\gamma_{\ca}(\bz)\right)=\sum_{(\na)\in \N^{\sg}} a_{(\na)}\prod \acc{\bz}{\roa}^{\na}.
\end{equation}
\end{lemme}
\begin{demo}
Compte tenu de \eqref{eq:P(bz)=},
il suffit de remarquer que pour $\alpha\in \sg$ on a
\begin{equation}
\acc{\gamma_{\ca}(\bz)}{\Da^{\vee}}
=\acc{\bz}{\gamma^{\vee}(\Da^{\vee})}
\end{equation}
et que, d'apr\`es \eqref{eq:gammaveeDaroa}, on a 
$\gamma^{\vee}(\Da^{\vee})=\roa$.
\end{demo}

\begin{defi}\label{defi:compatible}
Soit $M$ un sous-groupe de $\xtgd$.
\begin{enumerate}
\item
Soit $(\na)\in \Nsg$. 
Le mon\^ome $\produ{\alpha\in \sg}\za^{\na}$ 
est dit $M$-compatible si on a 
\begin{equation}
\gamma^{\vee}\left(\sum \na\Da^{\vee}\right)\in M,
\end{equation}
ce qui \'equivaut d'apr\`es \eqref{eq:gammaveeDaroa} \`a la condition
\begin{equation}
\sum \na\,\roa \in M.
\end{equation}
\item
Une s\'erie formelle est dite $M$-compatible si les mon\^omes
qui apparaissent dans son \'ecriture sont $M$-compatibles.
\end{enumerate}
\end{defi}
On a le lemme \'el\'ementaire suivant.
\begin{lemme}\label{lm:sgmon:comp}
Soit $M$ un sous-groupe de $\xtgd$. 
\begin{enumerate}
\item
Un mon\^ome $f$ est $M$-compatible si et seulement s'il existe
un \'el\'ement $m$ de $M$ v\'erifiant
\begin{equation}
\forall \bz\in \xtgca,\quad 
f\left(\gamma_{\ca}(\bz)\right)=\acc{\bz}{m}
\end{equation}
\item
Tout produit ou somme de s\'erie formelles $M$-compatibles est $M$-compatible.
\end{enumerate}
\end{lemme}
\begin{lemme}\label{lm:za^daDT0comp}
Pour tout $\alpha\in \sg$, le mon\^ome $\za^{\,\da}$ est $\DT^0$-compatible.
\end{lemme}
\begin{demo}
C'est imm\'ediat d'apr\`es la d\'efinition et le fait, donn\'e par le lemme
\ref{lm:carac:fonc}, que $\da\,\roa$ est dans $\DT^0$.
\end{demo}

\begin{defi}\label{defi:convergence}
Le rayon de convergence d'une s\'erie formelle $P$ est le plus grand r\'eel positif $R$ tel que 
pout tout $(\za)\in \csg$ v\'erifiant $\Max_{\alpha} \abs{\za}<R$, la s\'erie d\'efinissant $P(\za)$
est absolument convergente.
\end{defi}

\begin{lemme}\label{lm:rel:hauteur:carac:fonc}
Soit $v$ une place de $K$ v\'erifiant l'hypoth\`ese \ref{hyp:hinz}. 
\begin{enumerate}
\item
On a 
\begin{equation}\label{eq:chibz}
\forall\,\bz\in  \xtgu,
\quad
\forall t\in T(K_v),
\quad
\left(\chi_{\bz}\right)_v(t)
=
\frH_{v}
\left(\gamma_{\ca}(\bz),t\right).
\end{equation}
\item
Pour tout $t\in T(K_v)$, la fonction $\frH_v(\,.\,,t)$ s'\'etend en 
un mon\^ome qui est $\DT$-compatible. Ce mon\^ome est une constante
si et seulement si $t$ est dans $T(\Ov)$.
\end{enumerate}
\end{lemme}

\begin{demo}
Pour $\bz\in \xtgca$,
et $t\in T(K_v)$, on a,
compte tenu du lemme \ref{lm:degTv:degTLV:fonc}
et de la d\'efinition \eqref{eq:hauteur:locale:log:bis} de $\frH_v$
\begin{align}
\acc{\bz}{\deg_{T,v}(t)}^{f_v}&
=\left[\accsv{\gamma_{\ca} (\bz)}{t}\right]^{\,f_v}\\
&=\frH_v(\gamma_{\ca} (\bz),t).
\end{align}
Or, pour $\bz\in\xtgu$, on a d'apr\`es \eqref{eq:chizv}
\begin{equation}
(\chi_{\bz})_v(t)=\acc{\bz}{\deg_{T,v}(t)}^{\,f_v}
\end{equation}
d'o\`u la relation \eqref{eq:chibz}.

Soit $t\in T(K_v)$.
D'apr\`es le lemme \ref{lm:xtgv}
il existe un \'el\'ement $\sigma\in \Sigma^{\,G_v}$ tel que $\deg_{T,L,\V}(t)$ s'\'ecrit
\begin{equation}
\deg_{T,L,\V}(t)=\sum_{\beta\in \sigma(1)/G_v} \nb\,\taub 
\end{equation}
o\`u les $\nb$ sont des entiers strictement positifs.
Pour tout \'el\'ement $\varphi$ de $\plsg$ (identifi\'e \`a $\psg$), on a alors
\begin{equation}\label{eq:fvaccs=}
f_v\,\acc{\varphi}{t}=
\frac{f_v}{e_v}\,\accsv{\varphi}{\deg_{T,L,\V}(t)}
=\frac{f_v}{e_v}\sum_{\beta\in \sigma(1)/G_v} \nb\,\lb\,\acc{\varphi}{\Db^{\vee}}.
\end{equation}
Ainsi la condition $v$ v\'erifie l'hypoth\`ese \ref{hyp:hinz} entra\^\i ne qu'on a
\begin{equation}
\forall \beta\in \sigma(1)/G_v,\quad \frac{f_v\,\nb\,\lb}{e_v}\in \Z.
\end{equation}
On a alors, d'apr\`es \eqref{eq:fvaccs=} et la d\'efinition \eqref{eq:hauteur:locale:log:bis} de $\frH_v$, 
\begin{equation}
\forall \bz\in \psgca,\quad \frH_v(\bz,t)
=
\acc
{\bz}
{\sum_{\beta\in \sigma(1)/G_v}
\frac{f_v\,\nb\,\lb}{e_v}\,\Db^{\vee}
}
\end{equation}
ce qui montre, d'apr\`es le lemme \ref{lm:sgmon}, que $\frH_v(\,.\,,t)$ s'\'etend en un mon\^ome
sur $\csg$. Comme les $\nb$ sont strictement positifs, 
ce mon\^ome est une constante si et seulement si $\sigma(1)/G_v$ est vide,
c'est-\`a-dire si et seulement si $\deg_{T,L,\V}(t)$ est nul, soit encore
si et seulement si $t\in T(\Ov)$.

La relation
\begin{equation}
\forall \bz\in \xtgca,\quad \frH_v(\gamma_{\ca} (\bz),t)=\acc{\bz}{f_v\,\deg_{T,v}(t)},
\end{equation}
le fait que $\deg_{T,v}(t)$ soit dans $\DT$, et le lemme \ref{lm:sgmon:comp} montrent 
que $\frH_v(\,.\,,t)$ est un mon\^ome $\DT$-compatible.
\end{demo}

\begin{cor}\label{cor:hvsaxchibz}
Soit $v$ une place de $K$ v\'erifiant l'hypoth\`ese \ref{hyp:hinz}. 
Pour tout $t\in T(K_v)$, on a
\begin{equation}\label{eq:cor:hvsaxchibz:bis}
\forall\,\bz\in  \xtgu,
\quad
\forall\,\bs\in \plsgc,
\quad
H_v(\bs,t)
\,
\left(\chi_{\bz}\right)_v(t)
=
\frH_{v}
\left(\gamma_{\unit}(\bz)\,q^{\,\bs},t\right)
\end{equation}
\end{cor}

\begin{prop}\label{prop:frv:fv}
Soit $v$ une place de $K$ v\'erifiant l'hypoth\`ese \ref{hyp:hinz}
et $\chi$ un \'el\'ement de $\dualtop{T(\ak)/\K(T)}$.
La s\'erie formelle
\nindex{$\frf_v(\chi,\,.\,)$}
\begin{equation}\label{eq:demo:prop:frv:fv}
\frf_v(\chi,\,.\,)
=
\left(\int\limits_{T(\Ov)}\!d\mu_v\right)
\,
\sum_{t\in T(K_v)/T(\Ov)}
\chi_v(t)\,\frH_{v}(\,.\,,t)
\end{equation}
a un rayon de convergence sup\'erieur \`a $1$, 
est $\DT$-compatible, et v\'erifie
\begin{equation}\label{eq:prop:frv:fv}
\forall \bs\in \tube{\R^{\,\Sigma(1)/G}_{>0}},
\quad
\frf_v\left(\chi,q^{\,-\bs}\right)=f_v(\chi,\bs)
\end{equation}
\end{prop}
\begin{demo}
Comme le morphisme $\deg_{T,L,\V}$ induit une injection de $T(K_v)/T(\Ov)$ dans $\left(X(T)^{\vee}\right)^{\,G_v}$, 
$T(K_v)/T(\Ov)$ est un groupe ab\'elien de type fini, donc un mono\"\i de commutatif de type fini.
Par ailleurs, pour tout $t\in T(K_v)/T(\Ov)$ et tout $n\in \N$, on a
\begin{equation}
\frH_{v}(\,.\,,t^n)=\frH_{v}(\,.\,,t)^n,
\end{equation}
et d'apr\`es le lemme \ref{lm:rel:hauteur:carac:fonc}, $\frH_{v}(\,.\,,t)$ est une constante si et seulement si $t=0$.
Ce qui pr\'ec\`ede montre que pour tout entier $d$ il n'existe qu'un nombre fini de $t\in T(K_v)/T(\Ov)$
tel que le mon\^ome $\frH_v(\,.\,,t)$ soit de degr\'e total major\'e par $d$. Ainsi $\frf_v(\chi,\,.\,)$
est bien une s\'erie formelle.

Le fait que le rayon de convergence soit sup\'erieur \`a $1$ 
et la relation \eqref{eq:prop:frv:fv} viennent
de la remarque \ref{rem:fv=} et de \eqref{eq:hvphi=frhvqphi}.
La $\DT$-compatibilit\'e d\'ecoule aussit\^ot de la 
$\DT$-compatibilit\'e des $\frH_v(\,.\,,t)$
(lemme \ref{lm:rel:hauteur:carac:fonc}).
\end{demo}

\begin{lemme}
\label{lm:transfo:ram:fonc}
Soit $v$ une place de $K$ v\'erifiant l'hypoth\`ese \ref{hyp:hinz}.   
Soit $\chi$ un \'el\'ement de $\dualtop{T(\ak)/\K(T)}$.
Alors on a 
\begin{multline}\label{eq:lm:transfo:ram:fonc:bis}
\forall\,\bz\in \xtgu,
\quad
\forall\,\bs
\in 
\tube{\R^{\,\Sigma(1)/G}_{>0}},\\
f_v(\chi.\chi_{\bz},\bs)
=
\frf_v
\left(
\chi,
\gamma_{\unit}(\bz)
\,
q^{\,-\bs} 
\right).
\end{multline}
\end{lemme}
\begin{demo}
Par \eqref{eq:prop:frv:fv}, le membre de gauche de \eqref{eq:lm:transfo:ram:fonc:bis}
vaut $\frf_v
\left(
\chi.\chi_{\bz},
\,
q^{\,-\bs} 
\right)$,
et on a d'apr\`es \eqref{eq:demo:prop:frv:fv}
et \eqref{eq:chibz}
\begin{align}
\frf_v
\left(
\chi.\chi_{\bz},
\,
q^{\,-\bs} 
\right)
&
=
\left(\int\limits_{T(\Ov)}d\mu_v\right)
\,
\sum_{t\in T(K_v)/T(\Ov)}
\chi_v(t)\,(\chi_{\bz})_v(t)\,\frh_{v}(q^{\,-\bs},t)
\\
&
=
\left(\int\limits_{T(\Ov)}d\mu_v\right)
\,
\sum_{t\in T(K_v)/T(\Ov)}
\chi_v(t)\,\frh_{v}(\gamma_{\unit}(\bz)\,q^{\,-\bs},t)
\\
&
=
\frf_v
\left(
\chi,
\gamma_{\unit}(\bz)\,q^{\,-\bs}
\right)
\end{align}
\end{demo}

On consid\`ere d\'esormais $S$ un ensemble 
fini de places de $K$ contenant toutes les places ramifi\'ees dans $L$,
et tel que toutes les places de $S$ v\'erifient l'hypoth\`ese \ref{hyp:hinz}.
On note  
\begin{equation}
\K(T)^S=\produ{v\in S}T(\Ov)
\end{equation}
et
\begin{equation}
C_S=\int\limits_{\adh{T(K)}^{\,\,S}\cap \K(T)^S}
\left(\otimesu{v\in S}d\mu_v \right).
\end{equation}
Pour $t\in T(\ak)$, on note \nindex{$\frH_S(\,.\,,t)$}$\frH_S(\,.\,,t)$ le mon\^ome d\'efini par 
\begin{equation}
\frH_S(\,.\,,t)=\prod_{v\in S} \frH_v(\,.\,,t_v).
\end{equation}
\begin{lemme}\label{lm:frfS}
Soit $S$ un ensemble 
fini de places de $K$ contenant toutes les places ramifi\'ees dans $L$,
et tel que toutes les places de $S$ v\'erifient l'hypoth\`ese \ref{hyp:hinz}.
La 
s\'erie formelle 
$\frf_S$
d\'efinie par
\begin{equation}
\frf_S
=
C_S
\,
\sum_{t\in \adh{T(K)}^{\,\,S}/\left(\adh{T(K)}^{\,\,S}\cap \K(T)^S\right)}
\,\frH_{S}(\,.\,,t)
\end{equation}
a un rayon de convergence sup\'erieur \`a $1$, est $\DT^0$-compatible
et v\'erifie, pour tout $\bs \in \tube{\R^{\Sigma(1)/G}_{>0}}$
\begin{equation}
\int\limits_{\adh{T(K)}^{\,\,S}}\,
H(-\bs,t)
\chi_{\bz}(t)
\left(\otimesu{v\in S}d\mu_v \right)(t)
=
\frf_S\left(\gamma_{\unit}(\bz)\,q^{\,-\bs}\right)
\end{equation}
\end{lemme}
\begin{demo} 
Comme $\adh{T(K)}^{\,\,S}/\left(\adh{T(K)}^{\,\,S}\cap \,\K(T)^S\right)$
s'injecte dans $\produ{v\in S}T(K_v)/T(\Ov)$, un argument similaire \`a celui 
du d\'ebut de la proposition \ref{prop:frv:fv} montre que $\frf_S$ est bien
une s\'erie formelle.

Pour $\bs \in \tube{\R^{\Sigma(1)/G}_{>0}}$, la convergence 
de la s\'erie
\begin{equation}\label{eq:laserie}
\sum_{t\in \adh{T(K)}^{\,\,S}/\left(\adh{T(K)}^{\,\,S}\cap \,\K(T)^S\right)}
\,\prod_{v\in S} \frH_v\left(q^{-\Re(\bs)},t_v\right)
\end{equation}
se montre en injectant  $\adh{T(K)}^{\,\,S}/\left(\adh{T(K)}^{\,\,S}\cap \,\K(T)^S\right)$
dans $\produ{v\in S}T(K_v)/T(\Ov)$, ce qui permet de majorer \eqref{eq:laserie}
par 
\begin{equation}
\prod_{v\in S} \sum_{t\in T(K_v)/T(\Ov)} H_v\left(-\Re(\bs),t_v\right),
\end{equation}
qui est un produit de s\'eries convergentes (cf. remarque \ref{rem:fv=}).

On a ainsi, toujours pour  $\bs \in \tube{\R^{\Sigma(1)/G}_{>0}}$
\begin{equation}
\int\limits_{\adh{T(K)}^{\,\,S}}\,
H(-\bs,t)\chi_{ \bz }(t)
\left(\otimesu{v\in S}d\mu_v \right)(t)\\
\end{equation}
\begin{align}
&=
C_S
\,
\sum_{t\in \adh{T(K)}^{\,\,S}/\left(\adh{T(K)}^{\,\,S}\cap \,\K(T)^S\right)}
\,H(-\bs,t)\chi_{ \bz }(t)
\\
&=
C_S
\,
\sum_{t\in \adh{T(K)}^{\,\,S}/\left(\adh{T(K)}^{\,\,S}\cap \,\K(T)^S\right)}
\,H(-\bs,t)\acc{\bz}{\deg_T(t)}
\\
&=
C_S
\,
\sum_{t\in \adh{T(K)}^{\,\,S}/\left(\adh{T(K)}^{\,\,S}\cap \,\K(T)^S\right)}
\,\prod_{v\in S} \frH_v(\gamma_{\unit} (\bz)\,q^{-\bs},t_v)
\end{align}

Pour terminer la d\'emonstration du lemme, il suffit de montrer que pour 
$t\in \adh{T(K)}^{\,\,S}$, le mon\^ome $\frH_S(\,.\,,t)$
est $\DT^0$ compatible.
Soit donc $t=(t_v)_{v\in S}\in \adh{T(K)}^{\,\,S}$.
D'apr\`es la proposition \ref{lm:appfaible}, 
il existe $u\in T(K)$ et $t'\in T(\ak)$ un \'el\'ement de l'image de $T_{\Ps}(\ak)$ par 
$\gamma_{\,\ak}$ 
tels que $t=u\,t'$. 
Comme $\deg_T$ est trivial sur $T(K)$, on a $\deg_T(t)=\deg_{T}(t')$.
Comme  $t'$ est dans $\gamma_{\ak}(T_{\Ps}(\ak))$,  $\deg_T(t')$ est dans $\DT^0$.
Ainsi $\deg_T(t)$ est dans $\DT^0$.
Comme on a, pour tout $\bz\in \xtgc$,
\begin{equation}
\frH_S(\bz,t)=\acc{\bz}{\deg_{T}(t)}
\end{equation}
on voit, d'apr\`es le lemme \ref{lm:sgmon:comp},  que $\frH_S(\,.\,,t)$ est $\DT^0$-compatible.
\end{demo}

\subsubsection{Cas des places non ramifi\'ees}
\label{subsubsec:nonram}

Soit $v$ une place finie de $K$, $\V$ une place de $L$ divisant $v$ et $G_v$
son groupe de d\'ecomposition.
Les auteurs de \cite{BaTs:anis} d\'efinissent\footnote{
la d\'efinition utilis\'ee dans \cite{BaTs:anis} est en fait l\'eg\`erement diff\'erente 
car on y utilise le polyn\^ome $Q_{\Sigma,v}\left(X_{\beta}^{\,\lb }\right)$, 
ceci \'etant on obtient bien la m\^eme formule pour la transform\'ee de Fourier locale.
}
alors un polyn\^ome \`a coefficients entiers $Q_{\Sigma,v}$ 
en les $\card{\Sigma(1)/G_v}$ ind\'etermin\'ees 
$\left(X_{\beta}\right)_{\beta\in \Sigma(1)/G_v}$ 
par la formule
\begin{equation}
\left(\int\limits_{T(\Ov)}d\mu_v\right)
\sum_{\sigma\in \Sigma^{\,G_v}}
\prod_{\beta\in \sigma(1)/G_v}
\frac
{X_{\beta}}
{1-X_{\beta}}
=\frac
{Q_{\Sigma,v}(X_{\beta})}
{\underset{\beta\in \Sigma(1)/G_v}{\prod}\left(1-X_{\beta}\right)}.
\end{equation}

Un r\'esultat facile mais important pour les r\'esultats de convergence 
et de prolongement de la fonction z\^eta des hauteurs est (cf. \cite[Proposition 2.2.3]{BaTs:anis})

\begin{prop}[Batyrev, Tschinkel]\label{prop:convergence}
Le polyn\^ome $Q_{\Sigma,v}\left(X_{\beta}^{\lb }\right)-1$ 
ne contient que des mon\^omes de degr\'e total sup\'erieur ou \'egal \`a 2.
\end{prop}

D\'esormais nous fixons un sous-ensemble fini $S$
de $\placesde{K}$ contenant les places archim\'ediennes, 
et les places ramifi\'ees dans $L/K$.
Le r\'esultat suivant (\cite[Theorem 2.2.6]{BaTs:anis}) donne alors une expression explicite des transform\'ees 
de Fourier locales aux places $v\notin S$.

\begin{thm}[Batyrev, Tschinkel]
\label{thm:transfo:nonram}
Soit $v\notin S$ une place de $K$. 
Soit $\bs\in \tube{\R^{\,\Sigma(1)/G}_{>0}}$.
On a pour tout $\chi\in \dualtop{T(\ak)/\K(T)}$
\begin{multline}
\int\limits_{T(K_v)}\,H_{v}(-\bs,t)\chi_v(t)d\mu_v(t)\\
=\left(
\prod_{\beta\in \Sigma(1)/G_v}\,
\frac
{1}
{1-
\chib\left(\pi_{\wb}\right)
\,q_v^{\,-\lb \,\sb}}
\right)
\,Q_{\Sigma,v}
\left(
\chib\left(\pi_{\wb}\right)
\,q_v^{-\lb \sb}
\right)_{\beta\in\Sigma(1)/G_v}.
\end{multline}
\end{thm}

\begin{demo}
Nous rappelons la preuve de ce th\'eor\`eme, 
donn\'ee dans \cite{BaTs:anis}.
\`A l'aide du morphisme $\deg_{T,L,\V}$, on identifie $T(K_v)/T(\Ov)$ \`a $(X(T)^{\,\vee})^{\,G_v}$ 
(cf. proposition \ref{prop:draxl}). 
Comme les fonctions $H_{v}(-\bs,\,.\,)$ et $\chi_v$ sont $T(\Ov)$ invariante, 
elles induisent des fonctions sur $(X(T)^{\,\vee})^{\,G_v}$, 
qui seront not\'ees de fa\c con identique.
On a alors
\begin{equation}
\int\limits_{T(K_v)}\,H_{v}(-\bs,t)\chi_v(t)d\mu_v(t)
=
\int\limits_{T(\Ov)}d\mu_v
\,\sum_{t\in T(K_v)/T(\Ov)}\,H_{v}(-\bs,t)\chi_v(t)  
\end{equation}

Le lemme \ref{lm:xtgv} permet alors d'\'ecrire
\begin{align} 
&\quad\sum_{t\in T(K_v)/T(\Ov)}\,H_{v}(-\bs,t)\chi_v(t)\\
&
=
\sum_{\sigma\in\Sigma^{G_v}}\,\,
\sum_{t\,\in \,\text{intrel}(\sigma)\cap (T(T)^{\,\vee})^{\,G_v}}
H_{v}(-\bs,t)\,\chi_v(t)\\
&=\sum_{\sigma\in\Sigma^{G_v}}\,\,
\sum_{(\nb)\in \Z_{>0}^{\sigma(1)/G_v}}
H_{v}\left(-\bs,\sum \nb \,\taub\right)\prod \chi_v\left(\taub\right)^{\nb}.
\end{align}
Mais on a 
\begin{align}
&\quad \sum_{\sigma\in\Sigma^{G_v}}\,\,
\sum_{(\nb)\in \Z_{>0}^{\sigma(1)/G_v}}
H_{v}\left(-\bs,\sum \nb \,\taub\right)\prod \chi_v\left(\taub\right)^{\nb} 
\\
&=
\label{eq:gammavtaub}
\sum_{\sigma\in\Sigma^{G_v}}\,\,
\sum_{(\nb)\in \Z_{>0}^{\sigma(1)/G_v}}
H_{v}\left(-\bs,\sum \nb \,\taub\right)\prod \chib\left(\pi_{\wb}\right)^{\nb} 
\\
&=
\sum_{\sigma\in\Sigma^{G_v}}\,\,
\prod_{\beta\in \sigma(1)/G_v}\,
\label{eq:avantder}
\sum_{\nb \in \Z_{>0}} 
q_{v}^{\,-\nb \,\lb \,\sb}\,\chi_\beta\left(\pi_{\wb}\right)^{\,\nb }\\
&=
\sum_{\sigma\in\Sigma^{G_v}}\,\,
\prod_{\beta\in \sigma(1)/G_v}\,
\frac
{
\chi_\beta\left(\pi_{\wb}\right)
\,q_v^{\,-\lb \,\sb}}
{1-
\chi_\beta\left(\pi_{\wb}\right)
\,q_v^{\,-\lb \,\sb}}.
\end{align}
d'o\`u le r\'esultat par d\'efinition de $Q_{\Sigma,v}$. 

L'\'egalit\'e \eqref{eq:gammavtaub} vient du fait que si $v$ n'est pas ramifi\'ee dans $L$,
d'apr\`es le lemme \ref{lm:jvgammavpiwj}, on a
\begin{equation}
\deg_{T,L,\V}\left[\gamma_v\circ i_{\beta,v} (\pi_{\wb})\right]=\taub.
\end{equation}
L'\'egalit\'e \eqref{eq:avantder} vient 
de ce que par la d\'efinition \eqref{eq:hauteur:locale}
de la hauteur locale 
et le fait que $v$ est non ramifi\'ee on  a
\begin{align}
H_{v}
\left(
-\bs
,
\sum_{\beta\in \sigma(1)/G_v}
\nb \,\taub 
\right)&
=
\exp\left(\,
\log(q_v)\,
\,\,\underset{\beta\in \sigma(1)/G_v}{\sum}
\,\nb \,\accsv{\varphi}{\taub}
\right)
\\
&
=
\exp\left(\,
\log(q_v)\,
\,\,
\underset{\beta\in \sigma(1)/G_v}{\sum}\,\nb \,\lb \,\sb
\right).
\end{align}
\end{demo}

\begin{lemme}
\label{lm:transfo:nonram}
\begin{enumerate}
\item
Pour tout $\alpha\in \Sigma(1)/G$, on a 
\label{item:1:lm:transfo:nonram}
\begin{equation}
\prod_{\beta\in \alpha/G_v}\,
\frac
{1}
{1-
\chi_{\beta}\left(\pi_{\wb}\right)
\,
q_v^{\,-\lb \,\sb}}
=
\prod_{\substack{w\in \placesde{\Ka},\\ w|v}}
\frac
{1}
{1-\chia\left(\pi_w\right)\,q_w^{-\,\sa}}.
\end{equation}
En particulier, dans le cas fonctionnel, on a 
\begin{equation}
\prod_{\beta\in \alpha/G_v}\,
\frac
{1}
{1-
\chib\left(\pi_{\wb}\right)
\,
q_v^{\,-\lb \,\sb}}
=
\prod_{\substack{w\in \placesde{\Ka},\\ w|v}}
\frac
{1}
{1-\chia\left(\pi_w\right)\,q^{-f_w\,\da\,\sa}}.
\end{equation}
\item
\label{item:2:lm:transfo:nonram}
Dans le cas arithm\'etique,
pour $\chi\in \dualtop{T(\ak)/\K(T)}$ et $y\in X(T)^G_{\R}$, on a
\begin{equation}
\left(\chi.\chi_y\right)_{\alpha}\left(\pi_w\right)\,q_w^{-\sa}
=
\chia\left(\pi_w\right)\,q_w^{-\sa+i\,\acc{y}{\roa}}.
\end{equation}
\item
\label{item:3:lm:transfo:nonram}
Dans le cas fonctionnel, 
pour $\chi\in \dualtop{T(\ak)/\K(T)}$ et $\bz\in\xtgu$, on a
\begin{equation}
\left(\chi.\chi_{\bz}\right)_{\alpha}\left(\pi_w\right)\,q^{-f_w\,\da\,\sa}
=
\chia
\left(\pi_w\right)\,
\acc{\bz}{\da\,\roa} q^{\,-f_w\,\da\,\sa}.
\end{equation}
\end{enumerate}
\end{lemme}
\begin{demo}
Le point \ref{item:1:lm:transfo:nonram} provient des formules \eqref{eq:form:qwb} 
et \eqref{eq:form:fvlb}. 
Le point \ref{item:2:lm:transfo:nonram} d\'ecoule du lemme \ref{lm:carac:arit}
et le point \ref{item:3:lm:transfo:nonram} du lemme \ref{lm:carac:fonc}.
\end{demo}

Le r\'esultat suivant sera utile lors du calcul
du terme principal de la fonction z\^eta des hauteurs.
\begin{lemme}\label{lm:inttkvhv-s}
Soit $v\notin S$. Alors pour tout $s\in \tube{\R_{>0}}$, on a
\begin{equation}
\int\limits_{T(K_v)}\,H_{v}(-s\,\varphi_0,t)\,d\mu_v(t)\\
=
L_v(s,\Ps)
\,
Q_{\Sigma,v}
\left(
q_v^{-\lb s}
\right)_{\beta\in\Sigma(1)/G_v}.
\end{equation}
\end{lemme}
\begin{demo}
Ceci d\'ecoule aussit\^ot 
d'une part du th\'eor\`eme \ref{thm:transfo:nonram} et du lemme \ref{lm:transfo:nonram},
d'autres part des lemmes \ref{lm:lvzgh}, \ref{lm:comp:L:exseq} et de la d\'ecomposition 
\eqref{eq:idps}.
\end{demo}

\subsection{Propri\'et\'es analytiques de la transform\'ee de Fourier globale}\label{subsec:ptes:global}
On va d\'eduire des r\'esultats de la section \ref{subsec:transfo:four:locales}
que pour tout \'el\'ement  
$\bs$ de $\tube{\R^{\,\sg}_{>1}}$, 
la fonction $H(-\bs,\,.\,)$ est int\'egrable sur $T(\ak)$, 
ainsi que des renseignements sur le comportement analytiques des transform\'ees de Fourier globales.

Pour $\alpha\in\sg$,
nous notons $S_{\alpha}$ les places de $\Ka$ au-dessus des places de $S$.

\begin{lemme}\label{lm:prodvQsv}
Pour tout $\bs$ dans 
$\tube{\R^{\,\Sigma(1)/G}_{>\frac{1}{2}}}$   
le produit eul\'erien
\begin{equation}
\prod_{v\notin S}
Q_{\Sigma,v}
\left(
q_v^{-\lb \sb}
\right)_{\beta\in \Sigma(1)/G_v}
\end{equation}
est absolument convergent. 
\end{lemme}
\begin{demo}
Ceci d\'ecoule aussit\^ot de la proposition \ref{prop:convergence}.
\end{demo}

\begin{lemme}\label{lm:H-bsint}
Pour tout
$\bs\in \tube{\R^{\,\Sigma(1)/G}_{>1}}$,
$H(-\bs,\,.\,)$ est int\'egrable sur $T(\ak)$.
\end{lemme}

\begin{demo}
On a pour tout $\bs\in \C^{\,\Sigma(1)/G}$
\begin{equation}
\int\limits_{T(\ak)}
\abs{H(-\bs,t)}\,
\omega_{T}(t)
=
\int\limits_{T(\ak)}H(-\Re(\bs),x)\,
\omega_{T}(t).
\end{equation}
D'apr\`es la d\'efinition \eqref{eq:defomegaT} de $\omega_{T}$
cette derni\`ere expression vaut
\begin{equation}\label{eq:exprtransfo:arit}
c_{K,\dim(T)}\,\prod_{v\in \placesde{K}}\,\,\,\int\limits_{T(K_v)}H_{v}(-\Re(\bs),t)\,d\mu_v(t).
\end{equation}
D'apr\`es le th\'eor\`eme \ref{thm:transfo:nonram} et le lemme \ref{lm:transfo:nonram}, 
on a donc
\begin{multline}
\int\limits_{T(\ak)}
\abs{H(-\bs,t)}\,
\omega_{T}(t)\\
=
c_{K,\dim(T)}
\left(\prod_{\alpha\in \Sigma(1)/G} 
\zeta_{\Ka}
\left(\Re(\sa)\right)
\right)
\,
\left(
\prod_{\alpha \in \Sigma(1)/G} 
\prod_{w\in S_{\alpha}} (1-q_w^{-\Re(\sa)})
\right)
\\
\,
\prod_{v\notin S}
Q_{\Sigma,v}
\left(
q_v^{-\lb \Re(\sb)}
\right)_{\beta\in \Sigma(1)/G_v} 
\,
\prod_{v\in S} \,\,\,\int\limits_{T(K_v)}H_{v}(-\Re(\bs),t)\,d\mu_v(t) 
\end{multline}
D'apr\`es la proposition \ref{prop:fonc:zeta}, le lemme \ref{lm:prodvQsv},
le point \ref{item:1:prop:fv:place:finie:quelc} de la proposition \ref{prop:fv:place:finie:quelc}
et, dans le cas arithm\'etique, 
le point \ref{item:1:prop:estim:transfo:locale:archi} de la proposition \ref{prop:estim:transfo:locale:archi},
cette derni\`ere expression est finie d\`es que
$\Re(\bs)\in \R^{\,\Sigma(1)/G}_{>1}$
d'o\`u l'int\'egrabilit\'e de $H(-\bs,\,.\,)$ pour 
$\bs\in \tube{\R^{\,\Sigma(1)/G}_{>1}}$.
\end{demo}

\subsubsection{Cas arithm\'etique}

\begin{lemme}\label{lm:fchihol:arit}
Soit $\chi\in \dualtop{T(\ak)/\K(T)}$.
La fonction
\begin{multline}
\bs
\mapsto
c_{K,\dim(T)}
\prod_{v\notin S}
Q_{\Sigma,v}
\left(
\chib(\pi_{\wb})
\,
q_v^{-\lb \sb}
\right)_{\beta\in \Sigma(1)/G_v}
\\
\prod_{v\in S\setminus \placesde{K,\infty}}
f_v(\chi,\bs)
\prod_{\alpha \in \Sigma(1)/G} 
\,
\prod_{w\in S_{\alpha}} (1-\chia(\pi_w)\,q_w^{-\sa})
\end{multline}
est holomorphe sur
$
\tube{\R^{\,\Sigma(1)/G}_{>\frac{1}{2}}}.
$
\nindex{$f(\chi,\,.\,)$}
On la note $f(\chi,\,.\,)$.
\end{lemme}
\begin{demo}
Ceci d\'ecoule du lemme \ref{lm:prodvQsv} et 
du point \ref{item:2:prop:fv:place:finie:quelc} de la proposition \ref{prop:fv:place:finie:quelc}.
\end{demo}

Rappelons qu'on  a \'egalement not\'e
\begin{equation}
f_{\infty}(\chi,\,.\,)
=
\produ{v\in\placesde{K,\infty}}
f_{v}(\chi,\,.\,).
\end{equation}
\begin{lemme}\label{lm:finftychihol}
$f_{\infty}(\chi,\,.\,)$
est holomorphe sur
$
\tube{\R^{\,\Sigma(1)/G}_{>0}}.
$
\end{lemme}
\begin{demo}
D'apr\`es le point \ref{item:2:prop:estim:transfo:locale:archi} 
de la proposition \ref{prop:estim:transfo:locale:archi},
pour tout $v\in\placesde{K,\infty}$
la fonction $f_{v}(\chi,\,.\,)$ est holomorphe sur
$
\tube{\R^{\,\Sigma(1)/G}_{>0}}.
$
\end{demo}

\begin{lemme}
Soit $\chi\in \dualtop{T(\ak)/\K(T).T(K)}$.
Pour tout $\bs\in \tube{\R^{\,\Sigma(1)/G}_{>1}}$,
et tout $y\in X(T)^G_{\R}$, on a 
\begin{multline}\label{eq:expr:fourierHchisa:fchisa}
\fourier H\left(\chi.\chi_{y}, -\bs\right)\\
=
\left(\prod_{\alpha\in \Sigma(1)/G}
L_{\Ka}\left(\sa-i\,\acc{y}{\roa} ,\chia\right)\right)
f\left(\chi,\bs-i\,\gamma_{\R}(y) \right)
\,
f_{\infty}\left(\chi,\bs-i\,\gamma_{\R}(y)\right).
\end{multline} 
\end{lemme}
\begin{demo}
D'apr\`es la d\'efinition \eqref{eq:defomegaT} de $\omega_{T}$, on a
\begin{equation}\label{eq:intlimitst(ak)}
\int\limits_{T(\ak)}
(\chi\,\chi_y)(t) H(-\bs,t)\,\omega_{T}(t)
=
c_{K,\dim(T)}\prod_{v\in \placesde{K}}\,\,\int\limits_{T(K_v)}(\chi\,\chi_y)_v(t) H_{v}(-\bs,t)\,d\mu_v(t).
\end{equation}
D'apr\`es le th\'eor\`eme \ref{thm:transfo:nonram}, le lemme \ref{lm:transfo:nonram}
et le corollaire \ref{cor:transfo:ram:arit}, le membre de gauche de \eqref{eq:intlimitst(ak)}
est donc \'egal \`a 
\begin{gather}
c_{K,\dim(T)}
\prod_{\alpha\in \sg} 
L_{\Ka}(\chia,\sa-i\acc{y}{\roa})
\notag
\\
\times
\prod_{\alpha \in \Sigma(1)/G} 
\prod_{w\in S_{\alpha}} (1-\,\chia(\pi_w)\,q_w^{-\sa+i\,\acc{y}{\roa}})
\notag
\\
\times
\prod_{v\notin S}
Q_{\Sigma,v}
\left(
\chib(\pi_{\wb})\,q_v^{-\lb (\sb-i\,\acc{y}{\rob})}
\right)_{\beta\in \Sigma(1)/G_v} 
\notag
\\
\times
\prod_{v\in S} f_{v}(\chi,\bs-i\,\gamma_{\R}(y)),
\end{gather}
d'o\`u le r\'esultat.
\end{demo}

\begin{lemme}\label{lm:maj:fchifvchivfin}
Soit $\K$ un compact de $\R^{\Sigma(1)/G}_{>\frac{1}{2}}$.
Il existe une constante $C>0$ telle qu'on ait, pour
tout $\chi\in \dualtop{T(\ak)/\K(T)}$ et tout $\bs\in \tube{\K}$
la majoration
\begin{equation}
\abs{f(\chi,\bs)
}
\leq 
C.
\end{equation}
\end{lemme}
\begin{demo}
Ceci d\'ecoule imm\'ediatement du point \ref{item:lm:majfvchi:vnonar} 
de la proposition \ref{prop:fv:place:finie:quelc},
de la majoration
\begin{equation}
\abs{
Q_{\Sigma,v}
\left(
\chi_{\beta}\left(\pi_{\wb}\right)
\,
q_v^{-\lb \sb}
\right)_{\beta\in \Sigma(1)/G_v}
}
\leq
Q_{\Sigma,v}
\left(
\,
q_v^{-\lb \Re(\sb)}
\right)_{\beta\in \Sigma(1)/G_v}
\end{equation}
et du fait que 
\begin{equation}
\bs\mapsto 
\prod_{v\notin S}
Q_{\Sigma,v}
\left(
q_v^{-\lb \sb}
\right)_{\beta\in \Sigma(1)/G_v}
\end{equation}
est holomorphe sur $\R^{\Sigma(1)/G}_{>\frac{1}{2}}$.
\end{demo}

\subsubsection{Cas fonctionnel}\label{subsubsec:casfonc:frQ}
\begin{lemme}\label{lm:frQ}
Soit $\chi\in\dualtop{T(\ak)/\K(T)}$.
La s\'erie formelle 
\nindex{$\frQ(\chi,\,.\,)$}
$\frQ(\chi,\,.\,)$
d\'efinie par 
\begin{equation}
\frQ(\chi,(\za))
\eqdef
c_{K,\dim(T)}\,
\prod_{v\notin S}
{
Q_{\Sigma,v}
\left(
\chi_{\beta}\left(\pi_{\wb}\right)\,\zb^{\,f_{\wb}\,\db}\right)
}
\end{equation}
a un rayon de convergence sup\'erieur \`a $q^{-\frac{1}{2}}$,
est $\DT^0$-compatible et v\'erifie 
pour tout $\bs\in \tube{\R^{\Sigma(1)/G}_{>\frac{1}{2}}}$
\begin{equation}\label{eq:frQ=prodQv:bis}
\frQ(\chi,\,q^{-\bs})
=
c_{K,\dim(T)}\,
\prod_{v\notin S}
Q_{\Sigma,v}
\left(
\chib\left(\pi_{\wb}\right)
\,
q_v^{-\lb \sb}
\right).
\end{equation}
\end{lemme}
\begin{demo}
Pour $v\notin S$, $\beta\in \Sigma(1)/G_v$ et $s\in \C$, on a 
\begin{equation}
q_v^{-\lb\,s}=\left(q^{\,-s}\right)^{\,f_{\wb}\,\db}.
\end{equation}
Ceci joint au lemme \ref{lm:prodvQsv} montre aussit\^ot
que le rayon de convergence de $\frQ(\chi,\,.\,)$  est sup\'erieur \`a $q^{-\frac{1}{2}}$,
ainsi que la relation \eqref{eq:frQ=prodQv:bis}.

La $\DT^0$-compatibilit\'e provient du lemme \ref{lm:za^daDT0comp}.
\end{demo}

\begin{lemme}\label{lm:frf:chi:rayconv}
On suppose que toutes les places de $K$ v\'erifient l'hypoth\`ese
\ref{hyp:hinz}.
Soit $\chi\in \dualtop{T(\ak)/\K(T)}$.
La s\'erie formelle
\nindex{$\frf\left(\chi,.\right)$}
\begin{equation}\label{eq:deffrfchi}
\frf\left(\chi,\bz\right)
\eqdef
\frQ(\chi,\bz)
\,
\prod_{v\in S}
\frf_v(\chi,\bz)
\,
\prod_{\alpha\in \Sigma(1)/G}
\prod_{w\in S_{\alpha}}
\left(1-\chia(\pi_w)\,\za^{\,f_w\,\da}\right)
\end{equation}
a un rayon de convergence sup\'erieur \`a $q^{-\frac{1}{2}}$
et est $\DT$-compatible.
\end{lemme}
\begin{demo}
Ceci d\'ecoule du th\'eor\`eme \ref{thm:transfo:nonram} 
et les lemmes \ref{lm:transfo:nonram} et \ref{lm:transfo:ram:fonc}.
\end{demo}

\begin{lemme}\label{lm:fourierHchichibz}
On suppose que toutes les places de $K$ v\'erifient l'hypoth\`ese
\ref{hyp:hinz}.
Soit $\chi\in \dualtop{T(\ak)/\K(T)}$.
On a pour tout $\bs\in \tube{\R^{\,\Sigma(1)/G}_{>1}}$
et tout 
$\bz\in\xtgu$
\begin{equation}
\fourier H\left(\chi.\chi_{\bz},-\bs\right)
=
\left(
\prod_{\alpha\in \Sigma(1)/G} 
\ecL_{\Ka}
\left(
\chia,
\acc{\bz}{\da\,\roa}\,\,q^{\,-\da\,\sa}
\right)
\right)
\, 
\frf\left(\chi,\gamma_{\unit}(\bz)\,q^{\,-\bs }\right)
\end{equation}
\end{lemme}
\begin{demo}
D'apr\`es la d\'efinition \eqref{eq:defomegaT} de $\omega_{T}$, on a
\begin{equation}\label{eq:intlimitstak:fonc}
\int\limits_{T(\ak)}
(\chi\,\chi_{\bz})(t) H(-\bs,t)\,\omega_{T}(t)
=
c_{K,\dim(T)}\,\prod_{v\in \placesde{K}}\,\,\,\int\limits_{T(K_v)}(\chi\,\chi_{\bz})_v(t) H_{v}(-\bs,t)\,d\mu_v(t).
\end{equation}
D'apr\`es le th\'eor\`eme \ref{thm:transfo:nonram}, le lemme \ref{lm:transfo:nonram}
et le lemme \ref{lm:transfo:ram:fonc}, le membre de gauche de \eqref{eq:intlimitstak:fonc}
est donc \'egal \`a 
\begin{gather}
c_{K,\dim(T)}
\prod_{\alpha\in \sg} 
\ecL_{\Ka}(\chia,\acc{\bz}{\da\,\roa}\,q^{-\da\,\sa})
\notag
\\
\times
\prod_{\alpha \in \Sigma(1)/G} 
\prod_{w\in S_{\alpha}} (1-\,\chia(\pi_w)\,\acc{\bz}{\da\,\roa}q^{-f_w\,\da\,\sa})
\notag
\\
\times
\prod_{v\notin S}
Q_{\Sigma,v}
\left(
\chib(\pi_{\wb})\,\acc{y}{\db\,\rob} q^{-f_w\,\db \sb}
\right)_{\beta\in \Sigma(1)/G_v} 
\notag
\\
\times
\prod_{v\in S} \frf_{v}(\chi,\gamma_{\unit}(\bz)\,q^{-\bs}),
\end{gather}
d'o\`u le r\'esultat.
\end{demo}

\begin{lemme}\label{lm:prop:frf}
On suppose que toutes les places de $K$ v\'erifient l'hypoth\`ese
\ref{hyp:hinz}.
La s\'erie formelle
\nindex{$\frf(\,.\,)$}
\begin{equation}
\frf(\bz)
\eqdef
\frQ(\ind,\bz)\,\frf_{S}(\bz)\,
\prod_{\alpha\in \Sigma(1)/G}
\prod_{w\in S_{\alpha}}
\left(1-\za^{\,f_w\,\da}\right)
\end{equation}
a un rayon de convergence sup\'erieur \`a 
$q^{-\frac{1}{2}}$
et est $\DT^0$-compatible.
\end{lemme}
\begin{demo}
Ceci d\'ecoule aussit\^ot du lemme \ref{lm:frfS}
et du lemme \ref{lm:za^daDT0comp}.
\end{demo}

\begin{lemme}\label{lm:prop:frf:bis}
On suppose que toutes les places de $K$ v\'erifient l'hypoth\`ese
\ref{hyp:hinz}.
On a, pour tout $\bs\in \tube{\R^{\,\Sigma(1)/G}_{>1}}$
et tout $\bz\in\xtgu$
\begin{multline}\label{eq:prop:frf}
\int\limits_{\adh{T(K)}\,\cap\,T(\ak)}
H(-\bs,t)\,\chi_{\bz}(t)\,
\omega_{T}(t)
\\
=
\left(
\prod_{\alpha\in \Sigma(1)/G} 
Z_{\Ka}\left(\acc{\bz}{\da\,\roa}\,\,q^{\,-\da\,\sa}\right)
\right)
\frf\left(
\gamma _{\unit}(\bz)\,q^{\,-\bs}
\right).
\end{multline}
\end{lemme}
\begin{demo}
D'apr\`es le th\'eor\`eme \ref{thm:transfo:nonram} 
et le lemme \ref{lm:frfS}
on a 
\begin{align}
&
\quad
\int\limits_{\adh{T(K)}\,\cap\,T(\ak)}
H(-\bs,t)\,\chi_{\bz}(t)\,
\omega_{T}(t)\\
&
=
c_{K,\dim(T)}
\left(
\prod_{v\notin S}\,\,\,\int\limits_{T(K_v)}\,H_{v}(-\bs,t)\left(\chi_{ \bz }\right)_v(t)d\mu_v(t)
\right)
\,\notag\\
&\quad\quad\quad
\quad\quad\quad\quad
\times 
\int\limits_{\adh{T(K)}^{\,\,S}}\,
H(-\bs,t)
\left(
\chi_{\bz}
\right)(t)
\left(\otimesu{v\in S}d\mu_v\right)(t)
\\
&
=
\prod_{\alpha\in \Sigma(1)/G} 
Z_{\Ka}\left(\acc{\bz}{\da\,\roa}\,q^{\,-\da\,\sa}\right)
\notag\\
&
\quad\quad
\times
\prod_{\alpha\in \Sigma(1)/G}
\prod_{w\in S_{\alpha}}
\left(1-\bz^{\,f_w\,\da\,\roa} q^{\,-f_w\,\da\,\sa} \right)
\notag\\
&
\quad\quad
\times \frQ(\ind,\gamma _{\unit}(\bz)\,q^{\,-\bs })
\,
\frf_S\left(
\gamma _{\unit}(\bz)\,q^{\,-\bs }
\right)
\\
&
=
\left(
\prod_{\alpha\in \Sigma(1)/G} 
Z_{\Ka}\left(\acc{\bz}{\da\,\roa}\,q^{\,-\da\,\sa}\right)
\right)
\frf\left(
\gamma_{\unit}(\bz)\,q^{\,-\bs }
\right).
\end{align}
\end{demo}

\subsection{L'expression int\'egrale de la fonction z\^eta des hauteurs}
\label{subsec:expr:int}

Notons 
\nindex{$\gamma^{\ast}$}
$\gamma^{\ast}$
le morphisme de caract\`eres
\begin{equation}
\dualtop{T(\ak)/T(K)}
\longto 
\dualtop{T_{\Ps}(\ak)/T_{\Ps}(K)},
\end{equation}
induit par la composition avec le  morphisme $\gamma_{\ak}$.
Dans le cas fonctionnel comme dans le cas arithm\'etique,
nous choisissons arbitrairement et une fois pour toutes un scindage 
du morphisme quotient 
\begin{equation}
T(\ak)/T(K)\to T(\ak)/T(\ak)^1.
\end{equation}
Ce scindage induit par dualit\'e un isomorphisme
de $\dualtop{T(\ak)^1/T(K)}$ sur un facteur
direct de $\dualtop{T(\ak)/T(\ak)^1}$ dans
$\dualtop{T(\ak)/T(K)}$, lequel facteur direct sera not\'e $\wt{\ecU_T}$. 
On note \nindex{$\ecU_T$}$\ecU_T$ l'image de $\dualtop{T(\ak)^1/\K(T).T(K)}$ dans 
$\wt{\ecU(T)}$ par cet isomorphisme. En d'autre termes, 
$\ecU_T$ est le sous-groupe de $\wt{\ecU_T}$ constitu\'e des caract\`eres
qui sont triviaux sur $\K(T)$, et on a 
\begin{equation}\label{eq:scind:T:dualtop}
\dualtop{T(\ak)/\K(T).T(K)}=\dualtop{T(\ak)/T(\ak)^1}\times \ecU_T. 
\end{equation}

\begin{lemme}\label{lm:kergammast}
\begin{enumerate}
\item
Le noyau 
\label{item:1:lm:kergammast}
de $\gamma^{\ast}$ est isomorphe \`a $\dualt{A(T)}$. En particulier
$\Ker(\gamma^{\ast})$ est fini.
\item
Dans le cas arithm\'etique, $\Ker(\gamma^{\ast})$ est inclus dans $\wt{\ecU_T}$.
\item
Dans le cas 
\label{item:3:lm:kergammast}
fonctionnel, $\Ker(\gamma^{\ast})\cap \dualtop{T(\ak)/T(\ak)^1}$
est de cardinal $\KT$.
\end{enumerate}
\end{lemme}
\begin{demo}
D'apr\`es 
\eqref{eq:exsqtpaktakat}
et le lemme \ref{lm:picxsflasque},
on a une suite exacte
\begin{equation}
T_{\Ps}(\ak)/T_{\Ps }(K)
\longto 
T(\ak)/T(K)
\longto A(T)
\longto 0
\end{equation}
d'o\`u par dualit\'e une suite exacte
\begin{equation}
0
\longto 
A(T)^{\ast}
\longto 
\dualtop{T(\ak)/T(K)}
\longto 
\dualtop{T_{\Ps }(\ak)/T_{\Ps }(K)}.
\end{equation}
d'o\`u le premier point.

Dans le cas arithm\'etique, le fait que $\Ker(\gamma^{\ast})$
soit fini et que $\dualtop{T(\ak)/T(\ak)^1}$ soit sans torsion
donne aussit\^ot le second point.

Dans le cas fonctionnel, 
on note que le groupe $\Ker(\gamma^{\ast})\cap \dualtop{T(\ak)/T(\ak)^1}$
est le noyau du morphisme 
\begin{equation}
\dualtop{T(\ak)/T(\ak)^1}\longto \dualtop{T_{\Ps}(\ak)/T_{\Ps}(\ak)^1}
\end{equation}
induit par la composition avec $\gamma_{\ak}$,
et qu'on a un diagramme commutatif 
\begin{equation}
\xymatrix{
\DT^{\ast}
\ar[r]\ar[d]^{\deg_{T}^{\ast}}&
\DTqd^{\ast}
\ar[d]^{\deg^{\ast}_{T_{\Ps}}}\\
\dualtop{T(\ak)/T(\ak)^1}\ar[r]&
\dualtop{T_{\Ps}(\ak)/T_{\Ps}(\ak)^1}
}
\end{equation}
dont les fl\`eches verticales sont des isomorphismes.
Comme $\KT$ est le cardinal du conoyau du dual du morphisme
$\DT^{\ast}\longto \DTqd^{\ast}$, on a le troisi\`eme point.

\end{demo}

\subsubsection{Cas arithm\'etique}

Rappelons que $d\chi$ est par d\'efinition la mesure duale de $\wt{\omega_{T}}$, 
o\`u $\wt{\omega_{T}}$ est la mesure quotient de $\omega_{T}$ sur $T(\ak)/T(K)$,
$T(K)$ \'etant muni de la mesure discr\`ete.

Comme $\dualtop{T(\ak)^1/T(K)}$ est discret, $\dualtop{T(\ak)/T(\ak)^1}$
est un sous-groupe ouvert de $\dualtop{T(\ak)/T(K)}$, et on a pour toute
fonction $\varphi$ int\'egrable sur $\dualtop{T(\ak)/T(K)}$
\begin{equation}
\int\limits_{\dualtop{T(\ak)/T(K)}}\!\!\!\!\varphi(\chi)\,d\chi
=
\int\limits_{\chi\in \dualtop{T(\ak)/T(\ak)^1}}
\left(\sum_{\chi'\in \wt{\ecU_T}}\varphi(\chi.\chi')\right)d\chi
\end{equation}
En vue d'obtenir une formule int\'egrale explicite, 
il nous faut d\'eterminer la restriction de la mesure de Haar $d\chi$ au sous-groupe ouvert
$\dualtop{T(\ak)/T(\ak)^1}$.

Rappelons qu'on a d\'efini la mesure $\omega^1_{T}$
sur $T(\ak)^1/T(K)$ par la relation 
\begin{equation}
\wt{\omega_{T}}=\omega^1_{T}\,\left(\deg_T^{-1}\right)_{\ast}(dt),
\end{equation}
o\`u 
$dt$ est la mesure de Lebesgue sur
$\left(X(T)^G\right)^{\vee}_{\R}$
normalis\'ee par le r\'eseau $(X(T)^G)^{\vee}$.

Si on note $\dualtop{\omega^1_{T}}$ et $\dualt{dt}$ les mesures duales de
$\omega^1_{T}$ et $dt$ respectivement, 
on a donc
\begin{equation}
d\chi=\dualtop{\omega^1_{T}}.\,\left(\left[\deg_T^{\ast}\right]^{-1}\right)_{\ast}\dualt{dt}.
\end{equation}

Comme $T(\ak)^1/T(K)$ est compact,  d'apr\`es le lemme \ref{lm:ex:mesure:duale}
$\dualtop{\omega^1_{T}}$
est la mesure de Haar 
sur le groupe discret $\dualtop{T(\ak)^1/T(K)}$ 
pour laquelle chaque point 
a pour masse 
\begin{equation}
\frac{1}{\int\limits_{T(\ak)^1/T(K)}\!\!\!\!\omega^1_{T}}=\frac{1}{b(T)}.
\end{equation}

Par ailleurs, 
on note d\'esormais $dy$ la mesure de Lebesgue sur $X(T)^G_{\R}$
(identifi\'e au dual topologique de $\left(X(T)^G\right)^{\vee}_{\R}$ 
de la mani\`ere d\'ecrite au lemme \ref{lm:ex:mesure:duale}),
normalis\'ee par le r\'eseau $X(T)^G$.
D'apr\`es le lemme \ref{lm:ex:mesure:duale},
on a donc
\begin{equation}
dy=\frac{\dualt{dt}}{(2\,\pi)^{\rg\left(X(T)^G\right)}}.
\end{equation}
On obtient finalement le lemme suivant.
\begin{lemme}\label{lm:dec:int:arit}
Soit $\varphi$ une fonction int\'egrable sur $\dualtop{T(\ak)/T(K)}$.
On a alors
\begin{equation}
\int\limits_{\dualtop{T(\ak)/T(K)}}\!\!\!\!\varphi(\chi)\,d\chi
=
\frac{1}{(2\,\pi)^{\rg\left(X(T)^G\right)}\,b(T)}
\int\limits_{y\in X(T)^G_{\R}}
\left(\sum_{\chi\in \wt{\ecU_T}}\varphi(\chi_y.\chi)\right)dy
\end{equation}
\end{lemme}

Pour appliquer le lemme \ref{lm:formule:poisson}, 
on a besoin d'abord de prouver le caract\`ere $\L^1$ de 
$\fourier H(-\bs,\,.\,)$ 
sur $\dualtop{T(\ak)/\K(T).T(K)}$ 
(rappelons que cette int\'egrabilit\'e sera automatique dans le cas fonctionnel,
gr\^ace \`a la compacit\'e de $\dualtop{T(\ak)/\K(T).T(K)}$).
Pour ce faire, on utilise de mani\`ere cruciale
la majoration pour les transform\'ees de Fourier locales 
aux places archim\'ediennes obtenue au corollaire
\ref{cor:maj:fchiv:varch:2},
ainsi que la majoration uniforme des fonctions $L$
donn\'ee par le r\'esultat suivant,
cons\'equence du r\'esultat principal
de \cite{Rad:phrag}.

\begin{prop}\label{prop:maj:unif:L}
Soit $E$ un corps de nombres et $\delta>0$. 
Il existe un r\'eel $\eps>0$ (v\'erifiant $0<\eps<\frac{1}{3}$) 
et une constante $C>0$
telle qu'on ait pour tout $s$ v\'erifiant $\Re(s)>1-\eps$
et tout caract\`ere $\chi$ de $\G_m(\ade{E})$ trivial sur $\K(\G_m)$
la majoration
\begin{equation}
\abs{\frac{s-1}{s}\,L_E(\chi,s)}\leq C\,(1+\abs{\Im(s)})^{\delta}(1+\norm{\chi_{\infty}})^{\delta}
\end{equation}
\end{prop}

La proposition suivante est la proposition B.3 de \cite{CLTs:fibres}.
\begin{prop}[Chambert-Loir, Tschinkel]\label{prop:propB3} 
Soit $V$ un $\R$-espace -vectoriel de dimension finie 
muni d'une norme $\norm{\,.\,}$, 
$(\ell_j)$ une base du dual de $V$, 
$M$ un sous-espace vectoriel de $V$, $dm$ une mesure de Lebesgue sur $M$,
$V'$ un suppl\'ementaire de $M$ dans $V$, et $\delta$ un r\'eel v\'erifiant $0<\delta<1$.
Pout tout $\delta'>\delta$, il existe une constante $C>0$ et un ensemble fini
$(\ell_{i,j})_{i\in I,j\in J}$ d'\'el\'ements du dual de $V'$ tels que :
\begin{itemize}
\item
pour tout $i\in I$, la famille $(\ell_{i,j})_{j\in J}$ 
restreinte \`a $V$ est une base du dual de $V'$ ;
\item
pour tous 
$v_1$ et $v_2$ dans $V$ on a la majoration
\begin{multline}
\int\limits_M
\frac{1}{(1+\norm{v_1+m})^{1-\delta}}
\frac{dm}
{\prod(1+\abs{\ell_j(v_2+m)})}
\\
\leq 
\frac{C}{(1+\norm{v_1})^{1-\delta'}}
\sum_{i\in I}
\frac{1}
{\produ{j\in J}
(1+\abs{\ell_{i,j}(v_2)})}.
\end{multline}
\end{itemize}
\end{prop}
Nous aurons \'egalement besoin des r\'esultats suivants sur le morphisme
<<type \`a l'infini>>.
\begin{lemme}\label{lm:chi:infty}
\begin{enumerate}
\item
La restriction 
\label{item:1:lm:chi:infty}
du morphisme <<type \`a l'infini>>
\`a $\dualtop{T(\ak)/\K(T).T(K)}$ est de noyau fini.
\item
La composition 
\label{item:2:lm:chi:infty}
de l'isomorphisme 
\begin{equation}
\dualt{\deg_T}\,:\,
X(T)^G_{\R}
\longisom
\dualtop{T(\ak)/T(\ak)^1}
\end{equation}
avec le morphisme de $\dualtop{T(\ak)/T(\ak)^1}$ vers $X(T)_{\R,\infty}$
induit par le morphisme <<type \`a l'infini>>
co\"\i ncide avec l'injection
diagonale $X(T)^G_{\R}\to X(T)_{\R,\infty}$.
\item
On note  $\ecU_{T,{\infty}}$
l'image 
\label{item:3:lm:chi:infty}
dans $X(T)_{\R,\infty}$ de $\ecU_T$ par le morphisme <<type \`a l'infini>>. 
Alors $\ecU_{T,{\infty}}$ est un r\'eseau d'un sous-espace vectoriel suppl\'ementaire
de $X(T)^G_{\R}$ dans $X(T)_{\R,\infty}$.
\end{enumerate}
\end{lemme}
\begin{demo}
D'apr\`es la proposition \ref{prop:compacite}, le groupe 
\begin{equation}
T(\ak)/\K(T)\,T(K)\,T(\ak)_{\placesde{K,\infty}}
\end{equation}
est fini.
Par dualit\'e, ceci montre le point \ref{item:1:lm:chi:infty}.

Le point \ref{item:2:lm:chi:infty} 
d\'ecoule 
par dualit\'e
de la commutativit\'e du diagramme
\begin{equation}\label{eq:diag:eps}
\xymatrix{
\prod_{v\in \placesde{K,\infty}}
T(K_v)/T(\Ov)
\ar[r]\ar[d]^{\sum_{v} \deg_{T,v}}&T(\ak)/T(\ak)^1\ar[d]^{\deg_{T}} \\
\bigoplus_{v\in \placesde{K,\infty}} \Hom\left(X(T)^{\,G_v},\R\right)\ar[r] & \Hom(X(T)^G,\R)
}
\end{equation}
o\`u la fl\`eche horizontale du bas est la somme des restrictions.

Notons 
\begin{equation}
\eps\,:\,\prod_{v\in \placesde{K,\infty}}T(K_v)/T(\Ov)\longto \Hom(X(T)^G,\R)
\end{equation}
la composition de la fl\`eche verticale de gauche et de la fl\`eche horizontale du bas du 
diagramme \eqref{eq:diag:eps}.

Le groupe $T(\ak)^1_{\placesde{K,\infty}}/\K(T).T(K)\,\cap T(\ak)^1_{\placesde{K,\infty}}$
est un sous-groupe ouvert du groupe $T(\ak)^1/\K(T).T(K)$, et donc un sous-groupe
ferm\'e. D'apr\`es la proposition \ref{prop:compacite}, 
$T(\ak)^1_{\placesde{K,\infty}}/\K(T).T(K)\,\cap T(\ak)^1_{\placesde{K,\infty}}$ est donc compact.

Notons
\begin{equation}
T(\OK)=T(K)\cap \prod_{v\in \placesde{K,f}}T(\Ov)
\end{equation}
et $T(\OK)_{\infty}$ l'image de $T(\OK)$ dans $\produ{v\in \placesde{K,\infty}} T(K_v)/T(\Ov)$.
Ainsi $T(\OK)_{\infty}$ est un sous-groupe discret de $\produ{v\in \placesde{K,\infty}} T(K_v)/T(\Ov)$
dont l'image est contenue dans $\Ker(\eps)$.

On a 
\begin{equation}
\K(T).T(K)\,\cap T(\ak)^1_{\placesde{K,\infty}}=T(\OK).\prod_{v\in \placesde{K,\infty}} T(\Ov)
\end{equation}
et un isomorphisme
\begin{equation}
T(\ak)^1_{\placesde{K,\infty}}/\K(T).T(K)\,\cap T(\ak)^1_{\placesde{K,\infty}}
\longisom \Ker(\eps)/T(\OK)_{\infty}
\end{equation}
Ainsi $\Ker(\eps)/T(\OK)_{\infty}$ est compact et $T(\OK)_{\infty}$ est un r\'eseau de $\Ker(\eps)$
(lorsque $K=\G_m$, ce r\'esultat est le th\'eor\`eme des unit\'es de Dirichlet).

Le diagramme commutatif suivant
\begin{equation}
\xymatrix{
1\ar[d]&&0\ar[d]\\
\dualtop{T(\ak)/T(\ak)^1}\ar[rr]\ar[d]&&X(T)_{\R}^G\ar[d]\\
\dualtop{T(\ak)/\K(T).T(K)}\ar[rr]^{\sumu{v\in \placesde{K,\infty}}\dualt{\deg_{T,v}}}\ar[d]
&&\dualtop{\produ{v\in \placesde{K,\infty}} T(K_v)/T(\Ov)}\ar[d]\\
\dualtop{T(\ak)^1/\K(T).T(K)}\ar@{->>}[r]\ar[d]&\dualtop{\Ker(\eps)/T(\OK)_{\infty}}\ar@{^{(}->}[r]&\dualt{\Ker(\eps)}
\ar[d]\\
1&&0
}
\end{equation}
o\`u les deux colonnes sont exactes, montre alors le point \ref{item:3:lm:chi:infty}.
\end{demo}

Pour all\'eger les notations, on identifie dans l'\'enonc\'e et la d\'emonstration
de la proposition \ref{prop:contr:F} l'espace vectoriel $X(T)^G_{\R}$ \`a son image par
$\gamma_{\R}$ dans $\pls^G_{\R}$.
\begin{prop}\label{prop:contr:F}
\begin{enumerate}
\item
La s\'erie \label{item1:prop:contr:F}
de fonctions holomorphes
\begin{equation}
\bs
\mapsto 
 \sum_{\chi\in \ecU_T}
\left(\prod_{\alpha\in \Sigma(1)/G} 
L_{\Ka}
\left(\chia,\sa\right)\right)
\,f\left(\chi,\bs\right)
\,f_{\infty}\left(\chi,\bs\right)
\end{equation}
converge uniform\'ement sur tout compact de $\tube{\R^{\,\Sigma(1)/G}_{>1}}$
et d\'efinit donc une fonction holomorphe sur ce domaine, not\'ee $G(\bs)$.
\item
On note $F$ la fonction d\'ecal\'ee \label{item2:prop:contr:F}
\begin{equation}
\bs\mapsto \frac{1}{\card{A(T)}}\,G(\bs-\varphi_0).
\end{equation}
Il existe alors un r\'eel $\eps>0$ tel que
la fonction
\begin{equation}
\bs
\mapsto 
\left(\prod_{\alpha\in \Sigma(1)}
\sa
\right)
F(\bs)
\end{equation}
se prolonge en une fonction holomorphe sur $\tube{\R^{\,\Sigma(1)/G}_{>-\eps}}$,
dont la valeur en $\bs=0$ est 
\begin{equation}\label{eq:item2:prop:contr:F}
\lim_{s\to 1}\, 
(s-1)^{\,\card{\Sigma(1)/G}}
\int\limits_{\adh{T(K)} \cap T(\ak)}\!\!\!H(-s\,\varphi_0,t)\omega_{T}(t).
\end{equation}
\item
Soit $\delta>0$.
Il existe un r\'eel $\eps>0$ tel que la propri\'et\'e  \label{item3:prop:contr:F}
suivante est v\'erifi\'ee~:
pour tout $\delta'$ v\'erifiant $0<\delta'<1$ 
et tout compact $\K$ de $\R^{\Sigma(1)/G}_{>-\eps}$,
il existe une constante $C>0$ 
et un ensemble fini $(\ell_{ij})_{i\in I,j\in J}$
de formes lin\'eaires sur $\R^{\Sigma(1)/G}$ tels que :
\begin{itemize}
\item
pour tout $i\in I$, 
la famille $(\ell_{ij})_{j\in J}$ restreinte \`a $X(T)^G_{\R}$
est une base du dual de $X(T)^G_{\R}$ ;
\item
pour tout $\bs$ de $\tube{\K}$ et tout $y\in X(T)^G_{\R}$,
on a la majoration
\begin{multline}\label{eq:maj:Fsaya}
\abs{\left(\prod_{\alpha\in \sg} \frac{\sa+i\,\ya}{1+\sa+i\,\ya}\right)
\,F(\bs+i\,y)}
\\
\leq
C
\,
\frac{(1+\norm{\Im(\bs)})^{1+\delta}}{(1+\norm{y})^{1-\delta'}}
\,
\sum_{i\in I}
\,
\prod_{j\in J}
\,
\frac{1}{1+\abs{\ell_{i,j}(\Im(\bs)+y)}}.
\end{multline}
\end{itemize}
\end{enumerate}
\end{prop}
\begin{rems}\label{rem:prop:contr:F}
\begin{enumerate}
\item
La d\'emonstration qui suit est fortement inspir\'ee de 
\cite{CLTs:fibres}.
\item
Pour $\bs\in\tube{\R^{\,\Sigma(1)/G}_{>0}}$,
soit
\begin{equation}
\wt{G}(\bs)=\sum_{\chi\in \ecU_T}
\abs{
\left(\prod_{\alpha\in \sg} 
L_{\Ka}
\left(\chia,\sa\right)
\right)
\,f\left(\chi,\bs\right)
\,f_{\infty}\left(\chi,\bs\right)
}
\end{equation}
et $\wt{F}(\bs)=\frac{1}{\card{A(T)}}\,\wt{G}(\bs-\varphi_0)$.
La d\'emonstration de la proposition \ref{prop:contr:F}
montrera en fait que $\wt{F}$ v\'erifie la majoration
\eqref{eq:maj:Fsaya}.
\item
Le point 
\ref{item3:prop:contr:F}
montre que la fonction 
\begin{equation}
\bs \mapsto \left(\prod_{\alpha\in \sg} \frac{\sa}{1+\sa}\right)\,F(\bs)
\end{equation}
est $X(T)^G_{\R}$-contr\^ol\'ee 
au sens de \cite[D\'efinition 3.13]{CLTs:fibres}.
\end{enumerate}
\end{rems}

\begin{demo}

D'apr\`es la proposition \ref{prop:fonc:zeta}
et les lemmes \ref{lm:chi:non:triv}, 
\ref{lm:fchihol:arit}
et \ref{lm:finftychihol},
pour tout $\chi$, la fonction 
\begin{equation}\label{eq:def:gchisa}
\bs
\mapsto 
\left(
\prod_{\alpha\in \sg}
\frac{\sa-1}{\sa}
\,
L_{\Ka}
\left(\chia,\sa\right)
\right)
\,f\left(\chi,\bs\right)
\,f_{\infty}\left(\chi,\bs\right)
\end{equation}
se prolonge en une une fonction holomorphe sur 
$\tube{\R^{\,\Sigma(1)/G}_{>\frac{1}{2}}}$,
que l'on note $g(\chi,\bs)$.

Montrons alors qu'il existe un r\'eel $\eps$ v\'erifiant $0<\eps<\frac{1}{2}$ 
tel que la s\'erie de fonctions
\begin{equation}\label{eq:sumgchi}
\sum_{\chi\in 
\ecU_T
}
g(\chi,\bs)
\end{equation}
converge uniform\'ement sur tout compact de $\tube{\R^{\,\Sigma(1)/G}_{>1-\eps}}$,
ce qui donnera les points \ref{item1:prop:contr:F} et \ref{item2:prop:contr:F},
\`a l'exception de \eqref{eq:item2:prop:contr:F}.

Soit $\delta>0$.
D'apr\`es la proposition \ref{prop:maj:unif:L}, il existe un r\'eel  
$\eps>0$ et une constante $C_1>0$
telle qu'on ait, pour tout $\bs\in \tube{\R^{\Sigma(1)/G}_{>1-\eps}}$
et tout $\chi\in \dualtop{T(\ak)/\K(T)}$
\begin{equation}
\abs{
\prod_{\alpha\in \sg} 
\left(\frac{\sa-1}{\sa}\right)
\,
L_{\Ka}
\left(\sa,\chi_{\alpha}\right)
}
\leq
C_1
\prod_{\alpha\in \sg}
(1+\abs{\Im(\sa)})^{\delta}
(1+\norm{(\chi_{\alpha})_{\infty}})^{\delta}
\end{equation}
soit
\begin{equation}
\abs{
\prod_{\alpha\in \sg} 
\left(\frac{\sa-1}{\sa}\right)
\,
L_{\Ka}
\left(\sa,\chi_{\alpha}\right)
}
\leq
C_1
\,
{\left(1+\norm{\Im(\bs)}\right)^{\delta}}
{\left(1+\norm{\chi_{\infty}}\right)^{\delta}}.
\end{equation}
Soit $\K$ un compact de $\tube{\R^{\Sigma(1)/G}_{>1-\eps}}$.
D'apr\`es le corollaire \ref{cor:maj:fchiv:varch:2}, 
il existe une constante $C_2>0$ telle que 
pour tout $\bs\in \K$
et tout $\chi\in 
\ecU_T
$
on a 
\begin{equation}
\abs{f_{\infty}\left(\chi,\bs\right)}
\leq
\frac{C_2}
{1+\norm{\chi_{\infty}}}
\,
\sum_{\wt{\sigma}\in \Sigma^{\infty}_{\max}}
\frac
{1}
{
\produ{i\in \wt{\sigma}(1)}
(1+\abs{\acc{\roi}{\chi_\infty}+\Im(s_i)})
}
\end{equation}

D'apr\`es le lemme \ref{lm:maj:fchifvchivfin}, 
il existe une constante $C_3$
telle que 
pour tout $\bs\in \K$
et tout $\chi\in \ecU_T$
on a 
\begin{equation}
\abs{f\left(\chi,\bs\right)}
\leq C_3.
\end{equation}

Finalement, on a montr\'e qu'il existe une constante $C>0$
telle que pour tout $\bs\in \K$ 
et tout $\chi\in \ecU_T$
on a la majoration
\begin{equation}
\abs{g(\chi,\bs)}
\leq
C
\,
\frac
{1}
{\left(1+\norm{\chi_{\infty}}\right)^{1-\delta}}
\,
\sum_{\wt{\sigma}\in \Sigma^{\infty}_{\max}}
\frac
{
1
}
{
\produ{i\in \wt{\sigma}(1)}
(1+\abs{\acc{\roi}{\chi_\infty}+\Im(s_i)})
}
\end{equation}

Pour montrer que la s\'erie de fonctions \eqref{eq:sumgchi}
converge uniform\'ement sur $\K$ 
il suffit donc de montrer que pour tout $\wt{\sigma}\in \Sigma^{\infty}_{\max}$
la s\'erie
\begin{equation}
\sum_{\chi\in 
\ecU_T
} 
\frac
{1}
{\left(1+\norm{\chi_{\infty}}\right)^{1-\delta}}
\,
\frac
{
1
}
{
\produ{i\in \wt{\sigma}(1)}
(1+\abs{\acc{\roi}{\chi_\infty}+\Im(s_i)})
}
\end{equation}
est convergente uniform\'ement localement en $\bs$.

D'apr\`es le lemme \ref{lm:chi:infty}, il suffit de montrer que la s\'erie
\begin{equation}
\sum_{m\in 
\ecU_{T,\infty}
} 
\frac
{1}
{\left(1+\norm{m}\right)^{1-\delta}}
\,
\frac
{
1
}
{
\produ{i\in \wt{\sigma}(1)}
(1+\abs{\acc{\roi}{m}+\Im(s_i)})
}
\end{equation}
est convergente uniform\'ement localement en $\bs$.

Mais cette s\'erie est major\'ee par l'int\'egrale
\begin{equation}
\int\limits_{m\in 
\ecU_{T,\infty}
\otimes \R} 
\frac
{1}
{\left(1+\norm{m}\right)^{1-\delta}}
\,
\frac
{
dm
}
{
\produ{i\in \wt{\sigma}(1)}
(1+\abs{\acc{\roi}{m}+\Im(s_i)})
},
\end{equation}
o\`u $dm$ est la mesure de Lebesgue sur 
$\ecU_{T,\infty}\otimes \R$
normalis\'ee par le r\'eseau 
$\ecU_{T,\infty}$.
Pour conclure, on applique alors la proposition \ref{prop:propB3} 
avec $V=X(T)_{\R,\infty}$, 
$M=\ecU_{T,\infty}\otimes \R$,
$V'=X(T)^G_{\R}$ (qui est bien un suppl\'ementaire de $M$
dans $V$ d'apr\`es le point \ref{item:3:lm:chi:infty} 
du lemme \ref{lm:chi:infty}) 
$v_1=0$, $v_2=(\Im(s_i))_{i\in \wt{\sigma}(1)}$,
et en prenant pour base du dual de $V$
la famille $(\acc{\rho_i}{\,.\,})_{i\in \wt{\sigma}(1)}$.

Montrons \`a pr\'esent \eqref{eq:item2:prop:contr:F}.
Il s'agit de montrer qu'on a 
\begin{equation}\label{eq:item2:prop:contr:F:2}
\sum_{\chi\in 
\ecU_{T}
}
g(\chi,\varphi_0)
=
\card{A(T)}\,
\lim_{s\to 1}\, 
(s-1)^{\,\card{\Sigma(1)/G}}
\!\!\!\!\!\!\int\limits_{\adh{T(K)} \cap T(\ak)}\!\!\!\!H(-s\,\varphi_0,t)\omega_{T}(t).
\end{equation}
Compte tenu de la d\'efinition de $g(\chi,\bs)$
et du fait que $L_{\Ka}(\chia,\,.\,)$ 
est holomorphe sur $\C$ si $\chia$ est non trivial,
$g(\chi,\varphi_0)$ est nul d\`es qu'il existe un 
$\alpha$ tel que $\chia$ est non trivial,
en d'autre termes d\`es que $\chi$ n'est pas dans $\Ker(\gamma^{\ast})$.
Le premier membre de \eqref{eq:item2:prop:contr:F:2}
s'\'ecrit donc
\begin{equation}\label{eq:item2:prop:contr:F:3}
\sum_{
\chi\in 
\ecU_{T}
\cap \Ker(\gamma^{\ast})
}
\left(\prod_{\alpha\in \Sigma(1)/G} 
\underset{s=1}{\Res}
\,\,\,\zeta_{\Ka}
\left(s\right)
\right)
\,
f\left(\chi,\varphi_0\right)
\,
f_{\infty}\left(\chi,\varphi_0\right)
\end{equation}

Par ailleurs, pour $s$ complexe tel que $\Re(s)>1$,
appliquons la formule de Poisson 
\ref{thm:formule:poisson} 
avec $\ecG=T(\ak)$, 
$\cH=\adh{T(K)}\cap T(\ak)$, $dh=dg=\omega_{T}$ et $F=H(-s\,\varphi_0,\,.\,)$.
On obtient d'apr\`es le lemme \ref{lm:kergammast}
et le fait que $\fourier H\left(\chi,-s\,\varphi_0\right)$
est nulle d\`es que $\chi$ n'est pas trivial sur $\K(T)$
\begin{align}
\int\limits_{\adh{T(K)}\cap T(\ak)}H(-s\,\varphi_0,t)\,\omega_{T}(t)
&
=
\frac{1}
{\card{\Ker(\gamma^{\ast})}}
\sum_{\chi\in \Ker\left(\gamma^{\ast}\right)}
\fourier H\left(\chi,-s\,\varphi_0\right)
\\
&
=\label{eq:item2:prop:contr:F:4}
\frac{1}
{\card{A(T)}}
\sum_{\chi\in \Ker\left(\gamma^{\ast}\right)\cap 
\ecU_{T}
}
\fourier H\left(\chi,-s\,\varphi_0\right)
\end{align}

Mais d'apr\`es \eqref{eq:expr:fourierHchisa:fchisa} on a
pour tout 
$\chi\in \Ker(\gamma^{\ast})\cap \ecU_{T}$
\begin{equation}
\label{eq:item2:prop:contr:F:5}
\fourier H(\chi,-s\,\varphi_0)=
\left(
\prod_{\alpha\in \Sigma(1)/G}
\zeta_{\Ka}\left(s\right)\right)
f(\chi,s\,\varphi_0)\,
f_{\infty}(\chi,s\,\varphi_0).
\end{equation}
De
\eqref{eq:item2:prop:contr:F:3},
\eqref{eq:item2:prop:contr:F:4}
et
\eqref{eq:item2:prop:contr:F:5},
on d\'eduit aussit\^ot la relation
\eqref{eq:item2:prop:contr:F:2}.

\vskip0.5cm 
Montrons \`a pr\'esent le point \ref{item3:prop:contr:F}. 
Soit $\delta>0$. D'apr\`es la proposition \ref{prop:maj:unif:L}, 
il existe un r\'eel $\eps>0$ et une constante $C_1>0$
telle qu'on ait, pour tout $\bs\in \tube{\R^{\Sigma(1)/G}_{>-\eps}}$
et tout $\chi\in\ecU_T$
\begin{multline}\label{eq:maj:lchisa+iya}
\abs{
\left(
\prod_{\alpha\in \sg} 
\frac{\sa+i\ya}{\sa+1+i\ya}
\right)
\,
L_{\Ka}
\left(\chia,\sa+1+i\ya\right)
}
\\
\leq
C_1
\,
{\left(1+\norm{\Im(\bs)+y}\right)^{\delta}}
{\left(1+\norm{\chi_{\infty}}\right)^{\delta}}.
\end{multline}
Soit $\K$ un compact de $\R^{\Sigma(1)/G}_{>-\eps}$.
D'apr\`es 
\eqref{eq:maj:lchisa+iya},
le lemme \ref{lm:maj:fchifvchivfin}
et le corollaire \ref{cor:maj:fchiv:varch:2},
il existe une constante $C_2>0$ telle que pour tout
$\bs\in \tube{\K}$, 
tout $y\in X(T)^G_{\R}$ 
et tout $\chi\in \ecU_T$ 
l'expression
\begin{equation}
\abs{
g(\chi,(\sa-1+i\,\ya))}
\\
=
\abs{
g(\chi.\chi_y,\bs-\varphi_0)}
\end{equation}
est major\'ee par 
\begin{equation}
C_2
\frac
{\left(1+\norm{\Im(\bs)+y}\right)^{\delta}
\left(1+\norm{\chi_{\infty}}\right)^{\delta}
}
{1+\norm{\chi_{\infty}+y}}
\,
\sum_{\wt{\sigma}\in \Sigma^{\infty}_{\max}}
\frac
{
1+\sumu{i\in \wt{\sigma}(1)}\abs{\Im(s_i)}
}
{
\produ{i\in \wt{\sigma}(1)}
(1+\abs{\acc{\roi}{\chi_\infty+y}+\Im(s_i)})
}
.
\end{equation}
On en d\'eduit que 
pour tout
$\bs\in \tube{\K}$ et tout $y\in X(T)^G_{\R}$,
l'expression
\begin{equation}\label{eq:expr:a:majorer}
\abs{\left(
\prod_{\alpha\in \sg} 
\frac{\sa+i\ya}{\sa+1+i\ya}
\right)\,
F(\bs+i\,y)}
\end{equation}
est major\'ee par 
\begin{equation}
C_2
\left(1+\norm{\Im(\bs)+y}\right)^{\delta}
\sum_{\wt{\sigma}\in \Sigma^{\infty}_{\max}}
\left(1+\sumu{i\in \wt{\sigma}(1)}\abs{\Im(s_i)}\right)
\phi_{\wt{\sigma}}(\bs,y)
\end{equation}
o\`u 
$
\phi_{\wt{\sigma}}(\bs,y)
$
est d\'efini par la s\'erie
\begin{equation}
\phi_{\wt{\sigma}}(\bs,y)
=\sum_{\chi\in 
\ecU_{T}
}
\frac
{
\left(1+\norm{\chi_{\infty}}\right)^{\delta}
}
{1+\norm{\chi_{\infty}+y}}
\produ{i\in \wt{\sigma}(1)}
\frac
{1}
{
1+\abs{\acc{\roi}{\chi_\infty}+\Im(s_i)+y_i}
}.
\end{equation}
On obtient finalement
que l'expression \eqref{eq:expr:a:majorer} est major\'ee par
\begin{equation}
C_2
\left(1+\norm{\Im(\bs)}\right)^{1+\delta}\,\left(1+\norm{y}\right)^{1+\delta}
\sum_{\wt{\sigma}\in \Sigma^{\infty}_{\max}}
\phi_{\wt{\sigma}}(\bs,y).
\end{equation}
Fixons $\wt{\sigma}\in \Sigma^{\infty}_{\max}$.
De la majoration
\begin{equation}
1+\norm{\chi_{\infty}}\leq
1+\norm{\chi_{\infty}+y}+\norm{y}
\leq 
(1+\norm{\chi_{\infty}+y})(1+\norm{y}),
\end{equation}
on d\'eduit qu'on a 
\begin{equation}
\phi_{\wt{\sigma}}(\bs,y)
\leq\sum_{\chi\in 
\ecU_{T}
}
\frac
{
\left(1+\norm{y}\right)^{\delta}
}
{(1+\norm{\chi_{\infty}+y})^{1-\delta}}
\produ{i\in \wt{\sigma}(1)}
\frac
{1}
{
1+\abs{\acc{\roi}{\chi_\infty}+\Im(s_i)+y_i}
}
.
\end{equation}
D'apr\`es le lemme \ref{lm:chi:infty}, il  existe alors une constante $C_3>0$
telle qu'on ait 
\begin{equation}
\phi_{\wt{\sigma}}(\bs,y)
\leq
C_3
\sum_{m\in 
\ecU_{T,\infty}
}
\frac
{
\left(1+\norm{y}\right)^{\delta}
}
{(1+\norm{m+y})^{1-\delta}}
\produ{i\in \wt{\sigma}(1)}
\frac
{1}
{
1+\abs{\acc{\roi}{m}+\Im(s_i)+y_i}
}
\end{equation}
soit
\begin{multline}
\phi_{\wt{\sigma}}(\bs,y)
\\
\leq
C_3\,\left(1+\norm{y}\right)^{\delta}
\!\!\!\!\!\!\!\!\!\!
\int\limits_{m\in 
\ecU_{T,\infty}\otimes \R
}
\frac
{
1
}
{(1+\norm{m+y})^{1-\delta}}
\frac
{dm}
{
\produ{i\in \wt{\sigma}(1)}
\left(1+\abs{\acc{\roi}{m}+\Im(s_i)+y_i}\right)
}
.
\end{multline}
On conclut comme ci-dessus en appliquant la proposition \ref{prop:propB3} 
avec $V=X(T)_{\R,\infty}$, 
$M=\ecU_{T,\infty}$,
$V'=X(T)^G_{\R}$,  
$v_1=y$, $v_2=(\Im(s_i)+y_i)_{i\in \wt{\sigma}(1)}$,
et en prenant pour base du dual de $V$
la famille $(\acc{\rho_i}{\,.\,})_{i\in \wt{\sigma}(1)}$.
\end{demo}

\begin{cor}\label{cor:rep:int:arit}
Dans le cas arithm\'etique, 
pour tout $\bs\in \tube{\R^{\,\Sigma(1)/G}_{>0}}$ 
on a la repr\'esentation int\'egrale 
\begin{equation}\label{eq:cor:rep:int:arit}
\sum_{t\in T(K)}
\,
H\left(t,-\bs -\varphi_0\right)
=\frac
{\card{A(T)}}
{\,(2\,\pi)^{\rg\left(\xtg\right)}\,b(T)}
\int\limits_{X(T)^G_{\R}}
F(\bs-i\gamma_{\R}(y))dy
\end{equation}
o\`u $F$ est la fonction 
d\'efinie dans l'\'enonc\'e de la proposition \ref{prop:contr:F}.
\end{cor}
\begin{demo}
On veut appliquer le lemme \ref{lm:formule:poisson}.
On sait d\'ej\`a (lemme \ref{lm:H-bsint}) que $H\left(\,.\,,-\bs-\varphi_0\right)$
est int\'egrable sur $T(\ak)$, 
et il s'agit de montrer que
$\chi\mapsto \fourier H(\chi,-\bs-\varphi_0)$ est int\'egrable
sur $\dualtop{T(\ak)/\K(T)}$.
Or on a, d'apr\`es le lemme \ref{lm:dec:int:arit},
\begin{align}
&
\int\limits_{\dualtop{T(\ak)/\K(T)}}
\abs{\fourier H(\chi,-\bs-\varphi_0)}d\chi
\\
&
=
\frac{1}{(2\,\pi)^{\rg\left(X(T)^G\right)}\,b(T)}
\int\limits_{X(T)^G_{\R}}
\left(
\sum_{\chi \in \ecU_T}
\abs{\fourier H(\chi_y.\chi,-\bs-\varphi_0)}
\right)
dy
\\
&
=
\frac{1}{(2\,\pi)^{\rg\left(X(T)^G\right)}\,b(T)}
\int\limits_{X(T)^G_{\R}}
\left(
\sum_{\chi \in \ecU_T}
\abs{\fourier H(\chi,-\bs+i\gamma_{\R}(y)-\varphi_0)}
\right)
dy
\\
&
=
\frac{\card{A(T)}}{(2\,\pi)^{\rg\left(X(T)^G\right)}\,b(T)}
\int\limits_{X(T)^G_{\R}}
\wt{F}(\bs-i\,\gamma_{\R}(y))dy
\end{align}
o\`u $\wt{F}$ est la fonction d\'efinie \`a la remarque \ref{rem:prop:contr:F}.
Cette m\^eme remarque montre en outre que la derni\`ere int\'egrale est finie.

On peut donc appliquer le lemme \ref{lm:formule:poisson}, ce qui donne la formule
\ref{eq:cor:rep:int:arit},
compte tenu du lemme \ref{lm:dec:int:arit}.
\end{demo}

La technique pour \'evaluer le comportement analytique de la fonction 
d\'efinie par l'int\'egrale apparaissant dans \eqref{eq:cor:rep:int:arit}
est d\'evelopp\'ee par Batyrev et Tschinkel
dans \cite{BaTs:manconj}, et raffin\'ee par Chambert-Loir et Tschinkel dans \cite{CLTs:fibres}.
Dans la section \ref{sec:eval:int:arit}, 
nous rappelons le lemme technique principal de \cite{CLTs:fibres}.
Puis nous expliquons dans la section \ref{subsec:application:lemme:technique:arit}
comment appliquer ce r\'esultat  \`a l'\'evaluation de cette int\'egrale et donc \`a la fonction 
z\^eta des hauteurs des vari\'et\'es toriques (ceci est \'egalement 
expliqu\'e dans \cite{CLTs:fibres}, \`a ceci pr\`es que seules des vari\'et\'es toriques d\'eploy\'ees
sont consid\'er\'ees~; l'adaptation au cas non d\'eploy\'e est cependant ais\'ee).

\subsubsection{Cas fonctionnel}\label{subsubsec:expr:int:fonc}

Rappelons que $d\chi$ est par d\'efinition la mesure duale de $\wt{\omega_{T}}$, 
o\`u $\wt{\omega_{T}}$ est la mesure quotient de $\omega_{T}$ sur $T(\ak)/T(K)$,
$T(K)$ \'etant muni de la mesure discr\`ete.
Comme dans le cas arithm\'etique, en vue d'obtenir une formule int\'egrale explicite, 
il nous faut d\'eterminer la restriction de la mesure de Haar $d\chi$ 
au sous-groupe ouvert $\dualtop{T(\ak)/T(\ak)^1}$.

Rappelons que $\dualtop{T(\ak)/T(\ak)^1}$ est un sous-groupe 
de $\dualtop{T(\ak)/\K(T).T(K)}$ 
d'indice fini \'egal au cardinal de $T(\ak)^1/\K(T).T(K)$. 
Si on note $\cla(T)$ ce cardinal,
on a donc
\begin{equation}\label{eq:vol:St:1}
\int\limits_{\dualtop{T(\ak)/T(\ak)^1}}\!\!\!\!\!\!\!\!d\chi
=\frac{1}{\cla(T)}\,
\int\limits_{\dualtop{T(\ak)/\K(T).T(K)}}\!\!\!\!\!\!\!\!d\chi.
\end{equation}

Appliquons \`a pr\'esent le corollaire \ref{cor:formule:poisson} 
de la formule de Poisson
avec $\ecG=T(\ak)$, 
$dg=\omega_{T}$, 
$\cH=T(K)$, 
$dh$ la mesure discr\`ete, 
$\ecK=\K(T)$ 
(rappelons que $\int\limits_{\K(T)}\!\!\!\omega_{T}$ est non nul, cf. section \ref{sec:tamagawa:tore})
et $F$ l'indicatrice de $\K(T)$, 
on obtient
\begin{equation}\label{eq:vol:St:2}
\card{\K(T)\cap T(K)}
=
\int\limits_{\K(T)}\!\!\omega_{T}
\int\limits_{\dualtop{T(\ak)/\K(T).T(K)}}\!\!\!\!d\chi
.
\end{equation}
Or on a 
\begin{align}
\int\limits_{\K(T)}
\!\!\omega_{T}&
=
\card{\K(T)\cap T(K)}\int\limits_{\K(T).T(K)/T(K)}\!\!\!\!\omega_{T}^1\\
&
=
\frac{\card{\K(T)\cap T(K)}}{\cla(T)}\int\limits_{T(\ak)^1/T(K)}\!\!\!\!\omega_{T}^1
\\
&
=
\label{eq:vol:St:3}
\frac{\card{\K(T)\cap T(K)}}{\cla(T)} \log(q)^{\rg\left(X(T)^G\right)}\,b(T).
\end{align}
En combinant
\eqref{eq:vol:St:1}, \eqref{eq:vol:St:2}, et 
\eqref{eq:vol:St:3}, on obtient
\begin{align}
\int\limits_{\dualtop{T(\ak)/T(\ak)^1}}
\!\!\!\!
d\chi
&
=
\frac{1}{\cla(T)}
\,
\card{\K(T)\cap T(K)}
\,
\frac{\cla(T)}{\card{\K(T)\cap T(K)}\,\log(q)^{\rg\left(X(T)^G\right)}\,b(T)}
\\
&
=
\frac
{1}
{\log(q)^{\rg\left(X(T)^G\right)}\,b(T)}.
\end{align}

Soit  $\wt{d\chi}$ la mesure de Haar sur
$\dualtop{T(\ak)/T(\ak)^1}$ de masse totale $1$.
On d\'eduit de ce qui pr\'ec\`ede le lemme suivant.
\begin{lemme}\label{lm:dec:int:fonc}
Pour toute fonction  $\varphi$ int\'egrable sur 
$\dualtop{T(\ak)/T(K)}$,
on a 
\begin{equation}
\int\limits_{\dualtop{T(\ak)/T(K)}}\!\!\!\!\!\!\varphi(\chi)\,d\chi
=
\frac
{1}
{\log(q)^{\rg\left(X(T)^G\right)}\,b(T)}
\sum_{\chi'\in \wt{\ecU_T}}\,\,\,
\int\limits_{\dualtop{T(\ak)/T(\ak)^1}}
\!\!\!\!\!\!\!\!\varphi\left(\chi.\chi'\right)\wt{d\chi}
\end{equation}
\end{lemme}
\begin{conv}
Soit $N$ un $\Z$-module libre de rang fini.
On munira toujours $N_{\unit}$ de la mesure de Haar $d\bz$ de volume total 1.
En d'autres termes, 
dans toute expression du type
\begin{equation}
\int\limits_{N_{\unit}}\dots d\bz
\end{equation}
il sera toujours sous-entendu que $d\bz$ est la mesure de Haar de volume total 1.

Pour tout \'el\'ement $n^{\vee}$ de $N^{\vee}$, on a donc
\begin{equation}
\int\limits_{N_{\unit}}\acc{\bz}{n^{\vee}}d\bz
=
\left\{
\begin{array}{rl}
0& \text{ si }n^{\vee}\neq 0\\
1& \text{ si }n^{\vee}=0.
\end{array}
\right.
\end{equation}
\end{conv}
\begin{lemme}
\label{lm:intn=intm:bis}
Soient $M$ un $\Z$-module libre de rang fini, $N$ un sous-groupe de $M$ d'indice fini,
et $\varphi$ une fonction int\'egrable sur $M_{\unit}$. 
On a alors
\begin{equation}\label{eq:intn=intm}
\int\limits_{N_{\unit}}\varphi(\bz)d\bz=\int\limits_{M_{\unit}}\varphi(\bz)d\bz.
\end{equation}
\end{lemme}
\begin{demo}
Soit $(e_i)_{i\in I}$ une base de $M$ adapt\'ee \`a $N$. 
Au moyen de cette base, on identifie $M_{\unit}$ \`a $\unit^I$.
Si $dz$ d\'esigne la mesure uniforme sur $\unit$, 
le membre de droite de \eqref{eq:intn=intm} s'\'ecrit alors
\begin{equation}
\int\limits_{\unit^I}
\varphi(z_i)_{i\in I} \otimesu{i\in I} dz_i.
\end{equation}
Soit $(d_i)\in \N_{>0}^I$ tel que $(d_i\,e_i)$ soit une base de $N$.
Au moyen de cette base, on identifie $N_{\unit}$ \`a $\unit^I$.
Le membre de gauhe de \eqref{eq:intn=intm} s'\'ecrit alors
\begin{equation}
\int\limits_{\unit^I}
\varphi\left(z_i^{\,d_i}\right)_{i\in I} \otimesu{i\in I} dz_i.
\end{equation}
Pour conclure, on remarque que si $d$ est un entier non nul 
on a pour toute fonction $\psi$ int\'egrable sur $\unit$
\begin{equation}
\int\limits_{\unit}\psi\left(z^d\right)\,dz =\int\limits_{\unit}\psi(z)\,dz.
\end{equation}
\end{demo}
\begin{cor}\label{cor:dec:int:bis}
Pour toute fonction $\varphi$ int\'egrable sur $\dualtop{T(\ak)/T(\ak)^1}$,
on a 
\begin{equation}
\int\limits_{\dualtop{T(\ak)/T(\ak)^1}}
\!\!\!\!\!\!\varphi(\chi)
\,
\wt{d\chi}
=
\int\limits_{\xtgu} 
\!\!
\varphi\left(\chi_{\bz}\right)
\,
d\bz
\end{equation}
\end{cor}
\begin{demo}
On applique le lemme \ref{lm:intn=intm:bis} avec $M=\left(T(\ak)/T(\ak)^1\right)^{\vee}$ et $N=\xtg$.
\end{demo}
\begin{lemme}\label{lm:sumchiinkergamma1ast}
Pour $\bs\in \tube{\R^{\Sigma(1)/G}_{>1}}$, on a l'\'egalit\'e
\begin{multline}
\sum_{
\substack{
\chi'\in \ecU_T
\\
\chi'\in \Ker(\gamma^{\ast}).\dualtop{T(\ak)/T(\ak)^1}
}
}
\,\,\,\,\int\limits_{\dualtop{T(\ak)/T(\ak)^1}}
\!\!\!\!\!\!\!\!\!\!\fourier H\left(\chi.\chi',-\bs \right)\wt{d\chi}
\\
=
\frac{\card{A(T)}}{\card{\KT}}
\,
\int\limits_{\dualtop{T(\ak)/T(\ak)^1}}
\left[\,\int\limits_{\adh{T(K)}\cap T(\ak)}\!\!\!\!\!\!H(-\bs,t)\,\chi(t)\,\omega_{T}(t)\right]
\wt{d\chi}.
\end{multline}
\end{lemme}
\begin{demo}
Soit $\chi\in \dualtop{T(\ak)/T(\ak)^1}$ et $\bs\in \tube{\R^{\Sigma(1)/G}_{>1}}$.
Appliquons le corollaire \ref{cor:formule:poisson} 
de la formule de Poisson
avec $\ecG=T(\ak)$, 
$\ecH=\adh{T(K)}\cap T(\ak)$, $dh=dg=\omega_{T}$, $\ecK=\K(T)$ et $F=H(-\bs ,\,.\,)\,\chi$. 

On obtient, compte tenu du point \ref{item:1:lm:kergammast} du lemme \ref{lm:kergammast},
et du fait que la transform\'ee de Fourier en un caract\`ere non trivial sur $\K(T)$ est nulle
\begin{equation}\label{eq:intchi0:bis}
\int\limits_{\adh{T(K)}\cap T(\ak)}\!\!\!\!\!\!H(-\bs,t)\,\chi(t)\,\omega_{T}(t)
=
\frac{1}{\card{A(T)}}\,
\sum_{
\substack{
\chi'\in \Ker\left(\gamma^{\ast}\right)
\\
\chi'_{|\K(T)}=1
}
} 
\fourier H(\chi'\,.\,\chi,-\bs )
\end{equation}
On int\`egre \`a pr\'esent les deux membre de \eqref{eq:intchi0:bis} sur $\dualtop{T(\ak)/T(\ak)^1}$ :
\begin{multline}
\int\limits_{\dualtop{T(\ak)/T(\ak)^1}}
\left[\int\limits_{\adh{T(K)}\cap T(\ak)}\!\!\!\!\!\!H(-\bs,t)\,\chi(t)\,\omega_{T}(t)\right] \wt{d\chi}
\\
=
\frac{1}{\card{A(T)}}\,
\sum_{
\substack{
\chi'\in \Ker\left(\gamma^{\ast}\right)
\\
\chi'_{|\K(T)}=1
}
}
\,\,\,\,
\int\limits_{\dualtop{T(\ak)/T(\ak)^1}}
\!\!\!\!\!\!\!\!\!\!\fourier H(\chi'\,.\,\chi,-\bs )\wt{d\chi}.
\end{multline}
Or on a d'apr\`es \eqref{eq:scind:T:dualtop}
\begin{multline}
\sum_{
\substack{
\chi'\in \Ker\left(\gamma^{\ast}\right)
\\
\chi'_{|\K(T)}=1
}
}
\,\,\,\,
\int\limits_{\dualtop{T(\ak)/T(\ak)^1}}
\!\!\!\!\!\!\!\!\!\!\fourier H(\chi'\,.\,\chi,-\bs )\wt{d\chi}
\\
=
\sum_{\chi_1\in \ecU_T}
\,\,
\sum_{
\substack{
\chi_0 \in \dualtop{T(\ak)/T(\ak)^1}
\\
\chi_0\,\chi_1 \in \Ker(\gamma^{\ast})
} 
}
\,\,\,\,
\int\limits_{\dualtop{T(\ak)/T(\ak)^1}}
\!\!\!\!\!\!\!\!\!\!
\fourier H(\chi.\chi_{0}.\chi_1,-\bs )
\wt{d\chi}.
\end{multline}
Comme $\wt{d\chi}$ est une mesure de Haar sur 
$\dualtop{T(\ak)/T(\ak)^1}$, cette derni\`ere expression se r\'e\'ecrit
\begin{equation}
\sum_{\chi_1\in \ecU_T}
\card{\left\{\chi_0 \in \dualtop{T(\ak)/T(\ak)^1},\,\,
\chi_0\,\chi_1 \in \Ker(\gamma^{\ast})\right\}}
\int\limits_{\dualtop{T(\ak)/T(\ak)^1}}
\!\!\!\!\!\!
\fourier H(\chi.\chi_1,-\bs )
\wt{d\chi}.
\end{equation}
Pour conclure, on remarque que pour $\chi_1\in \ecU_T$, on a,
si $\chi_1\notin \Ker(\gamma^{\ast}).\dualtop{T(\ak)/T(\ak)^1}$,
\begin{equation}
\card{\left\{\chi_0 \in \dualtop{T(\ak)/T(\ak)^1},\,\,
\chi_0\,\chi_1 \in \Ker(\gamma^{\ast})\right\}}=0
\end{equation}
et, si $\chi_1\in \Ker(\gamma^{\ast}).\dualtop{T(\ak)/T(\ak)^1}$,
\begin{multline}
\card{\left\{\chi_0 \in \dualtop{T(\ak)/T(\ak)^1},\,\,
\chi_0\,\chi_1 \in \Ker(\gamma^{\ast})\right\}}
=
\card{\dualtop{T(\ak)/T(\ak)^1}\cap \Ker(\gamma^{\ast})}
\end{multline}
soit d'apr\`es le point \ref{item:3:lm:kergammast} du lemme \ref{lm:kergammast}
\begin{equation}
\card{\chi_0 \in \dualtop{T(\ak)/T(\ak)^1},\,\,
\chi_0\,\chi_1 \in \Ker(\gamma^{\ast})}=\KT.
\end{equation}
\end{demo}

Soit $I_1$ la fonction d\'efinie pour $\bs\in \tube{\R^{\Sigma(1)/G}_{>1}}$ par l'int\'egrale
\begin{equation}\label{eq:i1:fonc}
I_1(\bs)
=
\int\limits_{\xtgu}\left[\,\,\int\limits_{\adh{T(K)}\cap T(\ak)}
\!\!\!\!\!\!H(-\bs,t)\,\chi_{\bz}(t)\,\omega_{T}(t)\right]d\bz.
\end{equation}

Pour 
$\chi\in  \ecU_{T} 
$,
soit $I_{\chi}$ 
la fonction d\'efinie pour
$\bs\in \tube{\R^{\Sigma(1)/G}_{>1}}$
par l'int\'egrale
\begin{equation}\label{eq:ichis_fonc}
I_{\chi}(\bs)
=
\int\limits_{\xtgu} 
\fourier H\left(\chi_{\bz}\,\chi,-\bs \right)
d\bz.
\end{equation}

\begin{cor}\label{cor:rep:int:fonc}
Dans le cas fonctionnel, 
pour tout 
 \'el\'ement $\bs$ de 
$\tube{\R^{\,\Sigma(1)/G}_{>1}}$,
on a la repr\'esentation int\'egrale suivante
\begin{multline}\label{eq:repr:int:fonc}
\sum_{t\in T(K)}\,H\left(t,-\bs \right)
\\
=\frac
{1}
{\log(q)^{\rg\left(X(T)^G\right)}\,b(T)}\left(\frac{\card{A(T)}}{\card{\KT}}\,
I_1(\bs)+
\!\!\!\!\!\!
\sum_{
\substack{
\chi\in  \ecU_T  
\\
\chi\notin
\Ker\left(\gamma^{\ast}\right).\dualtop{T(\ak)/T(\ak)^1}
}
}
\!\!\!\!\!\!I_{\chi}(\bs)\right)
\end{multline}
o\`u $I_1(\bs)$ est l'int\'egrale donn\'ee par l'expression \eqref{eq:i1:fonc}
et, pour $\chi\in  \ecU_T  \setminus \Ker\left(\gamma^{\ast}\right).\dualtop{T(\ak)/T(\ak)^1}$, 
$I_{\chi}(\bs)$ est l'int\'egrale donn\'ee par l'expression \eqref{eq:ichis_fonc}.
\end{cor}
\begin{demo}
Soit $\bs$ un \'el\'ement de $\tube{\R^{\,\Sigma(1)/G}_{>1}}$.
On a vu (lemme \ref{lm:H-bsint})
que $H(-\bs,\,.\,)$ est int\'egrable sur $T(\ak)$. 
D'apr\`es la remarque \ref{rem:lm:formule:poisson}, 
on peut alors appliquer le lemme \ref{lm:formule:poisson}.
On obtient donc, compte tenu 
du lemme \ref{lm:dec:int:fonc}, l'\'egalit\'e
\begin{equation}\label{eq:app:poisson:fonc}
\sum_{t\in T(K)}
H\left(t,-\bs \right)
=
\frac{1}{\log(q)^{\rg\left(X(T)^G\right)}\,b(T)}
\sum_{\chi'\in  \ecU_{T} }
\,\,\,
\int\limits_{\dualtop{T(\ak)/T(\ak)^1}}
\!\!\!\!\!\!\!\!\!\!\fourier H\left(\chi\,\chi',-\bs\right)
\wt{d\chi}.
\end{equation}
On d\'ecompose la somme $\sumu{\chi'\in  \ecU_{T} }$ apparaissant 
dans le membre de droite
de \eqref{eq:app:poisson:fonc} 
(qui est, rappelons-le,
une somme finie) suivant que $\chi'$ appartient ou non \`a 
$\Ker(\gamma^{\ast}).\dualtop{T(\ak)/T(\ak)^1}$.

D'apr\`es le lemme \ref{lm:sumchiinkergamma1ast},
on a 
\begin{multline}
\sum_{
\substack{
\chi\in \ecU_T
\\
\chi'\in \Ker(\gamma^{\ast}).\dualtop{T(\ak)/T(\ak)^1}
}
}
\,\,\,\,
\int\limits_{\dualtop{T(\ak)/T(\ak)^1}}
\!\!\!\!\!\!\!
\fourier H\left(\chi.\chi',-\bs \right)\wt{d\chi}
\\
=
\frac{\card{A(T)}}{\card{\KT}}
\,
\int\limits_{\dualtop{T(\ak)/T(\ak)^1}}
\left[\,\,\int\limits_{\adh{T(K)}\cap T(\ak)}\!\!\!\!\!\!\!H(-\bs,t)\,\chi(t)\,\omega_{T}(t)\right]
\wt{d\chi}
\end{multline}

D'apr\`es le corollaire \ref{cor:dec:int:bis}, on a 
\begin{multline}
\int\limits_{\dualtop{T(\ak)/T(\ak)^1}}
\left[\,\,\int\limits_{\adh{T(K)}\cap T(\ak)}\!\!\!\!\!\!\!H(-\bs,t)\,\chi(t)\,\omega_{T}(t)\right]
\wt{d\chi}
\\
=\int\limits_{\xtgu}
\left[\,\,\int\limits_{\adh{T(K)}\cap T(\ak)}\!\!\!\!\!\!\!H(-\bs,t)\,\chi_{\bz}(t)\,\omega_{T}(t)\right]
d\bz.
\end{multline}
Toujours d'apr\`es le corollaire \ref{cor:dec:int:bis},
pour $\chi'\in \ecU_T  \setminus \Ker\left(\gamma^{\ast}\right).\dualtop{T(\ak)/T(\ak)^1}$, on a
\begin{equation}
\int\limits_{\dualtop{T(\ak)/T(\ak)^1}}
\!\!\!\!\!\!\!
\fourier H\left(\chi\,\chi',-\bs\right)
\wt{d\chi}
=
\int\limits_{\xtgu}
\fourier H\left(\chi_{\bz}\,\chi',-\bs\right)
d\bz.
\end{equation}
On en d\'eduit le r\'esultat annonc\'e.
\end{demo}

\paragraph{Le cas o\`u l'hypoth\`ese \eqref{hyp:hinz} est v\'erifi\'ee pour toute place.}

Si toutes les places $v$ de $K$ v\'erifient l'hypoth\`ese \eqref{hyp:hinz},
d'apr\`es les lemmes \ref{lm:fourierHchichibz} et \ref{lm:prop:frf:bis}
on a pour tout $\bs\in \tube{\R^{\Sigma(1)/G}_{>1}}$ 
et tout
$\chi \in \dualtop{T(\ak)/\K(T)}$
\begin{equation}\label{eq:i1:fonc:bis}
I_1(\bs)
=
\int
\limits_{\xtgu}
\!\!
\left(
\prod_{\alpha\in \Sigma(1)/G} 
Z_{\Ka}
\left(\acc{\bz}{\da\,\roa}\,q^{\,-\da\,\sa}\right)
\right)
\frf\left(\gamma_{\unit}(\bz)\,q^{\,-\bs}\right)
d\bz.
\end{equation}
et
\begin{equation}\label{eq:ichis_fonc:bis}
I_{\chi}(\bs)
=
\int\limits_{\xtgu }
\!\!
\left(
\prod_{\alpha\in \Sigma(1)/G} 
\ecL_{\Ka}
\left(
\chia,
\acc{\bz}{\da\,\roa}\,q^{\,-\da\,\sa}
\right)
\right)
\, 
\frf\left(\chi,
\gamma_{\unit}(\bz)\,q^{\,-\bs}\right)
d\bz.
\end{equation}

La technique pour \'evaluer de telles int\'egrales est d\'evelopp\'ee \`a la section \ref{sec:eval:int:fonc}.

\paragraph{Le cas d'une extension de d\'eploiement non ramifi\'ee.}

Le but de ce paragraphe est de d\'ecrire une situation techniquement plus simple que la
situation g\'en\'erale dans le cas fonctionnel, pour laquelle un lemme 
technique simplifi\'e sera suffisant. La compr\'ehension pr\'ealable de ce
qui se passe dans ce cas peut aider \`a la compr\'ehension du traitement du cas g\'en\'eral.

\begin{cor}\label{cor:repr:int:fonc:nonram}
On se place dans le cas fonctionnel et on suppose que l'extension de d\'eploiement $L/K$
est non ramifi\'ee. 
Pour tout 
\'el\'ement de $\bs$ de
$\tube{\R^{\,\Sigma(1)/G}_{>1}}$,
on a la relation
\begin{equation}\label{eq:repr:int:fonc:nonram}
\sum_{t\in T(K)}\,H\left(t,-\bs \right)
=\frac
{1}
{\log(q)^{\rg\left(X(T)^G\right)}\,b(T)}
\sum_{\chi\in  \ecU_T } 
I_{\chi}(\bs)
\end{equation}
o\`u, pour $\chi\in  \ecU_T$,
$I_{\chi}(\bs)$ est l'int\'egrale donn\'ee par 
\begin{equation}
I_{\chi}(\bs)
=
\int\limits_{\xtgu }
\!\!
\left(
\prod_{\alpha\in \Sigma(1)/G} 
\ecL_{\Ka}
\left(
\chia,
\acc{\bz}{\da\,\roa}\,q^{\,-\da\,\sa}
\right)
\right)
\, 
\frQ\left(\chi,
\gamma_{\unit}(\bz)\,q^{\,-\bs}\right)
d\bz.
\end{equation}
\end{cor}
\begin{demo}
D'apr\`es \eqref{eq:app:poisson:fonc} et le 
corollaire \ref{cor:dec:int:bis} on 
a pour tout $\bs\in \tube{\R^{\,\Sigma(1)/G}_{>1}}$
\begin{equation}\label{eq:app:poisson:fonc:nonram}
\sum_{t\in T(K)}
H\left(t,-\bs \right)
=
\frac{1}{\log(q)^{\rg\left(X(T)^G\right)}\,b(T)}
\sum_{\chi\in  \ecU_{T} }
\,\,\,
\int\limits_{\xtgu}
\!\!\!\fourier H\left(\chi\,\chi_{\bz},-\bs\right)
d\bz.
\end{equation}
D'apr\`es le lemme \ref{lm:fourierHchichibz}, on a 
pour tout $\chi\in \ecU_T$
\begin{equation}\label{eq:fourierHchi:nonram}
\fourier H\left(\chi.\chi_{\bz},-\bs\right)
=
\left(
\prod_{\alpha\in \Sigma(1)/G} 
\ecL_{\Ka}
\left(
\chia,
\acc{\bz}{\da\,\roa}\,\,q^{\,-\da\,\sa}
\right)
\right)
\, 
\frf\left(\chi,\gamma_{\unit}(\bz)\,q^{\,-\bs }\right)
\end{equation}
Comme la vari\'et\'e $\xs$ est d\'eploy\'ee par une extension non ramifi\'ee,
d'apr\`es la d\'efinition \eqref{eq:deffrfchi} de $\frf\left(\chi,\,.\,\right)$
on a 
$
\frf(\chi,\,.\,)
=
\,
\frQ(\chi,\,.\,)
$
d'o\`u le r\'esultat.
\end{demo}

\section{\'Evaluation de l'int\'egrale : le cas arithm\'etique}
\label{sec:eval:int:arit}

Soit $N$ un $\Z$-module libre de rang fini et $\Lambda$ un c\^one strictement convexe de $N_{\R}$.
On d\'efinit, pour tout \'el\'ement $\bs$ de $\tube{\intrel \left(\Lambda\right)}$,
\begin{equation}
\indic_{N,\Lambda}(\bs)
=
\int\limits_{\Lambda^{\vee}}
e^{\,-\acc{y}{\bs}}\,dy,
\end{equation}
o\`u $dy$ est la mesure de Lebesgue sur $N^{\vee}_{\R}$, 
normalis\'ee de sorte que le r\'eseau $N^{\vee}$
soit de covolume $1$. 
Cette fonction sera appel\'ee 
\dindex{$\indic$-fonction}
\termin{
$\indic$-fonction} du c\^one $\Lambda$.
Notons que si $\bs_0$ est un \'el\'ement de $\tube{\intrel \left(\Lambda\right)}$
on a 
\begin{equation}
\indic_{N,\Lambda}(\bs_0)
=
\lim_{t\to 0}t^{\,\rg(N)}\,\indic_{N,\Lambda}(t\,\bs_0).
\end{equation}
et que si $N'\subset N$ est un sous-groupe d'indice fini on a 
\begin{equation}
\label{eq:formule:indic}
\indic_{N,\Lambda}(\bs)
=
[N:N']
\,
\indic_{N',\Lambda}(\bs).
\end{equation}
\begin{rem}
D'apr\`es \eqref{eq:defalpha}, on a en particulier, pour toute vari\'et\'e projective et lisse $V$
telle que la classe du faisceau anticanonique appartient \`a l'int\'erieur du c\^one effectif de $V$  
\begin{equation}
\label{eq:defalpha:bis}
\alpha^{\ast}(V)
=
\indic_{\Pic(\xs),\ceff(V)}(\omega_{V}^{-1})
\end{equation}
o\`u $\ceff(V)$ est le c\^one effectif de $V$.
\end{rem}

L'\'evaluation est bas\'ee sur le r\'esultat suivant, cas particulier
de \cite[Thm 3.1.14]{CLTs:fibres}.

\begin{thm}\label{thm:CLTs} 
Soit $n\geq 1$ un entier et $\Gamma$ un sous-groupe de $\Z^n$ tel que 
$N=\Z^n/\Gamma$ est sans torsion 
et $\Gamma\cap \R_{\geq 0}^n=\{0\}$. 
On note $j$ le morphisme quotient
$\Z^n\to N$.

Soit $f$ une fonction holomorphe sur $\tube{\R^n_{>0}}$.
On suppose qu'il existe un $\eps>0$
tel que la fonction
\begin{equation}
\bs\mapsto f(\bs)\,\prod_{1\leq i\leq n}\frac{s_i}{1+s_i}
\end{equation}
se prolonge en une fonction holomorphe sur $\tube{\R^n_{>-\eps}}$
qui est $\Gamma_{\R}$-contr\^ol\'ee 
au sens de \cite[D\'efinition 3.13]{CLTs:fibres}.
On note $C$ la valeur de ce prolongement en $0$.

Alors l'int\'egrale
\begin{equation}
\frac{1}{(2\,\pi)^{\rg(\Gamma)}}
\int\limits_{y\in \Gamma_{\R}}
f(\bs+i\,y)
dy
\end{equation}
converge absolument en tout $\bs$ de $\tube{\R^n_{>0}}$
et d\'efinit une fonction holomorphe 
$\Gamma_{\C}$-invariante sur $\tube{\R^n_{>0}}$, not\'ee $g$.

Il existe en outre un $\eps'>0$
tel que la fonction
\begin{equation}
R\,:\,\bs\mapsto g(\bs)-C\,\indic_{N,j\left(\R^n_{\geq 0}\right)}(j(\bs))
\end{equation}
se prolonge en une fonction m\'eromorphe sur $\tube{\R^n_{>-\eps'}}$

Pour tout $\bs_0\in \tube{\R^n_{>0}}$, on a
\begin{equation}
\lim_{s\to 0} s^{\rg(N)} R(s\,\bs_0)=0.
\end{equation}
\end{thm}

\begin{rems}
\begin{enumerate}
\item
La $\Gamma_{\R}$-contr\^olabilit\'e  de $f$
signifie en deux mots que l'on dispose d'un bon contr\^ole de la croissance de $f$ sur les bandes en $\Gamma_{\R}$.
Nous n'avons pas estim\'e utile de rappeler ici la d\'efinition pr\'ecise de cette notion.
Il suffit de savoir que, pour l'application du th\'eor\`eme \`a l'\'evaluation du
comportement asymptotique de la fonction z\^eta des hauteurs, cette hypoth\`ese
est v\'erifi\'ee gr\^ace au point \ref{item3:prop:contr:F} de la proposition
\ref{prop:contr:F}, comme d\'ej\`a signal\'e \`a la remarque \ref{rem:prop:contr:F}.
\item
Le r\'esultat obtenu en appliquant  \cite[Thm 3.1.14]{CLTs:fibres} est en fait 
plus pr\'ecis. On obtient une description des p\^oles de $g$
au voisinage de z\'ero, et un contr\^ole de $g$ dans les bandes verticales.
Une fois le r\'esultat appliqu\'e \`a la fonction z\^eta des hauteurs,
ce contr\^ole joint \`a un th\'eor\`eme taub\'erien ad\'equat  permet de donner un d\'eveloppement asymptotique 
du nombre de points de hauteur born\'ee plus pr\'ecis que celui qui d\'ecoule du th\'eor\`eme \ref{thm:theo:ba:ts}.
\item
L'id\'ee g\'en\'erale de la preuve de \cite[Thm 3.1.14]{CLTs:fibres} 
est de mettre en oeuvre une r\'ecurrence 
sur le rang de $\Gamma$ utilisant le th\'eor\`eme des r\'esidus.
La notion de contr\^olabilit\'e sert \`a assurer l'int\'egrabilit\'e
(et le contr\^ole dans les bandes) des fonctions obtenues par int\'egrations successives.
Cette preuve nous semble inadaptable telle quelle au type de fonctions que l'on a \`a traiter dans
le cas fonctionnel. Une des raisons est l'apparition de p\^oles suppl\'ementaires. 
En outre, ce que devrait \^etre la d\'efinition 
de la contr\^olabilit\'e dans ce cadre n'est pas clair a priori  ; notamment, 
au vu de la p\'eriodicit\'e des fonctions mises en jeu, la notion de majoration dans les bandes 
perd tout son int\'er\^et. 
\end{enumerate}
\end{rems}

\section{\'Evaluation de l'int\'egrale : le cas fonctionnel}
\label{sec:eval:int:fonc}

Le but de cette partie est d'obtenir 
un r\'esultat analogue au th\'eor\`eme \ref{thm:CLTs}, 
adapt\'e \`a la forme des fonctions obtenues dans le cas fonctionnel.

\subsection{D\'efinition d'une certaine classe de fonctions}

Soit $N$ un $\Z$-module libre de rang fini.
Pour toute partie $A$ de $N_{\R}$ et toute application $a\,:\,A\cap N\to \C$ on d\'efinit
la s\'erie formelle
\nindex{$\ecL_{N,A,a}$}
\begin{equation}
\ecL_{N,A,a}(T)=\sum_{y\in A\,\cap\,N}a_y\,T^{y}.
\end{equation}
On pose
\begin{equation}\label{eq:eclnaz}
\ecL_{N,A,a}(\bz)=\sum_{y\in A\,\cap\,N}a_y\,\acc{\bz}{y}
\end{equation}
pour tout \'el\'ement $\bz$ de $N^{\vee}_{\ca}$
telle que le membre de droite de \eqref{eq:eclnaz}
est une s\'erie absolument convergente.
On a ainsi
\begin{equation}\label{eq:eclnaqs}
\ecL_{N,A,a}(q^{-\bs})
=
\sum_{y\in A\,\cap\,N}a_y\,q^{\,-\acc{y}{\bs}}
\end{equation}
pour tout \'el\'ement $\bs$ de $N^{\vee}_{\C}$ tel que 
telle que le membre de droite de \eqref{eq:eclnaqs}
est une s\'erie absolument convergente. 
Pour un tel $\bs$ on a donc, pour tout \'el\'ement $\bz$ de $N^{\vee}_{\unit}$,
\begin{equation}
\ecL_{N,A,a}(\bz\,q^{-\bs})=\sum_{y\in A\,\cap\,N}a_y\,\acc{\bz}{\,y}\,q^{\,-\acc{y}{\bs}}.
\end{equation}
Si $a$ est la fonction constante \'egale \`a 1, 
\nindex{$\ecL_{N,A}$}
on notera $\ecL_{N,A}$ pour $\ecL_{N,A,a}$.

\subsection{Un premier exemple}

Soit $\Upsilon$ un c\^one de $N_{\R}$, tel que $\Upsilon^{\vee}$ soit d'int\'erieur non vide. 
La s\'erie d\'efinissant $\ecL_{N,\Upsilon}\left(q^{-\bs}\right)$ converge alors absolument pour tout $\bs$ de 
$\tube{\text{int}\left(\Upsilon^{\vee}\right)}$.
La fonction $\ecL_{N,\Upsilon}\left(q^{-\bs}\right)$ est utilis\'ee par Peyre 
dans \cite{Pey:prepu:drap} pour une d\'efinition alternative 
de l'invariant $\alpha^{\ast}(V)$ 
attach\'e \`a une vari\'et\'e $V$ 
(cf. la formule \eqref{eq:defalpha} pour la d\'efinition adopt\'ee dans ce texte).
On a en effet le r\'esultat suivant :
\begin{lemme}\label{lm:fonction_cone}
Soit $\lambda_0$ un \'el\'ement de l'int\'erieur de $\Upsilon^{\vee}$. 
L'application 
\begin{equation}
s\longmapsto \ecL_{N,\Upsilon}\left(q^{\,-s.\lambda_0}\right)
\end{equation}
est bien d\'efinie et holomorphe sur $\tube{\R_{>0}}$, 
et se prolonge en une fonction m\'eromorphe sur $\C$, 
avec un p\^ole d'ordre au plus la dimension de $\Upsilon$ en $s=0$.
Si de plus $\Upsilon^{\vee}$ est strictement convexe, 
l'ordre de ce p\^ole est exactement $\rg(N)$, et on a 
\begin{equation}
\lim_{s\to 0} 
\left(s^{\,\rg(N)}\,
\ecL_{N,\Upsilon}\left(q^{\,-s.\lambda_0}\right)\right)
=
\log(q)^{\,-\rg(N)}\,\indic_{N,\Upsilon^{\vee}}(\lambda_0).
\end{equation}
\end{lemme}
\begin{rem}
On a donc en particulier d'apr\`es \eqref{eq:defalpha:bis}
\begin{equation}\label{eq:alpha}
\alpha^{\ast}(\xs)
=
\log(q)^{\,\rg(\Pic(\xs))}\,
\lim_{s\to 0} 
\left[
s^{\rg(\Pic(\xs))}\,
\ecL_{\Pic(\xs)^{\vee},\ceff(\xs)^{\vee}}
\left(q^{\,-s\,\left[-\can_{\xs}\right]}\right)
\right].
\end{equation}
\end{rem}
\begin{rem}
Si $N'\subset N$ est un sous-groupe d'indice fini de $N$ on a 
\begin{equation}\label{eq:sous-reseau}
\lim_{s\to 0} 
\left[s^{\,\rg(N)}\,\ecL_{N,\Upsilon}\left(q^{-s\,\lambda_0}\right)\right]
=
[N:N']\,\lim_{s\to 0}\,
\left[
s^{\,\rg(N)}\,\ecL_{N',\Upsilon}
\left(q^{-s\,\lambda_0}\right)
\right].
\end{equation}
\end{rem}

\begin{demo}
On \'ecrit $\Upsilon$ comme le support d'un \'eventail r\'egulier $\Delta$.
Concernant cet \'eventail, on reprend les notations introduites \`a la section \ref{subsec:vte:tor}.
 
On a alors
\begin{equation}
\ecL_{N,\Upsilon}(T)=\sum_{\delta\in \Delta}\,L_{N,N\cap \intrel(\delta)}(T),
\end{equation}
soit encore en termes plus explicites
\begin{equation}
\ecL_{N,\Upsilon}(T)
=
\sum_{\delta\in \Delta}\,\,\prod_{l\in \delta(1)}\,\,\left(\frac{1}{1-T^{\,\rho_l}}-1\right).
\end{equation}
On a donc
\begin{equation}\label{eq:expr_expl_LNupsilon}
\ecL_{N,\Upsilon}\left(q^{\,-\bs}\right)
=
\sum_{\delta\in \Delta}\,\,\prod_{l\in \delta(1)}\,\,\left(\frac{1}{1-q^{\,-\acc{\rho_l}{\bs}}}-1\right),
\end{equation}
d'o\`u le r\'esultat, car le cardinal maximal des ensembles $\delta(1)$ est \'egal \`a la dimension de $\Upsilon$. 

Supposons \`a pr\'esent $\Upsilon$ de dimension $\rg(N)$, 
et montrons la derni\`ere assertion du lemme.
Comme les c\^ones de dimension maximale de $\Delta$ recouvrent 
$\Upsilon$ et que leurs intersections sont de mesure de Lebesgue nulle,
on peut \'ecrire
\begin{align}
\indic_{N,\Upsilon^{\vee}}(\lambda_0)
&
=\int_{\Upsilon}
e^{\,-\acc{y}{\lambda_0}}\,dy,
\\
&
=
\sum_{
\substack{
\delta\in \Delta
\\
\dim(\delta)=\rg(N)
}}
\int_{\delta}
e^{\,-\acc{y}{\lambda_0}}\,dy,
\\
&
=
\sum_{
\substack{
\delta\in \Delta
\\
\dim(\delta)=\rg(N)
}}
\,\,\,
\int\limits_{\R_{>0}^{\delta(1)}}
e^{\,-\acc{\sumu{l\in \delta(1)}x_l\,\rho_l}{\lambda_0}}\,\otimesu{l\in \delta(1)}dx_l
\\
&
=
\sum_{
\substack{
\delta\in \Delta
\\
\dim(\delta)=\rg(N)
}}
\prod_{l\in \delta(1)}
\frac{1}{\acc{\rho_l}{\lambda_0}}.
\end{align}
Or, d'apr\`es \eqref{eq:expr_expl_LNupsilon} et le fait que pour tout complexe non nul
$z$ on a 
\begin{equation}
\lim_{s\to 0} \frac{s}{1-q^{\,s\,z}}=-\frac{1}{\log(q)\,z},
\end{equation}
on obtient 
\begin{equation}
\lim_{s\to 0} s^{\,\rg(N)}\,
\ecL_{N,\Upsilon}
\left(q^{\,-s\,\lambda_0}\right)
=
\sum_{
\substack{
\delta\in \Delta
\\
\dim(\delta)=\rg(N)
}
}\,\,
\log(q)^{\,-\rg(N)}\,\prod_{l\in \delta(1)}\,\frac{1}{\acc{\rho_l}{\lambda_0}},
\end{equation}
d'o\`u le r\'esultat.

\end{demo}

\subsection{Encore quelques d\'efinitions}

On fixe d\'esormais 
{\em pour toute la suite de la section \ref{sec:eval:int:fonc}} 
une base $(\lambda_i)_{i\in I}$ 
de $N$. On note $(\lambda_i^{\vee})_{i\in I}$ la base duale d'une telle base, 
$\Lambda$ le c\^one simplicial 
de $N_{\R}$ engendr\'e par cette base et 
\begin{equation}
\lambda^{\vee}=\sum_{i\in I}\lambda_i^{\vee}.
\end{equation}

On d\'efinit, pour tout r\'eel $\eta$,
\begin{equation}
\Lambda^{\vee}_{>\eta}=\sum_{i\in I}\R_{>\eta}\,\lambda_i^{\vee}\,\subset\,N^{\vee}_{\R}.
\end{equation}
Ainsi $\Lambda^{\vee}_{>0}$ est l'int\'erieur de $\Lambda^{\vee}$.

Soit 
\begin{equation}
a\,:\,\Lambda\cap N \to \C
\end{equation}
une application et $\eps>0$ un r\'eel
tels que la s\'erie d\'efinissant
$
\ecL_{N,\Lambda,a}\left(q^{\,-\bs}\right)
$
converge absolument 
pour tout $\bs$ 
du domaine $\tube{\Lambda^{\vee}_{>-\eps}}$.
Ceci \'equivaut \`a demander 
la convergence de la s\'erie
\begin{equation}
\sum_{y\in \Lambda\cap N}
\abs{a_{y}}
\,
q^{\,\eta\,\acc{y}{\lambda^{\vee}}}
\end{equation}
pour tout $\eta<\eps$.

La fonction 
\begin{equation}
\bs
\mapsto
\ecL_{N,\Lambda,a}\left(q^{\,-\bs}\right).
\end{equation}
est donc holomorphe sur 
le domaine
$
\tube{\Lambda^{\vee}_{>-\eps}}.
$
Une telle fonction sera appel\'ee \termin{fonction admissible \'el\'ementaire de multiplicit\'e sup\'erieure \`a z\'ero}.

Si $r\geq 1$ est un entier, 
une fonction $f$ holomorphe sur  
$
\tube{\Lambda^{\vee}_{>0}}
$
sera appel\'ee \termin{fonction admissible \'el\'ementaire de multiplicit\'e sup\'erieure \`a $-r$} 
s'il existe une \'ecriture 
\begin{equation}
f(\bs)=g(\bs)\,\ecL_{N',\intrel\left(\Upsilon\right)}\left(q^{-\bs}\right)
\end{equation}
o\`u $\Upsilon \subset \Lambda$ est un c\^one de dimension inf\'erieure \`a $r$, 
$N'$ est un sous-groupe de $N$ et $g$ est admissible \'el\'ementaire de multiplicit\'e sup\'erieure \`a $0$.
Une telle fonction $f$ se prolonge donc en une fonction m\'eromorphe sur le domaine
$
\tube{\Lambda^{\vee}_{>-\eps}}
$
pour un certain $\eps>0$.

Si $r\geq 0$ est un entier, on appellera 
\dindex{fonction admissible de multiplicit\'e sup\'erieure \`a $-r$} 
\termin{fonction admissible de multiplicit\'e sup\'erieure \`a $-r$} 
une fonction $f$ holomorphe sur le domaine
$
\tube{\Lambda^{\vee}_{>0}}
$
qui s'\'ecrit comme une somme finie de fonction admissibles \'el\'ementaires 
de multiplicit\'e sup\'erieure \`a $-r$.

Une telle fonction $f$ se prolonge donc en une fonction m\'eromorphe sur le domaine
$
\tube{\Lambda^{\vee}_{>-\eps}}
$
pour un certain $\eps>0$.
Par ailleurs, pour tout $\lambda_0^{\vee}$ \'el\'ement de l'int\'erieur de $\Lambda^{\vee}$, 
la fonction
d'une variable complexe
\begin{equation}
s\mapsto f\left(s\,\lambda_0^{\vee}\right)
\end{equation}
est m\'eromorphe sur un voisinage de z\'ero, et a un p\^ole d'ordre au plus $r$ en z\'ero.

\subsection{Avertissement au lecteur}\label{subsec:avertissement} 
Nous allons donner ci-dessous trois versions du lemme technique d'int\'egration
destin\'e \`a \'evaluer le comportement analytique de la fonction z\^eta des hauteurs
\`a partir de la repr\'esentation int\'egrale obtenue \`a la sous-section \ref{subsubsec:expr:int:fonc}.
Ces trois versions seront de g\'en\'eralit\'e (et de difficult\'e technique) croissante.
Le parti pris de ne pas pr\'esenter directement la version la plus g\'en\'erale nous semble
utile pour la compr\'ehension de la technique employ\'ee.

La premi\`ere version (lemme \ref{lm:tech0}) est une version-jouet, tr\`es simple et transparente. 
Elle est destin\'ee \`a faire comprendre
l'id\'ee \'el\'ementaire de base qui sous-tend les versions plus \'elabor\'ees qui vont suivre,
mais ne nous servira pas pour l'\'evaluation de la fonction z\^eta des hauteurs.

La deuxi\`eme version (lemme \ref{lm:tech1}) est une g\'en\'eralisation naturelle de la premi\`ere, et est 
suffisante pour traiter le cas des vari\'et\'es toriques d\'eploy\'ees par une extension non ramifi\'ee.
Elle s'appuie sur un <<lemme de d\'ecomposition>> de certaines fonctions caract\'eristiques
<<tordues>> associ\'ees \`a des c\^ones (lemme \ref{lm:crucial}).

La troisi\`eme version (lemmes \ref{lm:tech1:bis} et \ref{lm:tech1:bis:tilde}), la plus g\'en\'erale, 
est n\'ecessaire pour traiter le cas d'une extension 
de d\'eploiement quelconque, et pr\'esente quelques complications techniques 
qui peuvent la rendre un peu obscure au premier abord. Elle s'appuie sur des g\'en\'eralisations
du lemme de d\'ecomposition \ref{lm:crucial} (lemmes \ref{lm:crucial:bis} et \ref{lm:crucial:bis:tilde})

Nous conseillons donc de ne pas aborder en premi\`ere lecture 
les d\'emonstrations des lemmes \ref{lm:crucial:bis}, \ref{lm:crucial:bis:tilde}
\ref{lm:tech1:bis} et \ref{lm:tech1:bis:tilde} (i.e. les sous-sections
\ref{subsubsec:lemmedecom:generale} et \ref{subsubsec:lm:techgene}), ainsi que 
l'application qui en est faite \`a l'\'evaluation du comportement analytique de
la fonction z\^eta des hauteurs des vari\'et\'es toriques dans le cas g\'en\'eral
(sous-section \ref{subsubsec:application:lemme:technique:fonc:gene}).

\subsection{Un lemme de d\'ecomposition}\label{subsec:lemmedecom}
\subsubsection{Version simple}\label{subsubsec:lemmedecom:simple}

Soit $\Upsilon$ un c\^one de $N_{\R}$ contenu dans $\Lambda$. 
On notera $\left\langle \Upsilon \right\rangle$ le sous-espace vectoriel de $N_{\R}$ engendr\'e par $\Upsilon$.

On \'ecrit $\Upsilon$ comme le support d'un \'eventail r\'egulier $\Delta$.
Concernant cet \'eventail, on reprend les notations introduites \`a la section \ref{subsec:vte:tor}.
Pour $l\in\Delta(1)$, $\rho_l$ d\'esigne donc le g\'en\'erateur du mono\"\i de
$N\cap l$. Pour toute partie $I$ de $\Delta(1)$ nous noterons 
$
\ecC(I)
$
le c\^one engendr\'e par les $(\rho_i)_{i\in I}$.

On fixe un \'el\'ement $z$ de $\Lambda\cap N$.
On pose
\begin{equation}
\Upsilon_{z}=\Upsilon\cap \{z+\Lambda\cap N\}.
\end{equation}
On veut \'etudier la s\'erie formelle 
\begin{equation}
\ecL_{N,\Upsilon_{z}}(T)
=
\sum_{\substack{y\in \Upsilon\cap N\\ \\ y\in z+\Lambda\cap N}} T^{\,y}.
\end{equation}

Nous utilisons la m\^eme technique que dans la section 4.3.3. de \cite{Bou:crelle}. 
Nous \'ecrivons d'abord
\begin{equation}
\ecL_{N,\Upsilon_{z}}(T)
=
\sum_{\delta\in\Delta}
\,
\ecL_{N,\Upsilon_{z}
\,
\cap
\, 
\intrel(\delta)}(T)
.
\end{equation}

Pour $\delta\in\Delta$, notons qu'on peut \'ecrire
\begin{equation}\label{eq:dec_lnupsilon:0}
\ecL_{N,\Upsilon_{z}\cap\, \intrel(\delta)}(T)
=
\sum_{
\substack
{
y\in \intrel(\delta)
\\ 
\\ 
y-z\in \Lambda
}
}
T^{\,y}
=
\sum_{
\substack
{
y\in \intrel(\delta)
\\ 
\\ 
\forall i\in I,\,\,\acc{\lambda_i^{\vee}}{y}\geq \acc{\lambda_i^{\vee}}{z} 
}
}
T^{\,y}
.
\end{equation}
Pour tout sous-ensemble $K$ de $I$ et tout c\^one $C$ de $N_{\R}$,  
on note $C(K,z)$ le sous-ensemble de $\intrel(C)$ 
form\'e des \'el\'ements $y$ v\'erifiant la condition
\begin{equation}
\forall i\in K,\, 
\acc{\lambda^{\vee}_i}{y}<\acc{\lambda^{\vee}_i}{z}.
\end{equation}
En utilisant \eqref{eq:dec_lnupsilon:0}, on voit alors qu'on a une d\'ecomposition
\begin{equation}
\ecL_{N,\Upsilon_{z}\,\cap\,\intrel(\delta)}(T)
=
\sum_{K\subset I}\,
(-1)^{\card{K}}\,
\ecL_{N,\delta(K,z)}(T).
\end{equation}

Posons
\begin{equation}
\ecM=\card{I}\,\underset{l\in \Delta(1)}{\Sup}\acc{\rho_l}{\lambda}.
\end{equation}

\begin{lemme}\label{lm:crucial}
Soient $\delta$ un c\^one de $\Delta$, 
et $K$ une partie  de $I$. Soit $\delta(1)_K$ le sous-ensemble  de $\delta(1)$
donn\'e par  
\begin{equation}\label{eq:def_delta1K}
\delta(1)_K=\{\,l\in \delta(1),\,\forall\,i\in K,\quad\acc{\rho_l}{\lambda_i} =0\}
\end{equation}
et $F(\delta,K,z)$ le sous-ensemble de $\delta(K,z)$ donn\'e par
\begin{equation}
F(\delta,K,z)=\ecC(\delta(1)\setminus \delta(1)_K)(K,z).
\end{equation}

Alors
$\delta(1)_K$
et
$F(\delta,K,z)$
v\'erifient les propri\'et\'es suivantes :
\begin{itemize}
\item
On a une \'ecriture
\begin{equation}
\ecL_{N,\delta(K,z)}(T)
=
\ecL_{N,\intrel(\ecC(\delta(1)_K))}(T)
\,
\ecL_{N,F(\delta,K,z)}(T)
\end{equation}
\item
$F(\delta,K,z)$ est fini. 
Plus pr\'ecis\'ement on a la majoration
\begin{equation}\label{eq:maj1}
\card{F(\delta,K,z)}\leq \acc{z}{\lambda}^{\,\rg(N)}.
\end{equation}
\item
Pour tout $y\in F(\delta,K,z)$ on a
\begin{equation}\label{eq:maj2}
\acc{y}{\lambda}\leq \ecM \acc{z}{\lambda}.
\end{equation}
\item
Si $K$ est vide, on a $\delta(1)_K=\delta(1)$ et $F(\delta,K,z)=\varnothing$. 
\item
On suppose en outre que $\left(N/\left\langle \Upsilon \right\rangle \right)^{\vee}\cap \Lambda^{\vee}=\{0\}$,
\label{item:lm:crucial:4}
que $\delta$ est de dimension maximale, et que $K$ est non vide. Alors $\delta(1)_K$ est un sous-ensemble strict de $\delta(1)$. 
\end{itemize}

\end{lemme}

\begin{demo} Cette d\'emonstration est tr\`es similaire \`a celle du lemme 3 de \cite{Bou:crelle}.

Soit $y$ un \'el\'ement de $\intrel(\delta)$. 
Il s'\'ecrit donc de mani\`ere unique 
$y_1+y_2$ avec $y_1\in \intrel(\ecC(\delta(1)_K))$ 
et $y_2\in \intrel(\ecC(\delta(1)\setminus \delta(1)_K))$.

Au vu de la d\'efinition \eqref{eq:def_delta1K} de $\delta(1)_K$, on a
\begin{equation}
\forall\,i\in K,\quad\acc{y}{\lambda_i}=\acc{y_2}{\lambda_i}.
\end{equation}
Supposons en outre que $y$ est dans $\delta(K,z)$, i.e. v\'erifie
\begin{equation}
\forall\,i\in K,\quad\acc{y}{\lambda_i^{\vee}} < \acc{z}{\lambda_i^{\vee}}.
\end{equation}
On a donc 
\begin{equation}
\forall\,i\in K,\quad\acc{y_2}{\lambda_i^{\vee}} < \acc{z}{\lambda_i^{\vee}},
\end{equation}
en d'autres termes $y_2$ est dans $\ecC(\delta(1)\setminus \delta(1)_K,K,z)=F(\delta,K,z)$.

R\'eciproquement, on constate que 
si on a 
\begin{equation}
y_1\in \intrel(\ecC(\delta(1)_K)) 
\end{equation}
et
\begin{equation}
y_2\in \ecC(\delta(1)\setminus \delta(1)_K,K,z)=F(\delta,K,z)
\end{equation}
alors $y_1+y_2$ est un \'el\'ement de $\delta(K,z)$.

Ceci montre qu'on a une d\'ecomposition
\begin{equation}
\ecL_{N,\delta(K,z)}(T)=\ecL_{N,\intrel(\ecC(\delta(1)_K))}(T)\,\ecL_{N,F(\delta,K,z)}(T).
\end{equation}

Montrons que $F(\delta,K,z)$ est fini et majorons son cardinal. 
Soit $y_2$ un \'el\'ement de $F(\delta,K,z)$, 
qu'on \'ecrit
\begin{equation}
y_2=\sum_{l\in \delta(1)\setminus \delta(1)_K}\,\mu_l\,\rho_l
\end{equation}
avec les $\mu_l$ dans $\N_{>0}$.

Par d\'efinition de $\delta(1)_K$, pour tout $l$ de $\delta(1)\setminus \delta(1)_K$, 
il existe $i$ dans $K$ v\'erifiant $\acc{\rho_l}{\lambda_i}\,\geq \,1$ 
(rappelons que pour tout $i$ on a $\acc{\rho_l}{\lambda_i} \geq 0$). 
Comme $y_2$ v\'erifie
\begin{equation}
\forall i\in K,\quad\acc{y_2}{\lambda_i} < \acc{z}{\lambda_i}
\end{equation}
on a
\begin{equation}
\mu_l < \underset{i\in K}{\Sup} \,\acc{z}{\lambda_i}\leq \acc{z}{\lambda}.
\end{equation}
Ainsi $F(\delta,K,z)$ est fini et son cardinal est major\'e par 
\begin{equation}
\acc{z}{\lambda}^{\,\card{\delta(1)\setminus \delta(1)_K}}
\leq \acc{z}{\lambda}^{\,\rg(N)}.
\end{equation}
Par ailleurs un $y_2$ de $F(\delta,K,z)$ v\'erifie 
\begin{equation}
0
\leq 
\acc{y_2}{\lambda} 
\leq 
\sum_{l\in \delta(1)\setminus\delta(1)_K}\,
\acc{z}{\lambda}\,
\acc{\rho_l}{\lambda}\,
\leq M\,\acc{z}{\lambda}
\end{equation}
o\`u on rappelle que
\begin{equation}
\ecM=\card{I}\,\underset{l\in \Delta(1)}{\Sup}\left(\acc{\rho_l}{\lambda} \right).
\end{equation}

Supposons 
$\left(N/\left\langle \Upsilon \right\rangle\right)^{\vee} \cap \Lambda^{\vee}=\{0\}$.
Ceci \'equivaut \`a dire que si $\lambda\in \Lambda^{\vee}$ v\'erifie $\acc{m}{\lambda}=0$ 
pour tout $m\in \left\langle \Upsilon \right\rangle$ alors $\lambda=0$.
Mais si $\delta$ est de dimension maximale, 
les $(\rho_l)_{l\in \delta(1)}$ engendrent $\left\langle \Upsilon \right\rangle$, 
et donc si $K$ n'est pas vide, on ne peut avoir $\delta(1)_K=\delta(1)$. 
\end{demo}

\subsubsection{Version g\'en\'erale}\label{subsubsec:lemmedecom:generale}

On consid\`ere toujours $\Upsilon$ un c\^one de $N_{\R}$ contenu dans $\Lambda$.
On se donne en outre $N'\subset N$ un sous-groupe d'indice fini.

On \'ecrit cette fois $\Upsilon$ comme le support d'un \'eventail $N'$-r\'egulier $\Delta$
(i.e les c\^ones de $\Delta$ sont engendr\'es par des parties de bases de $N'$), et 
on reprend \`a cet effet les notations introduites \`a la section \ref{subsec:vte:tor}.
Soulignons que pour $l\in\Delta(1)$, $\rho_l$ d\'esigne le g\'en\'erateur du mono\"\i de
$N'\cap l$.

On fixe un \'el\'ement $z$ de $\Lambda\cap N$.
On pose
\begin{equation}
\Upsilon_{z}=\Upsilon\cap \{z+\Lambda\cap N'\}.
\end{equation}
On veut \'etudier la s\'erie formelle
\begin{equation}
\ecL_{N,\Upsilon_{z}}(T)
=
\sum_{\substack{y\in \Upsilon\cap N\\ \\ y\in z+\Lambda\cap N'}} T^{\,y}.
\end{equation}
Nous \'ecrivons d'abord
\begin{equation}
\ecL_{N,\Upsilon_{z}}(T)
=
\sum_{\delta\in\Delta}
\,
\ecL_{N,\Upsilon_{z}
\,
\cap
\, 
\intrel(\delta)}(T)
.
\end{equation}
Pour $\delta\in\Delta$, notons qu'on peut \'ecrire
\begin{equation}\label{eq:dec_lnupsilon}
\ecL_{N,\Upsilon_{z}\cap\, \intrel(\delta)}(T)
=
\sum_{
\substack
{
y\in \intrel(\delta)
\\ 
\\ 
y-z\in \Lambda
\\
\\
y-z\in N'
\\
}
}
T^{\,y}
=
\sum_{
\substack
{
y\in \intrel(\delta)
\\ 
\\ 
\forall i\in I,\,\,\acc{\lambda_i^{\vee}}{y}\geq \acc{\lambda_i^{\vee}}{z} 
\\
\\
y-z\in N'
}
}
T^{\,y}
.
\end{equation}
Pour tout sous-ensemble $K$ de $I$ et tout c\^one $C$ de $N_{\R}$,  
on note $C(K,z)$ le sous-ensemble de $\intrel(C)$ 
form\'e des \'el\'ements $y$ v\'erifiant la condition
\begin{equation}
\forall i\in K,\, 
\acc{\lambda^{\vee}_i}{y}<\acc{\lambda^{\vee}_i}{z}.
\end{equation}
Pour tout sous-ensemble $A$ de $N_{\R}$, on note $A_z$
l'ensemble des \'el\'ements $y$ de $A$ v\'erifiant la condition
\begin{equation}
y-z\in N'.
\end{equation}

En utilisant \eqref{eq:dec_lnupsilon}, on voit alors qu'on a une d\'ecomposition
\begin{equation}
\ecL_{N,\Upsilon_{z}\,\cap\,\intrel(\delta)}(T)
=
\sum_{K\subset I}\,
(-1)^{\card{K}}\,
\ecL_{N,\delta(K,z)_z}(T).
\end{equation}

Posons
\begin{equation}
\ecM=\card{I}\,\underset{l\in \Delta(1)}{\Sup}\acc{\lambda^{\vee}}{\rho_l}.
\end{equation}

\begin{lemme}\label{lm:crucial:bis}
Soient $\delta$ un c\^one de $\Delta$, 
et $K$ une partie  de $I$. 
Soit $\delta(1)_K$ le sous-ensemble de $\delta(1)$
d\'efini par  
\begin{equation}\label{eq:def_delta1K:0} 
\delta(1)_K=\{\,l\in \delta(1),\,\forall\,i\in K,\quad\acc{\lambda_i^{\vee}}{\rho_l} =0\}.
\end{equation}

On compl\`ete $(\rho_l)_{l\in \delta(1)}$ en une base $(\rho_l)_{l\in L}$ de $N'$.
Soit $N'_1$ (respectivement $N'_2$, respectivement $N'_3$)
le sous-groupe de $N$ engendr\'e par les
$(\rho_l)_{l\in \delta(1)_K}$
(respectivement
$(\rho_l)_{l\in \delta(1)\setminus \delta(1)_K}$,
respectivement
$(\rho_l)_{l\in L\setminus \delta(1)}$).

On \'ecrit 
\begin{equation}
z=z_1+z_2+z_3
\end{equation} 
avec  pour $i=1,2,3$, $z_i\in (N'_i)_{\Q}$.

On note $z'_1$ l'unique \'el\'ement de $N'$ v\'erifiant
\begin{equation}
z'_1-z_1=\sum_{l\in \delta(1)_K}\mu_l\,\rho_l\text{ avec }0\leq \mu_l<1.
\end{equation}

Soit $F(\delta,K,z)$
l'ensemble des \'el\'ements de $N$ qui s'\'ecrivent 
\begin{equation}
y_2+z_1-z'_1
\end{equation} 
o\`u $y_2$ est un \'el\'ement de $\ecC(\delta(1)\setminus \delta(1)_K)(K,z)_{z_2}$.

Alors
$\delta(1)_K$,
$\delta(K,z)_z$
et
$F(\delta,K,z)$
v\'erifient les propri\'et\'es suivantes.
\begin{enumerate}
\item
Si $\delta(K,z)_z$ est non vide, alors $z_3$ est un \'el\'ement de $N'$.
\item
Si $z_3$ est un \'el\'ement de $N'$, 
on a une d\'ecomposition 
\begin{equation}
\ecL_{N,\delta(K,z)_z}(T)
=
L_{N',\intrel(\ecC(\delta(1)_K))}(T)
\,
\ecL_{N,F(\delta,K,z)}(T).
\end{equation}
\item
$F(\delta,K,z)$ est fini. 
Plus pr\'ecis\'ement, on a la majoration
\begin{equation}\label{eq:maj1:bis}
\card{F(\delta,K,z)}\leq \left([N:N']\,\acc{\lambda^{\vee}}{z}\right)^{\,\rg(N)}.
\end{equation}
\item
Pour tout $y\in F(\delta,K,z)$ on a
\begin{equation}\label{eq:maj2:bis}
-\ecM\leq 
\acc{\lambda^{\vee}}{y}\leq \ecM \acc{\lambda^{\vee}}{z}.
\end{equation}
\item
Si $K$ est vide, on a $\delta(1)_K=\delta(1)$ et $F(\delta,K,z)=\{z_1-z'_1\}$. 
\item
On suppose en outre que $\left(N/\left\langle \Upsilon \right\rangle \right)^{\vee}\cap \Lambda^{\vee}=\{0\}$,
que $\delta$ est de dimension maximale, et que $K$ est non vide. Alors $\delta(1)_K$ est un sous-ensemble strict de $\delta(1)$. 
\end{enumerate}

\end{lemme}

\begin{demo}

Soit $y$ un \'el\'ement de $\intrel(\delta)$. 
Il s'\'ecrit donc de mani\`ere unique 
$y_1+y_2$ avec $y_1\in \intrel(\ecC(\delta(1)_K))$ 
et $y_2\in \intrel(\ecC(\delta(1)\setminus \delta(1)_K))$.

Supposons en outre que $y$ est dans $\delta(K,z)_z$, 
i.e. v\'erifie d'une part
\begin{equation}\label{eq:cond:dkz}
\forall\,i\in K,\quad\acc{y}{\lambda_i^{\vee}} < \acc{z}{\lambda_i^{\vee}}
\end{equation}
et 
d'autre part
\begin{equation}
y-z\in N'.
\end{equation}

Au vu de la d\'efinition \eqref{eq:def_delta1K:0} de $\delta(1)_K$,
on a 
\begin{equation}
\forall i\in K,\quad
\acc{y}{\lambda_i^{\vee}}=\acc{y_2}{\lambda_i^{\vee}}
.
\end{equation}
Ainsi la condition \eqref{eq:cond:dkz} montre que $y_2$ v\'erifie
\begin{equation}
\forall\,i\in K,\quad\acc{y_2}{\lambda_i^{\vee}} < \acc{z}{\lambda_i^{\vee}}.
\end{equation}

Par ailleurs on a $y-z=(y_1-z_1)+(y_2-z_2)+z_3$.
En outre $y-z$ est dans $N'$, et $y_1-z_1$ (respectivement $y_2-z_2$, respectivement $z_3$) 
est dans $(N'_1)_{\Q}$ (respectivement $(N'_2)_{\Q}$, respectivement $(N'_3)_{\Q}$).   
On en d\'eduit que $y_1-z_1$, $y_2-z_2$ et $z_3$ sont dans $N'$.

Ainsi on a 
\begin{equation}
y_2\in\ecC(\delta(1)\setminus \delta(1)_K)(K,z)_{z_2}
\end{equation}
et 
\begin{equation}
y_1\in \intrel[\ecC(\delta(1)_K)]_{z_1}.
\end{equation}
En outre ceci montre que la non vacuit\'e de $\delta(K,z)_z$ entra\^\i ne que $z_3$ est dans $N'$.

R\'eciproquement, on constate que  
si on a
\begin{equation}
y_1\in \intrel[\ecC(\delta(1)_K)]_{z_1},
\end{equation}
\begin{equation}
y_2\in \ecC(\delta(1)\setminus \delta(1)_K)(K,z)_{z_2},
\end{equation}
et 
\begin{equation}
z_3\in N'
\end{equation}
alors $y_1+y_2$ est un \'el\'ement de $\delta(K,z)_z$.

Ceci montre que si 
$z_3$ est dans $N'$ 
on a une d\'ecomposition
\begin{equation}
\ecL_{N,\delta(K,z)}(T)
=
\left(
\sum_{
\substack{
y_2\in \ecC(\delta(1)\setminus \delta(1)_K)(K,z)_{z_2}
}
}
T^{y_2}
\right)
\left(
\sum_{
\substack{
y_1\in \intrel(\ecC(\delta(1)_K))_{z_1}
}
}
T^{y_1}
\right).
\end{equation}

Par ailleurs on a 
\begin{align}
\sum_{
\substack{
y_1\in \intrel(\ecC(\delta(1)_K))_{z_1}
}
}
T^{y_1}
=
\sum_{
\substack{
y_1\in \intrel(\ecC(\delta(1)_K))
\\
y_1-z_1\in N'
}
}
T^{\,y_1}
&
=
T^{\,z_1-z'_1}
\,
\sum_{
\substack{
y\in \intrel(\ecC(\delta(1)_K))+z'_1-z_1
\\
y\in N'
}
}
T^{\,y}
\end{align}
Or, au vu de la d\'efinition de $z'_1$,  on a 
\begin{equation}
\left[\intrel(\ecC(\delta(1)_K))+z'_1-z_1\right]
\cap
N'
=
\intrel\left[\ecC(\delta(1)_K)\right]
\,
\cap
\,
N'
\end{equation}
soit
\begin{equation}
\sum_{
\substack{
y_1\in \intrel(\ecC(\delta(1)_K))
\\
y_1-z_1\in N'
}
}
T^{y_1}
=
T^{\,z_1-z'_1}
\,
L_{N',\intrel(\ecC(\delta(1)_K))}(T)
\end{equation}
d'o\`u 
\begin{align}
\ecL_{N,\delta(K,z)}(T)
&
=
\left(
\sum_{
y_2\in \ecC(\delta(1)\setminus \delta(1)_K)(K,z)_{z_2} 
}
T^{y_2+z_1-z'_1}
\right)
\,
L_{N',\intrel(\ecC(\delta(1)_K))}(T)
\\
&
=
\ecL_{N,F(\delta,K,z)}(T)\,L_{N',\intrel(\ecC(\delta(1)_K))}(T).
\end{align}

Montrons que $F(\delta,K,z)$ est fini et estimons son cardinal. 
On commence par remarquer que $F(\delta,K,z)$ est en bijection avec 
$\ecC(\delta(1)\setminus \delta(1)_K)(K,z_2)$.
Soit $y_2$ un \'el\'ement de $\ecC(\delta(1)\setminus \delta(1)_K)(K,z_2)$, 
qu'on \'ecrit
\begin{equation}
y_2=\sum_{l\in \delta(1)\setminus \delta(1)_K}\,\mu_l\,\rho_l
\end{equation}
avec les $\mu_l$ dans $\R_{>0}$. 

Comme $z\in N$, on a $[N:N']\,z\in N'$, donc $[N:N']\,z_2\in N'$
Comme $y_2-z_2\in N'$, $[N:N']\,y_2$ est dans $N'$.
Ainsi les $\mu_l$ sont dans $\frac{1}{[N:N']}\,\Z_{>0}$.

Par d\'efinition de $\delta(1)_K$, pour tout $l$ de $\delta(1)\setminus \delta(1)_K$, 
il existe $i$ dans $K$ v\'erifiant $\acc{\rho_l}{\lambda_i^{\vee}}\,\geq \,1$ 
(rappelons que pour tout $i$ on a $\acc{\rho_l}{\lambda_i^{\vee}} \geq 0$). 
Comme $y_2$ v\'erifie
\begin{equation}
\forall i\in K,\quad\acc{y_2}{\lambda_i^{\vee}} < \acc{z}{\lambda_i^{\vee}}
\end{equation}
on a
\begin{equation}
\mu_l < \underset{i\in K}{\Sup} \,\acc{z}{\lambda_i^{\vee}}\leq \acc{z}{\lambda^{\vee}}.
\end{equation}
Ainsi $\ecC(\delta(1)\setminus \delta(1)_K)(K,z_2)$ est fini et son cardinal est major\'e par 
\begin{equation}
[N:N']^{\,\card{\delta(1)\setminus \delta(1)_K}}
\,
\acc{z}{\lambda^{\vee}}^{\,\card{\delta(1)\setminus \delta(1)_K}},
\end{equation}
cette quantit\'e \'etant elle-m\^eme major\'ee par
\begin{equation} 
([N:N']\,\acc{z}{\lambda^{\vee}})^{\,\rg(N)}.
\end{equation}
Par ailleurs un \'el\'ement $y_2$ de $\ecC(\delta(1)\setminus \delta(1)_K)(K,z_2)$ v\'erifie 
\begin{equation}
0
\leq 
\acc{y_2}{\lambda^{\vee}} 
\leq 
\sum_{l\in \delta(1)\setminus\delta(1)_K}\,
\acc{z}{\lambda^{\vee}}\,
\acc{\rho_l}{\lambda^{\vee}}\,
\leq 
\ecM\,\acc{z}{\lambda^{\vee}}
\end{equation}
o\`u on rappelle que
\begin{equation}
\ecM=\card{I}\,\underset{l\in \Delta(1)}{\Sup}\left(\acc{\rho_l}{\lambda^{\vee}} \right).
\end{equation}
En outre on a 
\begin{equation}
0\geq \acc{z_1-z'_1}{\lambda^{\vee}}\geq -\sum_{l\in\delta(1)_K}  \acc{\rho_l}{\lambda^{\vee}}\geq -\ecM
\end{equation}
et finalement si $y$ est un \'el\'ement de $F(\delta,K,z)$ on a 
\begin{equation}
-\ecM
\leq 
\acc{y}{\lambda^{\vee}} 
\leq 
\ecM\,\acc{z}{\lambda^{\vee}}.
\end{equation}

Supposons 
$\left(N/\left\langle \Upsilon \right\rangle\right)^{\vee} \cap \Lambda^{\vee}=\{0\}$.
Ceci \'equivaut \`a dire que si $\lambda\in \Lambda^{\vee}$ v\'erifie $\acc{m}{\lambda}=0$ 
pour tout $m\in \left\langle \Upsilon \right\rangle$ alors $\lambda=0$.
Mais si $\delta$ est de dimension maximale, 
les $(\rho_l)_{l\in \delta(1)}$ engendrent $\left\langle \Upsilon \right\rangle$, 
et donc si $K$ n'est pas vide, on ne peut avoir $\delta(1)_K=\delta(1)$. 
\end{demo}

Nous pouvons g\'en\'eraliser le lemme \ref{lm:crucial:bis} de
la mani\`ere suivante. On conserve les notations introduites avant l'\'enonc\'e 
du lemme \ref{lm:crucial:bis}.
On consid\`ere en outre $\wt{I}\subset I$ un sous-ensemble strict de $I$
et $\wt{\Lambda}\subset \Lambda$ le c\^one engendr\'e par les $(\lambda_{i})_{i\in \wt{I}}$. 

Soit
\begin{equation}
\wt{\Upsilon_{z}}=\Upsilon\cap \{z+\wt{\Lambda}\cap N'\}.
\end{equation}
On veut \'etudier la s\'erie
\begin{equation}
\ecL_{N,\wt{\Upsilon_{z}}}(T)
=
\sum_{\substack{y\in \Upsilon\cap N\\ \\ y\in z+\wt{\Lambda}\cap N'}} T^{\,y}.
\end{equation}
Nous \'ecrivons d'abord
\begin{equation}
\ecL_{N,\wt{\Upsilon_{z}}}(T)
=
\sum_{\delta\in\Delta}
\,
\ecL_{N,\wt{\Upsilon_{z}}
\,
\cap
\, 
\intrel(\delta)}(T)
.
\end{equation}
Pour $\delta\in\Delta$, notons qu'on peut \'ecrire
\begin{align}
\ecL_{N,\wt{\Upsilon_{z}}\cap\, \intrel(\delta)}(\bs)
&
=
\sum_{
\substack
{
y\in \intrel(\delta)
\\ 
\\ 
y-z\in \wt{\Lambda}
\\
\\
y-z\in N'
\\
}
}
T^{\,y}
&
=
\sum_{
\substack
{
y\in \intrel(\delta)
\\ 
\\ 
\forall i\in \wt{I},\,\,\acc{\lambda_i^{\vee}}{y}\geq \acc{\lambda_i^{\vee}}{z} 
\\
\\
\forall i\in I\setminus \wt{I},\,\,\acc{\lambda_i^{\vee}}{y}=\acc{\lambda_i^{\vee}}{z} 
\\
\\
y-z\in N'
}
}
T^{\,y}
.
\end{align}
Pour tout sous-ensemble $K$ de $\wt{I}$ et tout c\^one $C$ de $N_{\R}$,  
soit $\wt{C}(K,z)$ le sous-ensemble de $\intrel(C)$ 
form\'e des \'el\'ements $y$ v\'erifiant les conditions
\begin{equation}
\forall i\in K,\, 
\acc{\lambda_i^{\vee}}{y}<\acc{\lambda_i^{\vee}}{z}
\end{equation}
et
\begin{equation}
\forall i\in I\setminus \wt{I},\, 
\acc{\lambda_i^{\vee}}{y}=\acc{\lambda_i^{\vee}}{z}.
\end{equation}

On a alors une d\'ecomposition
\begin{equation}
\ecL_{N,\wt{\Upsilon}_{z}\,\cap\,\intrel(\delta)}(T)
=
\sum_{K\subset \wt{I}}\,
(-1)^{\card{K}}\,
\ecL_{N,\wt{\delta}(K,z)_z}(T).
\end{equation}

Une l\'eg\`ere adaptation de la preuve du lemme 
\ref{lm:crucial:bis} permet alors de montrer le lemme suivant.

\begin{lemme}\label{lm:crucial:bis:tilde}
Soient $\delta$ un c\^one de $\Delta$, 
et $K$ une partie  de $\wt{I}$. 

Soit $\wt{\delta(1)_K}$ le sous-ensemble de $\delta(1)$
d\'efini par  
\begin{equation}
\wt{\delta(1)_K}=\{\,l\in \delta(1),\,
\forall\,i\in K\,\cup\,(I\setminus \wt{I}),\quad\acc{\rho_l}{\lambda_i^{\vee}} =0\}.
\end{equation}

On compl\`ete $(\rho_l)_{l\in \delta(1)}$ en une base $(\rho_l)_{l\in L}$ de $N'$.
Soit $N'_1$ (respectivement $N'_2$, respectivement $N'_3$)
le sous-module de $N$ engendr\'e par les
$(\rho_l)_{l\in \wt{\delta(1)_K}}$
(respectivement
$(\rho_l)_{l\in \delta(1)\setminus \wt{\delta(1)_K}}$,
respectivement
$(\rho_l)_{l\in L\setminus \delta(1)}$).

On \'ecrit $z=z_1+z_2+z_3$ avec  pour $i=1,2,3$, $z_i\in (N'_i)_{\Q}$.

On note $z'_1$ l'unique \'el\'ement de $N'$ v\'erifiant
\begin{equation}
z'_1-z_1=\sum_{l\in \wt{\delta(1)_K}}\mu_l\,\rho_l\text{ avec }0\leq \mu_l<1.
\end{equation}

Soit $\wt{F}(\delta,K,z)$
l'ensemble des \'el\'ements qui s'\'ecrivent $y_2+z_1-z'_1$ 
o\`u $y_2$ est un \'el\'ement de $\wt{\ecC(\delta(1)\setminus \wt{\delta(1)_K})}(K,z)_{z_2}$.

Alors les ensembles
$\wt{\delta(1)_K}$,
$\wt{\delta}(K,z)_z$
et
$\wt{F}(\delta,K,z)$
v\'erifient les propri\'et\'es suivantes :
\begin{enumerate}
\item
Si $\wt{\delta}(K,z)_z$ est non vide, alors $z_3$ est un \'el\'ement de $N'$.
\item
Si $z_3$ est un \'el\'ement de $N'$, 
on a une d\'ecomposition 
\begin{equation}
\ecL_{N,\wt{\delta}(K,z)_z}(T)
=
\,
\ecL_{N',\intrel\left(\ecC\left(\wt{\delta(1)_K}\right)\right)}(T)
\,
\ecL_{N,\wt{F}(\delta,K,z)}(T)
\end{equation}
\item
$\wt{F}(\delta,K,z)$ est fini. 
Plus pr\'ecis\'ement, on a la majoration
\begin{equation}\label{eq:maj1:bis:tilde}
\card{\wt{F}(\delta,K,z)}\leq \left([N:N']\,\acc{z}{\lambda^{\vee}}\right)^{\,\rg(N)}.
\end{equation}
\item
Pour tout $y\in \wt{F(\delta,K,z)}$ on a
\begin{equation}\label{eq:maj2:bis:tilde}
-\ecM\leq 
\acc{y}{\lambda^{\vee}}\leq \ecM \acc{z}{\lambda^{\vee}}.
\end{equation}
\item\label{item:lm:crucial:bis:tilde:4}
On suppose en outre que 
$\left(N/\left\langle \Upsilon \right\rangle \right)^{\vee}\cap \Lambda^{\vee}=\{0\}$
et que
$\delta$ est de dimension maximale.
Alors $\wt{\delta(1)_K}$ est un sous-ensemble strict de $\delta(1)$. 
\end{enumerate}

\end{lemme}

\subsection{Comportement des fonctions \'etudi\'ees par int\'egration}

Rappelons que nous avons fix\'e un $\Z$-module libre de rang fini $N$, 
et $\Lambda$ un c\^one simplicial de $N_{\R}$.
On consid\`ere en outre d\'esormais  un sous groupe $M$ de $N$ tel que le quotient $\Gamma=N/M$ soit sans torsion. 
On notera $j$ l'application quotient $N\to \Gamma$ et 
$i$ le morphisme d'inclusion $M\to N$.
On a donc une suite exacte de $\Z$-module libres de rang fini 
\begin{equation}\label{eq:exsqMNGamma}
0\longto
M
\overset{i}{\longto}
N
\overset{j}{\longto}
\Gamma
\longto
0
\end{equation}
et la suite exacte duale
\begin{equation}\label{eq:exsqMNGamma:dual}
0\longto
\Gamma^{\vee}
\overset{j^{\vee}}{\longto}
N^{\vee}
\overset{i^{\vee}}{\longto}
M^{\vee}
\longto
0.
\end{equation}

Nous aurons en particulier par la suite \`a consid\'erer la situation o\`u 
$
\Gamma^{\vee}\cap \Lambda^{\vee}=\{0\}.
$  
Cette hypoth\`ese \'equivaut au fait que $i^{\vee}(\Lambda^{\vee})$ est strictement convexe, 
ou encore que $\Lambda\cap M_{\R}=\left[i^{\vee}(\Lambda^{\vee})\right]^{\vee}$ est d'int\'erieur non vide,
ou encore que l'int\'erieur de $\Lambda$  rencontre $M$.

\begin{lemme}
Pour tout $\bs\in \tube{\Lambda^{\vee}_{>0}}$ on a
\begin{equation}
\ecL_{N,\Lambda\cap M_\R}\left(q^{\,-\bs}\right)
=
\ecL_{M,\Lambda\cap M_{\R}}
\left(q^{\,-i^{\vee}(\bs)}\right).
\end{equation}
\end{lemme}
\begin{demo}
En effet on peut \'ecrire 
\begin{align}
\ecL_{N,\Lambda\cap M_\R}\left(q^{\,-\bs}\right)
&=\sum_{y\in \Lambda \cap M_{\R} \cap N }q^{\,-\acc{y}{\bs}}
\\
&=\sum_{y\in \Lambda \cap M}q^{\,-\acc{y}{\bs}}
\\
&=\sum_{y\in \Lambda \cap M}q^{\,-\acc{i(y)}{\bs}}
\\
&=\sum_{y\in \Lambda \cap M}q^{\,-\acc{y}{i^{\vee}(\bs)}}
\\
&=
\ecL_{M,\Lambda\cap M_{\R}}
\left(q^{\,-i^{\vee}(\bs)}\right).
\end{align}
\end{demo}

\subsubsection{Le lemme technique : forme d\'epouill\'ee}

\begin{lemme}\label{lm:tech0}
On a
\begin{equation}
\int\limits_{\Gamma^{\vee}_{\unit}}
\ecL_{N,\Lambda}(j^{\vee}_{\unit}(\bz)\,q^{\,-\bs})\,d\bz 
=\ecL_{N,\Lambda\,\cap\,M_{\R}}\left(q^{\,-\bs}\right).
\end{equation}
\end{lemme}

\begin{demo} Il suffit d'\'ecrire
\begin{align}
\int\limits_{\Gamma^{\vee}_{\unit}}
\sum_{y\in \Lambda\cap N}\,\acc{j^{\vee}_{\unit}(\bz)}{y}\,q^{\,-\acc{y}{\bs}}d\bz&
=\sum_{y\in \Lambda\cap N}q^{\,-\acc{y}{\bs}}\,\int\limits_{\Gamma^{\vee}_{\unit}} \acc{j^{\vee}_{\unit}(\bz)}{y} d\bz\\
&=\sum_{y\in \Lambda\cap N}q^{\,-\acc{y}{\bs}}\,\int\limits_{\Gamma^{\vee}_{\unit}} \acc{\bz}{j(y)} d\bz\\
&=\sum_{\substack{y\in \Lambda\cap N \\ j(y)=0}}q^{\,-\acc{y}{\bs}}\\
&=\sum_{y\in \Lambda \cap M}q^{\,-\acc{y}{\bs}}\\
&=\ecL_{N,\Lambda\cap M_{\R}}\left(q^{\,-\bs}\right)
\end{align}
d'o\`u le lemme. 
\end{demo}

\subsubsection{Le lemme technique : forme simple}

On se donne $a\,:\,\Lambda\cap N\to \C$ et $\varepsilon>0$ tels que
la s\'erie d\'efinissant
$
\ecL_{N,\Lambda,a}\left(q^{\,-\bs}\right)
$
converge absolument pour tout $\bs\in \tube{\Lambda^{\vee}_{>-\varepsilon}}$.
Comme d\'ej\`a indiqu\'e,  
ceci \'equivaut \`a la condition suivante : 
\begin{equation}\label{eq:cond:conv:0}
\text{pour tout }\eta<\varepsilon,\text{ la s\'erie
 }\sum_{y\in \Lambda\cap N}\left|a_y\right|\,q^{\,\eta\,\acc{y}{\lambda}}
\text{
 est convergente.}
\end{equation}

Pour
$
\bs\in \tube{\Lambda^{\vee}_{>0}},
$
posons
\begin{equation}\label{eq:def:f1}
f_1(\bs)
=
\int\limits_{\Gamma^{\vee}_{\unit}}
\ecL_{N,\Lambda,a}\left(j^{\vee}_{\unit}(\bz)\,q^{\,-\bs}\right)
\,
\ecL_{N,\Lambda}\left(j^{\vee}_{\unit}(\bz)\,q^{\,-\bs}\right)
\,d\bz.
\end{equation}

Cela d\'efinit une fonction $f_1$ holomorphe sur 
$
\tube{\Lambda^{\vee}_{>0}}
$
.

\begin{rem}\label{rm:lm:tech1}
D'apr\`es le lemme \ref{lm:intn=intm:bis}, on peut, dans la d\'efinition \eqref{eq:def:f1}
de $f_1$, remplacer 
l'int\'egrale sur $\Gamma^{\vee}_{\unit}$ par une int\'egrale sur $\Gamma'_{\unit}$
o\`u $\Gamma'$ est n'importe quel sous-groupe d'indice fini de $\Gamma^{\vee}$.
\end{rem}

\begin{lemme}\label{lm:tech1}
On fait l'hypoth\`ese que $\Gamma^{\vee}\cap \Lambda^{\vee} =\{0\}.$ 
Alors la fonction
\begin{equation}
\bs\mapsto f_1(\bs)-\ecL_{N,\Lambda,a}(1)\,\ecL_{N,\Lambda\,\cap \,M_{\R}}\left(q^{\,-\bs}\right)
\end{equation}
est admissible de multiplicit\'e sup\'erieure \`a $1-\rg(M)$.
\end{lemme}

\begin{demo}
On a pour tout $\bs$ de $\tube{\Lambda^{\vee}_{>0}}$,
\begin{multline}
\ecL_{N,\Lambda,a}(j^{\vee}_{\unit}(\bz)\,q^{-\bs})
\,
\ecL_{N,\Lambda}(j^{\vee}_{\unit}(\bz)\,q^{\,-\bs})
\\
=
\left(
\sum_{y\in \Lambda\cap N}
a_y\,\acc{j^{\vee}_{\unit}(\bz)}{y}\,q^{\,-\acc{y}{\bs}}
\right)
\,
\left(
\sum_{y\in \Lambda\cap N} \acc{j^{\vee}_{\unit}(\bz)}{y} \,q^{\,-\acc{y}{\bs} }
\right)
\\
=
\sum_{(y_0,y_1)\in \left(\Lambda\cap N\right)^2}
a_{y_1}\,\acc{j^{\vee}_{\unit}(\bz)}{y_0+y_1}\,q^{\,-\acc{y_0+y_1}{\bs}}.
\end{multline}
d'o\`u
\begin{align}
\int\limits_{\Gamma^{\vee}_{\unit}}
\ecL_{N,\Lambda,a}
\left(j^{\vee}_{\unit}(\bz)\,q^{\,-\bs}\right)
\,
\ecL_{N,\Lambda}(j^{\vee}_{\unit}(\bz)\,q^{\,-\bs})
d\bz
&=
\sum_{(y_0,y_1)\in \left(\Lambda\cap N\right)^2}a_{y_1}\,q^{\,-\acc{y_0+y_1}{\bs}} 
\int\limits_{\Gamma^{\vee}_{\unit}} \acc{\bz}{j(y_0+y_1)}d\bz\\
&=\sum_{y_1\in \Lambda\cap N}a_{y_1}\,\sum_{\substack{y_0\in \Lambda\cap N  \\ j(y_0+y_1)=0}} q^{\,-\acc{y_0+y_1}{\bs} } \\
&=\sum_{y_1\in \Lambda\cap N}a_{y_1}\sum_{\substack{y\in \Lambda\cap M \\ y\in y_1+\Lambda\cap N}} q^{\,-\acc{y}{\bs}}.
\end{align}
On applique alors le lemme \ref{lm:crucial} avec $\Upsilon=\Lambda\cap M_{\R}$.  
On obtient la d\'ecomposition
\begin{equation}
f_1(\bs)
=
\sum_{\delta \in \Delta} 
\,
\sum_{K\subset I}(-1)^{\card{K}}
\,
\ecL_{N,\intrel(\ecC(\delta(1)_K))}\left(q^{\,-\bs}\right) 
\,
\sum_{y_1\in \Lambda\cap N} a_{y_1} 
\,
\ecL_{N,F(\delta,K,y_1)}\left(q^{\,-\bs}\right).
\end{equation}
Soient $\delta$ et $K$ donn\'es.
D'apr\`es les majorations \eqref{eq:maj1} et \eqref{eq:maj2} on a pour tout $\eta$
\begin{align}
\sum_{y_1\in \Lambda\cap N} |a_{y_1}| \,\sum_{y\in F(\delta,K,y_1)} q^{\eta\acc{y}{\lambda}}
&\leq \sum_{y_1\in \Lambda\cap N} |a_{y_1}| \acc{y_1}{\lambda}^{\,\rg(N)}\,q^{\eta\,M\,\acc{y_1}{\lambda}}.
\end{align}
et cette derni\`ere s\'erie converge pour tout $\eta<\frac{\varepsilon}{M}$. 
Ainsi la s\'erie
\begin{equation}
\sum_{y\in \Lambda\cap N}
\left(\sum_{\substack{y_1\in \Lambda\cap N\\ y\in F(\delta,K,y_1)}} a_{y_1}\right) q^{\,-\acc{y}{\bs}}
\end{equation}
d\'efinit une fonction admissible de multiplicit\'e sup\'erieure \`a $0$.
Comme $\Gamma^{\vee}\cap \Lambda^{\vee}=\{0\}$, d'apr\`es 
le point \ref{item:lm:crucial:4}
du lemme \ref{lm:crucial} 
la fonction 
\begin{equation}
\ecL_{N,\intrel(\ecC(\delta(1)_K))}\left(q^{\,-\bs}\right) \,\sum_{y_1\in \Lambda\cap N} a_{y_1} \,
\ecL_{N,F(\delta,K,y_1)}\left(q^{\,-\bs}\right).
\end{equation}
est, si $K\neq\varnothing$, 
une fonction admissible de multiplicit\'e sup\'erieure \`a $1-\rg(M)$.

La contribution des termes correspondant \`a $K=\varnothing$ dans la d\'ecomposition ci-dessus est
\begin{equation}
\left(\sum_{y_1\in \Lambda\cap N} a_{y_1}\right)
\times
\left(\sum_{\delta \in \Delta} \ecL_{N,\intrel(\delta)}\left(q^{\,-\bs}\right)\right)
=
\ecL_{N,\Lambda,a}(1)
\times 
\ecL_{M,\Lambda\,\cap\,M_{\R}}\left(q^{\,-\bs}\right)
\end{equation}
d'o\`u le lemme.
\end{demo}

\subsubsection{Le lemme technique : forme g\'en\'erale}\label{subsubsec:lm:techgene}

Comme dans la sous-section pr\'ec\'edente, on commence par se donner $a\,:\,\Lambda\cap N\to \C$ et $\eps>0$ tels que
la s\'erie d\'efinissant
$
\ecL_{N,\Lambda,a}\left(q^{\,-\bs}\right)
$
converge absolument pour tout $\bs\in \tube{\Lambda^{\vee}_{>-\varepsilon}}$.

On se donne en outre $N'$ un sous-groupe d'indice fini de $N$.

Pour
$
\bs\in \tube{\Lambda^{\vee}_{>0}},
$
posons
\begin{equation}\label{eq:def:f2}
f_2(\bs)
=
\int\limits_{\bz\in \Gamma^{\vee}_{\unit}}
\ecL_{N,\Lambda,a}
\left(j^{\vee}_{\unit}(\bz)\,q^{\,-\bs}\right)
\,
\ecL_{N',\Lambda}(j^{\vee}_{\unit}(\bz)\,q^{\,-\bs})
\,d\bz.
\end{equation}

Cela d\'efinit une fonction $f_2$ holomorphe sur 
$
\tube{\Lambda^{\vee}_{>0}}
$
.

\begin{rem}\label{rm:lm:tech1:bis}
D'apr\`es le lemme \ref{lm:intn=intm:bis}, on peut, dans la d\'efinition \eqref{eq:def:f2}
de $f_2$, remplacer 
l'int\'egrale sur $\Gamma^{\vee}_{\unit}$ par une int\'egrale sur $\Gamma'_{\unit}$
o\`u $\Gamma'$ est n'importe quel sous-groupe d'indice fini de $\Gamma^{\vee}$.
\end{rem}
\begin{lemme}\label{lm:tech1:bis}
On suppose que $\Gamma^{\vee}\cap \Lambda^{\vee} =\{0\}.$ 
\begin{enumerate}
\item
Il existe alors une fonction $g$ admissible de multiplicit\'e sup\'erieure \`a z\'ero
\label{item:1:lm:tech1:bis}
telle que la fonction
\begin{equation}
\bs\mapsto f_2(\bs)-g(\bs)\,\ecL_{N',\Lambda\,\cap \,M_{\R}}\left(q^{\,-\bs}\right)
\end{equation}
est admissible de multiplicit\'e sup\'erieure \`a $1-\rg(M)$.
\item
On suppose en outre que la condition suivante est v\'erifi\'ee : 
\label{item:2:lm:tech1:bis}
pour tout $y\in \Lambda\cap N$, si $a_y$
est non nul alors $j(y)$ est un \'el\'ement de $j(N')$.
Alors on peut choisir $g$ de sorte qu'on ait 
\begin{equation}
 g(0)=\ecL_{N,\Lambda,a}(1).
\end{equation}
\end{enumerate}
\end{lemme}

\begin{demo}
On a pour tout $\bs$ de $\tube{\Lambda^{\vee}_{>0}}$,
et tout $\bz$ de $\Gamma^{\vee}_{\unit}$ 
\begin{multline}
\ecL_{N,\Lambda,a}\left(j^{\vee}_{\unit}(\bz)\,q^{-\bs}\right)
\,
\ecL_{N',\Lambda}\left(j^{\vee}_{\unit}(\bz)\,q^{\,-\bs}\right)
\\
=
\left(
\sum_{y\in \Lambda\cap N}
a_y\,\acc{\bz}{j(y)}\,q^{\,-\acc{y}{\bs}}
\right)
\,
\left(
\sum_{y\in \Lambda\cap N'} \acc{\bz}{j(y)}\,q^{\,-\acc{y}{\bs} }
\right)
\end{multline}
d'o\`u
\begin{align}
&
\phantom{=}
\int\limits_{\bz\in \Gamma^{\vee}_{\unit}}
\ecL_{N,\Lambda,a}
\left(j^{\vee}_{\unit}(\bz)\,q^{\,-\bs}\right)
\,
\ecL_{N',\Lambda}
\left(j^{\vee}_{\unit}(\bz)\,q^{\,-\bs}\right)
d\bz
\\
&
=
\sum_{(y,z)\in 
\left(\Lambda\cap N'\right)\times \left(\Lambda\cap N\right)}
a_{z}\,q^{\,-\acc{y+z}{\bs}} 
\int\limits_{\Gamma^{\vee}_{\unit} } \acc{\bz}{\,j(y+z)}d\bz\\
&=\sum_{z\in \Lambda\cap N}
a_{z}
\,
\sum_{\substack{y\in \Lambda\cap N'  \\ j(y+z)=0}} 
q^{\,-\acc{y+z}{\bs}} 
\\
&=\sum_{z\in \Lambda\cap N}a_{z}
\,
\sum_{\substack{y\in \Lambda\cap N'  \\ y+z\in M}} 
q^{\,-\acc{y+z}{\bs}} 
\\
&
=
\sum_{z\in \Lambda\cap N}
a_{z}
\sum_{\substack{y\in \Lambda\cap M \\ y\in z+\Lambda\cap N'}} 
q^{\,-\acc{y}{\bs}}.
\end{align}

On applique alors le lemme \ref{lm:crucial:bis} avec $\Upsilon=\Lambda\cap M_{\R}$.  
On reprend \`a cet effet les notations du lemme \ref{lm:crucial:bis}.
On obtient la d\'ecomposition
\begin{equation}\label{eq:decompo}
f_2(\bs)
=
\sum_{\delta \in \Delta} 
\,
\sum_{K\subset I}(-1)^{\card{K}}
\,
\ecL_{N',\intrel(\ecC(\delta(1)_K))}\left(q^{\,-\bs}\right)
\,
\sum_{
\substack{
z\in \Lambda\cap N
\\
z_3\in N'
}
} 
a_{z} 
\,
\ecL_{N,F(\delta,K,z)}\left(q^{\,-\bs}\right).
\end{equation}
Soient $\delta$ et $K$ donn\'es.
D'apr\`es les majorations \eqref{eq:maj1:bis} et \eqref{eq:maj2:bis} 
on a pour tout $\eta$ la majoration
\begin{equation}
\sum_{
\substack{
z\in \Lambda\cap N
\\
z_3\in N'
}
} \abs{a_{z}} \,\sum_{y\in F(\delta,K,z)} q^{\eta\acc{y}{\lambda}}
 \leq 
\sum_{z\in \Lambda\cap N} 
|a_{z}| \left([N:N']\,\acc{z}{\lambda}\right)^{\,\rg(N)}\,q^{\eta\,\ecM\,\acc{z}{\lambda}}.
\end{equation}
D'apr\`es la condition \eqref{eq:cond:conv:0},
cette derni\`ere s\'erie converge pour tout $\eta<\frac{\varepsilon}{\ecM}$. 
Ainsi la s\'erie 
\begin{equation}
\sum_{
\substack{
z\in \Lambda\cap N
\\
z_3\in N'
}
} 
a_{z} 
\,
\ecL_{N,F(\delta,K,z)}\left(q^{\,-\bs}\right)
\end{equation}
d\'efinit une fonction admissible de multiplicit\'e sup\'erieure \`a $0$.

Comme $\Gamma^{\vee}\cap \Lambda^{\vee}=\{0\}$, 
l'assertion 6
du lemme \ref{lm:crucial:bis} 
montre que si $K\neq\varnothing$ ou $\dim(\delta)<\rg(M)$,  
la fonction 
\begin{equation}
\ecL_{N',\intrel(\ecC(\delta(1)_K))}\left(q^{\,-\bs}\right) \,
\sum_{z\in \Lambda\cap N} a_{z} \,L_{N,F(\delta,K,z)}\left(q^{\,-\bs}\right).
\end{equation}
est
une fonction admissible de multiplicit\'e sup\'erieure \`a $1-\rg(M)$.

\'Ecrivons la contribution des termes correspondant \`a $K=\varnothing$ et $\dim(\delta)=\rg(M)$ 
dans la d\'ecomposition \eqref{eq:decompo}.
Pour cela, on commence par remarquer la chose suivante : soit $\delta$ tel que $\dim(\delta)=\rg(M)$
et $K=\varnothing$. Avec les notations du lemme \ref{lm:crucial:bis}, on a
$N'_1=M \cap N'$, $N'_2=0$ et $N'=(N'\cap M)\oplus N'_3$.
Ainsi, si $z\in N$ s'\'ecrit $z=z_1+z_3$ avec $z_1\in  (N'\cap M)_{\Q}$
et $z_3\in (N'_3)_{\Q}$, la condition $z_3\in N'$ est \'equivalente \`a la condition
$j(z)\in j(N')$.
Ainsi la contribution consid\'er\'ee s'\'ecrit
\begin{equation}
\left(
\sum_{
\substack{
z\in \Lambda\cap N, 
\\
j(z)\in j(N')
}
} a_{z}\,q^{-\acc{z_1-z'_1}{\bs}}
\right)
\times
\left(
\sum_{
\substack{
\delta \in \Delta
\\
\dim(\delta)=\rg(M)
} 
}
\ecL_{N',\intrel(\delta)}\left(q^{\,-\bs}\right).
\right)
\end{equation}

Pour tout $z\in \Lambda\cap N$, on a $-\ecM\leq \acc{z_1-z'_1}{\lambda}\leq 0$.
Ainsi, pour $\eta\geq 0$ on a  
\begin{equation}
\sum_{
\substack{
z\in \Lambda\cap N, 
\\
j(z)\in j(N')
}
} 
\abs{a_{z}\,q^{\eta\acc{z_1-z'_1}{\lambda}}}
\leq 
\sum_{
z\in \Lambda\cap N, 
} 
\abs{a_{z}}<+\infty
\end{equation}
et pour $\eta\leq 0$ on a 
\begin{equation}
\sum_{
\substack{
z\in \Lambda\cap N, 
\\
j(z)\in j(N')
}
} 
\abs{a_{z}\,q^{\eta\acc{z_1-z'_1}{\lambda}}}
\leq 
q^{-\eta\,\ecM}
\sum_{
z\in \Lambda\cap N, 
}
\abs{a_{z}}
<+\infty.
\end{equation}

Ceci montre que la fonction $g$ d\'efinie par  
\begin{equation}
g(\bs)
=
\sum_{
\substack{
z\in \Lambda\cap N, 
\\
j(z)\in j(N')
}
} a_{z}\,q^{\acc{z_1-z'_1}{\bs}}
\end{equation}
est admissible de multiplicit\'e positive. 

Comme on a la d\'ecomposition
\begin{equation}
\ecL_{N',\Lambda\,\cap\,M_{\R}}\left(q^{\,-\bs}\right)
=\sum_{\delta \in \Delta} 
\,
\sum_{K\subset I}(-1)^{\card{K}}
\,
\ecL_{N',\intrel(\ecC(\delta(1)_K))}\left(q^{\,-\bs}\right)
\end{equation}
et que, pour $K\neq \vide$ ou $\dim(\delta)<\rg(M)$,
la fonction $\bs\mapsto \ecL_{N',\intrel(\ecC(\delta(1)_K))}\left(q^{\,-\bs}\right)$ est admissible
de multiplicit\'e sup\'erieure \`a $1-\rg(M)$,
on d\'eduit de ce qui pr\'ec\`ede que
\begin{equation}
\bs\mapsto f_2(\bs)-g(\bs)\,\ecL_{N',\Lambda\,\cap\,M_{\R}}\left(q^{\,-\bs}\right)
\end{equation}
est admissible de multiplicit\'e sup\'erieure \`a $1-\rg(M)$.

Supposons \`a pr\'esent que pour $z\in \Lambda\cap N$ la condition
$a_z\neq 0$ entra\^\i ne $j(z)\in j(N')$. 
On a alors l'\'egalit\'e
\begin{equation}
\sum_{
\substack{
z\in \Lambda\cap N, 
\\
j(z)\in j(N')
}
} a_{z}
=
\sum_{
z\in \Lambda\cap N
} a_{z}
\end{equation}
ce qui signifie exactement que $g(0)=\ecL_{N,\Lambda,a}(1)$.
\end{demo}
\begin{rem}
S'il existe $y$ tel que $a_y\neq 0$ et  $j(y)\notin j(N')$, la conclusion du point \ref{item:2:lm:tech1:bis}
du lemme \ref{lm:tech1:bis} peut \^etre mise en d\'efaut : il se peut que $f_2(\bs)$ soit admissible de 
multiplicit\'e sup\'erieure \`a $1-\rg(M)$, m\^eme si $\ecL_{N,\Lambda,a}(1)$ est non nul.
En d'autres termes, dans ce cas de figure, on <<perd>> le contr\^ole du terme principal (ce qu'il 
faudra bien s\^ur \'eviter absolument lors de l'application au calcul de la fonction z\^eta des hauteurs).

\`A titre d'exemple, consid\'erons le cas o\`u $N=\Z^2$, $\Lambda=\R_{\geq 0}^2$, $M=\Z\,(1,1)$, $N'=2\,N$,
et $a\,:\,\N^2\to \C$ est donn\'ee par $a_{(2,1)}=1$ et $a_{y}=0$ si $y\neq (2,1)$. 
On a alors $j(1,0)=-j(0,1)$ et on identifie $\Gamma$ \`a $\Z$ au moyen de la base $j(1,0)$.

Pour $\bz\in \Gamma^{\vee}_{\unit}\isom \unit$, on a alors
\begin{align}
\ecL_{N,\Lambda,a}
\left(j^{\vee}_{\unit}(\bz)\,q^{\,-\bs}\right)
\,
\ecL_{N',\Lambda}(j^{\vee}_{\unit}(\bz)\,q^{\,-\bs})
&
=
\frac{\bz^{j(2,1)}\,q^{-2\,s_1-s_2}}
{(1-\bz^{j(2,0)}\,q^{-2\,s_1})
(1-\bz^{j(0,2)}\,q^{-2\,s_2})
}
\\
&
=
\frac{\bz\,q^{-2\,s_1-s_2}}
{(1-\bz^2\,q^{-2\,s_1})
(1-\bz^{-2}\,q^{-2\,s_2})
}
\\
&
=
\sum_{(n_1,n_2)\in \N^2}
\bz^{2\,n_1-2\,n_2+1}\,
q^{-2\,(n_1+1)\,s_1-(2\,n_2+1)\,s_2}
\end{align}
et de cette derni\`ere expression, on d\'eduit aussit\^ot qu'on a
\begin{equation}
\int\limits_{\bz\in \Gamma^{\vee}_{\unit}}
\ecL_{N,\Lambda,a}
\left(j^{\vee}_{\unit}(\bz)\,q^{\,-\bs}\right)
\,
\ecL_{N',\Lambda}(j^{\vee}_{\unit}(\bz)\,q^{\,-\bs})
=
0.
\end{equation}

On pourra constater par contre que si on d\'efinit $a$ par $a_{(3,1)}=1$
et $a_{y}=0$ si $y\neq (3,1)$, les hypoth\`eses du point \ref{item:2:lm:tech1:bis} du lemme \ref{lm:tech1:bis} 
sont v\'erifi\'ees. De fait, un calcul similaire montre qu'on a alors
\begin{align}
\int\limits_{\bz\in \Gamma^{\vee}_{\unit}}
\ecL_{N,\Lambda,a}
\left(j^{\vee}_{\unit}(\bz)\,q^{\,-\bs}\right)
\,
\ecL_{N',\Lambda}(j^{\vee}_{\unit}(\bz)\,q^{\,-\bs})
&
=
\frac{q^{-3\,s_1-3\,s_2}}{1-q^{-2\,s_1-2\,s_2}}
\\
&
=
q^{-3\,s_1-3\,s_2}\, \ecL_{N',\Lambda\,\cap \,M_{\R}}\left(q^{\,-(s_1,s_2)}\right)
.
\end{align}
\end{rem}

Le lemme \ref{lm:tech1:bis} va en fait nous servir \`a traiter le terme principal
de la fonction z\^eta des hauteurs, qui correspond (\`a une constante multiplicative pr\`es)
au terme $I_1$ dans la formule \eqref{eq:repr:int:fonc}. Nous avons besoin d'un autre lemme pour
traiter les termes restants, \`a savoir les $I_{\chi}$, 
pour $\chi\in  \ecU_T  \setminus \Ker\left(\gamma^{\ast}\right)$,
apparasissant dans cette m\^eme formule.
On conserve les notations introduites avant l'\'enonc\'e du lemme \ref{lm:tech1:bis}.
On consid\`ere en outre $\wt{I}\subset I$ un sous-ensemble {\em strict} de $I$,
et $\wt{\Lambda}\subset \Lambda$ le c\^one engendr\'e par les $(\lambda_{i})_{i\in \wt{I}}$. 

Si on suppose que $\Gamma$ v\'erifie $\Gamma^{\vee} \cap \Lambda^{\vee} =\{0\}$,
$\wt{\Lambda}\cap M_{\R}$ est d'int\'erieur vide dans $M_{\R}$. 
En effet ce c\^one est le dual du c\^one $i(\wt{\Lambda})$, qui contient
la droite engendr\'ee par l'\'el\'ement non nul $i(\lambda_{i_0})$, o\`u $i_0\notin I$.

Pour
$
\bs\in \tube{\Lambda^{\vee}_{>0}},
$
posons
\begin{equation}
f_3(\bs)
=
\int\limits_{\Gamma^{\vee}_{\unit}}
\ecL_{N,\Lambda,a}\left(j^{\vee}_{\unit}(\bz)\,\,q^{\,-\bs}\right)
\,
\ecL_{N',\wt{\Lambda}}(j^{\vee}_{\unit}(\bz)\,\,q^{\,-\bs})
d\bz.
\end{equation}
Cela d\'efinit une fonction $f_3$ holomorphe sur 
$
\tube{\Lambda^{\vee}_{>0}}
$
.

\begin{lemme}\label{lm:tech1:bis:tilde}
On fait l'hypoth\`ese que $\Gamma^{\vee}\cap \Lambda^{\vee} =\{0\}.$ 
Alors la fonction $f_3$ est admissible de multiplicit\'e sup\'erieure \`a $1-\rg(M)$.
\end{lemme}

\begin{demo}
On a pour tout $\bs$ de $\tube{\Lambda^{\vee}_{>0}}$,
et tout $\bz$ 
de $\Gamma^{\vee}_{\unit}$ 
\begin{multline}
\ecL_{N,\Lambda,a}
\left(j^{\vee}_{\unit}(\bz)\,q^{\,-\bs}\right)
\,
\ecL_{N',\wt{\Lambda}}
\left(j^{\vee}_{\unit}(\bz)\,q^{\,-\bs}\right)
\\
=
\left(
\sum_{y\in \Lambda\cap N}
a_y\,\acc{\bz}{j(y)}\,q^{\,-\acc{y}{\bs}}
\right)
\,
\left(
\sum_{y\in \wt{\Lambda}\cap N'} \acc{\bz}{j(y)}\,q^{\,-\acc{y}{\bs} }
\right)
\end{multline}
d'o\`u
\begin{align}
&
\phantom{=}
\int\limits_{\Gamma^{\vee}_{\unit} }
\ecL_{N,\Lambda,a}
\left(j^{\vee}_{\unit}(\bz)\,q^{\,-\bs}\right)
\,
\ecL_{N',\wt{\Lambda}}
\left(j^{\vee}_{\unit}(\bz)\,q^{\,-\bs}\right)
d\bz
\\
&
=
\sum_{(y,z)\in 
\left(\wt{\Lambda}\cap N'\right)\times \left(\Lambda\cap N\right)}
a_{z}\,q^{\,-\acc{y+z}{\bs}} 
\int\limits_{\Gamma^{\vee}_{\unit}} \bz^{\,j(y+z)}d\bz\\
&=\sum_{z\in \Lambda\cap N}
a_{z}
\,
\sum_{\substack{y\in \wt{\Lambda}\cap N'  \\ j(y+z)=0}} 
q^{\,-\acc{y+z}{\bs}} 
\\
&=\sum_{z\in \Lambda\cap N}a_{z}
\,
\sum_{\substack{y\in \wt{\Lambda}\cap N'  \\ y+z\in M}} 
q^{\,-\acc{y+z}{\bs}} 
\\
&
=
\sum_{z\in \Lambda\cap N}
a_{z}
\sum_{\substack{y\in \Lambda\cap M \\ y\in z+\wt{\Lambda}\cap N'}} 
q^{\,-\acc{y}{\bs}}.
\end{align}

On applique alors le lemme \ref{lm:crucial:bis:tilde} avec $\Upsilon=\Lambda\cap M_{\R}$.  
On reprend \`a cet effet les notations du lemme \ref{lm:crucial:bis:tilde}.
On obtient la d\'ecomposition
\begin{equation}\label{eq:decompo:tilde}
f_3(\bs)
=
\sum_{\delta \in \Delta} 
\,
\sum_{K\subset \wt{I}}(-1)^{\card{K}}
\,
\ecL_{N',\intrel(\ecC(\wt{\delta(1)_K}))}\left(q^{\,-\bs}\right) 
\,
\sum_{
\substack{
z\in \Lambda\cap N
\\
z_3\in N'
}
} 
a_{z} 
\,
\ecL_{N,\wt{F}(\delta,K,z)}\left(q^{\,-\bs}\right).
\end{equation}
Soient $\delta$ et $K$ donn\'es.
D'apr\`es les majorations \eqref{eq:maj1:bis:tilde} et \eqref{eq:maj2:bis:tilde} 
on a pour tout $\eta$ la majoration
\begin{align*}
\sum_{
\substack{
z\in \Lambda\cap N
\\
z_3\in N'
}
} |a_{z}| \,\sum_{y\in \wt{F}(\delta,K,z)} q^{\eta\acc{y}{\lambda}}
&\leq 
\sum_{z\in \Lambda\cap N} 
|a_{z}| \left(d\,\acc{z}{\lambda}\right)^{\,\rg(N)}\,q^{\eta\,\ecM\,\acc{z}{\lambda}}.
\end{align*}
D'apr\`es la condition \eqref{eq:cond:conv:0},
cette derni\`ere s\'erie converge pour tout $\eta<\frac{\varepsilon}{\ecM}$. 
Ainsi la s\'erie
\begin{equation}
\sum_{
\substack{
z\in \Lambda\cap N
\\
z_3\in N'
}
} 
a_{z} 
\,
L_{N,F(\delta,K,z)}\left(q^{\,-\bs}\right)
\end{equation}
d\'efinit une fonction admissible de multiplicit\'e sup\'erieure \`a $0$.

Comme $\Gamma^{\vee}\cap \Lambda^{\vee}=\{0\}$, 
le point \ref{item:lm:crucial:bis:tilde:4} 
du lemme \ref{lm:crucial:bis:tilde} 
montre que 
la fonction 
\begin{equation}
\bs\mapsto \ecL_{N',\intrel(\ecC(\wt{\delta(1)_K}))}\left(q^{\,-\bs}\right) \,
\sum_{z\in \Lambda\cap N} a_{z} \,\ecL_{N,\wt{F}(\delta,K,z)}\left(q^{\,-\bs}\right)
\end{equation}
est
une fonction admissible de multiplicit\'e sup\'erieure \`a $1-\rg(M)$.
\end{demo}

\section{\hspace{-0.1cm}Application aux fonctions z\^eta des hauteurs}
\label{app_fonction_zeta}

\subsection{Pr\'eliminaires}

Rappelons que $X(T)^{\,G}$ s'identifie \`a un sous-groupe 
de $\Z^{\,\Sigma(1)/G}$ via la suite exacte
\begin{equation}
\label{eq:suite:ex}
0\longto X(T)^{\,G} \overset{\gamma}\longto  \Z^{\Sigma(1)/G}\overset{\pi}\longto \Pic(\xs)\longto H^1(G,X(T))\longto 0
\end{equation}
tir\'ee de la suite exacte \eqref{eq:resflT} en prenant les $G$-invariants.

Soit $m\in X(T)^{\,G}$ v\'erifiant
\begin{equation}
\forall\, \alpha \in \Sigma(1)/G,\quad\acc{m}{\roa} \geq 0.
\end{equation}
Alors on a
\begin{equation}
\forall\,  \alpha \in \Sigma(1),\quad\forall\,  \rho\in\alpha,\quad\acc{m}{\rho}\geq 0.
\end{equation}
L'\'eventail $\Sigma$ \'etant complet, 
ses rayons engendrent $X(T)^{\,\vee}$, 
ce qui entra\^\i ne $m=0$. 
Nous avons donc 
\begin{equation}\label{eq:xtg}
X(T)^{\,G}\cap \Z^{\Sigma(1)/G}_{\geq 0}=\{0\}
\end{equation}
(on aurait pu aussi invoquer directement le fait g\'en\'eral que $\ceff(\xs)$ est strictement convexe).

Dans le cas fonctionnel, 
on tire  par ailleurs de \eqref{eq:resflT} un complexe
\begin{equation}
\label{eq:complexe_dt}
\DTNS \longto  \DTqd \longto \DT.
\end{equation}
Notons $\DTNS^0$ le noyau du morphisme $\DTqd\to \DT$.
Rappelons que nous avons not\'e  $\DT^0$ l'image de ce morphisme 
dans $\DT$. 
On a donc une suite exacte de $\Z$-modules libre de rang fini.
\begin{equation}
0
\longto 
\DTNS^0
\overset{\pi^{\vee}}
\longto  
\DTqd
\overset{\gamma^{\vee}}
\longto
\DT^0
\longto 
0
\end{equation}
et la suite exacte duale
\begin{equation}
0
\longto 
\left(\DT^0\right)^{\vee}
\overset{\gamma}
\longto  
\left(\DTqd\right)^{\vee}
\overset{\pi}
\longto
\left(\DTNS^0\right)^{\vee}
\longto 
0.
\end{equation}

Comme \eqref{eq:resflT} est une r\'esolution flasque de $X(T)$
(lemme \ref{lm:picxsflasque}),
d'apr\`es la proposition \ref{prop:relationkt},
l'indice dans $\DTNS^0$
de 
l'image 
de $\DTNS$ dans $\DTqd$
est \'egal \`a 
\begin{equation}\label{eq:indice}
\frac
{\card{H^1(G,X(T))}\,\left[\CTNS\right]\,\card{\CT}\,\KT}
{\left[\CTqd\right]}.
\end{equation}

On pose 
\begin{equation}\label{eq:def:c0}
C_0
=
\lim_{s\to 1}\, (s-1)^{\,\card{\Sigma(1)/G}}\,
\int\limits_{\adh{T(K)}\cap T(\ak)}H(-s\,\varphi_0,t)\omega_{T}(t).
\end{equation}
On montrera \`a la sous-section \ref{subsubsec:calcul:c0} que $C_0$
est non nulle.

\subsection{Application du lemme technique dans le cas arithm\'etique}
\label{subsec:application:lemme:technique:arit}

Rappelons que d'apr\`es 
le corollaire \ref{cor:rep:int:arit},
on a
pour tout  \'el\'ement $\bs$ de $\tube{\R^{\,\Sigma(1)/G}_{>0}}$ 
la repr\'esentation int\'egrale 
\begin{align}
\zeta_{H}(\varphi_0+\bs)
&=
\sum_{t\in T(K)}\,H\left(t,-\bs -\varphi_0\right)
\\
&=
\frac
{\card{A(T)}}
{\,(2\,\pi)^{\rg\left(X(T)^G\right)}\,b(T)}
\int\limits_{y\in X(T)^G_{\R}}
\!\!\!\!F(\bs-i\,\gamma_{\R}(y))dy.
\end{align}
o\`u $F$ est la fonction d\'efinie dans l'\'enonc\'e de 
la proposition \ref{prop:contr:F}.
En particulier $F$ est une fonction holomorphe sur $\tube{\R^{\Sigma(1)/G}_{>0}}$,
et
il existe un $\eps>0$ tel que la fonction 
\begin{equation}
\bs\mapsto F(\bs)\,\prod_{\alpha\in \Sigma(1)/G} \frac{\sa}{1+\sa} 
\end{equation}
se prolonge en une fonction 
holomorphe sur $\tube{\R^{\,\Sigma(1)/G}_{>-\eps}}$
qui d'apr\`es le point 
\ref{item3:prop:contr:F}
de la proposition \ref{prop:contr:F}
est $X(T)^G_{\R}$-contr\^ol\'ee 
au sens de \cite[D\'efinition 3.13]{CLTs:fibres}.

D'apr\`es \eqref{eq:item2:prop:contr:F}, on a
\begin{equation}\label{eq:def:c0:bis}
C_0
=
\lim_{s\to 0}
s^{\card{\Sigma(1)/G}}
F(s\,\varphi_0)
.
\end{equation}

On d\'eduit alors du th\'eor\`eme \ref{thm:CLTs} qu'il existe un $\eps'>0$
tel que la fonction 
\begin{equation}\label{eq:zetah1+sa}
R(\bs)
=
\zeta_H(\varphi_0+\bs)
-
\frac
{\card{A(T)}\,C_0}
{b(T)}
\indic_{\pi\left(\Z^{\Sigma(1)/G}\right),\pi\left(\R^{\Sigma(1)/G}_{\geq 0}\right)}
(\pi(\bs)),
\end{equation}
se prolonge en une fonction m\'eromorphe sur
$\tube{\R_{>-\eps'}^{\Sigma(1)/G}}$,
v\'erifiant
\begin{equation}
\forall \bs\in \tube{\R^{\Sigma(1)/G}_{>0}},\quad
\lim_{s\to 0} s^{\rg\Pic(\xs)}\,R(s\,\varphi_0)=0
\end{equation}

La formule \eqref{eq:formule:indic} 
et la suite exacte \eqref{eq:suite:ex} montrent par ailleurs qu'on a 
\begin{equation}
\indic_{\pi\left(\Z^{\Sigma(1)/G}\right),\pi\left(\R^{\Sigma(1)/G}_{\geq 0}\right)}
(\pi(\bs))
=
\,
\card{H^1(G,X(T))}
\,
\indic_{\Pic(\xs),\ceff(\xs)}(\pi(\bs))
+R(\bs).
\end{equation}

On en d\'eduit que  
la fonction
\begin{equation}
s\longmapsto \zeta_{\varphi_0}(s)=\zeta_H(s\,\varphi_0)
\end{equation}
prolonge en une fonction m\'eromorphe 
sur $\tube{\R_{>1-\eps'}}$ avec 
un  p\^ole d'ordre $\rg(\Pic(\xs))$ en $s=1$, 
de terme principal en $s=1$ \'egal \`a 
\begin{equation}\label{eq:terme:princ:arit}
\card{H^1(G,X(T))}\,\alpha^{\ast}(\xs)\,\frac{\card{A(T)}\,C_0}{b(T)}.
\end{equation}

\subsection{Application du lemme technique dans le cas fonctionnel}
\label{subsec:application:lemme:technique:fonc}

\subsubsection{Le cas d'une extension de d\'eploiement non ramifi\'ee}\label{subsubsec:cas:part:nonram}

{\em On suppose dans toute la sous-section \ref{subsubsec:cas:part:nonram} que 
la vari\'et\'e torique $\xs$ est d\'eploy\'ee par une extension 
non ramifi\'ee.}

Rappelons que 
d'apr\`es le corollaire 
\ref{cor:repr:int:fonc:nonram}
on a alors pour tout  
\'el\'ement 
$
\bs
$
de 
$
\tube{\R^{\,\Sigma(1)/G}_{>1}}
$ 
la repr\'esentation int\'egrale 
\begin{equation}
\zeta_H
\left(\bs\right)
=
\frac
{1}
{\log(q)^{\rg\left(X(T)^G\right)}\,b(T)}
\sum_{\chi\in   \ecU_T } 
I_{\chi}(\bs)
\end{equation}
o\`u, pour $\chi\in \ecU_T$, $I_{\chi}(\bs)$ est la fonction donn\'ee par l'int\'egrale
\begin{equation}
\int\limits_{\bz\in \xtgu}
\left(
\prod_{\alpha\in \Sigma(1)/G} 
\cL_{\Ka}
\left(
\acc{\bz}{\da\,\roa}\,q^{\,-\da\,\sa},\chi_{\alpha}
\right)
\right)
\, 
\frQ\left(\chi,\gamma_{\unit}(\bz)\,q^{\,-\bs}\right)
d\bz.
\end{equation}

Soit $P$ la s\'erie formelle
\begin{multline}
P((\za))=
\prod_{\alpha\in\sg}
(1-\za)
\\
\times
\sum_{\chi\in \ecU_T}
\prod_{\alpha\in \sg}
\cL_{\Ka}
\left(
q^{-\,\da}\,\za,\chi_{\alpha}
\right)
\, 
\prod_{v\in\placesde{K}}
Q_{\Sigma,v}
\left(
\chi_{\beta}\left(\pi_{\wb}\right)\,q^{-\db}\,\zb^{\,f_{\wb}}
\right).
\end{multline}
D'apr\`es \eqref{eq:fct:zeta:fonc}, le lemme \ref{lm:chi:non:triv}
et le lemme \ref{lm:frQ}, $P$ a un rayon de convergence sup\'erieur \`a $q^{\frac{1}{2}}$ ;
en outre, pour tout $\bs\in \tube{\R^{\Sigma(1)/G}_{>-\frac{1}{2}}}$ 
et tout $\bz\in \xtgu$, on a d'apr\`es le lemme \ref{lm:frQ} 
et le lemme \ref{lm:accgammaubzza=accbzroa}
\begin{multline}\label{eq:Pqdassa}
\frac{P\left(\acc{\bz}{\da\,\roa} q^{\,-(\da\,\sa)}\right)}
{\produ{\alpha\in\sg}
\left(1-\acc{\bz}{\da\,\roa}\,q^{-\da\,\sa}\right)}
\\
=
\sum_{\chi\in \ecU_T}
\prod_{\alpha\in\sg}
\cL_{\Ka}
\left(
\chia,
\acc{\bz}{\da\,\roa}q^{\,-\da\,(\sa+1)}
\right)
\, 
\frQ\left(\chi,\gamma_{\unit}(\bz)\,q^{\,-\bs-\varphi_0} \right).
\end{multline}
On peut donc \'ecrire 
\begin{equation}\label{eq:zetaHbs+phi0:nonram}
\zeta_{H}(\bs+\varphi_0)
=
\frac
{1}
{\log(q)^{\rg\left(X(T)^G\right)}\,b(T)}
\int
\limits_{\bz \in \xtgu }
\frac
{P
\left(
\acc{\bz}{\da\,\roa}q^{\,-\da\,\sa}
\right)_{\alpha\in\Sigma(1)/G}
}
{
\produ{\alpha\in \sg}
\left(1- \acc{\bz}{\da\,\roa}\,q^{\,-\da\,\sa}\right)
}
\,d\bz.
\end{equation}
\begin{lemme}\label{lm:C0:casnonram}
On a 
\begin{equation}
P(1)
=
\log(q)^{\card{\Sigma(1)/G}}\,
\left(\prod_{\alpha} \da\right)
\,C_0.
\end{equation}
\end{lemme}
\begin{demo}
Comme $\xs$ est d\'eploy\'ee par une extension non ramifi\'ee, 
d'apr\`es le corollaire \ref{cor:nonramAT=0}
on a $\adh{T(K)}\cap T(\ak)=T(\ak)$.
Ainsi on a 
\begin{equation}
C_0
=
\lim_{s\to 1}\, (s-1)^{\,\card{\Sigma(1)/G}}\,
\fourier H(\ind,-s\,\varphi_0).
\end{equation}
D'apr\`es \eqref{eq:fourierHchi:nonram} et la remarque qui suit, on a
\begin{equation}
\fourier H(\ind,-s\,\varphi_0)=
\left(
\prod_{\alpha\in \Sigma(1)/G} 
Z_{\Ka}
\left(
\chia,
q^{\,-\da\,s}
\right)
\right)
\, 
\frQ\left(\ind,q^{\,-s\,\varphi_0 }\right)
\end{equation}
d'o\`u
\begin{equation}
C_0
=
\lim_{s \to 0} 
\,\,
s^{\,\card{\Sigma(1)/G}}\,
\prod_{\alpha\in \sg}
Z_{\Ka}\left(q^{-\da\,(s+1)}\right)
\,
\frQ(\ind,q^{-(s+1)\,\varphi_0})
\end{equation}
D'apr\`es le lemme \ref{lm:kergammast}, $\Ker(\gamma^{\ast})$
s'identifie \`a $A(T)^{\ast}$. Or, 
comme l'extension de d\'eploiement $L/K$ est non ramifi\'ee,
d'apr\`es le corollaire \ref{cor:nonramAT=0} $A(T)$ est trivial,
donc $\Ker\left(\gamma^{\ast}\right)$ \'egalement.
En particulier, si $\chi\in \ecU_T\setminus \{\ind\}$, 
$\chi$ n'est pas dans $\Ker(\gamma^{\ast})$,
et donc l'un au moins des caract\`eres $\chia$ est non trivial.
D'apr\`es le lemme \ref{lm:chi:non:triv}, ceci montre l'\'egalit\'e
\begin{multline}
\lim_{s \to 0} 
\,\,
s^{\,\card{\Sigma(1)/G}}\,
\prod_{\alpha\in \sg}
Z_{\Ka}\left(q^{-\da\,(s+1)}\right)
\,
\frQ(\ind,q^{-(s+1)\,\varphi_0})
\\
=
\lim_{s \to 0} 
\,\,
s^{\,\card{\Sigma(1)/G}}\,
\sum_{\chi\in \ecU_T}
\prod_{\alpha\in\sg}
\ecL_{\Ka}
\left(
\chia,
q^{\,-\da\,(\sa+1)}
\right)
\, 
\frQ\left(\chi,\,q^{\-(s+1)\,\varphi_0} \right).
\end{multline}
Mais d'apr\`es \eqref{eq:Pqdassa}, on a 
\begin{multline}
\lim_{s \to 0} 
\,\,
s^{\,\card{\Sigma(1)/G}}\,
\sum_{\chi\in \ecU_T}
\prod_{\alpha\in\sg}
\ecL_{\Ka}
\left(
\chia,
q^{\,-\da\,(s+1)}
\right)
\, 
\frQ\left(\chi,\,q^{\-(s+1)\,\varphi_0} \right).
\\
=
\lim_{s \to 0} 
\,\,
s^{\,\card{\Sigma(1)/G}}\, \frac{P(q^{-\da\,s})}{\produ{\alpha\in \sg} (1-q^{-\da\,s})}
=
P(1)\,\lim_{s \to 0} \frac{s^{\,\card{\Sigma(1)/G}}}{\produ{\alpha\in \sg} (1-q^{-\da\,s})}
\end{multline}
d'o\`u le r\'esultat.
\end{demo}

Reprenons \`a pr\'esent les notations de la section \ref{sec:eval:int:fonc}.
Le r\^ole de la suite exacte \eqref{eq:exsqMNGamma} est ici jou\'e par la suite
exacte
\begin{equation}
0
\longto 
\DTNS^0
\overset{\pi^{\vee}}
\longto  
\DTqd
\overset{\gamma^{\vee}}
\longto
\DT^0
\longto 
0.
\end{equation}
On pose $I=\Sigma(1)/G$,
et on prend pour base $(\lambda_i)$ la base $\left(\da\,\Da^{\vee}\right)$
de $\DTqd$. On prend pour $a$ la fonction donn\'ee par les coefficients de
la s\'erie formelle $P$. On a donc, d'apr\`es \eqref{eq:zetaHbs+phi0:nonram} et 
le lemme \ref{lm:accgammaubzza=accbzroa}
\begin{equation}
\zeta_H
\left(\varphi_0+\bs\right)
=
\frac
{1}
{\log(q)^{\rg\left(X(T)^G\right)}\,b(T)}
\int\limits_{\bz \in \xtgu}
\ecL_{N,\Lambda,a}
\left(\gamma_{\unit}(\bz)
\,q^{\,-\bs}\right)
\,
\ecL_{N,\Lambda}(\gamma_{\unit}(\bz)\,q^{\,-\bs})d\bz.
\end{equation}
D'apr\`es le lemme \ref{lm:C0:casnonram} on a
\begin{equation}
\ecL_{N,\Lambda,a}(1)
=
\log(q)^{\card{\Sigma(1)/G}}\,
\left(\prod_{\alpha} \da\right)
\,C_0
\neq 0.
\end{equation}
Nous appliquons alors le lemme \ref{lm:tech1} 
(ce qui est licite d'apr\`es \eqref{eq:xtg}), en tenant compte de la
remarque \ref{rm:lm:tech1} et du fait que $\xtg$ est un sous-groupe d'indice
fini de $\left(\DT^0\right)^{\vee}$.

Nous obtenons ainsi le lemme suivant.
\begin{lemme}\label{lm:lmprelimnonram}
La fonction 
\begin{equation}
\bs
\mapsto 
\zeta_{H}(\bs+\varphi_0)
-
\,
\frac
{
\log(q)^{\card{\Sigma(1)/G}}\,
\left(\produ{\alpha} \,\da\right)
\,C_0
}
{\log(q)^{\rg\left(X(T)^G\right)}\,b(T)}
\,
\ecL_{\DTNS^0,\Lambda\,\cap\,\DTNS^0}\left(q^{-\bs}\right)
\end{equation}
est admissible de multiplicit\'e sup\'erieure \`a $1-\rg(\Pic(\xs))$.
\end{lemme}

\begin{lemme}\label{lm:partprinc:fonczeta:nonram}
Il existe un $\eps>0$ tel que
la fonction  
\begin{equation}
s\mapsto \zeta_H(s\,\phi_0)
\end{equation}
se prolonge 
en une fonction m\'eromorphe 
sur
$\tube{\R_{>1-\eps}}$ 
avec un p\^ole d'ordre $\rg(\Pic(\xs))$ en $s=0$,
et v\'erifiant 
\begin{equation}\label{eq:terme:princ:fonc:nonram}
\lim_{s\to 0}s^{\,\rg(\Pic(\xs))}\,\zeta_H(s\,\varphi_0)
=
\card{H^1(G,X(T))}\,\card{\CT}\,\frac{\card{A(T)}\,C_0}{b(T)}\,\alpha^{\ast}(\xs).
\end{equation}
\end{lemme}
\begin{demo}
Du lemme \ref{lm:lmprelimnonram}, on d\'eduit aussit\^ot
qu'il existe un $\eps>0$ tel que
la fonction  
$
s\mapsto \zeta_H(s\,\varphi_0)
$
se prolonge 
en une fonction m\'eromorphe 
sur
$\tube{\R_{>1-\eps}}$ 
avec un p\^ole d'ordre $\rg(\Pic(\xs))$ en $s=0$,
et v\'erifiant 
\begin{multline}\label{eq:limsto0:nonram}
\lim_{s\to 0}s^{\,\rg(\Pic(\xs))}\,\zeta_H(s\,\varphi_0)
=
\log(q)^{\rg(\Pic(\xs))}\, \left(\produ{\alpha} \,\da\right)
\,C_0\,
\frac{1}{b(T)}
\\
\times \lim_{s\to 0}\left[s^{\,\rg(\Pic(\xs))}\,\ecL_{\DTNS^0,\Lambda\,\cap\,\DTNS^0}\left(q^{-s\,\varphi_0}\right)\right].
\end{multline}
Mais on a
\begin{multline}\label{eq:limsto0:eclDTNS0}
\lim_{s\to 0}s^{\,\rg(\Pic(\xs))}\,
\ecL_{\DTNS^0,\Lambda\,\cap\, \DTNS^0}\left(q^{-s\,\varphi_0}\right)\\
\\
=
\frac
{\card{H^1(G,X(T))}\,\card{\CT}}
{\produ{\alpha} \da}
\,\log(q)^{-\rg(\Pic(\xs))}\,\alpha^{\ast}(\xs).
\end{multline}
La relation 
\eqref{eq:limsto0:eclDTNS0}
sera d\'emontr\'e \`a la sous-section suivante (d\'emonstration du lemme \ref{lm:partprinc:fonczeta})
dans le cas o\`u l'extension de d\'eploiement n'est plus suppos\'ee non ramifi\'ee. Supposer l'extension 
de d\'eploiement non ramifi\'ee ne simplifie pas notoirement le calcul (on y gagne simplement le fait que certains
termes apparaissant dans le calcul, comme $\card{A(T)}$ et $\KT$, sont triviaux). 
\end{demo}

\subsubsection{Un cas plus g\'en\'eral}\label{subsubsec:application:lemme:technique:fonc:gene}

{\em On suppose dans toute la sous-section  \ref{subsubsec:application:lemme:technique:fonc:gene}, 
que le $G$-\'eventail  $\Sigma$ est tel que toutes les places de $K$ v\'erifient l'hypoth\`ese \ref{hyp:hinz}.}
Nous expliquerons en appendice comment adapter le raisonnement au cas g\'en\'eral.

Rappelons que d'apr\`es le corollaire \ref{cor:rep:int:fonc}
on a alors pour tout  
\'el\'ement $\bs$ de 
$
\tube{\R^{\,\Sigma(1)/G}_{>1}}
$ 
la repr\'esentation int\'egrale 
\begin{equation}\label{eq:reprintzetafonc}
\zeta_H
\left(\bs\right)
=
\frac
{1}
{\log(q)^{\rg\left(X(T)^G\right)}\,b(T)}
\left(
\frac{\card{A(T)}}{\card{\KT}}\,I_1(\bs)
+
\sum_{\chi\in   \ecU_T  \setminus \Ker\left(\gamma^{\ast}\right)} 
I_{\chi}(\bs)
\right)
\end{equation}
o\`u $I_1(\bs)$ est la fonction donn\'ee par l'int\'egrale
\begin{equation}
\int
\limits_
{\bz\in \xtgu}
\left(
\prod_{\alpha\in \Sigma(1)/G} 
Z_{\Ka}\left(\acc{\bz}{\da\,\roa}\,q^{\,-\da\,\sa}\right)
\right)
\, 
\frf
\left(\gamma_{\unit}(\bz)\,q^{\,-\bs}\right)
\,
d\bz
\end{equation}
et $I_{\chi}(\bs)$ est la fonction donn\'ee par l'int\'egrale
\begin{equation}
\int\limits_{\bz\in \xtgu}
\left(
\prod_{\alpha\in \Sigma(1)/G} 
\ecL_{\Ka}
\left(
\chia,
\acc{\bz}{\da\,\roa}\,q^{\,-\da\,\sa}
\right)
\right)
\, 
\frf\left(\chi,
\gamma_{\unit}(\bz)\,q^{\,-\bs}\right)
\,
d\bz.
\end{equation}

\paragraph{\'Etude de $\boldsymbol{I_{\chi}}$ pour $\boldsymbol{\chi\in \, \ecU_T  
\setminus \Ker\left(\gamma^{\ast}\right)}$.} 
Soit $\chi\in \, \ecU_T  \setminus \Ker\left(\gamma^{\ast}\right)$. 
\'Etudions
l'int\'egrale
\begin{equation}
I_{\chi}(\bs)
=
\int\limits_{\bz \in \xtgu }
\left(
\prod_{\alpha\in \Sigma(1)/G} 
\ecL_{\Ka}
\left(
\chia,
\acc{\bz}{\da\,\roa}\,q^{\,-\da\,\sa}
\right)
\right)
\, 
\frf\left(\chi,\gamma_{\unit}(\bz)\,q^{\,-\bs}\right) 
\,
d\bz.
\end{equation}

Rappelons que d'apr\`es le lemme \ref{lm:frf:chi:rayconv} 
la s\'erie $\frf\left(\chi,\,.\,\right)$ est une s\'erie formelle 
de rayon de convergence sup\'erieur \`a $q^{-\frac{1}{2}}$
et $\DT$-compatible.

Comme $\chi$ n'est pas dans $\Ker\left(\gamma^{\ast}\right)$
l'ensemble $(\Sigma(1)/G)_{\chi}$ des $\alpha\in \Sigma(1)/G$ tel que $\chi_{\alpha}$ 
est non trivial est non vide, 
et, par le lemme \ref{lm:chi:non:triv}, 
pour de tels $\alpha$ la fonction $\ecL_{\Ka}\left(\chia,\,.\,\right)$ 
est un polyn\^ome.
Pour les autres $\alpha$ la fonction $\ecL$ 
correspondante est une fonction z\^eta.

Ceci montre, compte tenu \'egalement de la proposition \ref{prop:fonc:zeta},
que la s\'erie
\begin{multline}
P(\chi,(\za))\eqdef
\left(
\prod_{\alpha\notin (\Sigma(1)/G)_{\chi}}
(1-\za^{\da})
\,
Z_{\Ka}
\left(
q^{\,-\da}\,\za^{\da}
\right)
\right)
\\
\times
\,
\left(
\prod_{\alpha\in(\Sigma(1)/G)_{\chi}}
\ecL_{\Ka}
\left(
\chia,
q^{\,-\da}\,\za^{\da}
\right)
\right)
\, 
\frf\left(\chi,\left(q^{\,-1}\,\za\right)\right)
\end{multline}
a un rayon de convergence sup\'erieur \`a $q^{\,\frac{1}{2}}$
et est $\DT$-compatible.
On a pour tout $\bs\in \tube{\R^{\Sigma(1)/G}_{>0}}$
et tout $\bz\in \xtgu$
\begin{multline}
P\left(\chi,\gamma_{\unit}(\bz)\,q^{\,-\bs}\right)
\\
=
\left(
\prod_{\alpha\notin (\Sigma(1)/G)_{\chi}}
(1-\acc{\bz}{\da\,\roa}\,q^{-\da\,\sa})
\,
Z_{\Ka}
\left(
\acc{\bz}{\da\,\roa}\,q^{\,-\da\,(\sa+1)}
\right)
\right)
\\
\,\,
\times
\,
\left(
\prod_{\alpha\in(\Sigma(1)/G)_{\chi}}
\ecL_{\Ka}
\left(
\chia,
\acc{\bz}{\da\,\roa}
\,
q^{\,-\da\,(\sa+1)}
\right)
\right)
\, 
\frf\left(\chi,\gamma_{\unit}(\bz)\,\,q^{\,-(\bs+\varphi_0)}\right)
\end{multline}

On a donc 
\begin{equation}
I_{\chi}(\bs+\varphi_0)
=
\int
\limits_{\bz \in \xtgu }
\frac
{P
\left(
\chi,
\gamma_{\unit}(\bz)\,q^{\,-\bs} 
\right)
_{\alpha\in\Sigma(1)/G}
}
{
\produ{\alpha\notin (\Sigma(1)/G)_{\chi}}
\left(1-\acc{\bz}{\da\,\roa}\,q^{\,-\da\,\sa}\right)
}
\,d\bz.
\end{equation}

Reprenons les notations de la section \ref{sec:eval:int:fonc}.
Le r\^ole de la suite exacte \eqref{eq:exsqMNGamma} est ici jou\'e par la suite
exacte
\begin{equation}
0
\longto 
\left(\pi(\Ps^G)\right)^{\vee}
\overset{\pi^{\vee}}
\longto  
\left(\Ps^G\right)^{\vee}
\overset{\gamma^{\vee}}
\longto
\left(X(T)^G\right)^{\vee}
\longto 
0.
\end{equation}
On pose 
$
N'=\DTqd
$,
$
I=\Sigma(1)/G
$,
$\wt{I}=(\Sigma(1)/G)\setminus (\Sigma(1)/G)_{\chi}$,
et on prend pour base
$(\lambda_i)$  de $N$ la base $\left(\Da^{\vee}\right)_{\alpha\in \sg}$.
On prend pour $a$ la fonction donn\'ee par les coefficients de la s\'erie formelle $P(\chi,\,.\,)$.
On peut alors \'ecrire
\begin{equation}
I_{\chi}(\bs+\varphi_0)
=
\int\limits_{\bz \in \xtgu}
\ecL_{N,\Lambda,a}
\left(\gamma_{\unit}(\bz)
\,q^{\,-\bs}\right)
\,
\ecL_{N',\wt{\Lambda}}(\gamma_{\unit}(\bz)\,q^{\,-\bs})d\bz
\end{equation}

On applique alors le lemme 
\ref{lm:tech1:bis:tilde}
et on obtient ainsi le lemme suivant.
\begin{lemme}\label{lm:comp:ichi}
La fonction 
\begin{equation}
\bs\mapsto I_{\chi}(\bs+\varphi_0)
\end{equation}
est une fonction admissible de multiplicit\'e au plus 
\begin{equation}
\card{\Sigma(1)/G}-\rg\left(X(T)^G\right)-1=\rg(\Pic(\xs))-1.
\end{equation}
En particulier, il existe un $\eps>0$,
tel que la fonction
$
s\mapsto I_{\chi}((s+1)\varphi_0)
$
se prolonge en une fonction m\'eromorphe
sur $\tube{\R_{>-\eps}}$, 
avec un p\^ole d'ordre au plus 
$
\rg(\Pic(\xs))-1
$ 
en z\'ero.
\end{lemme}

\paragraph{\'Etude de $\boldsymbol{I_1}$.}
Pour m\'emoire, $I_1(\bs)$ est donn\'ee par l'int\'egrale 
\begin{equation}\label{eq:int:I1} 
\int
\limits_{\bz\in \xtgu}
\left(
\prod_{\alpha\in \Sigma(1)/G} 
Z_{\Ka}\left(\acc{\bz}{\da\,\roa}\,q^{\,-\da\,\sa}\right)
\right)
\, 
\frf\left( 
\gamma_{\unit}(\bz)\,q^{\,-\bs}
\right)
\,
d\bz.
\end{equation}

Soit $P$ la s\'erie formelle 
\begin{equation}
P(\bz)\eqdef
\left(
\prod_{\alpha\in \Sigma(1)/G}
(1-q^{-\da}\,\za^{\da})
\,
Z_{\Ka}
\left(
q^{-\da}\,\za^{\da}
\right)
\right)
\times
\, 
\frf\left(\left(q^{-1}\,\za\right)\right).
\end{equation} 
D'apr\`es \eqref{eq:fct:zeta:fonc}, 
le lemme \ref{lm:prop:frf} et le lemme \ref{lm:za^daDT0comp}
$P$ est une
s\'erie formelle de rayon de convergence 
sup\'erieur \`a $q^{\frac{1}{2}}$
et $\DT^0$-compatible.
On a, pour tout $\bs\in \tube{\R_{>-\frac{1}{2}}^{\Sigma(1)/G}}$
et tout $\bz\in \xtgu$
\begin{multline}\label{eq:Pgammazq-s}
\frac
{P\left(\gamma_{\unit}(\bz)\,q^{\,-\bs}\right)}
{
\produ{\alpha\in \Sigma(1)/G} 
\left(1-\acc{\bz}{\,\da\,\roa}\,q^{\,-\da\,\sa}\right)
}
\\
=
\left(
\prod_{\alpha\in \Sigma(1)/G} 
Z_{\Ka}\left(\acc{\bz}{\da\,\roa}\,q^{\,-\da\,(\sa+1)}\right)
\right)
\,
\frf\left( 
\gamma_{\unit}(\bz)\,q^{\,-(\bs+\varphi_0)}
\right).
\end{multline}
On a ainsi, pour tout $\bs\in \tube{\R_{>-\frac{1}{2}}^{\Sigma(1)/G}}$,
\begin{equation}
I_1(\bs+\varphi_0)
=
\int\limits_{\bz \in \xtgu }
\frac
{P\left( \gamma_{\unit}(\bz)\,q^{\,-\bs}\right)
_{
\alpha\in\sg
}
}
{
\produ{\alpha\in \sg}
\left(1-\acc{\bz}{\,\da\,\roa}\,q^{\,-\da\,\sa}\right)
}
d\bz.
\end{equation}

\begin{lemme}\label{lm:C0}
On a
\begin{equation}
P(1)
=
\log(q)^{\card{\Sigma(1)/G}}\,
\left(\prod_{\alpha} \da\right)
\,C_0.
\end{equation}
\end{lemme}
\begin{demo}
D'apr\`es le lemme \ref{lm:prop:frf:bis}, on a pour tout $s\in \tube{\R_{>0}}$
\begin{equation}
\int\limits_{\adh{T(K)}\cap T(\ak)}\!\!H(-(s+1)\,\varphi_0),t)\omega_{T}(t)
=
\left(
\prod_{\alpha\in \sg} 
Z_{\Ka}\left(q^{\,-\da\,(s+1)}\right)
\right)
\frf\left(
q^{\,-(s+1)\,\varphi_0}
\right).
\end{equation}
Par ailleurs on a, d'apr\`es \eqref{eq:Pgammazq-s},
\begin{equation}
\left(
\prod_{\alpha\in \sg} 
Z_{\Ka}\left(q^{\,-\da\,(s+1)}\right)
\right)
\frf\left(
q^{\,-(s+1)\,\varphi_0}
\right)
=
\frac
{P\left(q^{\,-s\,\varphi_0}\right)
}
{
\produ{\alpha\in \sg}
\left(1-\,q^{\,-\da\,s}\right)
}
\end{equation}
On en d\'eduit qu'on a 
\begin{align}
C_0&=\lim_{s\to 0} s^{\card{\sg}}\,\int\limits_{\adh{T(K)}\cap T(\ak)}\!\!H(-(s+1)\,\varphi_0),t)\omega_{T}(t)
\\
&=\lim_{s\to 0} s^{\card{\sg}}\,\frac
{P\left(q^{\,-s\,\varphi_0}\right)
}
{
\produ{\alpha\in \sg}
\left(1-\,q^{\,-\da\,s}\right)
}
\\
&=P(1)\,
\lim_{s\to 0}
\frac
{
s^{\card{\sg}}
}
{
\produ{\alpha\in \sg}
\left(1-\,q^{\,-\da\,s}\right)
}
\end{align}
d'o\`u le r\'esultat annonc\'e.
\end{demo}

Reprenons les notations de la section \ref{sec:eval:int:fonc}.
Le r\^ole de la suite exacte \eqref{eq:exsqMNGamma} est ici jou\'e par la suite
exacte
\begin{equation}
0
\longto 
\DTNS^0
\overset{\pi^{\vee}}
\longto  
\DTqd
\overset{\gamma^{\vee}}
\longto
\DT^0
\longto 
0.
\end{equation}
On pose $I=\Sigma(1)/G$,
et on prend pour base $(\lambda_i)$ la base $\left(\da\,\Da^{\vee}\right)$
de $\DTqd$. On prend pour $a$ la fonction donn\'ee par les coefficients de
la s\'erie formelle $P$. On peut alors \'ecrire 
\begin{equation}
I_1(\bs+\varphi_0)
=
\int\limits_{\bz \in \xtgu}
\ecL_{N,\Lambda,a}\left(\gamma_{\unit}(\bz)\,q^{\,-\bs}\right)\,
\ecL_{N',\Lambda}(\gamma_{\unit}(\bz)\,q^{\,-\bs})d\bz.
\end{equation}
D'apr\`es le lemme \ref{lm:C0:casnonram} on a
\begin{equation}
\ecL_{N,\Lambda,a}(1)
=
\log(q)^{\card{\Sigma(1)/G}}\,
\left(\prod_{\alpha} \da\right)
\,C_0
\neq 0.
\end{equation}
Comme on a $j(N')=\gamma^{\vee}(\DTqd)=\DT^0$, d'apr\`es la d\'efinition \ref{defi:compatible}, 
le fait que $P$ soit $\DT^0$-compatible
signifie que $a$ v\'erifie la condition suivante :
si $a_y$ est non nul alors $j(y)$ est un \'el\'ement de $j(N')$.

Nous appliquons alors le lemme \ref{lm:tech1:bis} 
(ce qui est licite d'apr\`es \eqref{eq:xtg}) 
en tenant compte de la remarque \ref{rm:lm:tech1:bis}
et du fait que $\xtg$ est un sous-groupe d'indice
fini de $\left(\DT^0\right)^{\vee}$.
Nous obtenons ainsi le lemme suivant.
\begin{lemme}\label{lm:comp:i1}
Il existe une fonction admissible $g$
de multiplicit\'e sup\'erieure \`a z\'ero
et v\'erifiant
\begin{equation}
g(0)
=
\log q^{\card{\Sigma(1)/G}}\,
\left(\prod_{\alpha} \da\right)
\,
C_0
\end{equation}
telle
que la fonction  
\begin{equation}
\bs
\mapsto 
I_1(\bs+\varphi_0)-
g(\bs)
\,
\ecL_{\DTNS^0,\Lambda\,\cap\,\DTNS^0}\left(q^{-\bs}\right)
\end{equation}
est admissible de multiplicit\'e sup\'erieure \`a $1-\rg(\Pic(\xs))$.
\end{lemme}
\paragraph{Conclusion} 
\begin{lemme}\label{lm:partprinc:fonczeta}
Il existe un $\eps>0$ tel que
la fonction  
\begin{equation}
s\mapsto \zeta_H(s\,\varphi_0)
\end{equation}
se prolonge 
en une fonction m\'eromorphe 
sur
$\tube{\R_{>1-\eps}}$ 
avec un p\^ole d'ordre $\rg(\Pic(\xs))$ en $s=0$,
et v\'erifiant 
\begin{equation}\label{eq:terme:princ:fonc}
\lim_{s\to 0}s^{\,\rg(\Pic(\xs))}\,\zeta_H(s\,\varphi_0)
=
\card{H^1(G,X(T))}\,\card{\CT}\,\frac{\card{A(T)}\,C_0}{b(T)}\,\alpha^{\ast}(\xs).
\end{equation}
\end{lemme}
\begin{demo}
Des lemmes \ref{lm:comp:ichi}, \ref{lm:comp:i1}, 
et de \eqref{eq:reprintzetafonc}, 
on d\'eduit aussit\^ot
qu'il existe un $\eps>0$ tel que
la fonction  
$
s\mapsto \zeta_H(s\,\varphi_0)
$
se prolonge 
en une fonction m\'eromorphe 
sur
$\tube{\R_{>1-\eps}}$ 
avec un p\^ole d'ordre $\rg(\Pic(\xs))$ en $s=0$,
et v\'erifiant 
\begin{multline}\label{eq:limsto0}
\lim_{s\to 0}s^{\,\rg(\Pic(\xs))}\,\zeta_H(s\,\varphi_0)
\\
=
\frac{\card{A(T)}\,\log (q)^{\card{\Sigma(1)/G}}\,
\left(
\produ{\alpha} \da
\right)
\,C_0}{\KT\,\log(q)^{\rg\left(X(T)^G\right)}\,b(T)}
\,
\lim_{s\to 0}\left[s^{\,\rg(\Pic(\xs))}\,\ecL_{\DTNS^0,\Lambda\,\cap\,\DTNS^0}\left(q^{-s\,\varphi_0}\right)\right].
\end{multline}
Mais on a
\begin{align}
&\quad \lim_{s\to 0}s^{\,\rg(\Pic(\xs))}\,
\ecL_{\DTNS^0,\Lambda\,\cap\, \DTNS^0}\left(q^{-s\,\varphi_0}\right)\\
&=
\frac
{\card{H^1(G,X(T))}\,\card{\CTNS}\,\card{\CT}\,\KT}
{\card{\CTqd}}
\,
\lim_{s\to 0}
s^{\,\rg(\Pic(\xs))}\,
\ecL_{\gamma^{\vee}\left(\DTNS\right),\Lambda\,\cap\, \gamma^{\vee}
\left(\DTNS\right)}
\left(q^{-s\,\varphi_0}\right)\\
&=
\frac
{\card{H^1(G,X(T))}\,\card{\CT}\,\KT}
{\card{\CTqd}}
\lim_{s\to 0}s^{\,\rg(\Pic(\xs))}\,
\ecL_{\Pic(\xs)^{\vee},\ceff(\xs)^{\vee}}\left(q^{-s\,\gamma[\varphi_0]}\right)\\
&=\frac
{\card{H^1(G,X(T))}\,\card{\CT}\,\KT}
{\card{\CTqd}}
\,\log(q)^{-\rg(\Pic(\xs))}\,\alpha^{\ast}(\xs).
\end{align}
La premi\`ere \'egalit\'e vient de \eqref{eq:sous-reseau} et \eqref{eq:indice},
la seconde de \eqref{eq:sous-reseau} et du fait que, par d\'efinition,
$\DTNS$ est d'indice $\card{\CTNS}$ dans $\Pic(\xs)^{\vee}$, 
et la derni\`ere \'egalit\'e vient de \eqref{eq:alpha}.
Comme on a $\card{\CTqd}=\produ{\alpha}\da$ et 
\begin{equation}
\rg(\Pic(\xs))+\rg\left(X(T)^G\right)=\card{\sg},
\end{equation}
on en d\'eduit le r\'esultat annonc\'e.
\end{demo}

\subsection{Calcul du terme principal de la fonction z\^eta des hauteurs}
\label{subsec:calcul:constante}

\subsubsection{Calcul de $C_0$}
\label{subsubsec:calcul:c0}

Nous calculons \`a pr\'esent la constante $C_0$ 
d\'efinie en \eqref{eq:def:c0} 
(et montrons en particulier qu'elle est non nulle).

\begin{lemme}\label{lm:limsto1int=lxtgammaxs}
On a 
\begin{multline}
C_0
=
\lim_{s\to 1}\quad (s-1)^{\,\card{\Sigma(1)/G}}\,
\int\limits_{\adh{T(K)}\cap T(\ak)}H(-s\,\varphi_0,t)\omega_{T}(t)
\\
=
\ell\left(X(T)\right)\,\gamma_H(\xs).
\end{multline}
\end{lemme}
\begin{demo}
La premi\`ere \'egalit\'e est la d\'efinition de $C_0$. 
Montrons la deuxi\`eme \'egalit\'e.
D'apr\`es le scindage
\begin{equation}\label{eq:famous:splitting}
\adh{T(K)}\cap T(\ak)=\adh{T(K)}^{\,\,S}\times T(\ak)^S
\end{equation}
donn\'e par le lemme \ref{lm:scindage}
et la d\'efinition 
\eqref{eq:defomegaT}
de $\omega_T$
on peut \'ecrire pour tout $s\in \tube{\R_{>1}}$
\begin{multline}
\int\limits_{\adh{T(K)}\cap T(\ak)}H(-(s,\dots,s),t)\omega_{T}
\\
=c_{K,\dim(\xs)}\,
\int\limits_{\adh{T(K)}^{\,\,S}} 
\left(\prod_{v \in S}H_v(-(s,\dots,s),t)\right) \otimesu{v\in S}d\mu_v
\\
\times\prod_{v\notin S}H_v(-(s,\dots,s),t) d\mu_v.
\end{multline}
D'apr\`es le lemme \ref{lm:inttkvhv-s}, on en d\'eduit
\begin{multline}
\int\limits_{\adh{T(K)}\cap T(\ak)}H(-(s,\dots,s),t)\omega_{T}
\\
=
c_{K,\dim(\xs)}\,
\int\limits_{\adh{T(K)}^{\,\,S}} 
\left(\prod_{v \in S}H_v(-(s,\dots,s),t)\right) \otimesu{v\in S}d\mu_v
\\
\times
L_S(s,\Ps)\,\prod_{v\notin S}Q_{\Sigma,v}(q_v^{-\lb\, s})_{\beta\in \Sigma(1)/G_v}.
\end{multline}
On en tire
\begin{multline}
\lim_{s\to 1}\quad (s-1)^{\,\card{\Sigma(1)/G}}\,
\int\limits_{\adh{T(K)}\cap T(\ak)}H(-(s,\dots,s),t)\omega_{T}
\\
=
\frac{\ell\left(\Ps\right)}
{\produ{v\in S} L_v(1,\Ps)}
c_{K,\dim(\xs)}
\int\limits_{\adh{T(K)}^{\,\,S}} 
\left(\produ{v \in S}H_v(-(1,\dots,1),t)\right)
\otimesu{v\in S}
d\mu_v
\\
\times\produ{v\notin S}Q_{\Sigma,v}(q_v^{-\lb})_{\beta\in \Sigma(1)/G_v}.
\end{multline}
Par ailleurs, toujours d'apr\`es le lemme \ref{lm:inttkvhv-s},
on a pour $v\notin S$
\begin{align}
Q_{\Sigma,v}(q_v^{-\lb})
&
=
L_v(1,\Ps)^{-1}\,\int\limits_{T(K_v)}H_v(-(1,\dots,1),t)d\mu_v
\\
&
=
\frac{L_v(1,X(T))}
{L_v(1,\Ps)}
\int\limits_{T(K_v)}H_v(-(1,\dots,1),t)\omega_{T,v},
\end{align}
la deuxi\`eme \'egalit\'e provenant de la d\'efinition
\eqref{eq:def:dmuv} de $d\mu_v$.
D'apr\`es le lemme \ref{lm:comp:L:exseq}
et la suite exacte \eqref{eq:resflT}, on a donc
\begin{equation}
Q_{\Sigma,v}(q_v^{-\lb})
=
\int\limits_{T(K_v)}H_v(-(1,\dots,1),t)
\frac{\omega_{T,v}}
{L_v(1,\Pic(\xsl))}
\end{equation}
On en d\'eduit que la limite
\begin{equation}
\lim_{s\to 1}\quad (s-1)^{\,\card{\Sigma(1)/G}}\,
\int\limits_{\adh{T(K)}\cap T(\ak)}H(-(s,\dots,s),t)\omega_{T}
\end{equation}
est \'egale \`a 
\begin{multline}
\ell\left(\Ps\right)
\,
\frac{
\produ{v\in S} L_v(1,X(T))\,L_v(1,\Pic(\xsl))
}
{\produ{v\in S} L_v(1,\Ps)}
\\
\times c_{K,\dim(\xs)}
\int\limits_{\adh{T(K)}^{\,\,S}} 
\left(\produ{v \in S}H_v(-(1,\dots,1),t)\right) \otimesu{v\in S} \frac{\omega_{T,v}}
{L_v(1,\Pic(\xsl))}
\\
\times \produ{v\notin S} \int\limits_{T(K_v)}H_v(-(1,\dots,1),t)
\frac{\omega_{T,v}}
{L_v(1,\Pic(\xsl))}
\end{multline}
D'apr\`es le lemme \ref{lm:restriction:mesure}, on a 
\begin{multline}
c_{K,\dim(\xs)}
\int\limits_{\adh{T(K)}^{\,\,S}} 
\left(\produ{v \in S}H_v(-(1,\dots,1),x)\right) \otimesu{v\in S} \frac{\omega_{T,v}}
{L_v(1,\Pic(\xsl))}
\\
\produ{v\notin S} \int\limits_{T(K_v)}H_v(-(1,\dots,1),x)
\frac{\omega_{T,v}}
{L_v(1,\Pic(\xsl))}
\\
=
c_{K,\dim(\xs)}
\int\limits_{\adh{T(K)}^{\,\,S}} 
 \otimesu{v\in S} \frac{\omega_{\xs,v}}
{L_v(1,\Pic(\xsl))}
\produ{v\notin S} \int\limits_{T(K_v)}
\frac{\omega_{\xs,v}}
{L_v(1,\Pic(\xsl))}.
\end{multline}
En utilisant encore une fois le scindage 
\eqref{eq:famous:splitting}
et la d\'efinition \eqref{eq:defomegaxs} de $\omega_{\xs}$
on trouve
\begin{multline}
c_{K,\dim(\xs)}\,
\int\limits_{\adh{T(K)}^{\,\,S}} 
\otimesu{v\in S}\frac{\omega_{\xs,v}}{L_v(1,\Pic(\xsl))}
\times\prod_{v\notin S}\int\limits_{T(K_v)}\frac{\omega_{\xs,v}}{L_v(1,\Pic(\xsl))}
\\
=
\ell(\Pic(\xsl))^{-1}
\,
\int\limits_{\adh{T(K)}\cap T(\ak)}\omega_{\xs}.
\end{multline}
D'apr\`es le lemme \ref{lm:int_T=int_X}
on a 
\begin{equation}
\int\limits_{\adh{T(K)}\cap T(\ak)}\omega_{\xs}
=
\int\limits_{\overline{\xs(K)}}\omega_{\xs}.
\end{equation}
Finalement, on trouve 
\begin{multline}
\lim_{s\to 1}\quad (s-1)^{\,\card{\Sigma(1)/G}}\,
\int\limits_{\adh{T(K)}\cap T(\ak)}H(-(s,\dots,s),x)\omega_{T}
\\
=
\frac{\ell\left(\Ps\right)}{\ell(\Pic(\xsl))}
\,
\frac{
\produ{v\in S} L_v(1,X(T))\,L_v(1,\Pic(\xsl))
}
{\produ{v\in S} L_v(1,\Ps)}
\int\limits_{\overline{\xs(K)}}\omega_{\xs}.
\end{multline}
D'apr\`es le lemme \ref{lm:comp:L:exseq}
et la suite exacte \eqref{eq:resflT}, on a
\begin{equation}
\ell\left(\Ps\right)=\ell(X(T))\,\ell(\Pic(\xsl))
\end{equation}
et
\begin{equation}
\frac{
\produ{v\in S} L_v(1,X(T))\,L_v(1,\Pic(\xsl))
}
{\produ{v\in S} L_v(1,\Ps)}
=1
\end{equation}
d'o\`u le r\'esultat annonc\'e, au vu de la d\'efinition
\eqref{eq:defgammahxs}
de $\gamma_H(\xs)$.

\end{demo}

\subsubsection{Cas arithm\'etique}
\label{subsubsec:calcul:constante:arit}

D'apr\`es la relation \eqref{eq:A=H1/cha}
du lemme \ref{lm:appfaible},
on a 
\begin{equation}
\card{A(T)}
=
\frac
{\card{H^1(G,\Pic(X_{\Sigma,L}))}}
{\card{\cha(T)}}
.
\end{equation}
De plus, 
par le th\'eor\`eme d'Ono (th\'eor\`eme \ref{thm:ono}) ,
et la d\'efinition \eqref{eq:def:taut:arit} de $\tau(T)$,
on a 
\begin{equation}
\ell(X(T))=\frac{b(T)}{\tau(T)}=\frac{b(T)\,\card{\cha(T)}}{\card{H^1(G,X(T))}}.
\end{equation}

On en d\'eduit, compte tenu du lemme \ref{lm:limsto1int=lxtgammaxs},
\begin{align} 
\card{A(T)}\,C_0
&
=
\card{A(T)}\,\ell\left(X(T)\right)\,\gamma_H(\xs)\\
&=
\frac{b(T)}{\card{H^1(G,X(T))}}
\,\,
\card{H^1(G,\Pic(\xsl))} 
\,
\gamma_H(\xs)
\\
&=
\frac
{b(T)}
{\card{H^1(G,X(T))}}
\,
\beta(\xs)
\,
\gamma_H(\xs).
\end{align}

Dans le cas arithm\'etique, le terme principal de la fonction z\^eta des hauteurs en $s=1$ est donc,
d'apr\`es \eqref{eq:terme:princ:arit},
\begin{multline}
\card{H^1(G,X(T))}\,
\,\alpha^{\ast}(\xs)\,\frac{1}{\card{H^1(G,X(T))}}\,\beta(\xs)\,\gamma_H(\xs)
\\
=\alpha^{\ast}(\xs)\,\beta(\xs)\,\gamma_H(\xs).
\end{multline}

\subsubsection{Cas fonctionnel}
\label{subsubsec:calcul:constante:fonc}

D'apr\`es la relation \eqref{eq:A=H1/cha}
du lemme \ref{lm:appfaible},
on a 
\begin{equation}
\card{A(T)}\,\card{\cha(T)}
=
H^1(G,\Pic(X_{\Sigma,L}))
.
\end{equation}

Par le th\'eor\`eme d'Oesterl\'e (th\'eor\`eme \ref{thm:ono}) 
et la d\'efinition \eqref{eq:def:taut:fonc} de $\tau(T)$,
on a
\begin{equation}
\ell(X(T))=\frac{b(T)}{\tau(T)\,\card{\CT}}=\frac{b(T)\card{\cha(T)}}{\card{\CT}\,\card{H^1(G,X(T))}}.
\end{equation}

On en d\'eduit, compte tenu du lemme \ref{lm:limsto1int=lxtgammaxs}, 
\begin{align} 
\card{A(T)}\,C_0
&=
\card{A(T)}\,\ell\left(X(T)\right)\,\gamma_H(\xs)
\\
&=
\frac{b(T)}
{\card{\CT}\,\card{H^1(G,X(T))}}
\,\,
\card{H^1(G,\Pic(\xsl))} 
\,
\gamma_H(\xs)\\
&=
\frac
{b(T)}
{\card{\CT}\,\card{H^1(G,X(T))}}
\,
\beta(\xs)
\,
\gamma_H(\xs).
\end{align}

Le terme principal de la fonction z\^eta des hauteurs en $s=1$ est donc,
d'apr\`es le lemme \ref{lm:partprinc:fonczeta},
\begin{multline}
\card{H^1(G,X(T))}\,\card{\CT}\,\frac{1}{b(T)}\,\alpha^{\ast}(\xs)\,
\frac
{b(T)}
{\card{\CT}\,\card{H^1(G,X(T))}}
\,
\beta(\xs)
\,
\gamma_H(\xs)
\\
=\alpha^{\ast}(\xs)\,\beta(\xs)\,\gamma_H(\xs).
\end{multline}
Ceci ach\`eve la d\'emonstration 
du th\'eor\`eme \ref{thm:theoprinc}.

\appendix

\section{Appendice : le cas o\`u l'hypoth\`ese \ref{hyp:hinz} n'est pas v\'erifi\'ee}\label{appendix}

Dans le cas fonctionnel, nous indiquons \`a pr\'esent comment adapter 
ce qui pr\'ec\`ede au cas o\`u certaines places de $v$ ne v\'erifient pas l'hypoth\`eses
\ref{hyp:hinz}.
Soit $e$ un entier strictement positif divisibles par tous les $e_v$.

Ainsi, pour tout $v\in \placesde{K}$, pour tout $\varphi\in \eplsg$
et pour tout $t\in T(K_v)$, $f_v\,\accsv{\varphi}{t}$ est entier 

On peut donc d\'efinir  un accouplement
\begin{equation}\label{e:eq:hauteur:locale:log:bis}
\frH_{e,v}\,:\,
\map{\eplsgca\times T(K_v)}{\C^{\ast}}{(\psi,t)}
{\accsv{\psi}{t}^{\,f_v}}
\end{equation}
v\'erifiant
\begin{equation}\label{e:eq:hvphi=frhvqphi}
\forall \varphi\in \eplsgc=\plsgc,
\quad
\forall t\in T(K_v),
\quad
H_v(\varphi,t)=\frH_{e,v}\left(q^{\,\varphi},t\right).
\end{equation}

Dans toute la suite, on identifiera $\plsgc$ \`a $\csg$ non plus au moyen de la base
$(\Da)$, mais de la base $(e\,\Da)$. On notera $(\sa')$ l'\'el\'ement de $\csg$ correspondant
\`a un \'el\'ement de $\bs$ de $\plsgc$ via cette identification. Ainsi on a $\sa=e\,\sa'$ pour tout $\alpha\in \sg$.

Soit 
\begin{equation}
P=\sum_{(\na)\in \N^{\sg}} a_{(\na)}\prod \za^{\na}
\end{equation}
une s\'erie formelle, et $\bz$ un \'el\'ement de $\casg$ tel que la s\'erie d\'efinissant $P(\bz)$ converge
absolument. Si on voit $\bz$ comme un \'el\'ement de $\epsgca$, on a donc
\begin{equation}\label{e:eq:P(bz)=}
P(\bz)=\sum_{(\na)\in \N^{\sg}}a_{(\na)} \acc{\bz}{\frac{1}{e}\sum \na \,\Da^{\vee}}
\end{equation}
Pour tout $\bz\in \extgca$ tel que la s\'erie
d\'efinissant $P(\gamma_{\ca}(\bz))$ converge absolument, on a donc
\begin{equation}
P\left(\gamma_{\ca}(\bz)\right)=\sum_{(\na)\in \N^{\sg}} a_{(\na)}\acc{\bz}{\frac{1}{e} \sum \na \roa}.
\end{equation}

La d\'efinition de la compatibilt\'e est modifi\'ee comme suit : pour tout  sous-groupe $M$ de $\frac{1}{e}\,\xtgd$
et tout $(\na)\in \Nsg$, 
le mon\^ome $\produ{\alpha\in \sg}\za^{\na}$ 
est dit $M$-compatible si on a 
\begin{equation}
\gamma^{\vee}\left(\frac{1}{e}\,\sum \na\Da^{\vee}\right)\in M,
\end{equation}
autrement dit si on a 
\begin{equation}
\frac{1}{e}\,\sum \na\,\roa \in M.
\end{equation}
On voit qu'un mon\^ome $f$ est $M$-compatible si et seulement s'il existe
un \'el\'ement $m$ de $M$ v\'erifiant
\begin{equation}
\forall \bz\in \extgca,\quad 
f\left(\gamma_{\ca}(\bz)\right)=\acc{\bz}{m}
\end{equation}
soit encore
\begin{equation}
\forall \bz\in \xtgca,\quad 
f\left(\gamma_{\ca}(\bz^e)\right)=\acc{\bz^e}{m}.
\end{equation}
Par ailleurs pour tout $\alpha\in \sg$, le mon\^ome $\za^{\,e\,\da}$ est $\DT^0$-compatible.

On peut v\'erifier alors que tous les r\'esultats de la sous-section 
\ref{subsubsec:transfolocales:casfonc} (\`a partir du lemme \ref{lm:rel:hauteur:carac:fonc})
et de la sous-section \ref{subsubsec:casfonc:frQ} restent valable, en rempla\`c;ant dans les \'enonc\'es 
et d\'efinitions toutes
les occurences de $\frH_{v}$ (respectivement de $\chi_{\bz}$ pour $\bz\in \xtgu$, respectivement
de $\gamma_{\unit}(\bz)$ pour $\bz\in \xtgu$) par $\frH_{e,v}$ (respectivement $\chi_{\bz^e}$, 
respectivement $\gamma_{\unit}(\bz)$) et toutes les occurences des expressions du type <<de rayon
de convergence sup\'erieur \`a $q^{\,a}$>> avec $a\in \R$ par <<de rayon
de convergence sup\'erieur \`a $q^{\,\frac{a}{e}}$>>.

On obtient pout tout $\bs\in \tube{\R^{\Sigma(1)/G}_{>\frac{1}{e}}}$ 
et tout $\chi \in \dualtop{T(\ak)/\K(T)}$
\begin{align}\label{e:eq:i1:fonc:bis}
I_1(\bs)
&=
\int\limits_{\xtgu}\left[\int\limits_{\adh{T(K)}\cap T(\ak)}H(-\bs,t)\,\chi_{\bz}(t)\,\omega_{T}(t)\right]d\bz.
\\
&=
\int\limits_{\xtgu}\left[\int\limits_{\adh{T(K)}\cap T(\ak)}H(-\bs,t)\,\chi_{\bz^e}(t)\,\omega_{T}(t)\right]d\bz.
\\
&
=
\int
\limits_{\xtgu}
\!\!
\left(
\prod_{\alpha\in \Sigma(1)/G} 
Z_{\Ka}
\left(\acc{\bz^e}{\da\,\roa}\,q^{\,-\da\,e\,\sa'}\right)
\right)
\frf\left(\gamma_{\unit}(\bz^e)\,q^{\,-\bs}\right)
d\bz
\end{align}
et
\begin{align}\label{e:eq:ichis_fonc:bis}
I_{\chi}(\bs)
&
=
\int\limits_{\xtgu} 
\fourier H\left(\chi_{\bz}\,\chi,-\bs \right)
d\bz
\\
&
=
\int\limits_{\xtgu} 
\fourier H\left(\chi_{\bz^e}\,\chi,-\bs \right)
d\bz
\\
&
=
\int\limits_{\xtgu }
\!\!
\left(
\prod_{\alpha\in \Sigma(1)/G} 
\ecL_{\Ka}
\left(
\chia,
\acc{\bz^e}{\da\,\roa}\,q^{\,-\da\,e\,\sa'}
\right)
\right)
\, 
\frf\left(\chi,
\gamma_{\unit}(\bz^e)\,q^{\,-\bs}\right)
d\bz.
\end{align}

On peut alors adapter le raisonnement de la sous-section \ref{subsubsec:application:lemme:technique:fonc:gene}.
Pour $I_{\chi}$, le r\^ole de la suite exacte \eqref{eq:exsqMNGamma} est jou\'e par la suite
exacte
\begin{equation}
0
\longto 
\left(\pi(e\,\Ps^G)\right)^{\vee}
\overset{\pi^{\vee}}
\longto  
\frac{1}{e}\,\left(\Ps^G\right)^{\vee}
\overset{\gamma^{\vee}}
\longto
\frac{1}{e}\,\left(X(T)^G\right)^{\vee}
\longto 
0
\end{equation}
et pour $I_1$,
par la suite
exacte
\begin{equation}
0
\longto 
\DTNS^0
\overset{\pi^{\vee}}
\longto  
\DTqd
\overset{\gamma^{\vee}}
\longto
\DT^0
\longto 
0.
\end{equation}

Comme dans le cas d\'ej\`a trait\'e, on en d\'eduit qu'il existe une fonction admissible $g$
de multiplicit\'e sup\'erieure \`a z\'ero
et v\'erifiant
\begin{equation}
g(0)
=
\log q^{\card{\Sigma(1)/G}}\,
\left(\prod_{\alpha} \da\right)
\,
C_0
\end{equation}
telle
que la fonction  
\begin{equation}
\bs
\mapsto 
\zeta_H(\bs+\varphi_0)-
g(\bs)
\,
\ecL_{\DTNS^0,\Lambda\,\cap\,\DTNS^0}\left(q^{-\bs}\right)
\end{equation}
est admissible de multiplicit\'e sup\'erieure \`a $1-\rg(\Pic(\xs))$.

\printindex[not]

\printindex[def]

\begin{flushleft}
I.R.M.A.R\\
Campus de Beaulieu\\
35042 Rennes cedex \\
France \\
\verb|david.bourqui@univ-rennes1.fr|
\end{flushleft}

\begin{thebibliography}{40}

\bibitem[Ar]{Artin:uber}
\textsc{E. Artin}.
\newblock 
\"Uber eine neue {A}rt von {L}-{R}eihen.
\newblock {\em Abh. Math. Sem. Univ. Hamburg}, \textbf{3}, 1924.
\newblock p. 89-108.

\bibitem[BaTs1]{BaTs:anis}
\textsc{V.V. Batyrev, Y. Tschinkel}.
\newblock Rational points of bounded height on compactifications of anisotropic
  tori.
\newblock {\em Int. Math. Res. Notices}, \textbf{12}, 1995.
\newblock p. 591-635.


\bibitem[BaTs2]{BaTs:manconj}
\textsc{V.V. Batyrev, Y. Tschinkel}.
\newblock Manin's conjecture for toric varieties.
\newblock {\em J. of Algebraic Geometry}, \textbf{7}, 1998.
\newblock p. 15-53.

\bibitem[BaTs3]{BaTs:conicbundles}
\textsc{V.V. Batyrev, Y. Tschinkel}.
\newblock Rational points on some Fano cubic bundles.  
\newblock {\em  C. R. Acad. Sci. Paris S\'er. I Math.}, \textbf{323}, 1996.
\newblock p. 41-46.

\bibitem[BaTs4]{BaTs:hzf}
\textsc{V.V. Batyrev, Y. Tschinkel}.
\newblock Height zeta functions of toric varieties.
\newblock {\em  J. Math. Sci.}, \textbf{82}, 1996.
\newblock p. 3220-3239.


\bibitem[BoSp]{BoSp:rationality}
\textsc{A. Borel, T.A. Springer}.
\newblock Rationality properties of linear algebraic groups. II. 
\newblock {\em  T\^ohoku Math. J. (2)}, \textbf{20}, 1968.
\newblock p. 443-497.

\bibitem[Bki1]{Bki:spec}
\textsc{N. Bourbaki}.
\newblock {\em \'El\'ements de math\'ematiques XXXII : Th\'eories spectrales. Chapitres 1 et 2}.
\newblock Actualit\'es scientifiques et industrielles \textbf{1332}, Hermann, Paris, 1967.

\bibitem[Bki2]{Bki:diff}
\textsc{N. Bourbaki}.
\newblock {\em \'El\'ements de math\'ematiques XXXVI : Vari\'et\'es diff\'erentielles et analytiques. 
Fascicule de r\'esultats (Paragraphes 8 \`a 15)}.
\newblock Actualit\'es scientifiques et industrielles \textbf{1347}, Hermann, Paris, 1971.


\bibitem[Bo1]{Bou:jnt}
\textsc{D. Bourqui}.
\newblock Fonction z\^eta des hauteurs des surfaces de Hirzebruch dans le cas fonctionnel.
\newblock {\em J. of Number Theory}, \textbf{94}, 2002.
\newblock  p. 343-358.

\bibitem[Bo2]{Bou:crelle}
\textsc{D. Bourqui}.
\newblock Fonction z\^eta des hauteurs des vari\'et\'es toriques d\'eploy\'ees dans le cas fonctionnel.
\newblock {\em J. reine angew. Math.}, {\textbf{562}}, 2003.
\newblock p.~171-199.

\bibitem[Bo3]{Bou:PhD}
\textsc{D. Bourqui}.
\newblock Fonctions z\^eta des hauteurs des vari\'et\'es toriques en caract\'eristique positive.
\newblock Th\`ese de doctorat, Grenoble, 2003.
\newblock Disponible \`a l'URL \verb|http://tel.ccsd.cnrs.fr/tel-00004008|

\bibitem[Bro]{Bro:manin_dim_2}
\textsc{T.D. Browning}.
\newblock 
The Manin conjecture in dimension 2.
\newblock Preliminary version of lecture notes for the "School and conference on analytic number theory", 
ICTP, Trieste, 23/04/07-11/05/07,
\newblock \verb|arXiv:0704.1217v1|.

\bibitem[CLTs]{CLTs:fibres}
\textsc{A. Chambert-Loir, Y. Tschinkel}.
{\em \textnormal{Fonctions z\^eta des hauteurs des espaces fibr\'es} 
in Rational points on algebraic varieties, 
\textnormal{E. Peyre \& Y. Tschinkel, eds.}}
\newblock Birkh\"auser (2001).
\newblock  p. 71-115.


\bibitem[CTHaSk]{CTHaSk:comp}
\textsc{J.-L. Colliot-Th\'el\`ene, D. Harari, N. Skorobogatov}.
\newblock  Compactification \'equivariante d'un tore (d'apr\`es {B}rylinski et {K}\"unnemann)
\newblock  {\em Expo. Math.},  \textbf{23}, no. 2, 2005. 
\newblock p. 161--170. 


\bibitem[CTSa]{CTS:Requiv}
\textsc{J.-L. Colliot-Th\'el\`ene, J.-J. Sansuc}.
\newblock La {R}-\'equivalence sur les tores.
\newblock {\em Ann. scient. Ec. Norm. Sup.}, \textbf{10}, 1977.
\newblock p. 175-230.

\bibitem[CTSu]{CTSu}
\textsc{J.-L. Colliot-Th\'el\`ene, V. Suresh}.
\newblock Quelques questions d'approximation faible pour les tores alg\'ebriques. 
\newblock {\`A} para\^\i tre aux {\em Annales de l'Institut Fourier}.

\bibitem[Dr]{Dr}
{\textsc{P.K.J. Draxl}}.
\newblock L-{F}unktionen {K}-algebraischer {T}ori.
\newblock {\em J. of Number Theory}, \textbf{3}, 1971.
\newblock  p. 444-467.

\bibitem[FrMaTs]{FMT}
\textsc{J. Franke, Y. Manin, Y. Tschinkel}.
\newblock Rational points of bounded height on {F}ano varieties.
\newblock {\em Invent. Math.}, {\textbf{95}}, 1989.
\newblock p. 421-435.

\bibitem[Fu]{Ful:toric}
{\textsc{W. Fulton}}.
\newblock {\em Introduction to toric varieties}.
\newblock Annals of Mathematics Studies {\textbf{31}}, Princeton University
  Press, 1993.

\bibitem[Ha]{Has}
\textsc{H. Hasse}.
\newblock {\em Number theory}.
\newblock Grundlehren der mathematischen Wissenschaften {\textbf{229}}, Springer, 1980.

\bibitem[HoNa]{HoNa:cohomology_cft}
\textsc{G. Hochschild, T. Nakayama}.
\newblock Cohomology in class field theory.
\newblock{\em Ann. of Math.}, \textbf{55,} 1952. 
\newblock p. 348-366.

\bibitem[LaYe]{LaiYeu}
\textsc{K.F. Lai, K.M. Yeung}.
\newblock Rational points in flag varieties over function fields.
\newblock {\em J. of Number Theory}, \textbf{95}, 2002.
\newblock p. 142-149.

\bibitem[Na]{Na}
\textsc{T. Nakayama}.
\newblock Cohomology of class field theory and tensor product modules {I}.
\newblock{\em Ann. of Math.}, \textbf{65}, 1957.
\newblock p. 265-267.

\bibitem[Od]{Oda:conv}
{\textsc{T. Oda}}.
\newblock {\em Convex bodies and algebraic geometry}.
\newblock Ergebnisse der Mathematik und ihrer Grenzgebiete {\textbf{15}},
  Springer-Verlag, 1988.

\bibitem[Oe]{Oes:invent}
\textsc{J. Oesterl\'e}.
\newblock Nombre de {T}amagawa et groupes unipotents en caract\'eristique $p$.
\newblock {\em Invent. Math.}, \textbf{78}, 1984.
\newblock p. 13-88.

\bibitem[On1]{Ono_aritag}
\textsc{T. Ono}.
\newblock On some arithmetic properties of linear algebraic groups.
\newblock {\em Ann. of Math.}, \textbf{70}, 1959.
\newblock p. 266-290.

\bibitem[On2]{Ono:algtor}
\textsc{T. Ono}.
\newblock Arithmetic of algebraic tori.
\newblock {\em Ann. of Math.}, \textbf{74}, 1961.
\newblock p. 101-139.

\bibitem[On3]{Ono:tamnum}
\textsc{T. Ono}.
\newblock On the {T}amagawa number of algebraic tori.
\newblock {\em Ann. of. Math.}, \textbf{78}, 1963.
\newblock p. 47-73.

\bibitem[Pe1]{Pey:duke}
\textsc{E. Peyre}.
\newblock Hauteurs et mesures de {T}amagawa sur les vari\'et\'es de {F}ano.
\newblock {\em Duke Math. Journal}, {\textbf{79}}, 1995.
\newblock p. 101-218.

\bibitem[Pe2]{Pey:circle}
\textsc{E. Peyre}.
{\em \textnormal{Torseurs universels et m\'ethode du cercle}
in Rational points on algebraic varieties}.
\newblock Progr. Math. {\textbf{79}},
\newblock Birkha\"user, 2001.
\newblock  p. 221-274.

\bibitem[Pe3]{Pey:prepu:drap}
\textsc{E. Peyre}.
\newblock Points de hauteur born\'ee sur les vari\'et\'es de drapeaux en
  caract\'eristique finie.
\newblock Pr\'epublication, \verb|math.NT/0303067|.

\bibitem[Pe4]{Pey:ecoledete}
\textsc{E. Peyre}.
\newblock Etude asymptotique des points de hauteur born\'ee.
\newblock Notes de l'\'ecole d'\'et\'e sur les vari\'et\'es toriques, Grenoble,
  juin 2000.

\bibitem[Pe5]{Pey:bki}
{\textsc{E. Peyre}}.
\newblock Points de hauteur born\'ee et g\'eom\'etrie des vari\'et\'es, d'apr\`es {Y}. Manin et al.
\newblock S\'eminaire Bourbaki \textbf{891}, juin 2001.

\bibitem[Pe6]{Pey:bordeaux}
\textsc{E. Peyre}.
\newblock Points de hauteur born\'ee, topologie ad\'elique et mesure de {T}amagawa.
\newblock {\em J. Th. Nombres de Bordeaux}, {\textbf{15}}, 2003.
\newblock p. 319-349. 

\bibitem[PlRa]{PlRa:alggp}
\textsc{V. Platonov, A. Rapinchuk}.
\newblock {\em Algebraic groups and number theory}.
\newblock Pure and Applied Mathematics, {\textbf{139}}, Academic Press, Boston, 1994.

\bibitem[Ra]{Rad:phrag}
\textsc{H. Rademacher}.
\newblock On the {P}hragm\'en-{L}indel\"of theorem and some applications.
\newblock {\em Math. Z.}, \textbf{72}, 1959.
\newblock p. 192-204.

\bibitem[Sa]{Sal:tammes}
\textsc{P. Salberger}.
\newblock Tamagawa measure on universal torsors and points of bounded height on {F}ano varieties.
\newblock {\em Ast\'erisque}, \textbf{251}, 1998.
\newblock p. 91-258.

\bibitem[Se1]{Ser:corps}
\textsc{J.-P. Serre}.
\newblock {\em Corps locaux}.
\newblock Hermann, 1962.

\bibitem[Se2]{Ser:LMWT}
{\textsc{J.-P. Serre}}.
\newblock {\em Lectures on the Mordell-Weil theorem, \textnormal{translated and edited by M. Brown}}. 
\newblock Aspect of Mathematics, Vieweg, 1989.

\bibitem[Vo]{Vos:proj}
\textsc{V.E. Voskresenskii}.
\newblock Projective invariant Demazure models.
\newblock {\em Math. USSR Izv.}, \textbf{20}, no. 2, 1983.
\newblock p. 189-202.

\bibitem[Vo]{Vos:bir}
\textsc{V.E. Voskresenskii}.
\newblock Invariants birationnels des tores alg\'ebriques (en russe).
\newblock {\em Uspehi Mat. Nauk}, \textbf{30}, 1975.
\newblock p. 207--208.
 
\bibitem[Wa]{Wan}
{\textsc{D. Wan}}.
{\em \textnormal{Heights and zeta functions in function fields} 
in The arithmetic of function fields, \textnormal{D. Goss, D.R. Hayes \& M.I. Rosen, eds.}}
\newblock Mathematical Research Institute Publications, Ohio State University,
\newblock de Gruyter, 1992.
\newblock  p. 455-463.

\bibitem[We1]{Wei:BNT}
\textsc{A. Weil}.
\newblock {\em Basic number theory}.
\newblock Springer Verlag, 1967.

\bibitem[We2]{Wei:AAG}
\textsc{A. Weil}.
\newblock {\em Adeles and algebraic groups}.
\newblock Progress in mathematics \textbf{23}, Birkha\"user, 1982.

\end{thebibliography}
\end{document}